# Twenty Female Mathematicians

Hollis Williams

*Acknowledgements*

The author would like to thank Alba Carballo González for support and encouragement.

# *Preface*

In this piece of writing, I make an attempt to briefly discuss the work of twenty female mathematicians. I am aware that in some sections I have gone off on tangents a bit more or not said as much about the actual work which the mathematician does or the exact details of their achievements. In most cases, it would have been difficult to give more than a few details and it is completely impossible to give a good survey of someone's work within the confines of 10 pages per person along with added context: what I have tried to do is really give a *flavour* of the kind of things that the person works on and hopefully inspire you to read on further. The difficulty was particularly obvious in algebraic geometry.

By the time one has given all the definitions which are needed, one has run out of space and not had the opportunity to actually apply them to some interesting problems, so I had to try my best to give concise, roughly correct definitions whilst also trying to give an idea of the type of problem which can be tackled with the theory, even if there was not space for more specific examples. In some cases I have skipped some of the definitions, trusting that the reader will look them up in the references if unsure. The work does not necessarily represent the people who I think are the twenty greatest female mathematicians of all time and reflects my own interests in mathematics and the things which I felt I might be to talk about semi-competently for a few pages.

Given the range of topics and problems covered, there are likely to be errors and conceptual mistakes. Please email me to discuss mistakes and I will be happy to upload a corrected version.

Email: Hollis.Williams@warwick.ac.uk

# *Sofia Kovalevskaya*

Kovalevskaya was a nineteenth-century Russian mathematician, known mainly for her achievements in mechanics and the theory of PDEs. She gained qualifications and held a number of posts which were traditionally reserved only for men and was amongst the first women to work as an editor of a scientific journal[1]. Probably her most substantial mathematical contribution was the Cauchy-Kovalevskaya theorem, which belongs to PDE theory. PDEs (partial differential equations) are equations which model the rate of change of a quantity with respect to two more other quantities which are also changing, whereas ODEs (ordinary differential equations) model the rate of change with respect to one quantity only (usually a time variable). All of science, engineering and applied mathematics depends on differential equations. PDEs were originally brought into being as tools to describe physical phenomena such as wave motion and heat flow (a role which they still play), but over time PDE theory has also developed in another direction into a rich branch of mathematical theory in its own right (without completely losing sight of the physical origin of the equations).

A famous and well-studied example is the one-dimensional wave equation:

$$\frac{\partial^2 f}{\partial t^2} = c^2 \frac{\partial^2 f}{\partial x^2}.$$

One sees that there are no mixed derivatives in this expression (ie. derivatives with respect to more than one variable), but it should also be pointed out that mixed derivatives can in fact occur in equations which model physical processes: there are versions of the diffusion equation which contain mixed derivative terms, for example. A very interesting development in the twentieth century was the notion that one might be able to 'weaken' the notion of a solution to obtain weak solutions which exist where our traditional idea of a solution does not exist: this is rather surprising from the viewpoint of classical PDE theory, since physical intuition would seem to tell us that the notion of a solution has to be a concrete one which cannot be weakened in this way.

Essentially, we take a PDE and rewrite it so that it is in 'weak' form: a solution to this equation is known as a weak solution to the original PDE. A strong solution is automatically a weak solution, but it is unfortunately not always obvious that a weak solution should also be strong. I should point out that the act of rewriting an equation is not necessarily by itself a profound one. The notion of having an equation or solution with particular types of operator in it and rewriting it so those operators do not appear is fairly common and trivial.

[1] Ann Koblitz, *A Convergence of Lives: Sofia Kovalevskaya: Scientist, Writer, Revolutionary* (New Jersey, USA: Rutgers University Press, 1993).

As an example, one can obtain a fundamental solution to a particular system of equations modelling dilute gas flow, where the solution represents the physical pressure due to a gas:

$$p = -\mathbf{f}.\nabla\phi + \mathbf{G}(\nabla\nabla\phi) - \frac{8\gamma_2^2 \mathrm{Kn} g}{15}\mu_{\gamma_2} + \mathbf{G}(\nabla\nabla\mu_{\gamma_2}),$$

where $\mathbf{G}$ is some arbitrary trace-free symmetric 2-tensor, $\phi$ is the fundamental solution to the PDE known as Laplace's equation (derived essentially by considering the symmetries which the solution must have), $\mu_{\gamma_2}$ is the fundamental solution to the Helmholtz equation (similar to the Laplace equation) and $\mathrm{Kn}$ is a dimensionless parameter similar to the Reynolds number[2]. One notes the repeated appearance of the $\nabla$ operator, known as 'nabla' or 'the gradient'. Using the definition of $\nabla$ as a differential operator, it is a trivial operation to rewrite the pressure term so that it no longer includes nabla:

$$\begin{aligned} p = \mathbf{f}.\left(-\frac{\mathbf{r}}{4\pi|\mathbf{r}|^3}\right) + G_{ij}r_i r_j\left(\frac{3}{4\pi|\mathbf{r}|^5}\right) - \frac{\frac{40}{6\mathrm{Kn}}g}{15}\frac{1}{4\pi|\mathbf{r}|}e^{-\gamma_2 r} \\ + G_{ij}r_i r_j \frac{e^{-\gamma_2 r}\left(5|\mathbf{r}| + 3\mathrm{Kn}\left(6\mathrm{Kn} + \sqrt{30}|\mathbf{r}|\right)\right)}{24\mathrm{Kn}^2\pi|\mathbf{r}|^5}. \end{aligned}$$

This looks a bit messier, but the differential operator no longer appears and so it might be faster to process the solution with a computer. However, in this case, there is something more powerful going on, since the reformulation we propose allows us to keep a notion of a solution to the equation, even if it cannot be differentiated a sufficient number of times to be counted as a classical solution!

At the bottom of things, a variational technique is being employed to find existence of solutions. As the simplest possible example in what is called the 'elliptic' theory, one takes a standard second-order elliptic PDE such as the Poisson equation defined for a boundary value problem on a three-dimensional ball $B$ and then rewrites it as:

$$\int_B -\Delta u \cdot \phi \,\mathrm{d}x = \int_B f \cdot \phi \,\mathrm{d}x,$$

where $\phi$ belongs to $C^\infty(B)$ and $u$ solves the original PDE. The fact that this reformulation is equivalent to the PDE we had before is due to what is known as the fundamental lemma of calculus of variations. The next step is to integrate by parts to get

$$\int_B -\Delta u \cdot \phi = \int_B Du \cdot D\phi =: (u, \phi)_{H_0^1},$$

where $D$ is the gradient vector of derivatives and we have defined a new inner product in the final step. Similarly, in a different variational argument, one can find a function $\bar{u}$ which minimises an energy functional similar to the Dirichlet energy over the space of admissible functions:

---

[2] Rory Claydon, Abhay Shrestha, Anirudh Rana, James Sprittles and Duncan Lockerby, 'Fundamental solutions to the regularised 13-moment equations: efficient computation of three-dimensional kinetic effects'. Journal of Fluid Mechanics, 833, R4 (2017).

$$\varepsilon[u] := \int_{\Omega} |\nabla u(\mathbf{x})|^2 \ \mathrm{d}\mathbf{x}$$

and then deduce from the following result (the weak form of Laplace's equation)

$$\int_{\Omega} \nabla \bar{u} \cdot \nabla v \ \mathrm{d}\mathbf{x} = 0$$

and the fundamental lemma of calculus of variations that $\bar{u}$ satisfies the Laplace equation[3].

The point of this reformulation is that it makes sense even if $u$ is not in $C^2$ and cannot be differentiated twice, whereas the classical theory tells us that we must be able to differentiate the solution twice by definition. One then defines a linear functional $L$ which maps $C_0^\infty(B)$ to $\mathbb{R}$ as follows:

$$L(\phi) = \int_B f \cdot \phi \ \mathrm{d}x.$$

One can see that this is a bounded linear functional using a standard inequality (the Friedrichs-Poincaré inequality). A linear operator $T$ between two normed spaces $X$ and $Y$ is said to be bounded if the norm of $T(x)$ is less than or equal to a non-negative constant multiplied by the norm of $x$.

$$\|T(x)\|_Y \leq c\|x\|_X.$$

In our weak reformulation, we have now changed the problem of solving the Poisson equation to the problem of solving

$$L(\phi) = (u, \phi)_{H_0^1}.$$

One would like to apply the Riesz representation theorem to our bounded linear functionals, but they are bounded operators on the space $\left(C_0^\infty(B), (\cdot,\cdot)_{H_0^1}\right)$, which is not complete (ie. you can find a Cauchy sequence of functions in $C_0^\infty(B)$ which converges to a function which is not in $C_0^\infty(B)$).

The Riesz theorem for representation of functionals on Hilbert spaces states that every bounded linear functional $f$ on a Hilbert space $H$ can be represented in a simple way via the inner product:

$$f(x) = \langle x, z \rangle,$$

where $z$ depends on and is determined uniquely by $f$ and has the norm

$$\|f\| = \|z\|.$$

[3] Michael Renardy and Robert Rogers, *An Introduction to Partial Differential Equations* (New York: Springer, 2004).

This is a fundamental theorem and more details can be found in any book or set of lecture notes on functional analysis[4]. A more useful way of putting it is that for a Hilbert space $H$ and a linear functional in the dual space $H^*$, there exists a unique element $u$ in $H$ such that

$$L(\phi) = (u, \phi)$$

for any $\phi$ in $H$. $L$ is then a linear isometric isomorphism from the dual space to the original Hilbert space.

One example of an application of the theorem would be to gradient flows. Imagine that we start with a linear functional $f$ which maps a space $X$ to $\mathbb{R}$. If we complete that space to a new space $\hat{X}$ by equipping it with an inner product and taking the closure, then $f$ lives in the dual space $(\hat{X})^*$. Recall that Hilbert spaces are complete by definition. Riesz representation tells us that

$$f(x) = \langle \nabla f, x \rangle,$$

for any $x$ in $\hat{X}$ and any $\nabla f$ in $(\hat{X})^*$. Obviously, that inner product could be the $L^2$ inner product, depending on the Hilbert space. $f$ can always be written as $\langle \nabla f, x \rangle$ and f having gradient just means that $\nabla f \in X$. An $L^2$ gradient flow is then a flow associated with the $L^2$ gradient vector field of a functional (it is important to be clear that the gradient flow depends on the functional). For an example from differential geometry, one can take a surjective continuous map $\pi$ from a topological space $E$ to a compact orientated manifold $M$. This type of map is known as a bundle map and we equip the manifold with an associated bundle metric $\langle \cdot, \cdot \rangle_E$. In this case, the $L^2$ inner product has to take $M$ to $\mathrm{Sym}^2(T^*M)$, so it takes the form

$$\langle s, t \rangle_{L_2} = \int_M \langle s(x), t(x) \rangle_E \, \mathrm{d}v(x),$$

where $s$ and $t$ live in the space of sections of $E$. One can then define gradient flows in an analogous way.

As we hinted at above, one can always complete a space, so the way out is to take the completion of $\left(C_0^\infty(B), (\cdot,\cdot)_{H_0^1}\right)$ and view it as a subspace of the well-understood space $L^2$. This completion is denoted by $H_0^1$ and since the space is complete, one can apply the Riesz theorem to find a weak solution $u$ in $H_0^1(B)$ using the fact that there always exists a $u$ such that

$$(u, \phi) = L(\phi),$$

where the inner product can be defined in terms of weak derivatives and assuming that we can prove that the weak solution $u$ is also smooth (which is possible).

Also note that the exact definition for a weak solution might vary depending on the type of equation and how one formulates the problem. For example, in rigorous fluid dynamics, a vector field $v$ in $L^2(T^2)$ is defined to be a weak solution of the incompressible Euler equations in $T^3 \times [0, T]$

$$\partial_t v + \operatorname{div}(v \otimes v) + \nabla\phi = 0,$$

$$\operatorname{div} v = 0,$$

if for any $\phi$ and any $\psi$ in $C_c^\infty([0,T] \times T^2)$ such that $\operatorname{div} \psi = 0$, we have

$$\int v \, \partial_t\psi + (v \otimes v)\nabla\psi = 0,$$

$$\int v \cdot \nabla\phi = 0,$$

where $T^2$ is the 2-torus and the subscript $c$ for $C^\infty$ refers to the fact that the functions have compact support[5]. $C_c^\infty(U)$ is a particularly nice function space to work with because all the derivatives of the functions living in this space have to vanish within a neighbourhood of the boundary of the set $U$ and they must also vanish in some sense as one gets closer and closer to infinity. These weak solutions we have defined are quite different from the corresponding classical solutions which hold pointwise $C^1(T^3 \times [0, T])$.

[4] Erwin Kreyszig, *Introductory Functional Analysis with Applications* (USA: John Wiley & Sons, 1978).

Although not as famous as their cousins, the Navier-Stokes equations which govern viscous fluid flow, there are still many deep and interesting mathematical questions regarding the Euler equations. One of these (proved recently) is Onsager's conjecture. This a technical conjecture which states that if $v$ is a $C^\beta$ weak solution of the Euler equations such that $\beta \in (0,1)$, then when $\beta > 1/3$ we must have that $E(t) = c$ for any $t \in [0, T]$, where $E$ denotes the Euler equations and $c$ is a constant. If, on the other hand, $\beta < 1/3$, then there must exist $v$ such that $E$ is non-increasing and not a constant. So, in some sense, the numerical value $1/3$ in the exponent offers a kind of turning point. In more physical terms, Hölder continuous solutions to the incompressible Euler equations with exponent greater than $1/3$ must conserve the total kinetic energy, whereas Hölder continuous solutions with exponent less than $1/3$ must dissipate the total kinetic energy.

If you are not familiar, in analysis there are various notions of continuity, some of them stronger than others. If $U$ is a subset of $\mathbb{R}^n$ and $\beta \in (0,1)$, then a function $f$ mapping $U$ to $\mathbb{R}$ is uniformly Hölder continuous with exponent $\beta$ if there exists a constant $c$ such that $|f(x) - f(y)| \leq c|x - y|^{\beta}$ for any $x$ and $y$ in $U$. One can also define a Hölder continuous function with exponent $\beta$ as a function which returns a finite value when we act on it with the semi-norm defined for the Hölder space $C^{0,\beta}$.

$$\|f\|_{C^{0,\beta}(U)} \coloneqq \|f\|_\infty + \sup \frac{|f(x) - f(y)|}{|x - y|^\beta}$$

[5] Alexander Shnirelman, 'On the nonuniqueness of weak solution of the Euler equation', Communications on Pure and Applied Mathematics, Volume 50, Issue 12, 1261-1286 (1997).

In the special case where $\beta = 1$, the function is Lipschitz continuous[6]. A function $f$ is locally Hölder continuous if it is uniformly Hölder continuous on every compact subset of $U$.

The conjecture is based on Kolmogorov's theory of turbulence, so a proof of the conjecture gives strong indirect evidence in favour of that theory. This is probably the most famous theory of isotropic turbulence and the details should be in a good book on turbulence. In words, a turbulent flow is formed of eddies: the energy flux cascades down from the larger eddies to the smaller ones, where it is dissipated by viscosity[7]. As with many difficult conjectures, the proof of the Onsager conjecture is spread over several decades and many papers and intermediate lemmas. Some progress was made using mollification and commutator estimates, and it was found that convex integration could be used to produce an infinite number of bounded weak dissipative solutions with non-constant energy. Convex integration was a technique originally used by John Nash for his famous embedding theorem, which states that every strictly short embedding of a Riemannian manifold into some Euclidean space can be uniformly approximated by $C^1$ isometric embeddings. An eventual proof was finalised recently using convex integration once more and building on the work of many other mathematicians who had contributed to the problem[8]. The overall strategy is to start with what is roughly speaking a 'subsolution' of the problem ie. a solution of a relaxed version of the problem which is nicer and easier to handle. One then iteratively adds perturbations to the subsolution to obtain at each step another subsolution which is closer to being an actual solution. Careful estimates are required which imply convergence to a solution in the required topology. The convex integration parts of the proof rely on 'Mikado flows' (so-called because the components of the flow look like straws and so resemble the Japanese 'pocky' biscuits which are sold under the name Mikado in Europe). The proof also uses a technique of 'gluing approximation' at certain points: for example, to approximate Reynolds stresses with disjoint time steps when mollifying with the integration kernel.

Within the area of PDE theory, the Cauchy-Kovalevskaya theorem relates to a particular type of problem called a Cauchy problem. The Cauchy problem starts with initial data along a smooth curve in $\mathbb{R}^2$ and looks for a solution of a PDE which assumes the initial data we have specified on that curve. If the PDE is of order $n$, then the Cauchy data are the values which the dependent variable can take along with all of its possible PDEs up to and including $n - 1$ on the curve. If the PDE, the smooth curve and the Cauchy data can be written in terms of functions which can all be represented locally in terms of convergent power series, then the Cauchy-Kovalevskaya theorem guarantees the existence of a unique solution to the Cauchy problem close to a point on the smooth curve, and that solution can also be represented locally in terms of convergent power series (ie. it is analytic). If you have studied complex analysis, you might recall that the word 'analytic' is generally replaced with holomorphic when working over the complex numbers. This theorem is a rigorous version of our intuitive feeling that there should be a unique solution to a PDE which achieves the values specified by the initial data along the smooth curve and it was proved in the general form by Kovalevskaya. If you not sure what to imagine by a convergent power series, you could think of a Taylor series, but generalise it and remember that a Taylor series need not converge. Although we have stated things for a curve, in the general setting for PDEs of more than 2 variables, this would be replaced with a hypersurface of the appropriate dimension[9].

[6] Filip Rindler, *Calculus of Variations* (Cham, Switzerland: Springer, 2018).
[7] Marcel Lesieur, *Turbulence in Fluids* (Dordrecht, The Netherlands: Kluwer Academic Publishers, 1997).
[8] Philip Isett, 'A proof of Onsager's conjecture', Annals of Mathematics, Vol. 188, Issue 3, 871 – 963 (2018).

The method of proving the Cauchy-Kovalevskaya theorem demonstrates that Kovalevskaya's work blends into analysis, since it relies on the method of convergent power series. If you have studied ODEs, you will be familiar with the idea of obtaining solutions to differential equations in terms of power series. In a similar way, one can find the functions involved in the general solution of a PDE by expanding the solution as a power series and making substitutions to find the coefficients which are involved. As mentioned before, this would obviously restrict us to analytic PDEs and analytic solutions which happen to have representations as convergent power series. The strategy which one takes for such an 'analytic' problem with analytic initial conditions should also be familiar if you have done some study of basic analysis. One calculates all the possible partial derivatives of the solution to the PDE at the origin and uses these to find the Taylor series expansion. Whether or not this series is a reasonable solution depends on whether or not it converges.

At this stage, it is useful to take the PDE (which might be second or higher order) and convert it to a canonical system of first order. I will skip the details here, but essentially one writes the PDE as a function of particular variables, and then views these variables in the argument as new unknown functions of two new independent variables. After manipulations with integrals and derivatives, one then ends up with a canonical system of first-order quasi-linear PDEs for new unknowns $u_k$.

$$\frac{\partial u_j}{\partial \xi} = \sum_{k=1}^{m} a_{jk}(u_1, \dots, u_m)\frac{\partial u_k}{\partial \eta}, \qquad j = 1, \dots, m$$

A linear differential equation is one which is linear (ie. of degree 1) in the dependent variable and its derivatives, whereas 'quasi-linear' merely means that the partial derivatives appear linearly. It is described as being 'canonical' because we have found $\partial u_k/\partial \xi$ in terms of $u_k$ and $\partial u_k/\partial \eta$. The coefficients in this system have convergent power series representations in a neighbourhood of the origin[10].

In fact, one can prove the Cauchy-Kovalevskaya theorem without switching to a canonical system, but if we assume that the canonical system is sufficient for our purposes, we must now find Taylor series expansions about the origin for the unknown functions $u_k$. The coefficients in such a series are found by mechanically computing the partial derivatives of $u_k$ at the origin. To do this, one has to find recursion formulae for partial derivatives with respect to the variable $\xi$ (again, you might already be familiar with recursion formulae for finding coefficients in power series solutions to ODEs). The PDEs themselves are already recursion formulas of this type and one can find the derivatives of the $u_k$ with respect to $\eta$ by differentiating the initial conditions. After some more manipulations, we find all the partial derivatives of the functions $u_k$ and realise that they are always positive. This eventually allows us to prove the convergence of the series expansion

$$u_k(\xi,\eta) = \sum_{\mu,\nu=0}^{\infty} \frac{1}{\mu!\ \nu!} \frac{\partial^{\mu+\nu} u_k(0,0)}{\partial\xi^{\mu}\partial\eta^{\nu}} \xi^{\mu}\eta^{\nu}$$

for the unknowns $u_k$. One does this essentially by creating a majorant for the above series by solving the canonical initial value problem which we have already defined. The majorant is a convergent series of positive terms, each of which is at least as large in absolute value as the corresponding term in the Taylor series solution above. Since this solution solves the initial value problem by design, then the expansion above must always converge for sufficiently small $\xi$ and $\eta$.

[9] Peter Olver, *Introduction to Partial Differential Equations* (Cham, Switzerland: Springer, 2016).
[10] Paul Garabedian, *Partial Differential Equations* (New York: Chelsea Publishing Company, 1986).

This proves that the canonical problem has a unique analytic solution within a neighbourhood of the origin and leads us to the Cauchy-Kovalevskaya theorem, which can be stated in the following way: in the neighbourhood of a point at which the coefficients $a_{jk}$ and the functions $h_j$ have power series representations, one can find a unique vector with components $u_k$ which solves the following initial value problem[11]:

$$u_{\xi} = A(u)u_{\eta}, \quad u(0,\eta) = h(\eta).$$

The functions $h_j$ are analytic data defined as

$$h_j^{(n)}(0) = \frac{n!\,M}{\rho^n},$$

where $M$ is the upper bound on the absolute value of each term in the Taylor series and $\rho$ is a sufficiently small positive number. This theorem is more limited than it appears (even if we assume that we are working with problems where all the functions are analytic), but it does show that if are willing to focus on analytic functions, the number of functions needed for a general solution of the PDE is always the same as the order of the PDE and that the PDE has one more variable in comparison to the number of variables in the arguments of the arbitrary functions.

The Cauchy problem for the gravitational field plays a key role in general relativity. One starts with a spacelike hypersurface $M^3$ equipped with a metric $g_{ij}$ and a symmetric 2-tensor $k_{ij}$ known as the extrinsic curvature. The word 'extrinsic' refers to the fact that we are studying how an immersion of a space into a smaller space bends and distorts the larger space inside the smaller one. Intrinsic geometry studies operations on a manifold which do not depend on how we immerse it in another manifold[12]. As an example, Gauss's *Theorema Egregium* states that the Gaussian curvature of a 2-dimensional submanifold $(M, g)$ of $\mathbb{R}^3$ (where $g$ is the metric 'induced' or inherited from the ambient space into which the manifold is embedded) is a local isometry invariant of $(M, g)$. In different words, a surface in $\mathbb{R}^3$ has a type of curvature called Gaussian curvature which is defined by the fact that the surface is embedded into $\mathbb{R}^3$, but that definition coincides with the definition given by the Riemann tensor, which only depends on the metric with which the surface is equipped, and the metric does not care about the ambient space into which the surface is embedded and only measures intrinsic bending. Therefore, the Gaussian curvature does not actually depend on the embedding of the surface into $\mathbb{R}^3$, and is an intrinsic invariant. This is a rather beautiful and remarkable result of classical intrinsic geometry[13].

[11] *Ibid.*

[12] Peter Petersen, *Riemannian Geometry* (New York: Springer, 2006).

In GR, the 3-manifold along with the metric and the 2-tensor is known as an 'initial data' for the Einstein equations (similar to the idea of initial data which you might already be familiar with from elementary PDE theory, where you would need to know the initial temperature distribution of the body before being able to deduce its final temperature from the heat equation). For the Cauchy problem, one starts with a 3-manifold $N$ with some more complicated initial data, and must then find a 4-manifold $M$, plus an embedding $E$ which maps $N$ to $M$ and a metric on $M$ which satisfies the Einstein equations. The metric must also agree with the initial values of $E(N)$ and $E(N)$ must be a Cauchy surface for $M$, where a Cauchy surface is a spacelike hypersurface such that every curve in the surface which is not spacelike intersects only once. The final manifold is said to be a development of the initial manifold with its initial data and there always exist developments of the initial 3-manifold as long as the initial data satisfy well-defined constraint equations (we will talk about this in more depth later when we discuss the work of Choquet-Bruhat)[14].

Another area where Kovalevskaya made a key contribution is that of classical mechanics. I have discussed some of the basics of classical mechanics in another work, so I will not go over these again here[15]. In essence, Kovalevskaya described a case of a rigid spinning top, whose motion due to free precession under gravity is simple enough to be solvable analytically. This is not possible generally due to the complexity of precession. To get some insight about the Kovalevskaya top, it might be helpful to consider the case of a symmetric top with one point fixed in space. The top is then acted on by a uniform gravitational field ie. the centre of mass of the top is acted upon by a constant vertical force. The fixed point is a point on the symmetry axis of the object (recall that we are assuming the spinning top to be a symmetric rigid body). This axis can be chosen to be the $z$-axis of the coordinate frame which is embedded in the body and which spins with it. As there is a base point which is fixed, we can specify the motion of the top with the Euler angles (again, I have discussed these elsewhere).

The rates of change of the Euler angles give us three motions which have nice physical interpretations. $\dot{\psi}$ is the ordinary rotation of the top or gyroscope about the $z$-axis inside the top, $\dot{\phi}$ is the rotation of this axis about the 'regular' $z$-axis which is existing outside the top (known as 'precession') and $\dot{\theta}$ is the bobbing of the embedded axis relative to the regular vertical axis (think of a float bobbing up and down, but it is the axis itself which is doing the motion in this case). This bobbing motion is known as nutation. In many cases, the size of the rate of the change of the precession is much larger than the rate of change for the bobbing motion, and the rate of change for the ordinary rotation is much larger than either of them. In this case, and given that we are studying a symmetric top where the two of the moments of inertia are the same by definition, the Euler equations simplify to[16]

$$I_1\dot{\omega}_1 + \omega_2\omega_3(I_3 - I_2) = N_1,$$

$$I_2\dot{\omega}_2 + \omega_1\omega_3(I_1 - I_3) = N_2,$$

$$I_3\dot{\omega}_3 = N_3.$$

[13] Jürgen Jost, *Riemannian Geometry and Geometric Analysis* (Heidelberg: Springer, 2008).
[14] Stephen Hawking and George Ellis, *The Large Scale Structure of the Universe* (Cambridge: Cambridge University Press, 2006).
[15] Isaac Williams, *Worlds of Motion: Why and How Things Move* (London: Austin Macauley Publishers, 2018).

One can see the types of motion I have described in the photos of a gyroscope below.

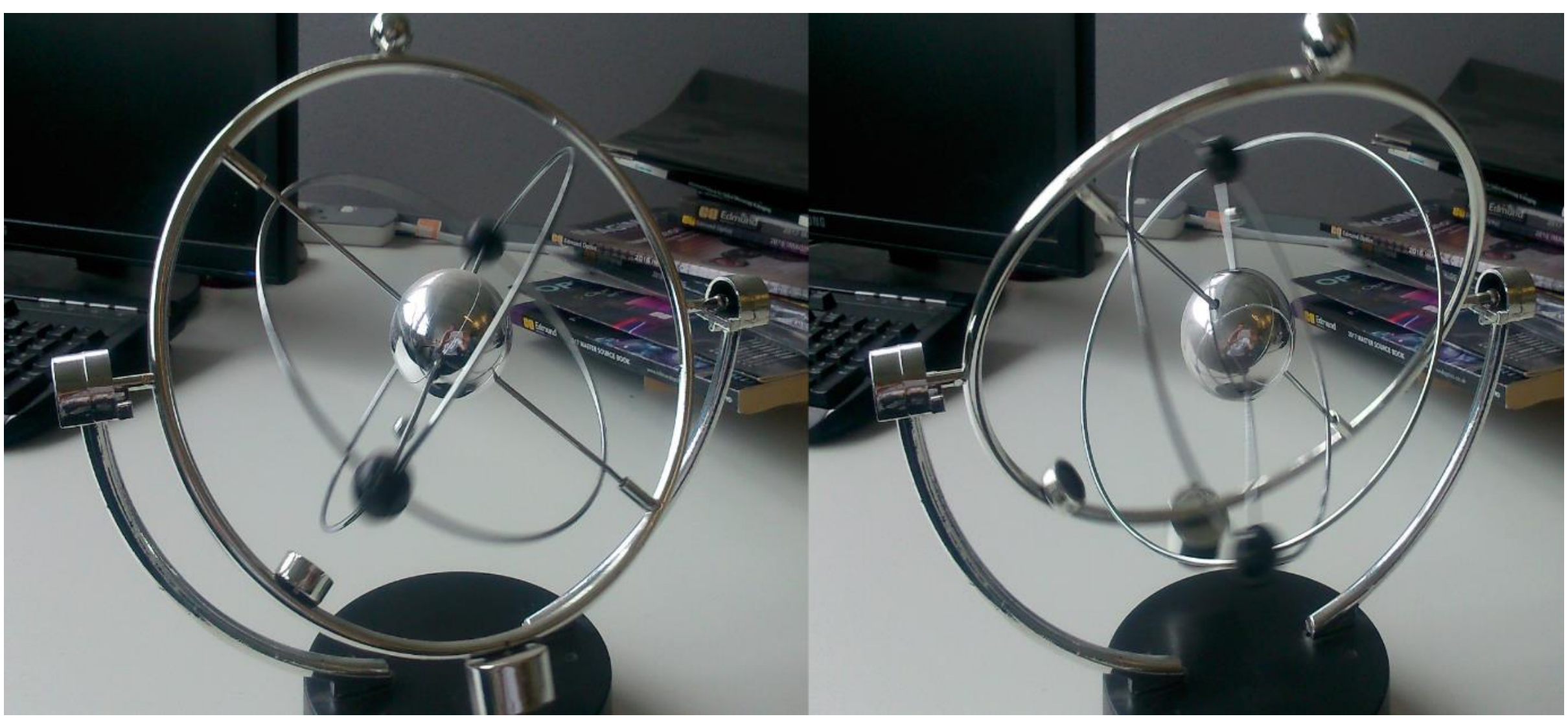

Figure 1: Motion of a gyroscope.

The Euler equations look appealing, but they are not particularly useful to solve in this particular situation, so it is instead easier to use the Lagrangian formulation. You might see elsewhere the Hamiltonian formulation for tops. The kinetic energy of the top is

$$T = \tfrac{1}{2}I_1(\omega_1^2 + \omega_2^2) + \tfrac{1}{2}I_3\omega_3^2.$$

The Lagrangian of a physical system is

$$L = T - V,$$

where $V$ is the potential energy. In this case, we are working with a uniform, constant gravitational field acting on the centre of mass of the top, so the potential energy is fairly easy to write down in terms of the Euler angles. The kinetic energy can also be rewritten in terms of these angles, so we end up with

$$L = \frac{I_1}{2}(\dot{\theta}^2 + \dot{\phi}^2\sin^2\theta) + \frac{I_3}{2}\left(\dot{\psi} + \dot{\phi}\cos\theta\right)^2 - Mgl\cos\theta.$$

[16] Herbert Goldstein, Charles Poole and John Safko, *Classical Mechanics* (USA: Addison Wesley, 2002).

For this system, the components of the angular momentum along the spatial vertical axis and along the $z$- axis embedded in the top are constant, so we have two constants of the motion of the top.

$$p_\psi = \frac{\partial L}{\partial \dot{\psi}} = I_3(\dot{\psi} + \dot{\phi}\cos\theta),$$

$$p_\phi = \frac{\partial L}{\partial \dot{\phi}} = (I_1 \sin^2\theta + I_3\cos^2\theta)\dot{\phi} + I_3\dot{\psi}\cos\theta.$$

The system does not lose energy over time, so the kinetic energy added to the potential energy must give us a third constant of the motion:

$$E = \frac{I_1}{2}(\dot{\theta}^2 + \dot{\phi}^2\sin^2\theta) + \frac{I_3}{2}\left(\dot{\psi} + \dot{\phi}\cos\theta\right)^2 + Mgl\cos\theta.$$

Some manipulations with these first integrals of the motion give us differential equations for the rates of change of the three Euler angles, which effectively tells us everything we need to know about the motion after some more work. For example, we can obtain the frequency of the nutation and the angular velocity of the precession. The average precession frequency of the top is found to be

$$\bar{\dot{\phi}} = \frac{Mgl}{I_3}\omega_3.$$

Hence the rate of precession decreases as one increase the initial rotational velocity of the top[17].

Another related example is the Kovalevskaya top[18]. This top is also symmetric, but one now assumes that there is a ratio amongst the principal moments of inertia such that two of the moments are equal and double the size of the remaining moment.

$$I_1 = I_2 = 2I_3$$

In our previous example, the centre of gravity was located on the symmetry axis, but for the Kovalevskaya top it is located on a plane which is perpendicular to the symmetry axis. One can write down the Hamiltonian for this top and find that as before the first integrals of motion include the energy (energy will always be a constant of the motion as long as the system is conservative) and the component of the angular momentum in the direction of the embedded $z$-axis, but that there is also a new conserved quantity called the Kovalevskaya invariant. It is quite surprising that a mathematician of the nineteenth century was able to find another classical top whose description is sufficiently simple as to allow for an analytic solution for the motion of the precession. One might have thought naively that Euler and Lagrange had already dealt with all the cases which are simple enough that finding their precession is an integrable problem.  In fact, another example of a solvable top was found in the twentieth century, but it is very difficult to integrate and there was a long period of time between its introduction and its eventual integration.

---

[17] *Ibid.*

[18] Sofia Kovalevskaya, 'Sur le problème de la rotation d'un corps solide autour d'un point fixe', Acta Mathematica, 12: 177 – 232 (1899).

# *Emmy Noether*

Noether was one of the most important and influential mathematicians of the nineteenth century, especially known for her contributions to algebra and mathematical physics. Despite her brilliance, she faced many obstacles in her academic life due to the fact that she was a woman, and had to work without pay for several years. She did much to shape the face of abstract algebra when it was relatively under-developed, including the basic theories of rings and fields which are now so well-known to undergraduate mathematics students. If you consider the state which abstract algebra was in at the time, the difficulty of this task cannot be overestimated. Workers in other fields such as fluid dynamics, mathematical physics, or differential geometry, usually have hints to work with, or analogies and simple pictures to fall back on. Even with something like the fundamental group in algebraic topology, there is a simple geometric idea at the bottom of the abstract concept: that of counting the number of holes encircled by a closed loop. These hints tend to vanish when one studies pure abstract algebra, especially for the beginner.

Noether laid the foundations of much of the theory of commutative rings, known as commutative algebra. This is still a basic and important branch of algebra, and one has to know a fair amount of commutative algebra before even beginning to consider studying problems in algebraic geometry. A ring $R$ in general is a set equipped with addition and multiplication such that $R$ forms an abelian group with respect to addition ($R$ contains a zero element and every element of the ring has an additive inverse which can be added to it to obtain the zero element). Furthermore, multiplication is associative

$$(xy)z = x(yz),$$

and distributive over addition

$$x(y+z) = xy + xz, \qquad (y+z)x = yx + zx,$$

for any $x, y, z \in R$. One can also require that the ring has an identity element denoted by 1 and that $xy = yx$. The latter requirement means that the ring is commutative, and the former is also often assumed for simplicity. A ring homomorphism is a map $f$ taking a ring $A$ into another ring $B$ such that

$$f(x+y) = f(x) + f(y),$$

$$f(xy) = f(x)f(y),$$

$$f(1) = 1.$$

The first criterion means $f$ is a group homomorphism of abelian groups. As you might expect, you can form a subset of a ring which is then a subring if the product of sum of two elements is also an element of the set and if the set contains the identity element for the original ring $R$. A composition of ring homomorphisms is also a ring homomorphism.

Another basic definition which one needs is that of an ideal: these are simple objects, but they were employed to great effect by Noether. An ideal $\mathfrak{r}$ of a ring $R$ is a subset of the ring such that $\mathfrak{r}$ is an additive subgroup (ie. a linear sum of elements of $\mathfrak{r}$ with coefficients taken from $R$ is still in $\mathfrak{r}$). One might write $(x, y)$ for the ideal generated in this way by two elements of $R$. An ideal is also defined such that

$$R\mathfrak{r} \subseteq \mathfrak{r}.$$

If $\mathfrak{r} \subset R$ is an ideal, then the quotient group $R/\mathfrak{r}$ inherits multiplication from $R$ such it becomes a ring in itself, known as the quotient ring. There exists a surjective ring homomorphism $\phi$ from $R$ to $R/\mathfrak{r}$ which takes an element of $R$ to its coset and the quotient ring along with the quotient homomorphism are defined uniquely modulo isomorphisms. It is an important and useful fact that the map $\phi^{-1}$ which takes ideals of $R/\mathfrak{r}$ to ideals of $R$ which contain $\mathfrak{r}$ is a one-to-one order-preserving correspondence. As an example of a fact which uses rings and ideals, take a ring $R$ which is not the zero ring. This is the ring whose only element is $0$, which can only occur if we allow the identity element to equal $0$, since this means that

$$x = x1 = x0 = 0.$$

Excluding this possibility, the following three statements are equivalent: $R$ is a field (that is, every non-zero element of the ring is a unit), the only ideals which the ring contains are the ideals generated by $0$ and $1$, and every homomorphism of $R$ into a non-zero ring $A$ is an injection. Bear in mind that whereas the ideal generated by $1$ is denoted by $(1)$, the ideal generated by $0$ is usually denoted just by $0$ (as is the zero ring). Perhaps the hardest part to prove in this chain is that the third statement implies the first statement. By definition, every ideal $\mathfrak{r}$ is the kernel of the quotient homomorphism from $R$ to $R/\mathfrak{r}$. We have assumed that the homomorphism is injective, but that implies that the only proper ideal is $(0)$, where a proper ideal is an ideal which is smaller than the ring. For any proper ideal $(x)$, $x = 0$. This implies that every non-zero $x$ is invertible, which means that the ring is a field. Alternatively, we could prove much the same thing by a contradiction argument. Assume that the third statement is true, but that $R$ is not a field. This means that there has to be an element of $R$ which is not a unit. If it is not a unit, the ideal $(x)$ does not contain $1$, so $(x)$ is a proper ideal of $R$. If we take the quotient homomorphism $\phi$, this homomorphism is not an injection by definition because $\phi(x) = 0$ and $x$ is not $0$. $(x)$ is not a proper ideal, so the quotient ring $R/(x)$ is not the zero ring (since the only way to get the zero ring from a quotient is via the quotient $R/R$ and we just said that $(x)$ is smaller than $R$). This is a contradiction, which proves the claim[19].

We are not quite finished with the preliminary definitions, as we still need to define modules. The definition of modules might seem a bit strange at first, but they are essential to modern commutative algebra and they help to unify many concepts which you may have seen before in linear algebra without really knowing what was going on or why the computations were producing the correct answers. Start with a commutative ring $R$. An $R$-module is an abelian group $M$ on which $R$ acts via linear transformations. More formally, it is an abelian group $M$ along with a map from $M \times R$ to $R$ such that for elements $x, y$ in $M$ and $a, b$ in $R$, the following axioms hold:

$$(x + y)a = xa + ya,$$

$$x(a + b) = xa + xb,$$

$$x(ab) = (xa)b,$$

$$x1 = x.$$

[19] Michael Atiyah and Ian MacDonald, *Introduction to Commutative Algebra* (USA: Westview Press, 1969).

You might see these axioms written such that the operators occur in the opposite order: this results in a left $R$-module, whereas we have stated the axioms for a right $R$-module. The classic example of a module which you might have seen is a vector space over a field: this vector space is an $R$-module, where $R$ is a field. As another example, if $R = \mathbb{Z}$, then the $R$-module is an abelian group.

In order to obtain more meaningful results about modules and rings, one has to make use of finiteness conditions for partially ordered sets called chain conditions[20]. If $S$ is a set equipped with a partial order for the relation $\leq$ such that $x \leq y$ and $y \leq x$ implies that $y = x$, then the following conditions are equivalent: firstly, for an increasing sequence in the set,

$$x_1 \leq x_2 \leq \cdots$$

there exists $n_0$ such that $x_m = x_n$ for all $n, m \geq n_o$. Secondly, every non-empty subset of $S$ contains a maximal element. The first condition is known as the ascending chain condition, and the second is the maximum condition. The particular application we have in mind is when $S$ is the set of submodules of a module and the partial ordering relation is given by inclusion of closed subsets. A module which satisfies either of the above conditions is Noetherian. If we order the set via the opposite inclusion, the conditions become the descending chain condition and the minimum condition and a module which satisfies either is Artinian. If we take a set in which the maximum condition holds, the principle of Noetherian induction says that if $X$ is a subset of $S$ which contains an element $s$ of $S$ whenever it contains all elements $x \in S$ such that $x > s$, then $X = S$. Technically, this is based on the fact that the set of all submodules forms a modular lattice under partial ordering by inclusion, and it is this lattice which then satisfies one of the conditions[21]. Noetherian modules are more useful than Artinian ones, since the Noetherian condition is the one which is required for many theorems and propositions. For example, $M$ being a Noetherian $R$-module implies that every submodule of $M$ is finitely generated.

In many cases, however, the basic properties of Noetherian and Artinian modules are the same. If we take a short exact sequence of $R$-modules,

$$0 \to M' \to M \to M'' \to 0$$

[20] Emmy Noether, 'Idealtheorie in Ringbereichen', *Mathematische Annalen*, 83 (1), 24 – 66 (1921).

[21] Paul Cohn, *Basic Algebra: Groups, Rings and Fields* (Trowbridge: Springer, 2003).

$M$ being Noetherian or Artinian implies that $M'$ and $M''$ are also of this type. If a set of modules $M_i$ are all Noetherian or Artinian, then so is the direct sum. This can be proved using a similar short exact sequence.

$$0 \to M_n \to \prod_{i=1}^{n} M_i \to \prod_{i=1}^{n-1} M_i \to 0.$$

We can take a ring $R$ and regard it as a right $R$-module in itself, and by the same terminology the ring is said to be a Noetherian or an Artinian ring, where the maximum and minimum conditions are defined with respect to right ideals. If $R$ is a right Noetherian ring, then a finitely generated $R$-module is also Noetherian (or Artinian, respectively). Most rings which one uses when applying commutative algebra are Noetherian: in particular, the rings which are used in algebraic geometry when one wishes to study the links between rings and geometric objects are Noetherian. It might be worth mentioning that the ascending chain condition also has applications in topology: a topological space $X$ is Noetherian if the open subsets of the space satisfy the ascending chain or maximum condition. Since closed subsets are complements of open subsets (and vice versa), it is an equivalent requirement that the closed subsets of the space satisfy the descending chain or minimum condition. One can show that every subspace of a Noetherian topological space is also Noetherian, and that the space is compact.

There are three equivalent conditions for a ring $A$ to be Noetherian: we already established their equivalence in the setting of modules. Firstly, every non-empty set of ideals in the ring contains a maximal element. Secondly, every ascending chain of ideals in the ring is stationary: that is, the chain of ideals

$$\mathfrak{a}_1 \subset \mathfrak{a}_2 \subset \cdots \subset \mathfrak{a}_k \subset \cdots$$

stops at some point, so we end up with $\mathfrak{a}_k = \mathfrak{a}_{k+1}$ for some $k$. Thirdly, every ideal $\mathfrak{a} \in A$ is finitely generated. If $A$ is Noetherian and $\phi$ is a ring homomorphism from $A$ onto another ring $B$, then $B$ is also Noetherian. If $A$ is a subring of $B$ such that $A$ is Noetherian and $B$ is finitely generated as an $A$-module, then $B$ is a Noetherian ring. A trivial example of a ring which is not Noetherian is the polynomial ring $k[X_1, \dots, X_n, \dots]$ formed from the set of polynomials in an infinite number of indeterminates, where the coefficients are drawn from the field $k$ (more generally, another ring).

Many examples of rings which *are* Noetherian can be found using the Hilbert basis theorem. This states that the polynomial ring $A[x]$ is Noetherian if $A$ is Noetherian. One can prove this by defining auxiliary sets which form an ascending chain of ideals and showing that any ideal $\mathfrak{a} \subset A[x]$ is finitely generated. (The word 'basis' was used historically to refer to a generating set). As a direct corollary by induction on $n$, if $A$ is Noetherian, so is $A[x_1, \dots, x_n]$. Also, if $A$ is a Noetherian ring and $B$ is a finitely generated $B$-algebra, then $B$ is Noetherian. Every finitely generated ring, and every finitely generated algebra over $\mathbb{Z}$ or a field $k$, is Noetherian. One can follow similar lines to the proof of the Hilbert basis theorem to show that $A$ being a Noetherian ring implies that the formal power series ring $A[\![x]\!]$ is Noetherian[22]. We can also state that if $k$ is a field and $E$ is a finitely generated $k$-algebra such that $E$ is also a field, then $E$ is a finite algebraic extension of $k$. A field extension is a monomorphism $i: k \to l$, where $k$ and $l$ are both subfields of a larger field (often taken to be $\mathbb{C}$). We might say that $k$ is the small field and $l$ is the large field. There are many interesting results for field extensions and one studies this type of extension in some detail in Galois theory: specifically, the group of field automorphisms of the extension. The famous Tower Law, for example, is a result about degrees of extensions where there is a chain of subfields (the degree of a field extension being the dimension of the large field considered as a vector space over the small field)[23].

A finite extension is one with finite degree and an extension (denoted as $l : k$) is algebraic if every element of $l$ is algebraic over $k$. Every finite extension is algebraic, but the converse is not always true. An element is algebraic over $k$ if it satisfies a dependence relation for an $n$-th degree polynomial in the polynomial ring over the field:

$$f(y) = a_n y^n + \cdots + a_1 y + a_0 = 0.$$

One can divide through by $a_n$ when working in a field, so we can assume that $a_n = 1$: that is, the polynomial is monic. However, a polynomial ring does not have to be over a field. In the case when it is over a ring, the relation

$$f(y) = a_n y^n + \cdots + a_1 y + a_0 = 0,$$

where $a_n = 1$ is known as an integral dependence relation for $y$. There is a useful relation between integral dependence and finiteness conditions in a ring extension: as an example, if one knows that two elements of a ring extension $B \subset A$ are integral over $A$, one can use finiteness conditions to show that the sum and the product of the two elements are integral over $A$ without needing to compute their integral dependence relations. One can start with an $A$-algebra $\varphi: A \to B$ and an element $y \in B$ and show that the following statements are all equivalent: $y$ is integral over $A$, the polynomial subring $A'[y] \subset B$ generated by $A' = \varphi(A)$ and $y$ is finite over $A$, and there exists an $A$-subalgebra $C \subset B$ such that $A'[y] \subset C$ and $C$ is finite over $A$. One can prove the final equivalence using the determinant trick[24].

These finiteness conditions lead us to an important lemma, the Noether normalization lemma (really a theorem). For this lemma one takes a field $k$ and a finitely generated $k$-algebra $A$. In this case, there exist elements $x_1, \dots, x_n$ in $A$ such that $x_1, \dots, x_n$ are algebraically independent over $k$ and $A$ is independent over $B = k[x_1, \dots, x_n]$. A $k$-algebra is said to be finitely generated over $k$ if $A = k[x_1, \dots, x_n]$ for some finite set $x_1, \dots, x_n$. In more detail for the case where $A$ is an integral domain, the result says that for an infinite field $k$ and $K = k(x_1, \dots, x_n)$ the field of fractions of $k$, there exists a nonnegative integer $d$ and $d$ linear combinations $y_1, \dots, y_d$ of $x_1, \dots, x_n$ with coefficients of $k$ such that $y_1, \dots, y_d$ are algebraically independent over $k$ and such that every element of $R$ is integral over $k[y_1, \dots, y_d]$. If $K$ is separably generated over $k$, then the $y_i$ can be chosen such that $K$ is a separable extension of $k(y_1, \dots, y_d)$[25]. The field of fractions of the ring $R$ is a field $K$ containing a subring isomorphic to $R$ such that every element of $K$ can be expressed as $r/s$, where $s$ is non-zero and both $r$ and $s$ belong to the subring. Noether introduced this lemma in the different context of a paper on the theory of invariants of finite groups[26].

[22] Miles Reid, *Undergraduate Commutative Algebra* (Cambridge: Cambridge University Press, 2010).
[23] Ian Stewart, *Galois Theory* (USA: Chapman & Hall, 2004).
[24] Miles Reid, *Undergraduate Commutative Algebra* (Cambridge: Cambridge University Press, 2010).

One important use of the lemma is to establish a result known as Hilbert's Nullstellensatz (or the weak form, at least). This states if $k$ is a field and $K$ is a finitely generated $k$-algebra, then $K$ is algebraic over $k$ such that $k$ is a finite field extension. Another way of stating this is that if one takes an algebraically closed field $k$, each maximal ideal in $R = k[x_1, \dots, x_n]$ has the form

$$\mathfrak{m} = (x_1 - \alpha_1, \dots, x_n - \alpha_n).$$

In any case, the proof starts off by using the Noether normalisation lemma, since one uses the fact that one immediately has a set of elements in $K$ which are algebraically independent, plus the fact that $K$ is finite over the polynomial ring in the indeterminates formed by those elements where the coefficients are drawn from $k$. The full form of the Nullstellensatz is somewhat drawn out, but it can be stated in a few lines using the concept of the radical of an ideal. If $I$ is an ideal contained in commutative ring with an identity, the radical of the ideal $\sqrt{I}$ is the set of elements in the ring such that the $k$-th power is in the ideal for some $k \geq 1$. Using this concept, start with $R = k[x_1, \dots, x_n]$, a polynomial ring over an algebraically closed field $k$. It follows that for any ideal $I \subset R$,

$$IV(I) = \sqrt{I}.$$

Although it might not be obvious from the statement using radicals, this theorem marks the beginning of modern algebraic geometry (the study of geometric properties defined by algebraic equations), since it states the basic result that a family of polynomials over an algebraically closed field has a common solution when the ideal they generate in the polynomial ring is proper. In general, a proper ideal does not necessarily need to have a non-empty zero locus, so the theorem in part is making the non-trivial statement that a proper ideal does have a non-empty zero locus if we assume that the field is algebraically closed. If we think of a system of polynomial equations which we are looking to solve simultaneously in the appropriate field $F_i(x_1, \dots, x_n)$, the set of solutions is the zero locus for the system. The key is that increasing the size of the system does not affect the zero locus, so we can take everything to be a member of an ideal in the polynomial ring $k[x_1, \dots, x_n]$. The Hilbert basis theorem already tells us that an ideal $I \subset k[x_1, \dots, x_n]$ is finitely generated, so studying the zero locus of an ideal is the same as studying the zero locus of a set of polynomials.

Another important area of commutative algebra which Noether contributed to is the theory of decompositions of ideals in Noetherian rings, where the decompositions are primary since they break ideals up into primary ideals. Recall that integers can be factorised into products of prime numbers and powers of prime numbers. This concept of factorisation generalises to ideals, since one can have prime ideals and ideals which are analogous to products of prime powers. A ideal $\mathfrak{p} \subset A$ is prime if $\mathfrak{p} \neq (1)$ and if $xy$ belonging to the ideal implies that either $x$ or $y$ by themselves belongs to the ideal. $\mathfrak{p}$ is prime if and only if $A/\mathfrak{p}$ is an integral domain. One might be able to guess this equivalence, since a ring $R$ is an integral domain if for two elements $x$ and $y$ of the ring, $xy = 0$ implies that $x = 0$ or $y = 0$. A classic example is the ring $\mathbb{Z}$ (also obviously a commutative ring), since the product of two integers coming out as zero implies that one of those integers has to be zero. Another example is the polynomial ring $\mathbb{Z}[t]$: one can show by a contradiction that one cannot have $fg = 0$ such that $f, g \neq 0$. Another way of stating this is that the domain has no zero-divisors, since one cannot find a non-zero element which can be multiplied with another non-zero element to get zero. As well as prime ideals, one can have primary ideals. An ideal $I \subset A$ is primary if $I \neq A$ and if $xy$ belonging to $A/I$ implies that either $x$ is in $I$ or $y^n$ is in $I$ for some positive $n$.

[25] Anthony Knapp, *Advanced Algebra* (New York: Birkhäuser, 2007).
[26] Emmy Noether, 'Der Endlichkeitsatz der Invarianten endlicher linearer Gruppen der Charakteristik $\rho$', *Nachr. Ges. Wissm*, DE: 28 – 35 (1926).

An alternative way of stating the definition is that $I$ is a primary ideal if and only if $\frac{A}{I} \neq 0$ and every zero-divisor in the quotient is nilpotent (a element of a ring $x$ being nilpotent if there exists some positive $n$ such that $x^n = 0$). For example, the primary ideals in $\mathbb{Z}$ are the ideals generated by $0$ and $p^n$, for prime $p$. One needs th second type of ideal to allow for ideals which are not radical. A primary decomposition of an ideal $I$ in $A$ is simply the expression of that ideal as a finite intersection of primary ideals:

$$I = \bigcap_{i=1}^{n} \mathfrak{q}_i.$$

Not all ideals have primary decompositions, but Noether proved the basic existence theorem that every ideal contained in a Noetherian ring has a primary decomposition. The proof is an application of an argument by Noetherian induction and also makes use of the concept of an indecomposable ideal. An ideal is indecomposable if it cannot be written as the intersection of two larger ideals, so a prime ideal is obviously indecomposable. There are also several options for uniqueness or partial uniqueness theorems. For example, take a decomposable ideal $\mathfrak{a}$ and write out a minimal primary decomposition for $\mathfrak{a}$:

$$\mathfrak{a} = \bigcap_{i=1}^{n} \mathfrak{q}_i.$$

A decomposition is said to be minimal if the radicals $r(\mathfrak{q}_i)$ are all distinct and if the intersection $\bigcap_{j \neq i} \mathfrak{q}_i$ is not properly contained in $\mathfrak{q}_i$. Any primary decomposition can be reduced to a minimal decomposition. Take $\mathfrak{p}_i = r(\mathfrak{q}_i)$, then the $\mathfrak{p}_i$ are the prime ideals in the set of ideals $r(\mathfrak{a}: x)$ and so do not depend on the decomposition which is chosen for $\mathfrak{a}$. Recall that a famous uniqueness result exists for the integers known as the fundamental theorem of arithmetic (first proved by Euclid): every integer is either a prime number or can be decomposed into a product of prime numbers which is unique modulo ordering. The fundamental theorem was used by Gauss to prove another basic result, the quadratic reciprocity law. This states that for any odd prime number $p$ and $q$, the following relation holds

$$\left(\frac{p}{q}\right)\left(\frac{q}{p}\right) = (-1)^{\frac{p-1}{2}\frac{q-1}{2}}$$

where the Legendre symbol $\left(\frac{q}{p}\right)$ is defined as $0$ if $q$ divides $p$, $1$ if the modular equation $x^2 \equiv q \pmod{p}$ has a solution, and $-1$ otherwise. The last two possibilities are equivalent to $q$ being a quadratic residue[27].

For another uniqueness theorem, take a decomposable ideal $\mathfrak{a}$. Let $\bigcap_{i=1}^{n} \mathfrak{q}_i$ be a minimal primary decomposition of $\mathfrak{a}$ and $\{\mathfrak{p}_{i_1}, \dots, \mathfrak{p}_{i_n}\}$ be a set of prime ideals of $\mathfrak{a}$. It follows that the intersection $\mathfrak{q}_{i_1} \cap \dots \cap q_{i_m}$ does not depend on the decomposition. As a corollary, the primary components $\mathfrak{q}_i$ corresponding to minimal prime ideals are uniquely determined by $\mathfrak{a}$. The term 'isolated' originates in algebraic geometry, as does the Spec notation, referring to a type of topological space known as a 'prime spectrum of a ring' where the topology is the Zariski topology. In algebraic geometry, one typically considers objects known as varieties, where a variety is a geometric object cut out by a set of homogenous polynomial equations in one or more variable

$$\{P_1(x_1, \dots, x_d) = \cdots = P_k(x_1, \dots, x_d) = 0\}.$$

If $A = k[x_1, \dots, x_n]$ for a field $k$, the ideal $\mathfrak{a}$ cuts out a variety $X$. The minimal prime components $\mathfrak{p}_i$ correspond to the components of the variety which are irreducible, whereas the embedded primes correspond to sub-varieties which are embedded in the irreducible components. Although we have said that the isolated primary components are determined uniquely by the ideal, this is not generally true of the embedded primary components[28].

As with the example of Noetherian topological spaces, there are other 'Noetherian' objects which follow the type of finiteness condition associated with Noether. For example, in algebraic geometry a scheme $X$ is said to be Noetherian if it has a finite cover by affine set

$$X = \bigcup_i U_i\,, \quad \text{with } U_i = \operatorname{Spec} A_i,$$

such that the $A_i$ are Noetherian rings. This is equivalent to a scheme having a kind of finiteness condition on its dimensionality. If the affine scheme $\operatorname{Spec} A$ is Noetherian, then it can be proved that $A$ is a Noetherian ring using one of those standard arguments with a chain of ideals which you are hopefully getting used to hearing about[29].

A famous theorem of Noether relates to invariants computed in variational problems: it is usually known as Noether's theorem and it is one of the most important results in the history of mathematics and theoretical physics (modern physics would certainly not have developed in the same direction without this theorem)[30]. In modern terminology, the theorem can be stated as follows: take $f$ mapping from $\Omega \times \mathbb{R}^m \times \mathbb{R}^{m\times d}$ to $\mathbb{R}$ such that $f$ is twice differentiable in $v$ and $A$ and satisfies growth bounds

$$|D_v f(x, v, A)|, |D_A f(x, v, A)| \leq C(1 + |v|^P + |A|^p),$$

for some positive $C$ and $p \in [1, \infty)$. Take an associated functional for $f$

$$\mathcal{F}[u] \coloneqq \int_\Omega f\big(x, u(x), \nabla u(x)\big)\, \mathrm{d}x$$

such that $u_*$ is a minimizer for this functional and let this functional be invariant under transformations specified by the following maps

$$g(x, 0) = x,$$

$$H(x, 0) = u_*(x).$$

[27] Henryk Iwaniec and Emmanuel Kowalski, *Analytic Number Theory* (USA: American Mathematical Society, 2004).
[28] Michael Atiyah and Ian MacDonald, *Introduction to Commutative Algebra* (USA: Westview Press, 1969).
[29] Igor Shafarevich, *Basic Algebraic Geometry: Schemes and Complex Manifolds* (Berlin: Springer, 1994).

The equations for $g$ and $H$ resemble those for a homotopy, and one can think of the transformation in this way. Assume that there exists a majorant $h$ in $L^p(\Omega)$ such that

$$|\partial_\tau H(x, \tau)|, |\partial_\tau g(x, \tau)| \leq h(x)$$

for almost every $x$ and real $\tau$. It follows that for any minimizer of the functional $\mathcal{F}[u]$, there is a conservation law

$$\mathrm{div}[\mu^T D_A f(x, u_*, \nabla u_*) - \nu f(x, u_*, \nabla u_*)] = 0$$

almost everywhere in $\Omega$, where

$$\mu(x) \coloneqq \partial_\tau H(x, 0),$$

$$\nu(x) \coloneqq \partial_\tau g(x, 0),$$

are known as the Noether multipliers. This essentially means that invariances of a functional which can be differentiated give rise to conservation laws. In more quantitative terms, the invariances of the functional give rise to differential equations which a minimizer of the functional must satisfy[31].

The theorem can also be stated in more geometric terms if preferred. Start with a one-parameter group of motions $\varphi_s\colon M \to M$ on $M$, which generates a vector field $X$ by

$$X_q f = \frac{\partial f\big(\varphi_s(q)\big)}{\partial s},$$

when the derivative is evaluated at $s = 0$. If this looks strange, recall that vector fields are operators on the space of smooth functions such that a vector field acting on a real-valued function is another real-valued function $Xf$. The tangent map $(\varphi_s)_*$ induces a local flow on the tangent bundle and the Lagrangian is invariant under the one-parameter group if

$$L((\varphi_s)_* u) = L(u),$$

for all $u$ in the tangent bundle $TM$. Noether's theorem states that $\langle \omega_u, X\rangle$ is a constant of the motion, where $\omega$ is a 1-form on $M$ defined as

$$\omega = \left(\frac{\partial L}{\partial \dot{q}^i}\right) \mathrm{d}q^i.$$

---

[30] Emmy Noether, 'Invariante Variationsprobleme', *Nachrichten von der Königlichen Gesellschaft der Wissenschaften zu Göttingen, Mathematisch-Physikalische Klasse*, 235 – 257 (1918).

[31] Filip Rindler, *Calculus of Variations* (Cham, Switzerland: Springer, 2018).

This can be proved in local coordinates on the tangent bundle by setting

$$\dot{\mathbf{q}}(s,t) = \varphi_s \mathbf{q}(t),$$

where $q^i(t)$ is a solution of the Euler-Lagrange equations and differentiating with respect to $t$. After a computation and assuming the invariance of the Lagrangian, we end up with

$$\frac{\partial}{\partial t}\left(p_i X^i\right) = 0.$$

However,

$$\langle \omega_u, X\rangle = g_{ij}\dot{q}^i X^j = p_i X^i.$$

The term on the far right must be a constant from the previous equation, so this establishes the result. The one-parameter group is usually known as a symmetry group for the system, hence there is a physical relationship between the conservation laws of a system and its mathematical symmetries[32].

I should emphasise that this has been an extremely brief and selective survey from Noether's body of work and given the space constraints I have chosen to focus on her most basic and important contributions to algebra and mathematical physics. Noether worked in many other areas besides commutative algebra and I have not touched upon her other contributions: see, for example, her work on early algebraic topology and representation theory of groups and algebras.

[32] Peter Szekeres, *A Course in Modern Mathematical Physics: Groups, Hilbert Space and Differential Geometry* (Cambridge: Cambridge University Press, 2004).

# *Mary Cartwright*

Cartwright was the first female mathematician to be elected as a Fellow of the Royal Society and the first woman to be elected President of the Mathematical Association and the London Mathematical Society. She was also the first woman to receive the Sylvester Medal, a prize which had previously been awarded to Poincaré, Cantor, Darboux and Levi-Civita. She was supervised in her doctoral studies by G. H. Hardy and E. C. Titchmarsh and went on to do important collaborative work with J. E. Littlewood, including foundational work on chaos theory and non-linear dynamics. Curiously, there is no mention of Cartwright or Kovalevskaya in the well-known popular book on chaos theory by Gleick[33]. Chaos theory did not start to gain in popularity until the more visual approach of Lorenz. The approach of Cartwright and Littlewood was really analytic in nature, and in fact, Cartwright is regarded primarily as an analyst, rather than a dynamicist.

Cartwright's first work on analysis concerned the zeroes of certain types of entire function. Recall that a complex-valued function is entire if it is holomorphic over the finite complex plane (in older terminology, an entire function may be referred to as an 'integral' function)[34]. Typical examples are the exponential and trigonometric functions. A function of this kind will have a Taylor series expansion

$$f(z) = a_0 + a_1 z + \cdots + a_n z^n + \cdots,$$

which converges for all finite $z$. If $f(z)$ has a regular point at infinity, it must be identically constant due to Liouville's theorem, a basic theorem which says that a holomorphic function $f(z)$ for which there exists a number $M > 0$ such that $|f(z)| \leq M$ for all $z \in \mathbb{C}$ is a constant function. The theorem allows for a convenient proof of the fundamental theorem of algebra. If $f(z)$ has an essential singularity at infinity, $f(z)$ is said to be entire transcendental.

We are often concerned mostly with the transcendental type when working with entire functions. Define the maximum modulus function:

$$M(r) = \max_{|z|=r} |f(z)|.$$

This is a strictly increasing function of $r$. If $f(z)$ is an entire transcendental function with maximum modulus function $M(r)$, it follows that

$$\lim_{r \to \infty} \frac{\log M(r)}{\log r} \longrightarrow \infty.$$

In words, an entire transcendental function grows faster than a positive power of $r$. The simplest rapidly growing entire transcendental function is the exponential function, which suggests that we could quantify the growth of $|f(z)|$ using the exponential function. Suppose there exists a positive number $M(r)$ such that

$$M(r) < e^{r^{\mu}}$$

for sufficiently large $r$. In that case, the entire function is of finite order and the greatest lower bound on the $\inf$ of $\mu$ for which the previous inequality holds is called the order of $f(z)$. The complex exponential function is of order $1$, for example. It is another basic theorem of entire functions that the order $\rho$ of an entire function is given by a formula

$$\rho = \overline{\lim_{r\to\infty}} \frac{\log\log M(r)}{\log r}.$$

[33] James Gleick, *Chaos: Making a New Science* (New York: Vintage, 1997).
[34] Mary Cartwright, *Integral Functions* (Cambridge: Cambridge University Press, 1956).

This can be proved by using the previous definitions and finding a sequence $\{r_n\}$ such that

$$M(r_n) > e_n^{r^{\rho-\varepsilon}}.$$

One can also show that

$$\rho = \overline{\lim_{n\to\infty}} \frac{\log n}{\log\left(\frac{1}{\sqrt[n]{|a_n|}}\right)}.$$

Furthermore, suppose for a function of finite order $\rho$ that there is a number $K > 0$ such that

$$M(r) < e^{Kr^{\rho}}$$

for sufficiently large $r$. The greatest lower bound on the $\inf$ of $K$ for which this inequality holds is called the type of $f(z)$. Another theorem states that the type $\sigma$ of an entire function of order $\rho$ is given by a formula

$$\rho = \overline{\lim_{r\to\infty}} \frac{\log M(r)}{r^{\rho}}.$$

For example, the complex sine function is defined by

$$\sin z = \frac{e^{iz} - e^{-iz}}{2i}.$$

It can be shown with some fairly trivial estimates that the order and type of this function are both 1. The coefficients of the Taylor series expansion of an entire function also depend on its order and type.

We earlier mentioned Liouville's theorem from elementary complex analysis. This result can be improved significantly with Picard's little theorem. This states that if $f(z)$ is an entire transcendental function of order $\rho$, with $A$ any finite complex number, then the set of $A$-points of $f$ is infinite unless $\rho \in \mathbb{Z}^{+}$ and $f$ is of the form

$$f(z) = A_0 + p(z)e^{P(z)},$$

where $p(z)$ and $P(z)$ are both polynomials[35]. Cartwright herself contributed too many results to the theory of entire functions to list here, but we can quote a few examples. Many of her results were concerned with entire functions of a form such that

$$\log|f(z)| < A|z|,$$

for large $z$ and some constant $A$, and the conclusions which could be drawn about the zeroes of this type of $f$ given certain assumptions. For example, Cartwright proved that if we assumed that

$$\int_{-\infty}^{\infty} \frac{\log^{+}|f(x)|}{1+x^2} dx < \infty,$$

then the number of zeroes of $f(z)$ in the disc $|z| < r$ satisfies the following relation for some constant $B$

$$n \sim r.$$

Also, if the zeroes of $f(z)$ take the form

$$z_n = r_n e^{i\theta_n},$$

then most of the zeroes lie on the real line[36]:

$$\sum \frac{|\sin\theta_n|}{r_n} < \infty.$$

The set of functions satisfying the first assumption is sometimes known as the Cartwright class.

Another result of Cartwright compares $f(x)$ with a function which changes as slowly as possible ie. a constant. Take $f(z)$ to be holomorphic over the interval $|\arg z| \leq \alpha \leq \pi/2$, and let

$$h(\theta) \leq a|\cos\theta| + b|\sin\theta|,$$

where $a, b < \infty$. If $b < \pi$ and $|f(n)| \leq M$, then $|f(x)|$ is bounded for $x \geq 0$. If $f(n)$ tends to zero as $n$ tends to infinity, then $f(x)$ also tends to zero as $x$ tends to infinity. A more specific version of this theorem states that if $f(z)$ is an entire function satisfying the above inequality for $h(\theta)$ such that $b < \pi$ and $|f(n)| \leq M$ where $n = 0, \mp 1, \mp 2, \ldots$, then $|f(x)| \leq K$ for all real $x$, where $K$ is a constant depending only on $b$[37].

Cartwright also studied holomorphic functions on the unit disc. A function is said to be univalent in the unit disc if the equation $f = w$ never has more than one root in the disc. For these functions, one can write down an inequality relating the first coefficient in the corresponding Taylor series to the zero coefficient. We can have a corresponding function which is $p$-valent if the equation $f = w$ never has more than $p$ roots. It was conjectured in this case that one would have an inequality

$$|f| < C(1-r)^{-2p}.$$

[35] Alexei Markushevich, *Theory of Functions of a Complex Variable* (New York: Chelsea Publishing Company, 1985).
[36] Mary Cartwright, 'On functions which are regular and of finite order in an angle', *Proc. London Math. Soc.* (2) 38, 158 – 179 (1935).
[37] Ralph Boas, *Entire Functions* (New York: Academic Press, 1954).

Cartwright proved this conjecture in full generality[38].

Continuing in the theme of entire functions, Cartwright demonstrated that an entire function of order $\rho$ can have at most $2\rho$ asymptotic values: the proof uses the method of conformal maps. Two copies of the complex plane can be related by a map. If this map preserves the angle between two curves which intersect in the first plane at coordinates $(x_0, y_0)$, then the map is conformal at the point $(x_0, y_0)$. If $f(z)$ is holomorphic and $f'(z) \neq 0$ at $z_0 = x_0 + iy_0$, then the map $f(z)$ is conformal at $z_0$. If $f'(z) \neq 0$ in a region, then the map is conformal at all points of the region. A conformal map always takes the lines for constant $x$ and $y$ into two families of curves which intersect each other at right angles. A standard example of a conformal map is the Schwarz-Christoffel mapping. This transformation maps a polygonal region in the first complex plane (the $w$ plane) to the upper half of the $z$-plane. The polygon boundary becomes the $x$-axis and the interior of the polygon becomes the upper half-plane. The existence of this map is guaranteed by the Riemann mapping theorem, which says that there always exists a biholomorphic map $f$ which maps a simply connected open subset of $\mathbb{C}$ onto the open unit disc. Explicitly, the Schwarz-Christoffel mapping is given by

$$\frac{dw}{dz} = A(z-x_1)^{\frac{\alpha_1}{\pi}-1}(z-x_2)^{\frac{\alpha_2}{\pi}-1} \dots (z-x_n)^{\frac{\alpha_n}{\pi}-1},$$

which can theoretically be integrated to get a map

$$w = A\int (z-x_1)^{\frac{\alpha_1}{\pi}-1}(z-x_2)^{\frac{\alpha_2}{\pi}-1} \dots (z-x_n)^{\frac{\alpha_n}{\pi}-1}\, dz + B,$$

where $A$ and $B$ are complex constants and the $\alpha_i$ are the internal angles of the polygonal region. This integration can usually only be carried out for special cases.

For example, consider an equilateral plane triangle whose angles are all $\pi/3$ radians. The first vertex of the triangle is at the origin $O$, the second vertex $A$ is at a distance $a$ along the $u$-axis ($u$ and $v$ being equivalent to $x$ and $y$ for this complex plane) and the third vertex $B$ is given by coordinates $\left(\frac{a}{2}, \frac{\sqrt{3}}{2}a\right)$. The magnitude of the length of the line joining $O$ to $B$ is also $a$. We wish to map to a section of the upper half-plane given by taking $B$ to $x = -1$, $O$ to $x = 1$, and leaving aside the point $A$ by sending it to infinity. This gives us some values which we can plug into our formula:

$$x_1 = -1, \quad x_2 = 1.$$

The interior angles which we are considering are both the same:

---

[38] Mary Cartwright, 'On analytic functions regular in the unit circle II', *Quart. J. Math. Oxford Ser.* (2) 6, 94 – 105 (1935).

$$\alpha_1 = \frac{\pi}{3}, \quad \alpha_2 = \frac{\pi}{3}.$$

This gives us

$$\begin{aligned}\frac{dw}{dz} &= C(z+1)^{-\frac{2}{3}}(z-1)^{-\frac{2}{3}}, \\ &= C(z^2-1)^{-\frac{2}{3}},\end{aligned}$$

Bear in mind that I am changing the value of the constant as I go along although the symbol for it does not change: this is common practice. Hence we have

$$w = C\int_{z_0}^{z} \frac{dt}{(t^2-1)^{\frac{2}{3}}} + D.$$

The values of $z_0$, $C$ and $D$ are constants to be found (although the lower limit of integration is not too important). Set $z_0 = 0$.

$$w = C\int_{0}^{z} \frac{dt}{(t^2-1)^{\frac{2}{3}}} + D.$$

We know that the point at $w = 0$ is mapped to $z = 1$, which implies that

$$D = -C\int_{0}^{1} \frac{dt}{(t^2-1)^{\frac{2}{3}}}.$$

By considering the limits of integration, we have

$$w = C\int_{1}^{z} \frac{dt}{(t^2-1)^{\frac{2}{3}}}.$$

We also know that the point $w = ae^{\pi i/3}$ is mapped to $z = -1$, so we must have

$$ae^{\pi i/3} = C\int_{1}^{-1} \frac{dt}{(t^2-1)^{\frac{2}{3}}}.$$

Also, because $A$ gets sent to infinity, we have

$$a = C\int_{1}^{\infty} \frac{dt}{(t^2-1)^{\frac{2}{3}}}.$$

Using some tricks, you can determine the value for $C$ and you then have the conformal map $w$ written down as explicitly as possible, noticing that the integral cannot be carried out using only elementary functions.

We earlier discussed examples of rigid symmetric tops. For example, the Kovalyevskaya top has equations of motion can be integrated to provide analytic solutions (although this is highly non-trivial for the case which Kovalevskaya studied and involves many pages of calculations). There are a lot of interesting problems where an exact solution cannot be found, so the next best thing which a physicist would like to study is a problem which has a potential with an integrable part and a small part which can be dealt with using standard perturbation theory. As an example of perturbation theory, take the following differential equation:

$$\frac{d^2y}{dx^2}+y=\frac{\cos 2x}{1+\varepsilon y}, \qquad y\left(-\frac{\pi}{4}\right)=y\left(\frac{\pi}{4}\right)=0, |\varepsilon| \ll 1.$$

It is crucial that the size of $\varepsilon$ is small. Let

$$y=y_0+\varepsilon y_1+\varepsilon^2 y_2+O(\varepsilon^3).$$

In that case, we can use the binomial expansion to re-write

$$\begin{aligned}\frac{1}{1+\varepsilon y}&=1-\varepsilon y+\varepsilon^2 y^2+O(\varepsilon^3),\\ &=1-\varepsilon\big(y_0+\varepsilon y_1+O(\varepsilon^2)\big)+\varepsilon^2\big(y_0^2+O(\varepsilon)\big)+O(\varepsilon^3),\\ &=1-\varepsilon y_0-\varepsilon^2 y_1+\varepsilon^2 y_0^2+O(\varepsilon^3).\end{aligned}$$

Truncate to second order and substitute everything in:

$$\frac{d^2}{dx^2}(y_0+\varepsilon y_1+\varepsilon^2 y_2)+(y_0+\varepsilon y_1+\varepsilon^2 y_2)=(1-\varepsilon y_0-\varepsilon^2 y_1+\varepsilon^2 y_0^2)\cos 2x.$$

$$y_0''+y_0+\varepsilon(y_1''+y_1)+\varepsilon^2(y_2''+y_2)=(1-\varepsilon y_0-\varepsilon^2 y_1+\varepsilon^2 y_0^2)\cos 2x.$$

Collect like terms and construct an ODE for each coefficient $\varepsilon^k$ for $k=0,1,2$.

$$y_0''+y_0-\cos 2x,$$

$$y_1''+y_1+y_0\cos 2x=0,$$

$$y_2''+y_2-(y_0^2-y_1)\cos 2x=0.$$

The boundary conditions for each ODE are

$$y_i\left(-\frac{\pi}{4}\right)=y_i\left(\frac{\pi}{4}\right)=0.$$

These equations can be solved by hand or using a computer if preferred. The first equation can be solved to get

$$y_0=-\frac{1}{3}\cos 2x.$$

Plug this into the second equation to get

$$y_1''+y_1=\frac{1}{3}\cos^2 2x=\frac{1}{6}(1+\cos 4x).$$

This can be solved by plugging in a trial function and using the method of undetermined coefficients. Combined with the boundary conditions, we find that

$$y_1=\frac{1}{6}-\frac{8\sqrt{2}}{45}\cos x-\frac{1}{90}\cos 4x.$$

Solving the third equation for $y_2$, we find that the approximate solution to the original ODE in second-order perturbation theory is

$$y = -\frac{1}{3}\cos 2x + \varepsilon\left(\frac{1}{6} - \frac{8\sqrt{2}}{45}\cos x - \frac{1}{90}\cos 4x\right) + \varepsilon^2\left(\frac{2\sqrt{2}x}{45}\sin x - \frac{\sqrt{2}}{90}(\pi+1)\cos x + \frac{7}{270}\cos 2x - \frac{\sqrt{2}}{90}\cos 3x - \frac{1}{1050}\cos 6x\right).$$

This type of working is usually not shown and the author will just quote the result and say that it is derived using perturbation theory (they may even have run it straight through a computer).

However, we might have a weak interaction term in a potential which couples two of the equations of motion together and so which cannot be easily treated with perturbation theory. If the term is no longer 'small' in size, the solutions to the uncoupled equations could be considerably different to the solutions for the coupled equations. These new solutions might be nice, but they could also behave badly: it could be that a tiny change to the initial conditions leads to a massive change in the motion. These solutions are chaotic, in the sense that they are sensitive to initial conditions. Two bounded regular solutions which begin close to each other in a small region of phase space will evolve over time to cover a region which is still quite small. This is due to a fundamental theorem known as Liouville's theorem, which states that the density of points representing systems which are close to a point representing another system in phase will remain constant over time. In equation form:

$$\frac{dD}{dt} = [D,H] + \frac{\partial D}{\partial t}.$$

A chaotic phase space trajectory is something between a regular trajectory from a solution of integrable equations of motion and the complete randomness which we would normally associate with stochastic processes. Cartwright and Littlewood were the first to realize that chaotic phase space trajectories arise from the motion due to non-linear systems[39]. This type of motion is non-periodic, but not impossible to predict.

To begin with, we need to go back to our old picture of the simple harmonic oscillator which follows closed periodic motion. We will not go into the details, but it should be believable to you that the motion of such an oscillator translates into a uniform circular trajectory in phase space. This allows for phase space representation of the motion of an uncoupled double oscillator: the circular motions of the lower-frequency and higher-frequency harmonic oscillators together generate a torus. In real physical situations, we are often concerned with the dynamics of a system which can be described by an integrable Hamiltonian plus a small interaction. Consider the orbit of a planet in the Solar System which is periodic, but very weakly perturbed by the presence of another planet. One could write the Hamiltonian as

$$H = H_0 + \Delta H.$$

[39] Mary Cartwright and John Littlewood, 'On non-linear differential equations of the second order', *J. London Math. Soc.,* 20, 180 (1945).

If $\Delta H$ is small relative to $H_0$, perturbation theory provides us with a solution to the equations of motion, but it is not clear whether the perturbed solution is regular or if the phase space orbits will remain close to the phase space paths for the unperturbed solutions. Roughly speaking, the Kolmogorov-Arnold-Moser theorem tells us that if an integrable Hamiltonian is perturbed by a small perturbation which rends it non-integrable, then the motion will remain confined to an $N$-torus similar to the one which we have already described, as long as the perturbation is small and if the frequencies $\omega_i$ are incommensurable. This is with the exception of a set of initial conditions which cause a wandering trajectory. In more topological terms, an integrable system is represented by a torus in phase space. This torus remains a torus when we perturb the system and try to deform it, but under some perturbations, the torus can be destroyed.

Stable orbits of the integrable Hamiltonian generally continue to be stable when modified to be orbits of the total perturbed Hamiltonian. Another important possibility is that the initial conditions might start the motion on a trajectory which evolves over time to some fixed point or stable orbit in phase space. These are both examples of attractors. An attractor is a set in phase space to which the solution evolves over a long time. The dimension of a regular attractor is one less than the dimension of the phase space. The attractors might be $N$-tori, but there can also be some so-called strange attractors with fractal dimension ie. the dimension is a fraction or an irrational number[40]. The Lorenz equations and Henón map are examples of strange attractors. An example of a stable phase space orbit towards which a motion could evolve is given by the van der Pol equation used to describe driven oscillations in mechanical and electrical systems:

$$m\frac{d^2x}{dt^2} - \varepsilon(1 - x^2)\frac{dx}{dt} + m\omega_0^2 x = F\cos\omega_D t.$$

This was the original oscillator which Cartwright and Littlewood were studying when they came upon chaotic motion. They noticed that for some parameter ranges one would obtain a stable motion, whereas in other ranges the Poincaré section has a variety of structures[41]. When $\varepsilon = 0$ this is a simple harmonic oscillator with forcing, where the resonant frequency is $\omega_0$ and the forcing frequency $\omega_D$. With the second term for damping, the motion will tend towards a phase orbit given by a circle of unit radius. For positive damping, the motion will spiral inwards towards the limit cycle, and for negative damping, it will spiral outwards towards the limit cycle. If $\varepsilon$ is sufficiently large, the damping term becomes large, and the trajectory is still drawn towards the limit cycle, but the cycle is now deformed from a circular shape. If the damping becomes very large, the shape starts to become more like a square.

---

[40] Herbert Goldstein, Charles Poole and John Safko, *Classical Mechanics* (USA: Addison Wesley, 2002).
[41] Frank Wang, 'Pioneer women in chaos theory', Contribution to the conference 'Mathematics, Information Technology, and Education' held at Orenburg State University, Russia (2008).

We mentioned earlier that that chaotic trajectories are sensitive to initial conditions, but there are in fact several other technical properties which this type of motion has. Another one is mixing. This means that if we take two arbitrarily small regions in the domain and follow an orbit with passes through the first region, then it must also pass through the second region (perhaps after a long time). Another property is that the motion has dense quasi-periodic orbits. An orbit is quasi-periodic if it repeatedly passes through the domain in a non-regular way without closing itself off. An orbit is dense if it passes through or close to every point of the domain. The example we gave earlier of the torus in phase space due to an uncoupled double oscillator is an example of a dense periodic orbit, since the trajectory never closes, covering the torus and coming arbitrarily close to every point on the surface. This is a bounded, non-closed orbit. However, at least in the literature, it is usually sensitivity to initial conditions which is regarded in some sense as the phenomenon which 'causes' chaos, since one has an instability in the phase space which causes trajectories to separate exponentially. We know that the Kolmogorov-Arnold-Moser theorem is valid for small perturbations, but if the perturbation is sufficiently large, the behaviour may become chaotic. The orbits of a chaotic sequence remain close to the 'correct', original value after a few iterations, but the rate of increase of separation can rise exponentially as the number of iterations is increased.

The exponential divergence of trajectories is quantified via the Liapunov exponent. For a wide class of systems, if two chaotic trajectories are initially separated by a small distance $s_0$, then their separation will given by

$$s(t) \sim s_0 e^{\lambda t}$$

for a later time $t$, where $\lambda$ is the Liapunov exponent. If $\lambda > 0$ the motion will be chaotic, and $\lambda$ then quantifies the growth of a deviation of the trajectory for a regular solution due to a perturbation. It also sets a time scale for the growth of exponential divergences caused by large perturbations. When the time becomes much larger than this time scale $\tau$ becomes obvious, until the separation $s(t)$ is comparable in size to the dimension of the phase space, at which point the separation can only vary randomly without increasing. If we have a system which evolves by an iterative process (this might be possible in electrical engineering, for example), we will instead have

$$s(n) \sim s_0 e^{\lambda n},$$

where $n$ is the number of iterations. The divergence of phase space orbits in a chaotic region becomes something like a Markov process, where the present state of the system only allows us to deduce something at most about the previous state, and nothing about states in the distant past.

The Liapunov exponent can also be negative: in this case, it will quantify the rate at which a system approaches a regular attractor. If $\lambda < 0$ the motion is non-chaotic and the distance from an attractor at time $t$ is given by

$$s(t) \sim s_0 e^{-|\lambda| t},$$

for some initial distance $s_0$. For an iterative process, there is a similar expression to the one we had before:

$$s(n) \sim s_0 e^{-|\lambda| n},$$

for $n$ iterations[42]. There have been some suggestions in the literature that the criterion of the Liapunov exponent should be refined, since one can have unstable trajectories which converge very strongly over extremely long periods of time. This appears to be the case, for example, for the following simple chaotic system, which models small bodies in a turbulent flow:

$$\dot{x} = v,$$

$$\dot{v} = \gamma(u(x,t) - v,$$

where $\gamma$ is a constant for the rate of damping of the motion of a small particle and $u(x,t)$ is a velocity field which fluctuates randomly to imitate turbulence[43]. We should mention that we have overall presented chaos theory as something which develops from the study of non-integrable Hamiltonian systems in classical mechanics, but there are other ways of arriving at similar ideas: for example, by studying macroscopic physical equations such as the Navier-Stokes equations[44]. It is also possible to study the quantum mechanical behaviour of classically chaotic systems, a subject known as 'quantum chaos'[45].

We will finish by mentioning that the collaboration between Littlewood and Cartwright extended as far as topology. With Littlewood, Cartwright proved the following fixed point theorem: if $\tau$ is a one-to-one continuous and orientation-preserving transformation of the plane to itself which leaves a bounded continuum interval $I$ invariant, and if the complement of $I$ is one single simply connected domain, then $I$ contains a fixed point[46].

[42] Herbert Goldstein, Charles Poole and John Safko, *Classical Mechanics* (USA: Addison Wesley, 2002).

[43] Marc Pradas, Alain Pumir, Greg Huber and Michael Wilkinson, 'Convergent chaos', *J. Phys. A: Math. Theor.* 50 275101 (2017).

[44] Hao Bai-Lin, *Chaos* (Singapore: World Scientific, 1984).

[45] Hans-Jürgen Stöckmann, *Quantum Chaos: An Introduction* (Cambridge: Cambridge University Press, 2006).

[46] Mary Cartwright and John Littlewood, 'Some fixed point theorems', *Ann. of Math.* (2) 54, 1 – 37 (1951).

# *Julia Robinson*

Robinson was the first female mathematician to be elected to the National Academy of Sciences and the first woman president of the American Mathematical Society. She made contributions to game theory and decision problems, but is especially renowned for her fundamental work on Hilbert's tenth problem. The influential mathematician David Hilbert published a list of problems at the beginning of the twentieth century which he thought should guide research in that century. Most of the problems are now at least partially resolved, although there are some famous ones which still remain (most obviously the eighth problem, known as the Riemann hypothesis). This simply states that the real part of a non-trivial zero of the Riemann zeta function is always $1/2$. The Riemann zeta function is a function of a complex variable $z = a + ib$ which in the case where $\mathrm{Re}(z) > 1$ can be written as a convergent infinite series:

$$\zeta(z) = \sum_{n=1}^{\infty} \frac{1}{n^z}.$$

Technically speaking, the above is a very special case of what is called a Dirichlet series, and the Riemann zeta function is the analytic continuation of this function to the entire complex plane. Any arithmetic function can be attached to an infinite Dirichlet series which is the generating series of that function:

$$D_f(z) = \sum_{n=1}^{\infty} \frac{f(n)}{n^z}.$$

In the case where $\mathrm{Re}(z) > 1$, the Riemann zeta function can be represented as an Euler product:

$$\zeta(z) = \prod_{p} (1 - p^{-z})^{-1}.$$

This is an expression of the fact that the natural numbers can be uniquely factorized into prime powers.

An arithmetic function is a function which is defined only on the positive integers, the two most famous examples being the Möbius function and the Euler totient function. The Möbius function is defined as:

$$\mu(1) = 1,$$

$$\mu(n) = (-1)^k \text{ if } a_1 = \cdots = a_k = 1,$$

$$\mu(n) = 0 \text{ otherwise},$$

where $n$ is decomposed into prime powers as $p_1^{a_1} \ldots p_k^{a_k}$. The definition of this function seems strange, but it has many applications in number theory. The Euler totient function is the number of positive integers not above $n$ which are relatively prime to $n$. As an equation, this could be written as

$$\varphi(n) = \sum_{k=1}^{n} 1,$$

where we only sum over $k$ if $k$ and $n$ are coprime. If $n$ is greater than or equal to $1$, we can also write a nice formula which relates both these functions[47].

$$\varphi(n) = \sum_{d|n} \mu(d) \frac{n}{d}.$$

Many approaches have been suggested to try to prove the hypothesis. For example, there is the idea of a spectral proof where one finds a self-adjoint operator whose eigenvalues magically match up with the non-trivial zeroes of the Riemann zeta function. This idea caused a lot of interest when it was first suggested, but has not come close to working so far. The spectral proof is distinct from the 'spectral theory' of the Riemann zeta function, where one tries to deduce things about the zeta function using summation formulae involving a different type of analytic function[48]. Another interesting angle of attack comes from random matrix theory, where one uses the fact that the Riemann zeta function can be modelled by the characteristic polynomial of a large complex Hermitian random matrix whose eigenvalue distribution is known as the Gaussian unitary ensemble and uses this to make predictions about the spacings of the zeroes. In fact, there is the pair correlation conjecture of Montgometry, which says that the pair correlation between pairs of zeroes of the zeta function is the same as the pair correlation between pairs of eigenvalues of random Hermitian matrices.

There is numerical evidence which shows that in the limit the distribution of the zeroes of the zeta function approaches the curve given by the eigenvalue distribution of a GUE random matrix, but ultimately the random matrix approach is a model and so it can only make interesting predictions which then have to be proved rigorously by other means. For that reason, it seems extremely unlikely that random matrix theory will provide a proof of the Riemann hypothesis. It has to be said that some of the suggested approaches to proving the hypothesis seem unconvincing apart from the fact that they draw on several different branches of mathematics in an exciting way, but there is no evidence that just having a proof which links in with multiple areas of mathematics somehow makes the proof more convincing, as one could easily construct a fallacious proof for anything which draws on insights from ten different branches of mathematics.

Another very significant problem on the list is Hilbert's nineteenth problem: are the solutions of regular problems in calculus of variations always analytic? This does not seem that difficult, but in the case of regular variational integrals, we end up with a non-linear Euler-Lagrange equation. The problem was solved by De Giorgi in 1957 and Nash independently in 1958. De Giorgi's regularity theorem states that if we have $S: \Omega \to \mathbb{R}^{d\times d}$ measurable and symmetric which satisfies the following estimates

$$\mu|v|^2 \leq v^T S(x) v \leq M|v|^2,$$

for non-negative constants $\mu$ and $M$, then if $u$ in the Sobolev space $W^{1,2}(\Omega)$ weakly solves

$$-\mathrm{div}(S\nabla u) = 0,$$

---

[47] Tom Apostol, *Introduction to Analytic Number Theory* (New York: Springer, 1976).
[48] Yoichi Motohashi, *Spectral Theory of the Riemann Zeta-Function* (Cambridge: Cambridge University Press, 1997).

$u$ is $\alpha$-Hölder continuous for some $\alpha$ between $0$ and $1$ depending on $d$ and $M/\mu$. The theorem is called the regularity theorem, as it shows 'Hölder regularity'. Combined with standard Schauder estimates, we can obtain the De Giorgi-Nash-Moser theorem: take a regular variational integral $\mathcal{F}$ with an integrand $f: \mathbb{R}^{d\times d} \to \mathbb{R}$ which can be continuously differentiated $n$ times, then if $u_*$ is a minimizer of $\mathcal{F}$, then if $u_* \in C^{n-1,\alpha}_{\mathrm{loc}}(\Omega)$ for some $\alpha$ between $0$ and $1$. If the integrand $f$ is analytic, then the minimizer is also analytic[49]. Nash's methods were different to De Giorgi's, and it is worth mentioning that he was motivated by physical intuition with one eye on the obvious applications of a regularity result to the types of problem which arise in physics and applied mathematics (this was almost certainly part of Hilbert's motivation for introducing the problem).

Hilbert's tenth problem is easier than the Riemann hypothesis, but still an immensely difficult problem in its own right: it asks if one can find an algorithm which tells you whether a polynomial Diophantine equation with coefficients in $\mathbb{Z}$ has a solution in $\mathbb{Z}$. It was eventually shown by Matiyasevich that this algorithm does not exist, making extensive use of techniques developed by Robinson over a period of two decades of work on the problem. A Diophantine equation is a polynomial equation which is solved over the integers. This type of equation has been studied since ancient times, but it was not until the twentieth century that a systematic theory was developed for them. An example might be:

$$y^2 = x^3 + k,$$

where $k \in \mathbb{Z}$. One is then looking to see if the equation has solutions $x, y \in \mathbb{Z}$ for a given $k$. For example, in the case we have described, our equation has no solutions when $k$ is of the form

$$k = (4n-1)^3 - 4m^2,$$

where $m, n \in \mathbb{Z}$ such that no prime $p \equiv -1 \pmod{4}$ divides $m$. This can be proved via a contradiction by assuming that such a solution does exist. Another theorem states that the Diophantine equation

$$y^2 = f(x)$$

has a finite number of solutions when $f(x)$ is a polynomial of degree greater than or equal to 3 with distinct values such that $f(x) = 0$ and coefficients in $\mathbb{Z}$.

The most famous Diophantine equation of them all is

$$x^n + y^n = z^n.$$

---

[49] Filip Rindler, *Calculus of Variations* (Cham, Switzerland: Springer, 2018).

As we know, you could construct any number of solutions when $n = 2$ ($x = 3, y = 4, z = 5$ might be the first one that comes to mind), but Fermat's Last Theorem says that there are no positive integer solutions for $n \geq 3$: this was eventually proved by Wiles in the twentieth century using techniques that Fermat could not have had access to. You might also be familiar with linear systems of Diophantine equations. In general, a system of linear congruences might have no solution even if the congruences themselves have solutions, but one can show that a system of solvable congruences can be solved simultaneously much as you would a simultaneous system of linear equation, but only if the moduli are pairwise coprime: this is an elementary theorem called the Chinese remainder theorem. Recall that $a$ is congruent to $b$ modulo $m$ if $m$ divides $a - b$. This is written in equation form as

$$a \equiv b \pmod{m},$$

where $m$ is the modulus of the congruence. If you are not familiar with this notation, you might want to think of a clock face. If you start at the midday point on the clock and go forward by 14 intervals around the clock, the first 12 intervals will just take you back to the start, so you will only end up moving 2 intervals, hence

$$14 \equiv 2 \pmod{12}.$$

The Chinese remainder theorem states that if $m_1, \ldots, m_r$ are positive integers which are pairwise coprime, then the system of congruences

$$x \equiv b_1 \pmod{m_1}, \ldots, x \equiv b_r \pmod{m_r}$$

has one solution modulo $m_1, \ldots, m_r$.

The system of congruences can also be rewritten as a linear system of Diophantine equations:

$$x = a_1 + x_1 n_1, \ldots, x = a_k + x_k n_k.$$

for $x, x_i \in \mathbb{Z}$. You could solve this system by writing it in the usual matrix form

$$Ax = B,$$

and determining the Smith normal form of $A$[50]. In the context of algebra, you will probably see the theorem stated in terms of isomorphisms of rings, but this is an equivalent formulation. Linear Diophantine equations are pretty easy to deal with, since an equation of this type

$$ax + by = c$$

has a solution if and only if the greatest common denominator of $a$ and $b$ divides $c$ (this follows from an elementary theorem about the GCD). Writing down explicit solutions to the above equation is then just a matter of finding solutions to the equation

[50] Tom Apostol, *Introduction to Analytic Number Theory* (New York: Springer, 1976).

$$ax + by = \gcd(a, b),$$

which boils down to some computations with the Euclidean algorithm[51].

Diophantine equations are distinct from Diophantine approximation, which is a separate type of problem which occurs in number theory. Most of us are familiar with rational numbers. These are fractions $r = p/q$, where $p$ and $q$ are positive or negative integers. Any non-repeating decimal number represents an element of $\mathbb{Q}$ whose denominator does not contain a factor apart from 2 or 5. $\mathbb{Q}$ forms a countably infinite subset of $\mathbb{R}$ and a real number which is not also rational is an irrational number (almost every element of $\mathbb{R}$ is irrational, although there are some famous irrational numbers like $\pi$ and $e$). The question which occurs to a number theorist is exactly how well can real numbers be approximated by rational numbers. In the case of dimension 1, there is a strong link with continued fraction expansions. If we start with an irrational number $\alpha$, then there will be an infinite number of relatively prime integers $m$ and $n$ such that

$$\left|\alpha - \frac{m}{n}\right| < \frac{1}{n^2}.$$

Another way of seeing it is that $\alpha n$ is dense $(\mathrm{mod}\ 1)$ when $\alpha$ is irrational. For an irrational number $\alpha$ and a positive integer $q$, the vector

$$(x_1, \dots, x_1) = (\alpha \bmod 1, \dots, q\alpha \bmod 1)$$

has an element which is closest to $0 \bmod 1$. The discrepancy of the list $\alpha n \ (\mathrm{mod}\ 1)$ is controlled by the continued fraction expansion of $\alpha$. One thing which we can say is that if $\alpha \in \mathbb{R}$ is between 0 and ½ and $p, q \in \mathbb{Z}$ such that $(0/1, p_1/q_1, p_2/q_2, \dots)$ are the convergents of $\alpha$ and furthermore if $q > q_1$ and $|\alpha - p/q| \leq 1/2q^2$, then $p/q$ is a continued fraction convergent of the real number $\alpha$[52].

Although an algorithm for determining whether a Diophantine equation is solvable or not does not exist, there are some useful results which emerge from the work of Robinson and others: for example, the fact that a Diophantine set can formally be almost anything. Recall that we can write out a system of $m$ Diophantine equations. This system will only have a solution if the sum of the squares of each of the individual equations

$$D_1(x_1, \dots, x_n)^2 + \cdots + D_m(x_1, \dots, x_n)^2 = 0$$

has a solution, so we can always convert a problem of solving a system of Diophantine equations to the problem of solving one single Diophantine equation. We say that a set $S$ of $m$-tuples is Diophantine if there exists a family of Diophantine equations with parameters $\alpha_1, \dots, \alpha_m$ and unknowns $x_1, \dots, x_n$

$$D(\alpha_1, \dots, \alpha_m, x_1, \dots, x_n) = 0$$

which has a solution $x_1, \dots, x_n$ if and only if the $m$-tuple $(\alpha_1, \dots, \alpha_m)$ is contained in $S$. In the case where $m = 1$ this reduces to asset of non-negative integers being a Diophantine set if and only if there is a Diophantine family $D = 0$ which has a solution if and only if the parameter $\alpha$ belongs to that set. A set of positive integers is Diophantine if and only if it is the set of positive values of a polynomial which has non-negative integer values of the variables. The set of all even numbers is a Diophantine set because there is a solvable Diophantine family

$$a - 2x = 0.$$

[51] Martin Erickson and Anthony Vazzana, *Introduction to Number Theory* (New York: Chapman & Hall, 2008).
[52] Doug Hensley, *Continued Fractions* (Singapore: World Scientific Publishing, 2006).

The set of all odd numbers is also Diophantine (just add 1 to the family).

$$a - (2x + 1) = 0.$$

The Fibonacci numbers form a Diophantine set because there exists a Diophantine representation

$$y^2 - xy - x^2 = \mp 1.$$

This representation is in fact used in the unsolvability proof for the tenth problem. The class of all Diophantine sets is closed under the operators of intersection and union.

The overall problem which we wish to answer is whether there is an algorithm which can tell us if a Diophantine equation is solvable, so we have to be clear about what an algorithm is and what it means to compute an algorithm (in fact, the proof of the tenth problem might be unique in the extent to which it involves intersections between logic and number theory). For our purposes, it is sufficient to think of an algorithm as a procedure which is followed to generate positive numbers. You can check that the Euclidean algorithm is in fact an algorithm when you follow this definition. The algorithm can be used by anyone, regardless of the amount of ingenuity required to discover it: you put in the input, and then obtain an output. Turing famously imagined a computer called the Turing machine which could compute any algorithm and there is a hypothesis which states (informally speaking) that any algorithm which is capable of being computed is computable by a Turing machine[53].

A set of natural numbers which can be produced via computations carried out with a Turing machine is called 'recursively enumerable'. Within this category, one can imagine a set such that the Turing machine can tell is a given natural number is a member of this set. If this can be done for all natural numbers, then the set is said to be recursive. A set of natural numbers is recursive if and only if both the set and its complement are recursively enumerable. It seems difficult to think of a set which is not recursively enumerable, but this occurs when you attempt to determine whether a program can keep running for an infinite amount of time (a program being a code $P$ along with its data $d$). It is an important result that every Diophantine set is recursively enumerable: you can simply write an algorithm which takes a Diophantine family and goes through all the $n$-tuples to see if they are solutions. One can represent a set of non-negative integers $S$ with an arithmetic formula $\mathfrak{F}$: we say that the set is represented by a formula with a variable $a$ if there is an equivalence $a \in S \Leftrightarrow \mathfrak{F}$. Gödel showed that any recursively enumerable set is represented by an arithmetic formula of some kind, but Robinson and Matiyasevich significantly improved this result by showing that all recursively enumerable sets of non-negative integers can be represented by two kinds of arithmetic formula which can be written down explicitly and which only contain 3 quantifiers[54]:

$$\exists b\ \exists c\ \&\ \exists d\ [P_\iota(a,b,c) < D_\iota(a,b,c)d < Q_\iota(a,b,c)],$$

$$\exists b\ \exists c\ \forall f\ [f \leq F(a,b,c) \Longrightarrow W(a,b,c,f) > 0].$$

[53] Martin Erickson and Anthony Vazzana, *Introduction to Number Theory* (New York: Chapman & Hall, 2008).

We are now in a position to sketch the reasons as to why there is a negative solution for the tenth problem. Robinson, Matiyasevich and others had shown that every recursively enumerable set is Diophantine. However, we have already said that every Diophantine set is recursively enumerable. Since the correspondence goes both ways, it must be the case that Diophantine sets and recursively enumerable sets are the same. The recursive sets form a sub-class of the recursively enumerable sets, so let's take a set of non-negative integers $S$ which is recursively enumerable but not recursive. Since $S$ is also Diophantine, it has a representation via a Diophantine family

$$D(a, x_1, \dots, x_n) = 0.$$

If the tenth problem could be solved, then the decision problem would be solved for this family. This is a contradiction, so the tenth problem cannot be solved. The key step, then, is to prove that all recursively enumerable sets are Diophantine. Davis, Putnam and Robinson began by first proving that all recursively enumerable sets are exponential Diophantine[55]. Any Diophantine set has an associated Diophantine function: this is just a function such that the graph $\{(m,n): n = f(m)\}$ is a Diophantine set. An example of a Diophantine function might be the factorial function. One can define an exponential Diophantine set such that the defining function is built from polynomials and an exponential function. It is not too hard to prove that all exponential Diophantine sets are Diophantine if the exponential function itself is Diophantine. It was already known that the exponential function is Diophantine as long as one could exhibit some Diophantine function which *grows* exponentially. Matiyasevich proved that such a Diophantine function exists in the form of a function $n = f_{2m}$, where $f_{2m}$ is the $2m$-th Fibonacci number (recall that the Fibonacci numbers form a sequence with exponential growth).

There are however simpler proofs which use solutions to a version of the Pell equation

$$x^2 - (a^2 - 1)y^2 = 1.$$

The original Pell equation is

$$x^2 - dy^2 = 1.$$

---

[54] Yuri Matiyasevich and Julia Robinson, 'Two universal 3-quantifier representations of recursively enumerable sets', *Teoriya Algorifmov i Matematicheskaya Logika* (Moscow: Vychislitel'nyi Tsentr Akademii Nauk SSR, 1974) 112 – 123.

[55] Martin Davis, Hilary Putnam and Julia Robinson, 'The decision problem for exponential Diophantine equations', *Annals of Mathematics*, 74(3), 425 – 436 (1961).

This is a Diophantine equation, so given some $d$, a solution is a pair of non-negative integers. In fact, another proof alters the Pell equation further and considers the equation

$$x^2 - \lambda xy + y^2 = 1.$$

Take $\lambda \geq 2$. It can be proved that a pair of integers satisfies the above equation if and only if $(x, y) = (a_n, a_{n+1})$ or $(x, y) = (a_{n+1}, a_n)$ for a sequence $\{a_n\}$ defined by a recurrence relation

$$a_n = \lambda a_{n-1} - a_{n-2}, \qquad n \geq 2,$$
$$a_0 = 0, a_1 = 1.$$

Recall that a recurrence relation defines a sequence once one or more initial values have been specified. In the case of the Fibonacci numbers, the recurrence relation is

$$f_n = f_{n-1} + f_{n-2},$$

so two initial values need to be specified. You can usually spot the relation by inspection. For example, if we have a sequence of terms

$$a_1 = -\frac{1}{4}, a_2 = \frac{1}{36}, a_3 = -\frac{1}{576}, a_4 = \frac{1}{14400}, \ldots,$$

then the recurrence relation is

$$a_n = -\frac{a_{n-1}}{(n+1)^2}.$$

We are particularly interested in the sequence $\{a_n(\lambda)\}$ which satisfies a recurrence relation

$$a_n(\lambda) = \lambda a_{n-1}(\lambda) - a_{n-2}(\lambda), \qquad n \geq 2,$$
$$a_0(\lambda) = 0, a_1(\lambda) = 1.$$

Consecutive terms of this sequence are coprime. One can show that the function $a = a_b(\lambda)$ is Diophantine when $\lambda \geq 4$ and from there we can prove that the exponential function $a = \lambda^b$ is Diophantine. We start with inequalities

$$(\lambda - 1)^{n-1} \leq a_n(\lambda) \leq \lambda^{n-1}, n \geq 1.$$

One can obtain an upper bound on the quotient

$$\frac{a_{b+1}(\lambda x)}{a_{b+1}(x+1)} \leq \frac{(\lambda x)^b}{x^b} = \lambda^b,$$

as well as a lower bound

$$\frac{a_{b+1}(\lambda x)}{a_{b+1}(x+1)} \geq \frac{(\lambda x - 1)^b}{(x+1)^b} > \lambda^b - 1.$$

We can then write the exponential function in terms of the ceiling function (ie. the function which gives the smallest integer which is greater than or equal to the quantity in the brackets) of this quotient:

$$\lambda^b = \left\lceil \frac{a_{b+1}(\lambda x)}{a_{b+1}(x+1)} \right\rceil, x > 2ba_{b+1}(\lambda + 1) - 1.$$

It can be proved that the ceiling function of a quotient is Diophantine. This proves that the exponential function is Diophantine, which proves that Hilbert's tenth problem is unsolvable[56].

We will finish by mentioning Robinson's contributions to game theory. Probably the two most well-known researchers in this field are Nash and von Neumann (in fact, modern game theory did not really exist until the work of von Neumann). In essence, game theory studies the phenomena which are observed during interactions between rational decision-makers (or rather, it studies simplified mathematical models of those interactions). Mathematical models and the application of mathematics to real-world situations are certainly not restricted to physics, and occur heavily in finance and economics. We assume not only that the decision-makers are rational, but that they are using some kind of strategy which takes into account the expected behaviour of any other decision-makers in the game. A game which we all know (or some variant) is the prisoner's dilemma. Imagine that two suspects of a crime are placed into separate prison cells with no way of communicating or influence each other. If they both confess to the crime, they will both receive a sentence of three years in prison. If one of them confesses, he will walk free, and the other person will receive a sentence of four years. If neither of them confesses, they will both spend only one year in prison.

This is an interesting scenario, because the most overall outcome for both of the prisoners would be if neither of them said anything, but there is an incentive to take a risk and confess instead. In fact, whatever the other prisoner does, it is always better to confess rather than say nothing. This means that this game only has one Nash equilibrium, where both players confess. In the case where we play a game where two people wish to go out together to a shopping mall but one wants to go to one shop and the other prefers another shop, there are two possible Nash equilibria. The Nash equilibrium is a kind of strategic steady state into which a game falls. Mathematically, a strategic game is a triple $\langle N, (A_i), (\succcurlyeq_i) \rangle$ where $N$ is the set of players, $A_i$ is the set of actions which are available to each player $i$ and $\succcurlyeq_i$ is a binary preference relation on the set $A = \times_{j \in N} A_j$ for each player. A collection of values of a variable, one for each player, is called a profile. A Nash equilibrium of a game is a profile $a^*$ in the set $A$ such that for every player the following relation holds:

$$(a^*_{-i}, a^*_i) \succcurlyeq_i (a^*_{-i}, a_i) \text{ for all } a_i \in A_i.$$

For the profile $a^*$ to be a Nash equilibrium it must be that no player has an action which gives an outcome which he or she prefers to the outcome generated when they choose $a^*_i$, given that the other players choose their respective equilibrium actions $a^*_j$. Not every strategic game has a Nash equilibrium, and it is of interest to study the conditions under which we have equilibria for a game. Of course, we cannot say much if we try to talk about all games in general and so we have to make a restriction (all competitive two-player games, for example)[57]. It is possible to have some kind of adjustment process which leads to a Nash equilibrium. An example of such a process is called fictitious play, where a player always chooses a best response to the frequency of the other players' past actions during the game. Robinson proved that in any strictly competitive game, the process of fictitious play always converges to a mixed strategy Nash equilibrium, where an equilibrium of this kind is a mixed strategy profile $\alpha^*$ with the property that for every player, every action in the support of $\alpha_i^*$ is a best response to $\alpha^*_{-i}$. Informally speaking, a mixed strategy is a strategy such that the choices are being made according to probabilistic rules. It can be proved that every finite strategic game has a mixed strategy Nash equilibrium[58].

[56] Martin Erickson and Anthony Vazzana, *Introduction to Number Theory* (New York: Chapman & Hall, 2008).

[57] Martin Osborne and Ariel Rubinstein, *A Course in Game Theory* (Cambridge, Massachusetts: MIT Press, 1994).
[58] Julia Robinson, 'An iterative method of solving a game', *Annals of Mathematics* 54, 296 – 301 (1951).

# *Olga Ladyzhenskaya*

Ladyzhenskaya was a mathematician known for her work on PDEs and mathematical physics and for obtaining rigorous convergence results for numerical solution of PDEs. She was shortlisted for the Fields Medal (the most prestigious prize in mathematics) in 1958, indicating that female mathematicians had achievements which were worthy of a Fields Medal long before Mirzakhani was born[59]. Some of her best-known results were proofs of convergence for FDM (the finite difference method) applied to the Navier-Stokes equations. This is significant, as many numerical methods for solving PDEs (although useful in industry and engineering) lack proofs for convergence. In some cases, engineers will employ a numerical method which has no known convergence result at all or a method which is known to have bad convergence properties and simply 'hope' that everything converges to give a good result[60].

Computational fluid dynamics is a complicated subject in its own right requiring the user to learn various pieces of computer software, with the software often depending on the exact problem in hand: a problem about spherical bubble collapse, for example, might require a particular piece of software which needs to be learned and then altered and validated for the specific physical situation[61]. Generally speaking, almost every branch of science will have a numerical or computational sub-discipline which is a science in its own right (numerical general relativity, for example), but the same can also occur in pure mathematics, where one sees things such as computational group theory and computational topology, and it is quite common to see papers on number theory and algebra where a piece of computer software is used at some point for some calculations which would otherwise be intractable[62].

There are several well-known methods for numerical solution of PDEs, including FEM (finite element method), the finite difference method and FVM (finite volume method). The basic idea of these 'finite' methods is to divide a complicated real geometry or problem up into a mesh. For example, in nuclear fusion one might consider the tokamak which is used to contain plasma. The geometry of a tokamak is particularly complicated and difficult to model due to the presence of holes and cuts which cause it to deviate away from the geometry of a torus. The domain is now discrete because we only calculate values at certain nodes on the mesh. The PDE or PDEs defining the problem are made discrete by converting them into algebraic equations, which can be solved iteratively or simultaneously. If the resulting system is linear, this is useful as there are a number of iterative methods for solving linear systems.

[59] Michael Barany, 'The Fields medal should return to its roots', *Nature*. 553 (7688), 271 – 273 (2018).
[60] Zi-Cai Li, Yimin Wei, Yunkun Chen and Hung-Tsai Huang, 'The method of fundamental solutions for the Helmholtz equation', *Applied Numerical Mathematics*. Vol. 135, 510 – 536 (2019).
[61] Max Koch, Christiane Lechner, Fabian Reuter, Karsten Köhler, Robert Mettin and Werner Lauterborn, 'Numerical modelling of laser generated cavitation bubbles with the finite volume and volume of fluid method, using OpenFOAM', *Computers and Fluids*. 126, 71 – 90 (2016).
[62] Jeremy Rouse, 'Explicit bounds for sums of squares', *Mathematical Research Letters*. 19(2), 359 – 376 (2012).

A well-known iterative method is the Newton-Raphson method, which looks as follows in equation form:

$$x_+ = x_c - F'(x_c)^{-1}F(x_c).$$

We could use this to numerically find a root for a cubic equation given that the initial iterate is close to the root. The idea that we start close to the solution we are trying to find might seem prohibitive, but this could occur in a physical situation where you are trying to find the height of a body of water in a canal and know that the height is not going to deviate too far from the height of the canal which you already know. In the example below, one sees that after 15 iterations we are already converging in on the solution with an accuracy which might be permissible in applications.

| $n$ | $x_n$ | $n$ | $x_n$ | $n$ | $x_n$ |
|---|---|---|---|---|---|
| 1 | 3 | 6 | 2.579 | 11 | 2.521 |
| 2 | 2.758 | 7 | 2.561 | 12 | 2.519 |
| 3 | 2.688 | 8 | 2.548 | 13 | 2.518 |
| 4 | 2.639 | 9 | 2.534 | 14 | 2.517 |
| 5 | 2.604 | 10 | 2.523 | 15 | 2.517 |

Table 1: Convergence to a root of a cubic equation using the Newton-Raphson method.

Another way of achieving an iterative solution is to rewrite the linear system as $A\mathrm{x} = b$ as a linear fixed-point iteration. For example, we write $A\mathrm{x} = b$ as

$$\mathrm{x} = (I - a) + b.$$

We then define the Richardson iteration as

$$\mathrm{x}_{k+1} = (I - A)\mathrm{x}_k + b.$$

In general we can also write

$$\mathrm{x}_{k+1} = M\mathrm{x}_k + c.$$

An alternative way of converting the system to a linear fixed-point iteration is to perform a splitting of $A$ of the form $A = A_1 + A_2$, where $A_1$ is a non-singular matrix. Examples of this type of method are the Jacobi and Gauss-Seidel methods. Both methods involve splitting $A$ and applying a scheme to it. In the Gauss-Seidel method, $A$ is split such that

$$A = D - L - U,$$

where $D$ is the diagonal of $A$ and $-L$ and $-U$ are the lower and upper triangular parts. Even with preconditioners and improvements, it can be shown that there is a more efficient method known as the conjugate gradient method. Unlike Gauss-Seidel, this is not a stationary iterative method, as the transition between iterates depends on the iteration history. The conjugate gradient method is ideal for large, sparse positive definite systems and it can be defined as

$$\mathrm{x}_{i+1} = \mathrm{x}_i + \alpha_i p_i,$$

where

$$\alpha_i = -\frac{(p_i, p_i)}{(p_i, Ap_i)}.$$

In FVM, the domain is discretized into non-overlapping elements known as finite volumes. The PDEs are integrated over each volume to obtain algebraic equations. The particular thing about FVM is that some of the terms in the conservation equation are volume integrals with a divergence term: these are converted to surface integrals using the divergence theorem

$$\int_V \nabla \cdot \mathbf{F}\, \mathrm{d}V = \oint_S \mathbf{F} \cdot \mathbf{n}\, \mathrm{d}S$$

and evaluated as fluxes at the volume faces. Since the flux which goes into a volume as the same as the flux leaving the neighbouring volume, the method is obviously conservative. FVM can also be formulated for polygonal meshes with inherent structure and it is fairly easy to implement different types of boundary conditions, since the unknown variables are not evaluated at boundary faces. All these attributes make FVM suitable for computational fluid dynamics[63].

VOF (volume of fluid) is a method which is particular to computational fluid dynamics. This method was historically used to model and follow the free surface of a flow of one fluid, although it is now common to apply it to the interface of a two-fluid flow (often a gas and a liquid in one flow). The first aim of the method is to reconstruct an approximation of the shape of the free surface or the interface by using knowledge of the volume fraction in each cell (the volume fraction being the fraction of a computational cell which is filled by the fluid assumed to be the reference phase). You may also see the terminology 'colour function' or 'marker function', where the word 'marker' refers to the fact that we have some function which can take different values in the different fluids in the system and the volume fraction is the discrete version of the marker function. Depending on the reconstruction method, it might be necessary to find the normal vector $\mathbf{m} = -\nabla C$ pointing away from the reference phase, where $C$ is a colour function which varies from taking the value 1 in a full cell, 0 in an empty cell and some value in between in a mixed transition cell. In the second step, the interface is advected in a velocity field: this is equivalent to exchanging volumes across boundaries of neighbouring cells.

Conversation of mass rather than flux is specific to VOF. For example, start with the advection equation for the marker function $H$ in an incompressible flow:

$$\frac{\partial H}{\partial t} + \nabla \cdot (\mathbf{u}H) = 0.$$

---

[63] Fadl Moukalled, Luca Mangani and Marwan Darwish, *The Finite Volume Method in Computational Fluid Dynamics: An Advanced Introduction with OpenFOAM and Matlab* (Switzerland: Springer, 2015).

Integrating this equation over a square cell of side length $h$ and using the definition of $C$, we obtain

$$h^2 \frac{\partial C_{i,j}(t)}{\partial t} + \int_\Gamma \mathbf{u} \cdot \mathbf{n}\, H(\mathbf{x}, t)\, \mathrm{d}l = 0,$$

where $\Gamma$ is the cell boundary line. If you integrate this over discrete time steps and sum the result over all grid cells with the correct boundary conditions, we obtain an equation for conservation of total area.

$$\sum_{i,j} C_{i,j}^{n+1} = \sum_{i,j} C_{i,j}^{n}.$$

In an interface reconstruction method where the normal vector is used, the normal in a cell is first determined using the colour function of that cell and the neighbouring cells, and one then obtains a geometric equation for the segment of the interface in question

$$\mathbf{m} \cdot \mathbf{xb} = \alpha,$$

where $\alpha$ is a parameter to be adjusted until the area under the interface equals the square of the length multiplied by the colour function of the cell[64].

As mentioned before, the particular method which Ladyzhenskaya contributed to was FDM. We will attempt to illustrate the principle of FDM using a one-dimensional example. One starts with a domain $\Omega$ which for the purposes of applications can be taken to be a $d$-dimensional unit cube. One defines a computational grid $G$, which is defined as

$$G = \left\{ x \in \mathbb{R} \text{ such that } x = x_j = jh, j = 0,1, \dots, n, \quad h = \frac{1}{n} \right\}.$$

One defines the forward and backward difference operators respectively as

$$\Delta\varphi_j = \frac{(\varphi_{j+1} - \varphi_j)}{h},$$

$$\nabla\varphi_j = \frac{(\varphi_j - \varphi_{j-1})}{h}.$$

Under certain assumptions, a general second-order elliptic PDE can be reduced to

$$-a \frac{\partial^2 \varphi}{\partial x_j^2} = q.$$

This PDE is used to study flows in porous media. If we take the one-dimensional case:

$$-a \frac{\mathrm{d}^2 \varphi}{\mathrm{d}x_1^2} = q,$$

---

[64] Grétar Tryggvason, Ruben Scardovelli and Stéphane Zaleski, *Direct Numerical Simulations of Gas-Liquid Multiphase Flows* (Cambridge: Cambridge University Press, 2011).

a finite difference approximation is obtained by replacing differential operators with the relevant difference operators. Probably the neatest possible formula would be

$$-\frac{1}{2}\big(\nabla(a\Delta)+\Delta(a\nabla)\big)\varphi_j, \qquad j=1,\dots,n-1.$$

This can be written out in full and certain terms eliminated or altered depending on the nature and form of each of the boundary conditions at $x=0$ and $x=1$.

In two dimensions, the computational grid becomes a vertex-centred grid, with the approximation of the second-order elliptic PDE in this case giving us a vertex-centred discretization. The difference operators are now defined in the analogous way:

$$\Delta_\alpha\varphi_j=\frac{\left(\varphi_{j+e_\alpha}-\varphi_j\right)}{h_\alpha},$$

$$\nabla_\alpha\varphi_j=\frac{\left(\varphi_j-\varphi_{j-e_\alpha}\right)}{h_\alpha},$$

where $e_1=(1,0)$ and $e_2=(1,0)$. As before, the finite difference approximation is obtained by replacing the relevant differential operators (partial derivatives this time) with some linear combination of the difference operators. Again, there is an analogous neat formula for the approximation:

$$-\frac{1}{2}\Big(\nabla_\beta\big(a_{\alpha\beta}\Delta_\alpha\big)+\Delta_\beta\big(a_{\alpha\beta}\nabla_\alpha\big)\Big)\varphi+\frac{1}{2}(\nabla_\alpha+\Delta_\alpha)(b_\alpha\varphi)+c\varphi=q.$$

This finite difference scheme relates $\varphi_j$ to $\varphi$ in the neighbouring grid points, with the set of neighbour points along with the grid points $\mathbf{x}_j$ being known as the finite difference stencil for the approximation. The stencil in this case is not symmetric but it is possible to have symmetric stencils, meaning that the discrete approximation inherits any symmetries of the solution[65].

Ladyzhenskaya's work on finite difference schemes began quite early when she proposed difference analogues for Fourier series and then used them to study finite difference schemes and prove rigorous results about their stability. In order to analyse a difference scheme with the difference analogue of the Fourier method, we do have to assume that the scheme corresponds to a problem whose solution can be represented using classical Fourier series. The usual results such as Parseval's theorem and boundedness of trigonometric polynomials with Fourier coefficients have analogues in the Fourier difference method. Other theorems then show that the stability of a scheme in fact depends on the convergence. As an example of a basic convergence result in difference schemes, start with $u$ in $L^2(\Omega)$ and let the difference quotients $u_{x_i}(x)$ converge weakly in $L^2(\overline{\Omega})$ to a function $u_i(x)$ as $h_i$ tends to 0 for all $\overline{\Omega}'\subset\Omega$. The difference quotients for the product of two grid functions $uv$ can be written as follows (where $u_{x_i}(x)$ and $u_{\bar{x}_i}(x)$ correspond to the forward and backward difference operators and $\mathbf{e}_i$ is a unit vector along the $x_i$-axis):

[65] Pieter Wesseling, *Principles of Computational Fluid Dynamics* (Berlin: Springer, 2001).

$$(uv)_{x_i}(x) = u_{x_i}(x)v(x) + u(x + h_i\mathbf{e}_i)v_{x_i}(x),$$

$$(uv)_{\bar{x}_i}(x) = u_{\bar{x}_i}(x)v(x) + u(x - h_i\mathbf{e}_i)v_{\bar{x}_i}(x).$$

If the above convergence occurs, the function $u_i(x)$ is just the generalized derivative $\partial u/\partial x_i$ of $u(x)$ in the domain $\Omega$. It follows from this that we can use any difference quotients which approximate the derivatives.

As another example, begin with the following inequality:

$$\Delta_h \sum_{\bar{\Omega}_h} u_h^2 \leq c, \qquad \Delta_h = h_1 \dots h_n,$$

where the $h_n$ are numbers ranging over sequences which have limit zero. If this inequality holds, one can prove that the fact that either one of the sequences $\{\tilde{u}_h\}$, $\{u_h'\}$ or $\{u_{(m)}\}$ converges weakly in $L^2(\Omega)$ to a function $u(x)$ as $h_1, \dots, h_n$ tend to zero implies that they all converge to $u(x)$ in the same way. $\tilde{u}_h(x)$, $u_h'(x)$ and $u_{(m)}(x)$ correspond to elementary types of interpolation of an arbitrary function $u_h$ for a grid. The simplest of these is the piecewise constant interpolation:

$$\tilde{u}_h(x) = u_h(kh), \qquad x \in \omega_{(kh)}$$

where $kh$ is a vertex such that $kh = (k_1h_1, \dots, k_nh_n)$ and $\omega_{(kh)} = \{x \colon k_ih_i < x_i < (k_i + 1)h_i\}$. There are many more convergence lemmas and theorems, plus the analogous results on compactness, pre-compactness and boundedness which one might expect or hope for.

All this theory can be put together and used to solve some illustrative PDEs: for example, the Dirichlet problem for a second-order elliptic PDE

$$\mathcal{L}u = f(x) + \frac{\partial f_i(x)}{\partial x_i},$$

where $u$ restricted to some boundary $S$ and the domain is bounded. One writes out the general form of a solution to the PDE in the appropriate Sobolev space and constructs an approximation to it by replacing integrals over the domain with integrals over cells. Derivatives are replaced with interpolations of approximating difference quotients. After considering the necessary identities, we end up with a system of difference equations to be satisfied at points of the grid $\Omega_h$:

$$\mathcal{L}_h u_h = \left(a_{ijh}u_{x_j} + a_{ih}u_h\right)_{\bar{x}_i} + b_{ih}u_{x_i} + a_h u_h = f_{ih\bar{x}_i} + f_h.$$

This is then coupled with the appropriate boundary condition to form a linear system whose unknowns are the values of $u_h$. One can prove stability of this scheme in the grid norm, plus a convergence result that the sequence of interpolations $\{u_h'(x)\}$ formed from the sequence of solutions to the difference system converges strongly in $L^2(\Omega)$ to the generalized solution of the original Dirichlet problem and that the derivatives $\left\{\frac{\partial u_h'}{\partial x_i}\right\}$ converge weakly in $L^2(\Omega)$ to $\partial u/\partial x_i$. This assumes that the coefficients of the original PDE are bounded, measurable functions which satisfy a particular condition[66].

Ladyzhenskaya is particularly well-known for publishing a large number of results relating to the Navier-Stokes equations for viscous fluid flow. The Navier-Stokes equations are certainly the most important PDEs in fluid mechanics and currently the most famous PDEs in the pure mathematics community. Proof of the existence and smoothness of solutions to the three-dimensional Navier-Stokes equations is considered to be one of the most important outstanding open problems across all branches of mathematics, with a very large cash prize for a solution of the problem. The non-linearity of the equations means they are notoriously difficult to solve in general. Interestingly, Stokes did not regard the Navier-Stokes equations as his greatest achievement in fluid dynamics, as he was aware that both Poisson and Navier already knew of the equations and merely provided an alternative derivation. (Cauchy and Saint-Venant also derived the equations before Stokes)[67]. Stokes was most pleased with his discovery of Stokes' law for the drag force on a small sphere moving through a viscous flow with low Reynolds number:

$$D = 6\pi\mu R v,$$

where $v$ is the speed of the sphere, $R$ is the radius and $\mu$ is the viscosity. The result is derived in various textbooks on fluid mechanics[68]. Since the result is analytic and exact, it is useful as a benchmark against which numerical methods can be compared: obviously, a method is no good if it cannot reproduce simple analytic solutions which have been known for many years. For example, one can linearize the Navier-Stokes equations to obtain the Stokes equations for Stokes flow (a type of flow where the forces due to inertia are small in comparison to the viscous forces). Stokes' law is an important result in the theory of Stokes flow. In numerical methods such as the boundary integral method which rely on Green's functions or other types of singular solution, one can obtain the Stokeslet (the Green's function for the Stokes equation) and then implement it in the method to model creeping flow around a small sphere and expect that the numerical result should match the analytic expression above[69].

A result which Ladyzhenskaya originally published in the context of two-dimensional Navier-Stokes boundary value problems is known as the Ladyzhenskaya inequality[70]. This is an example of an interpolation inequality and a particular case of the Gagliardo-Nirenberg inequality. Assume $u$ to be a function from $\mathbb{R}^2$ to $\mathbb{R}$ with compact support such that both $u$ and its gradient $\nabla u$ are in $L^2$. It follows that

$$\|u\|_{L^2}^4 \leq C\|u\|_{L^2}\|\nabla u\|_{L^2}.$$

---

[66] Olga Ladyzhenskaya, *The Boundary Value Problems of Mathematical Physics* (New York: Springer Science and Business, 1985).

[67] Mark McCartney, 'Fluids, fluorescence and a hat full of beetles'. *Mathematics Today*. Vol. 55 No. 4, 142 – 144 (2019).

[68] Lev Landau and Evgeny Lifschitz, *Fluid Mechanics* (Oxford: Pergamon, 2013).

[69] Constantine Pozrikidis, *Boundary Integral and Singularity Methods for Linearized Viscous Flow* (Cambridge: Cambridge University Press, 1992).

[70] Olga Ladyzhenskaya, 'Solution 'in the large' of the nonstationary boundary value problem for the Navier-Stokes system in two space variables'. *Comm. Pure Appl. Math.* 12, 427 – 433 (1959).

The inequality can be used to establish existence and uniqueness for weak solutions of the Navier-Stokes equations in two dimensions for a smooth bounded domain $\Omega$ and sufficiently smooth initial conditions[71]. If we use a general domain $\Omega$ (possibly bounded), then $u$ must belong to the Sobolev space $H_0^1(\Omega)$. Whether or not the inequality is true for a general domain $\Omega \subset \mathbb{R}^2$ will depend on the domain in question and the conditions we impose.

The existence proof is one which would never work in three dimensions, as it relies on the fact that there are 'nice' exponents in the Sobolev embedding theorem which allow us to obtain a regularity result for the time derivative of $u$:

$$\frac{\partial u}{\partial t} \in L^2(0,T;\ V^*),$$

for some long time $T$. This regularity allows one to use weak solutions of the Navier-Stokes equations as test functions whose uniqueness is relatively easy to prove. The first step in the proof is to say that a function $u$ is a weak solution of the Navier-Stokes equations if it satisfies

$$\int_0^s -\langle u, \partial_t \varphi\rangle + \int_0^s \langle \nabla u, \nabla \varphi\rangle + \int_0^s \langle (\mathrm{u}\cdot\nabla)u, \varphi\rangle = \langle u_0, \varphi(0)\rangle - \langle u(s), \varphi(s)\rangle,$$

where $u_0$ is an initial condition and $\varphi$ is a test function in a suitable space. One then sets $w$ to be the difference between two solutions with the same initial condition and substitutes it in as a test function[72]. The Sobolev embedding theorem is a basic result in the theory of Sobolev spaces. Start with a Sobolev space $W^{1,p}(\Omega)$ (that is, the space of all functions $u$ in $L^p(\Omega)$ whose first weak derivatives exist and are also in $L^p(\Omega)$). The theorem can be stated in separate parts depending on whether $p < d$, $p = d$ or $p > d$. If $p < d$, $u$ is in the space $L^{p^*}(\Omega)$ where

$$p^* \coloneqq \frac{dp}{d-p},$$

and there is a positive constant depending on $p$ and the domain such that

$$\|u\|_{L^{p^*}} \leq C\|u\|_{W^{1,p}}.$$

If $p = d$, then $u$ is in $L^q(\Omega)$ for $q$ between 1 and $\infty$ and

$$\|u\|_{L^q} \leq\ C\|u\|_{W^{1,p}}.$$

If $p > d$, then $u$ is in $C(\Omega)$ and

$$\|u\|_{\infty} \leq\ C\|u\|_{W^{1,p}}.$$

---

[71] Jacques-Louis Lions and Giovanni Prodi, 'Un théoréme d'existence et unicité dans les équations de Navier-Stokes en dimension 2'. *C. R. Acad. Sci. Paris* 248, 3519 – 3521 (1959).

[72] James Robinson, José Rodrigo and Witold Sadowski, *The Three-Dimensional Navier-Stokes Equations* (Cambridge: Cambridge University Press, 2016).

This is equivalent to the statement that for $1 \leq p < n$ and $p^*$ defined as above, there is a continuous embedding

$$W^{1,p}(\Omega) \subset L^{p^*}(\Omega).$$

There are also similar results for embeddings of Sobolev spaces known as Morrey's theorem and the Rellich-Kondrachov theorem[73]. In the study of finite difference schemes there are embedding theorems for grid functions which are the analogues of the theorems we have stated above.

Ladyzhenskaya was a key figure in the introduction of mathematical rigour into the study of fluid dynamics and was partly responsible for pushing mathematical fluid dynamics to its current status as a branch of mathematics, rather than physics or engineering. Experiment tells us that a steady, incompressible flow is stable as long as the Reynolds number is quite low, and that the flow becomes turbulent once the Reynolds number is increased past a certain point. It can be shown with more rigour that, given a stationary solution $v_S$ of the two-dimensional Navier-Stokes equations on a bounded domain, another solution $v$ with smooth initial data will approach the stationary solution with exponential speed

$$v - v_s = O(e^{-\alpha t}),$$

if what we define as the generalized Reynolds number is less than 1:

$$\mathrm{Re}^* = \frac{2v_*}{\nu\sqrt{\lambda_1(\Omega)}},$$

where $\nu$ is the viscosity, $\lambda_1$ is the first eigenvalue of the Dirichlet Laplacian on the domain, and $v_*$ is the characteristic velocity. Surprisingly, the parameter $\alpha$ which fixes the rate of convergence to the stationary solution can be stated in terms of the quantities above, one of which is traditionally physical and one of which is unfamiliar in classical, more experimental fluid dynamics[74]:

$$\alpha = \nu\lambda_1(\Omega)(1 - \mathrm{Re}^*).$$

This Reynolds number is to be distinguished from other 'starred' quantities in the literature, which usually refer to some non-dimensional quantity which has been scaled in some way: for example, one might multiply the Reynolds number by $\sin\beta$ for a problem related to a liquid film falling down a slope, where $\beta$ is the angle of inclination[75]:

$$\mathrm{Re}^* = \mathrm{Re}\ \sin\beta.$$

One final idea worth mentioning is the modelling of turbulent flow via a stochastic PDE known as the Ladyzhenskaya model: of course, the idea that one might model turbulence or thermal fluctuations in a flow using a PDE with stochastic terms in it has now become commonplace. The Ladyzhenskaya model models a higher-order variant of a type of non-Newtonian fluid known as a power-law fluid such that the stress tensor takes the form[76]

$$\tau(u) = 2\mu_0(1 + |e(u)|^2)^{\frac{p-2}{2}} e(u) - 2\mu_1 \Delta e(u).$$

---

[73] Filip Rindler, *Calculus of Variations* (Cham, Switzerland: Springer, 2018).

[74] Olga Ladyzhenskaya, *The Mathematical Theory of Viscous Incompressible Flow* (New York: Gordon and Breach, 1969).

[75] Fabian Denner, Marc Pradas, Alexandros Charogiannis, Christos Markides, Berend van Wachem and Serafim Kalliadasis, 'Self-similarity of solitary waves on inertia-dominated falling liquid films'. *Phys. Rev. E* 93, 033121 (2016).

As with any stochastic PDE or stochastic dynamical system, random attractors can be studied for this model, where the attractor of a random dynamical system is a random compact set to which the system 'tends over time' in some sense, although there are some subtle technical requirements in the definition which we cannot go into. Ladyzhenskaya herself did work on global attractors for various PDEs: the aim of much of this was to try to study and understand the set of all limit states of a system (this set being a type of global attractor), developing a rigorous theory of global stability for dissipation problems in mathematical physics[77].

---

[76] Olga Ladyzhenskaya, 'New equations for the description of the viscous incompressible fluids and solvability in the large of the boundary value problems for them', in *Boundary value problems of mathematical physics V*, (Providence, RI: Amer. Math. Soc., 1970).

[77] Olga Ladyzhenskaya, 'On finding the minimum global attractors for the Navier-Stokes equations and other PDE'. *Uspechi Math. Nauk.* 42, n. 6, 25 – 60 (1987).

## *Yvonne Choquet-Bruhat*

Choquet-Bruhat is known for her applications of PDE theory and differential geometry to mathematical physics and is responsible for much of the modern mathematical formulation of general relativity, with a unique simultaneous grasp of both the mathematics and the physics of the theory. Although she is a mathematician, Choquet-Bruhat emphasises that she has never lost sight of the fact that general relativity is a physical theory and takes an interest in observations and new experimental data, especially if it is relevant to GR[78]. She was the first woman to be elected as a full member of the French Academy of Sciences.

I will not recapitulate all the details of GR here, as this is really a book about mathematics, but essentially at the heart of it you have a set of PDEs called the Einstein equations, which tell you that the presence of matter causes the curvature of spacetime (roughly speaking).

$$G = 8\pi T,$$

where $G$ is the Einstein tensor for the curvature and $T$ is the energy-momentum tensor which describes classical flow and density of energy and momentum. Einstein was originally unsure of how he was going to deal with the geometric problem of curvature of spacetime manifolds, but realised by asking around that the basic mathematical theory of curvature of manifolds of arbitrary dimension had already been laid down many years before by Riemann. General relativity explains several effects which cannot be explained by Newtonian mechanics, and can also be applied to cosmology. There are many books and review articles which outline relativistic cosmology in more detail[79].

A famous correction from general relativity is to the precession of the perihelion of Mercury. The observed rate of precession of Mercury's perihelion disagrees with Newtonian mechanics, even after all other factors have been accounted for (the perihelion of an orbit of a body being the point where it comes closest to the Sun). Most of the observed secular precession of the perihelion of Mercury can be explained via precession of the equinoxes and perturbations of the orbit by other planets, but after these effects are dealt with, there is still an observed precession of approximately 43 arcseconds per century. If one takes the Kepler problem for a bound orbit perturbed by a central force and includes a perturbation Hamiltonian due to the Schwarzchild solution for weak gravitational fields, one obtains a value for the average rate of precession which is 42.98 arcseconds per century. Observation currently tells us that the precession due to relativistic effects is 43.1 $\mp$ 0.5 arcseconds per century, an insignificant deviation from theory[80].

To illustrate this visually, one can use any piece of computer algebra software to numerically solve an ODE for the orbit

---

[78] Yvonne Choquet-Bruhat, *General Relativity and the Einstein Equations* (New York: Oxford University Press, 2009).
[79] Hollis Williams, 'Applications of General Relativity to Cosmology', Conference Proceedings of Tomorrow's Mathematicians Today, Manchester Metropolitan University (2017).
[80] Herbert Goldstein, Charles Poole and John Safko, *Classical Mechanics* (USA: Addison Wesley, 2002).

$$\frac{2}{r^3}\left(\frac{\mathrm{d}r}{\mathrm{d}\phi}\right)^2 - \frac{1}{r^2}\frac{\mathrm{d}^2 r}{\mathrm{d}\phi^2} + \frac{1}{r} = 1,$$

with initial conditions

$$r(0) = \frac{2}{3}, r'(0) = 0.$$

We then create a plot of the orbit after several revolutions.

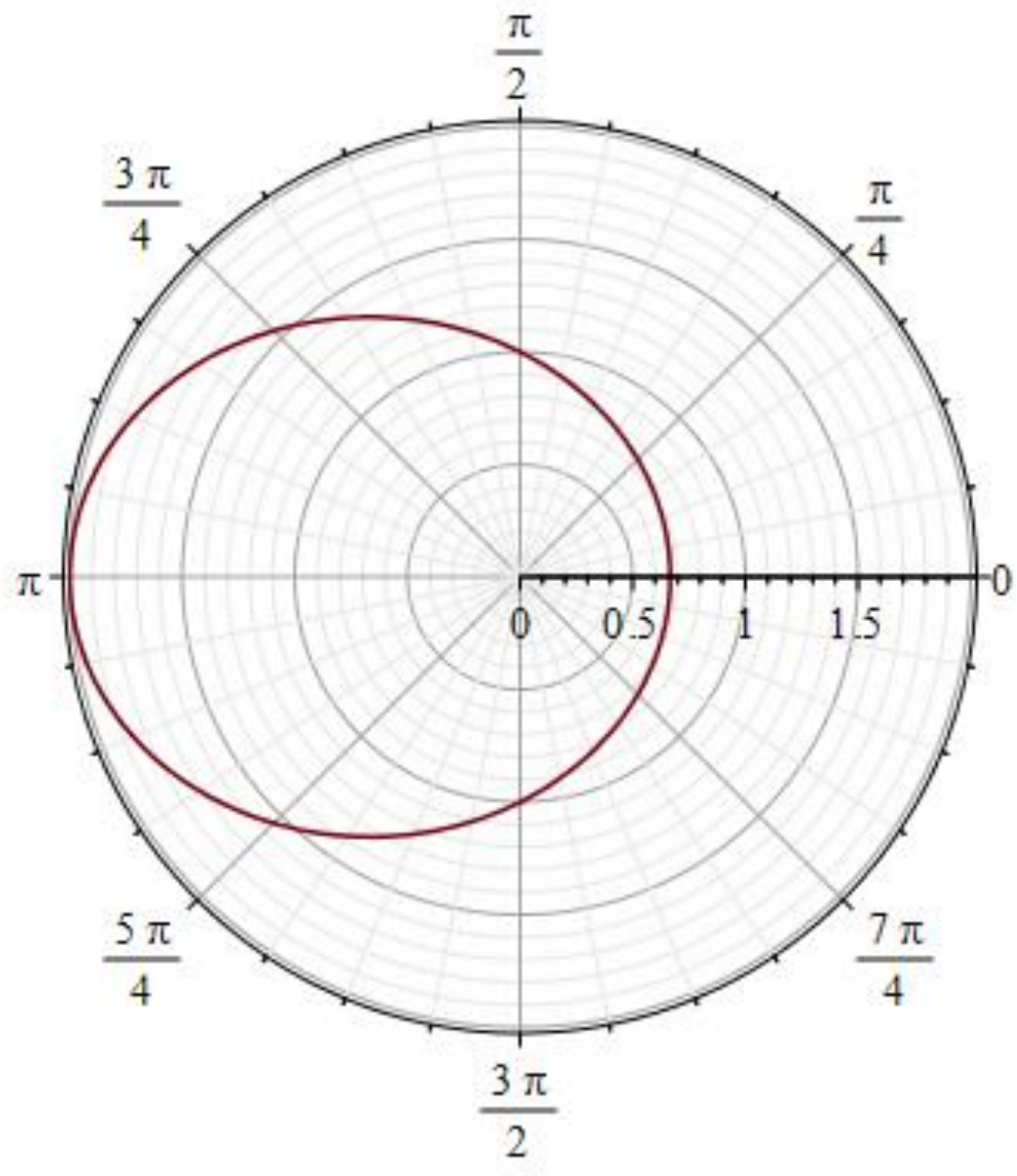

Figure 2: Plot of orbit from numerical solution of orbital equation.

We repeat the process for the orbital equation with an additional term added for the correction due to general relativity.

$$\frac{2}{r^3}\left(\frac{\mathrm{d}r}{\mathrm{d}\phi}\right)^2 - \frac{1}{r^2}\frac{\mathrm{d}^2 r}{\mathrm{d}\phi^2} + \frac{1}{r} = 1 + \frac{3}{64r^2}.$$

(In classical mechanics, it can be shown that adding a potential which is proportional to $1/r^2$ corresponds to an elliptical orbit in a rotating coordinate system). When we create the plot for the orbit after several revolutions, we see that there is a shift in the perihelion which was not present in the Newtonian case.

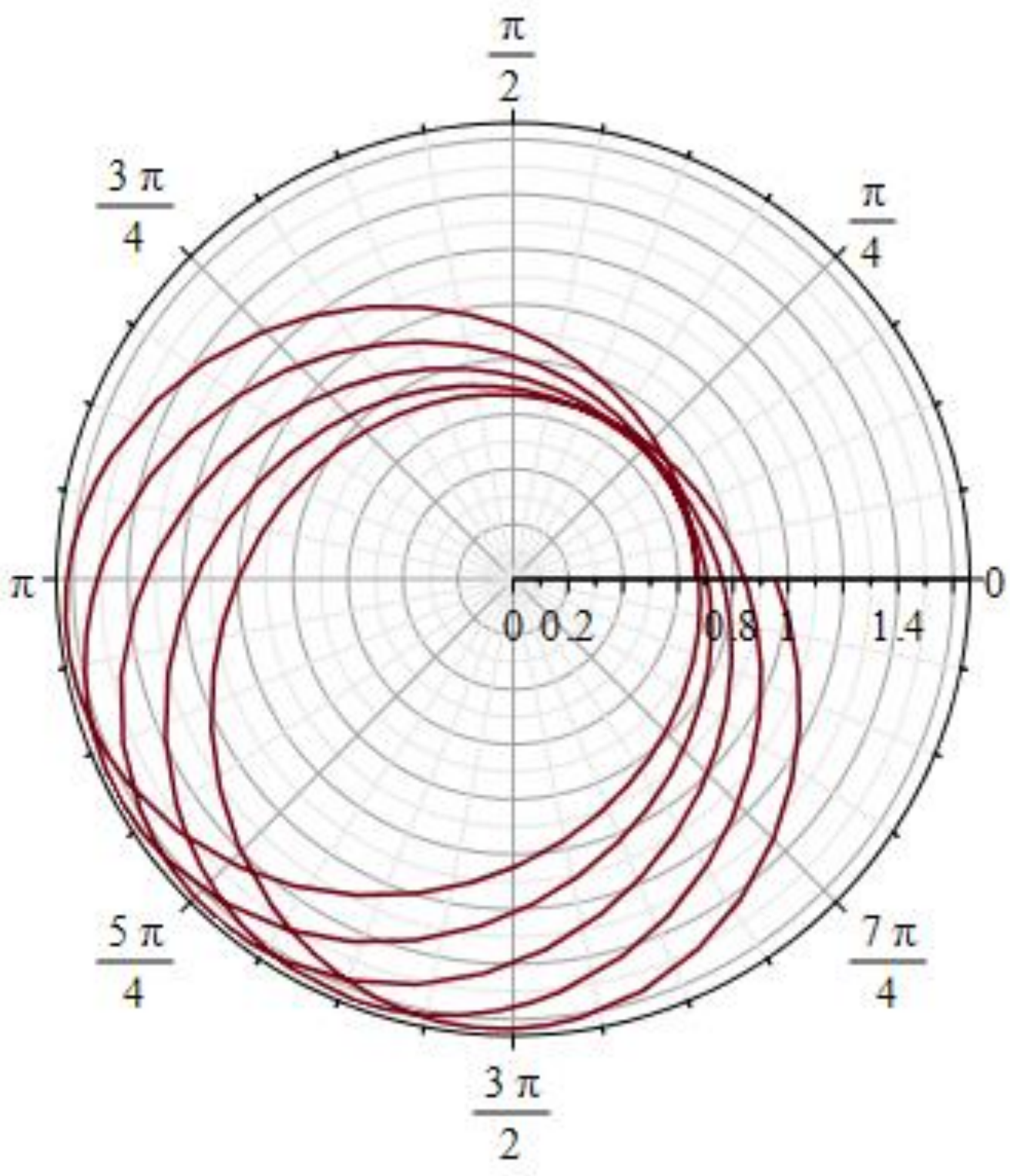

Figure 3: Plot of orbit from solution with relativistic correction.

One can then find the angular positions of the consecutive perihelia by taking the minimum distance of the orbit from the origin for every revolution (ie. wherever $r(\phi)$ is a minimum).

This is all well and good and there is no doubting the physical success of the theory, but one must still answer a number of mathematical questions and problems regarding the Einstein equations, a project which Choquet-Bruhat began to contribute to under the supervision of Einstein. For example, it is not immediately obvious that the Einstein equations are mathematically well posed as set of differential equations. A PDE problem is said to be well posed if there exists a unique solution which depends continuously on the initial data which you feed in. If you recall our earlier discussion of the Cauchy problem, this problem is not well posed for elliptic PDEs and solutions of the problem must be analytic. If the equations are hyperbolic, the problem is well posed for analytic and also for smooth functions. For linear second-order systems on a Lorentzian manifold $(N^{n+1}, g)$ which read as

$$g^{\alpha\beta}\nabla^2_{\alpha\beta} f^I + b_J^{\alpha,I}\nabla_\alpha f^J + c_J^I f^J = h^I,$$

the Cauchy problem with intitial data on a spacelike hypersurface of $N^{n+1}$ is globally well posed, if you assume the right function spaces and the right assumptions on the Cauchy data. The solution to the system is smooth if $(N^{n+1}, g)$ and the Cauchy data are smooth. In particular, the solution at a point only depends on the values for the initial data in the past of this point, as we would expect from the physical assumption of relativistic causality[81]. To see that a solution is global simply means that it exists on all of $\mathbb{R}^{n+1}$, although it might still grow. The Einstein equations for a vacuum form a quasi-linear, second-order system of PDEs for the $(n+1)(n+2)/2$ coefficients of the metric: the PDEs have both hyperbolic and elliptic character (not easy to deal with). As the system is geometric with an obvious geometric interpretation, the equations are invariant under diffeomorphisms of the manifold in question and the isometries of the metric.

As a consequence of diffeomorphism invariance, one can use the contracted Bianchi identities to reduce the number of independent equations. In the classical case where the dimension of the Lorentzian manifold is 4, one can use the Bianchi identities to reduce the number of equations for the metric coefficients from 10 to 6, implying 4 'loose' degrees of freedom, which have to be fixed by the user by choosing a 4-dimensional local coordinate system. The choice of local coordinates is called 'choosing the gauge' in theoretical physics lingo. A natural choice for many years was to make the requirement that the functions which define the local coordinates satisfy a system of wave equations:

$$F^{\alpha} \coloneqq \partial^{\mu}\partial_{\mu}x^{\alpha} = g^{\lambda\mu}\nabla_{\lambda}\partial_{\mu}x^{\alpha} = 0.$$

A computation shows that the Ricci tensor can now be written as:

$$R_{\alpha\beta} = R_{\alpha\beta}^{(h)} + L_{\alpha\beta},$$

where the $R_{\alpha\beta}^{(h)}$ form a system of quasi-linear, quasi-diagonal operators

$$R_{\alpha\beta}^{(h)} = -\frac{1}{2}g^{\lambda\mu}\partial_{\lambda\mu}^{2}g_{\alpha\beta} + P_{\alpha\beta}(g,\partial g).$$

The system of PDEs

$$R_{\alpha\beta}^{(h)} = \rho_{\alpha\beta}$$

are known as the harmonically reduced Einstein equations. The Cauchy problem becomes the problem of constructing a solution along with its first derivatives given initial values on a spacelike submanifold of the Lorentzian manifold. This leads to a theorem that this Cauchy problem (the Cauchy problem for the vacuum Einstein equations in wave gauge) is well posed and has a unique local solution which depends continuously on the initial data (modulo all the relevant theory where we take the functions to be living in the right Sobolev space with the right number of derivatives). Again, causality is implied, since the value at a point only depends on the initial in the past of that point[82].

A further result states that the Cauchy problem for the vacuum Einstein equations with initial data $g_{\alpha\beta}$ on a spacelike hypersurface and $\partial_t g_{\alpha\beta}$ satisfying the wave gauge constraints has a time-local solution which is a Lorentzian metric depending on the initial data in a continuous manner. This solution exists and is locally unique[83]. This is for an Einstein spacetime which is constructed in local coordinates in a neighbourhood of a hypersurface from initial data $g_{\alpha\beta}$ and a Lorentzian metric for which the hypersurface is spacelike. The wave gauge constraints are the following constraints on the initial data when restricted to the hypersurface:

$$G^{\alpha 0} - T^{\alpha 0} = 0.$$

---

[81] Yvonne Choquet-Bruhat, *Introduction to General Relativity, Black Holes, and Cosmology* (Oxford: Oxford University Press, 2015).

[82] *Ibid*.

These constraints must hold for a solution of the harmonically reduced Einstein equations to also be a solution to the original Einstein equations: these are called the wave gauge constraints since they fix the metric in wave gauge.

In order to obtain results on global uniqueness, we must give a geometric formulation of the Cauchy problem which we mentioned earlier, where an initial data set is a triple $(M, g, k)$ with $(M, g)$ a Riemannian $n$-manifold equipped with a metric and a symmetric 2-tensor. A vacuum development of the initial manifold is a Lorentzian $(n+1)$-manifold solution of the vacuum Einstein equations such that $M$ is an embedded submanifold in the larger Lorentzian manifold, the metric on the submanifold is induced by the metric on ambient manifold and the 2-tensor $k$ is the extrinsic curvature of the submanifold as it sits in the ambient space. Results about hyperbolic PDEs then tell us that an initial data set for the Einstein vacuum equations satisfying a pair of geometric constraints admits a vaccum development which is locally unique modulo isometries. The constraints are:

$$-2G_0^0 = \bar{R} - k_{ij}k^{ij} + \left(k_h^h\right)^2 = \rho,$$

$$N^{-1}R_{0j} = \partial_j k_h^h - \bar{\nabla}_h k_j^h = J_i,$$

where $N$ is a gauge variable known as the lapse. The geometric constraints correspond to the Gauss-Codazzi equations from differential geometry, if that helps. A fairly recent conjecture (known as the bounded $L^2$ curvature conjecture) is that sufficient conditions for the existence of a Lorentzian metric solution of the vacuum Einstein equations in a neighbourhood of a Riemannian 3-manifold $M$ equipped with a metric $\bar{g}$ and a 2-tensor $k$ are that

$$\mathrm{Ric}(\bar{g}) \in L^2$$

(ie. that the Ricci tensor as a function which picks up the metric and evaluates it lives in the space $L^2$), that $\bar{\nabla}k$ lives in $L^2$ locally on the 3-manifold and that the volume radius for the manifold equipped with that metric has to be strictly positive. The conjecture was recently proved by its originators[84].

From here, one would like to obtain information about global properties of solutions. An important definition is that of global hyperbolicity for general hyperbolic PDEs. This is interpreted as meaning that a path in the relevant function space between two points whose tangent lives in the cone which determines the domain of dependence of solutions is compact. We earlier defined a Cauchy surface, and it can be proved that the existence of a Cauchy surface implies global hyperbolicity. These definitions can be used to obtain the result that a vacuum Einstein development of an initial data triple satisfies the vacuum Einstein constraints always exists and that it is unique in the space of globally hyperbolic spacetime. The manifold in the initial data is a Cauchy hypersurface. A solution is global if it is a complete Lorentzian manifold. A discussion of the problem of global existence for solutions of the vacuum Einstein equations would require more advanced mathematics which is well beyond the scope of this small book[85].

[83] Yvonne Choquet-Bruhat, 'Théorème d'existence pour certains systèmes d'équations aux dérivées partielles non linéaires', *Acta Math*. Volume 88, 141 – 225 (1952).
[84] Sergiu Klainerman, Igor Rodnianski, and Jérémie Szeftel, 'The resolution of the bounded $L^2$ curvature conjecture in general relativity', *Séminaire Laurent Schwartz.* Volume 202 (1), 1 – 18 (2015).

One of the most important problems in general relativity is that of the positive mass theorem (also known as the positive energy theorem). This problem has a long and interesting history which is in a sense unfinished, since it leads to the general Penrose conjecture. Essentially speaking, there is no well-defined way of talking about the energy density of a gravitational field. The only satisfactory picture of the total energy density field is an asymptotic one, where the density is defined in terms of the behaviour of the field as one goes away from the spacetime out towards spatial infinity. It is not clear that this definition for the total gravitational energy has to be zero or positive. The positive energy theorem makes the claim that a gravitational system with non-negative local matter density has to have non-negative total energy, when that energy is measured asymptotically as one goes to infinity.

The theorem can also be stated in an equivalent way which is known as the positive mass theorem (note that physicists seem to prefer calling it the positive energy theorem, whereas mathematicians tend to call it positive mass). Most importantly, the theorem can be formulated in geometric terms, which is more satisfactory given that general relativity is written in the language of differential geometry. One starts with a three dimensional Riemannian manifold equipped with a metric $(N, g)$. We then take an initial data set consisting of $N$ and a symmetric 2-tensor $h$ (the extrinsic curvature) which satisfy equations of constraint for the local matter and momentum densities:

$$\mu = \frac{1}{2}\left(R - \sum_{a,b} h^{ab}h_{ab} + \left(\sum_{a} h_a^a\right)^2\right),$$

$$J^a = \nabla_b\left(h^{ab} - \left(\sum_{c} h_c^c\right)g^{ab}\right).$$

where $R$ is the scalar curvature. An initial data is asymptotically flat if for some compact set $C$, $N$ with $C$ cut out is made up of a finite number of components $N_i$ such that each component is diffeomorphic to the complement of a compact set in $\mathbb{R}^3$. Under these diffeomorphisms, the metric tensor is then written as

[85] Yvonne Choquet-Bruhat, *Introduction to General Relativity, Black Holes, and Cosmology* (Oxford: Oxford University Press, 2015).

$$\mathrm{d}s^2 = \left(1 + \frac{M}{2r}\right)^4 \left(\sum_i \left(\mathrm{d}x^i\right)^2\right) + \sum_{i,j} p_{ij}\mathrm{d}x^i\mathrm{d}x^j.$$

The quantity $M_i$ is the ADM mass of the end $N_i$. The theorem then says that for an asymptotically flat initial data set, each end has non-negative total mass and if only one of the ends has zero mass, then the initial data set vanishes and the Riemann tensor is zero[86]. This is like saying that an isolated gravitating system with non-negative local mass density has non-negative total mass, when the mass is measured asymptotically at infinity.

The overall strategy from the geometric point of view is to first prove what one might call the 'Riemannian' positive mass theorem, where the second fundamental form of the spacelike hypersurface in the spacetime vanishes to zero. This is equivalent to the hypersurface being 'time-symmetric' or totally geodesic, where a totally geodesic Riemannian submanifold is one such that for every $V$ in the tangent bundle, the geodesic $\gamma_V$ lies entirely in the larger manifold $M$ when viewed in terms of the induced metric associated with the submanifold. In reality, it is only quite special hypersurfaces which have a vanishing second fundamental form even from the mathematical point of view, and one would not expect a 'real' hypersurface embedded in a spacetime to have zero second fundamental form, so it is certainly not sufficient to prove the special case from a physical point of view. However, one can then employ a quasilinear elliptic PDE called the Jang equation to reduce the theorem for general Cauchy data back to the time-symmetric case, which then proves the theorem in general[87].

Although other mathematicians and physicists made progress on proving particular cases of the theorem (Choquet-Bruhat and Marsden showed that it was true if the data are close enough to flat data in a particular smooth norm, for example), it was actually two differential geometers Schoen and Yau who gave the first proof, using the arguments sketched above[88]. Another proof by Witten followed quite soon: this proof was completely different and used identities for Dirac spinors[89]. The idea is to take a Lorentzian 4-manifold (spacetimes in general relativity are always modelled by this type of manifold) which contains an embedded space-like hypersurface. We assume that the manifold satisfies the Einstein equations and that the manifold of the spacelike hypersurface is asymptotic to Euclidean 3-space as one goes to infinity in the 4-manifold. We then take the spinor bundle $S$ of the 4-manifold with its canonical connection and restrict the bundle equipped with this connection to the spacelike hypersurface. If we take a local orthonormal frame field in the usual way, we have a Dirac operator on the restriction of the spinor bundle:

[86] Richard Schoen and Shing-Tung Yau, 'On the Proof of the Positive Mass Conjecture in General Relativity', *Commun. Math. Phys*. 65, 45 – 76 (1979).

[87] Richard Schoen and Shing-Tung Yau, 'Proof of the Positive Mass Theorem. II', *Commun. Math. Phys*. 79, 231 – 260 (1981).

[88] Yvonne Choquet-Bruhat and Jerrold Marsden, 'Solution of the local mass problem in general relativity', *Commun. Math. Phys.* 51, 283 – 296 (1976).

[89] Edward Witten, 'A New Proof of the Positive Energy Theorem', *Commun. Math. Phys*. 80, 381 – 402 (1981).

$$D = \sum_{j=1}^{3} e_j \cdot \nabla_{e_j}.$$

One then writes the square of the Dirac operator

$$D^2 = \nabla^* \nabla + \Im,$$

where the second term only involves the second fundamental form of the spacelike hypersurface.

One can then demonstrate the positivity of $\Im$, which shows that the Dirac operator can be inverted on the right subspaces of the space of sections of $S$. As the hypersurface is asymptotically flat, we just have to take solutions of the equation

$$D\sigma = 0$$

and show that there exist unique solutions for every asymptotically constant value. Once we have such a solution $\sigma$, we take a region of the hypersurface and use our identity for the square of the Dirac operator to get

$$0 = \int_{\Omega} \langle \nabla^* \nabla \sigma, \sigma \rangle + \langle \Im(\sigma), \sigma \rangle = -E_{\Omega}(\sigma) + \int_{\Omega} \|\nabla \sigma\|^2 + \langle \Im(\sigma), \sigma \rangle,$$

where $E_{\Omega}(\sigma)$ is an integral over the boundary of the region. The fact that this integral is strictly positive allows us to show that the total energy is also strictly positive[90]. Witten's proof is appealing, since the mass appears as an integral of a positive quantity over the manifold. It is worth mentioning, however, that although Witten's proof appears simpler than the proof of Schoen and Yau, some work still needs to be done to establish the existence of a Green's function for the Dirac operator on a hypersurface. This is essential, since Witten assumes that it makes mathematical sense to work with a Dirac operator on a spacelike hypersurface[91]. There are now a variety of proofs up to arbitrary dimension and without the assumption of spin manifolds, including a proof of the three-dimensional case with Ricci flow[92]. Choquet-Bruhat herself uploaded a slicker, cleaner spinor proof which is valid for spin manifolds of any dimension and which does not require spacetime spinors[93].

The story does not end here, however, as one still needs to put black holes into the picture (probably the most famous object associated with general relativity). Although the derivation is obviously not this simple, one could intuitively think of this an extension of the positive mass theorem, but now instead of saying the mass is zero or more than zero, we are saying that the mass of the manifold is equal either to the mass of the black holes which it contains, or that it exceeds the mass of those black holes. Also, instead of retreating back to Minkowski spacetime as the flat space, we now go back to the Schwarzchild spacetime for a non-rotating spherically symmetric black hole. Using a physical argument based on various theorems from general relativity and the cosmic censorship hypothesis (which proposes that a naked singularity can never be observed in the Universe outside of a black hole), Penrose conjectured that the following inequality always holds:

$$m \leq \sqrt{\frac{A}{16\pi}},$$

where $m$ is the ADM mass of an asymptotically flat spacelike hypersurface in a spacetime and $A$ is the area of the event horizons of all the black holes in the spacetime.

---

[90] H. Blaine Lawson, Jr. and Maire-Louise Michelsohn, *Spin Geometry* (New Jersey, USA: Prince University Press, 1989).

[91] Thomas Parker and Clifford Henry Taubes, 'On Witten's Proof of the Positive Energy Theorem', *Commun. Math. Phys.* 84, 223 – 238 (1982).

[92] Yu Li, 'Ricci flow on asymptotically Euclidean manifolds' (2018) arXiv:1603.05336.

[93] Yvonne Choquet-Bruhat, 'Positive Gravitational Energy in Arbitrary Dimensions' (2011) arXiv:1107.4283.

If we assume that the second fundamental form vanishes, then the problem can be formulated so that we only need to prove what is known as the Riemannian Penrose inequality. This states that

$$m \leq \sqrt{\frac{A_0}{16\pi}},$$

where the hypersurface has non-negative scalar curvature and $A_0$ is the area of the outermost minimal surface $\Sigma_0$ of the submanifold $(M^3, g)$. Interestingly, the two proofs of the Riemannian Penrose inequality so far rely on geometric flows. A geometric flow is a gradient flow on a manifold with a simple geometric interpretation (a PDE with specific geometric meaning, if you like). This is different from the notion of a global flow which you may have encountered on other manifolds courses, which is a family of diffeomorphisms of a manifold determined by the collection of all the integral curves of some chosen vector field on the manifold. As an example, the coordinate vector field $V = \partial/\partial x$ on $\mathbb{R}^2$ generates a flow which shifts the plane to the right or the left (depending on whether $t$ is positive or negative). This flow is determined by the integral curves, which are straight lines parallel to the $x$-axis[94].

The two most famous examples of geometric flows are probably Ricci flow and mean curvature flow (curve-shortening flow is just mean curvature flow in one dimension). A map $F$ from a family of hypersurfaces $M^n \times [0, T)$ to $(N^{n+1}, \bar{g})$ solves mean curvature flow if

$$\frac{\mathrm{d}}{\mathrm{d}t} F(p,t) = \vec{H}\big(F(p,t)\big),$$

where $\vec{H}$ is the mean curvature vector. This is a quasi-linear, second-order system of parabolic PDEs. The equation is not strictly parabolic, as it is invariant under tangential diffeomorphisms (diffeomorphisms which arbitrarily slide the manifold around, in some sense), meaning that the principal symbol has some zero directions. However, a PDE can only be strictly parabolic if the associated principal symbol is always positive when evaluated for any pair $(x, \xi)$ in $T^*M$. In geometric terms, the flow evolves a family of surfaces in a Riemannian manifold by evolving the normal components of vectors at points on the surface by a speed which is equal to the mean curvature at that point. It is also possible to define a flow called inverse mean curvature flow, where the same is done, but now for the inverse of the mean curvature. This would seem like a difficult flow to use in applications, since it can develop singularities easily and looks to be undefined if the mean curvature at a point is zero.

[94] John M. Lee, *Introduction to Smooth Manifolds* (New York: Springer, 2003).

However, it is possible to use a level sets approach for the Penrose conjecture to show that the surface $\Sigma(t)$ which results from inverse mean curvature flow of a minimal surface after time $t$ can be defined as the level set of a function $u$ on $(M^3, g)$ such that

$$\operatorname{div}\left(\frac{\nabla u}{|\nabla u|}\right) = |\nabla u|.$$

One can define a suitable weak existence theory for this degenerate elliptic PDE and prove existence of solutions. Since the Hawking mass is monotone on the level sets of $u$ and is eventually bounded by the ADM mass as the surfaces evolve, this proves the Riemannian Penrose inequality and provides yet another independent proof of the positive mass theorem[95]. The Huisken-Ilmanen proof actually establishes a slightly weaker result: the Riemannian Penrose inequality for the *largest* black hole in the manifold. One could imagine more than one black hole in a spacetime and, indeed, multiple black holes have been observed in Nature in the form of black hole mergers. A different flow (conformal flow of metrics) was used by Bray to prove the inequality for any number of black holes. Essentially, one takes a flow of metrics which meets a certain number of criteria and which converges back to the Schwarzchild metric, although there are many other technical details and lemmas to deal with. Whereas Huisken and Ilmanen give a new proof of the positive mass theorem, Bray assumes it and applies it at several key points of his argument[96]. It is not entirely clear what, if any, is the physical interpretation of Bray's proof, although such an interpretation has been suggested in the physics literature[97].

The general case where the hypersurface in the spacetime is not totally geodesic is still an open question. One might hope that we could follow the example of Schoen and Yau and use something like the Jang equation to reduce the general case back to the Riemannian Penrose inequality, but there are technical reasons why this is difficult or impossible. The PDE approach to proving the conjecture where one tries to find existence theory for a set of PDEs which imply the Riemannian inequality is perhaps an intractable problem, since one must somehow establish an existence theorem with the correct boundary conditions compatible with the type of horizon being considered (whatever those conditions might be). A counterexample was found to a stronger version of the Penrose conjecture for objects called generalized apparent horizons: if an explicit counterexample can be found for a strong version, this suggests that it might be almost impossible to establish the necessary existence theory for the normal weaker version of the conjecture[98].

---

[95] Gerhard Huisken and Tom Ilmanen, 'The inverse mean curvature flow and the Riemannian Penrose inequality', *Journal of Differential Geometry*. 59 (3): 353 – 437 (2001).
[96] Hubert Bray, 'Proof of the Riemannian Penrose inequality using the Positive Mass Theorem', *Journal of Differential Geometry*. 59 (2): 177 – 267 (2001).
[97] Seiju Ohashi, Tetsuya Shiromizu and Sumio Yamada, 'The Riemannian Penrose inequality and a virtual gravitational collapse', *Physical Review D*. 80(4) (2009).

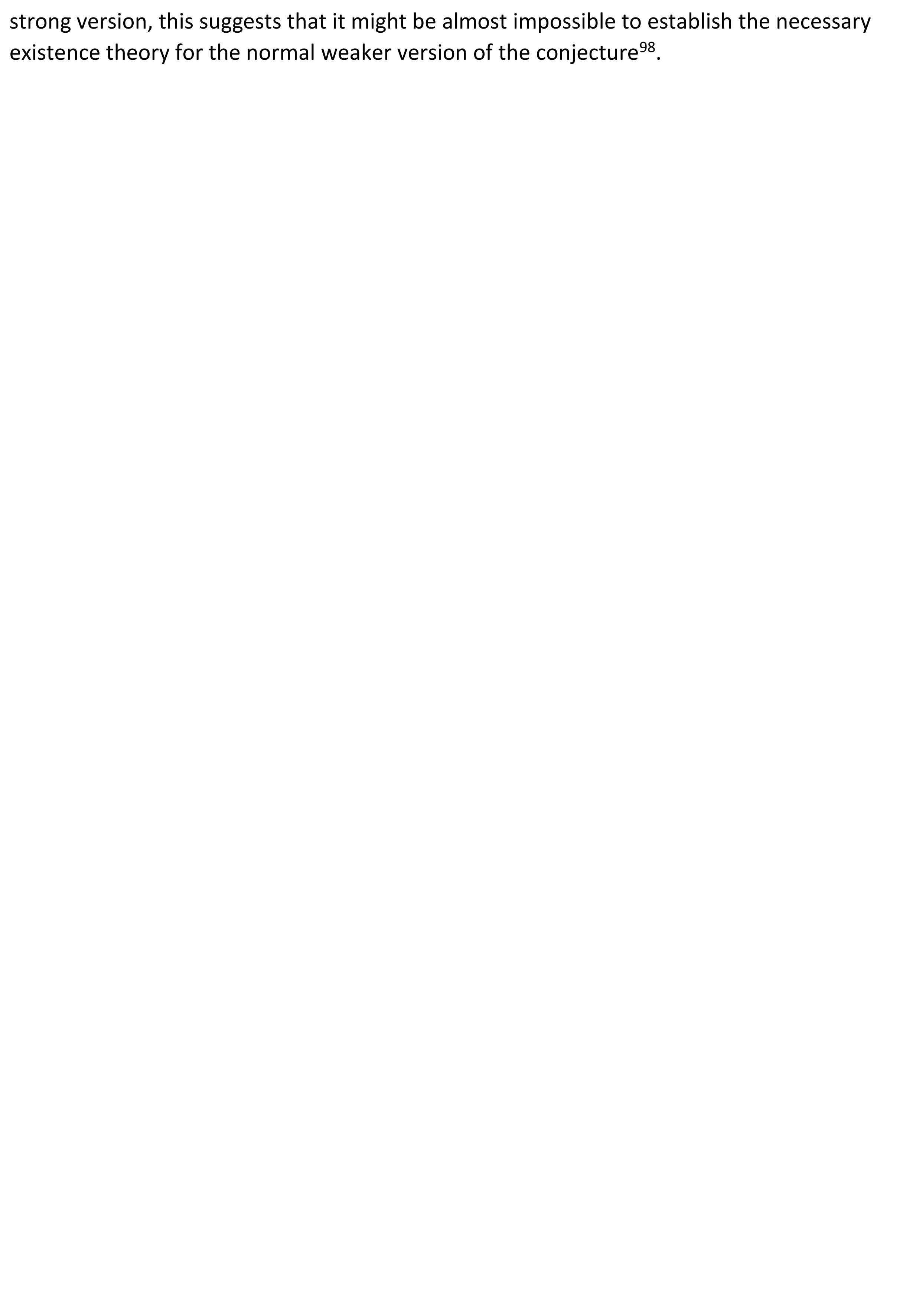

[98] Hubert Bray and Marcus Khuri, 'PDEs which imply the Penrose Conjecture', *Asian Journal of Mathematics*. 15 (2011).

# Olga Oleinik

Oleinik received many prizes and was known for her wide-ranging and extensive achievements in several areas of mathematics, including PDE theory, boundary layer theory and elasticity theory. The first equation of elasticity which many of us are likely to come into contact with is Hooke's law. This can be written in tensor form as

$$u_{ik} = \frac{\delta_{ik}\sigma_{ll}}{9K} + \frac{\left(\sigma_{ik} - \frac{1}{3}\delta_{ik}\sigma_{ll}\right)}{2\mu}.$$

This allows us to write the strain tensor in terms of the stress tensor, where the stress tensor is the rank two tensor which is defined such that the force on a volume can be written as

$$\int F_i \, dV = \int \frac{\partial \sigma_{ik}}{\partial x_k} dV = \oint \sigma_{ik} \, df_k.$$

$\sigma_{ik} df_k$ is the $i$th component of the force on the surface element $d\mathbf{f}$. The $ik$-th component of the stress tensor represents the the $i$th component of the force on a unit area which is perpendicular to the $x_k$-axis. In the case where the body is undergoing uniform compression from all sides, the stress tensor becomes

$$\sigma_{ik} = -p\delta_{ik},$$

such that its non-zero components are equal to the pressure.

The square of the distance between two points after deformation can be written as

$$\begin{aligned} \mathrm{d}l'^2 &= \mathrm{d}l^2 + 2\frac{\partial u_i}{\partial x_k}\mathrm{d}x_i\mathrm{d}x_k + \frac{\partial u_i}{\partial x_k}\frac{\partial u_i}{\partial x_l}\mathrm{d}x_k\mathrm{d}x_l, \\ &= \mathrm{d}l^2 + 2u_{ik}\mathrm{d}x_i\mathrm{d}x_k, \end{aligned}$$

where $\mathrm{d}l$ is the distance between the two points before the deformation and $u_i$ is the displacement vector for the displacement of a point due to a deformation. The symmetric rank two tensor $u_{ik}$ is known as the strain tensor and is defined as

$$u_{ik} = \frac{1}{2}\left(\frac{\partial u_i}{\partial x_k} + \frac{\partial u_k}{\partial x_i} + \frac{\partial u_l}{\partial x_i}\frac{\partial u_l}{\partial x_k}\right).$$

Assuming that the deformations are small, the modulus of compression $K$ can be written as

$$\frac{1}{K} = -\frac{1}{V}\left(\frac{\partial V}{\partial p}\right)_T.$$

The quantity $K^{-1}$ is known as the coefficient of compression. Hooke's law is applicable to most elastic deformations, since most deformations only stop being elastic at a point where Hooke's law still gives a surprisingly good approximation. There are, however, some well-known deformable materials such as rubber for which this is not the case.

Another basic property which you may be familiar with is the Young modulus. If we consider a simple extension or compression of a rod along the $z$-axis, it can be seen by considering the equation which relates the strain tensor to the stress tensor that only the diagonal components of the strain are non-zero. The remaining components are given by

$$u_{xx} = u_{yy} = -\frac{1}{3}\left(\frac{1}{2\mu} - \frac{1}{3K}\right)p, \qquad u_{zz} = \frac{1}{3}\left(\frac{1}{3K} + \frac{1}{\mu}\right)p,$$

where $\mu$ is one of the Lamé coefficients in the second order expansion of the free energy in powers of $u_{ik}$ (an appropriate level of approximation to describe the free energy of a deformed isotropic body):

$$F = F_0 + \frac{1}{2}\lambda u_{ii}^2 + \mu u_{ik}^2.$$

The longitudinal component $u_{zz}$ gives the relative lengthening of the rod as it is stretched in both directions. The reciprocal of the coefficient of $p$ is called Young's modulus $E$, so we can write

$$u_{zz} = \frac{p}{E},$$

where

$$E = \frac{9K\mu}{3K + \mu}.$$

The $u_{xx}$ and $u_{yy}$ components give the compressions of the rod in the transverse direction. The ratio of the transverse compression to the extension in the longitudinal direction (more generally, the transverse strain to the longitudinal strain) is known as the Poisson ratio:

$$u_{xx} = -\sigma u_{zz},$$

where

$$\sigma = \frac{\frac{1}{2}(3K - 2\mu)}{3K + \mu}.$$

This ratio is usually somewhere between $0$ and $0.5$ for real materials. Rubber has a Poisson ratio which is almost $0.5$, corresponding to a Young modulus which is small in comparison with the modulus of compression.

Deformations are often accompanied by a change of temperature, so one needs to take account of the fact that the temperature of a body undergoing deformation is not constant. We can, however, simplify things significantly by assuming that the heat transfer from one part of the body to another by conduction is very slow. In that case, the general form for the equations of motion can be written down using Newton's laws, and one can write down the equations for an isotropic elastic medium:

$$\rho\ddot{\mathbf{u}} = \frac{E}{2(1+\sigma)}\Delta\mathbf{u} + \frac{E}{2(1+\sigma)(1-2\sigma)}\nabla(\operatorname{div}\mathbf{u}).$$

Isotropic means that the deformation behaviour of the material does not depend on its orientation: the stress-strain response of a sample of the material is the same when measured in any direction. In the case of linear elasticity where the deformations are small, the motions are elastic waves. As a simple example, consider a plane wave in an infinite isotropic medium. From the above equations, we obtain equations of motion which you might recognise as one-dimensional wave equations

$$\frac{\partial^2 u_x}{\partial x^2} - \frac{1}{c_l^2}\frac{\partial^2 u_x}{\partial t^2} = 0, \qquad \frac{\partial^2 u_y}{\partial x^2} - \frac{1}{c_l^2}\frac{\partial^2 u_y}{\partial t^2} = 0,$$

where

$$c_l = \sqrt{\frac{E(1-\sigma)}{\rho(1+\sigma)(1-2\sigma)}}, \qquad c_t = \sqrt{\frac{E}{2\rho(1+\sigma)}}.$$

are the velocities of propagation. We can from this that an elastic wave (a sound wave, for example) is basically two waves which propagate independently. The wave whose displacement in the direction of propagation with velocity $c_l$ is the longitudinal wave, and the wave whose displacement is in a plane perpendicular to the displacement vector is the transverse wave with velocity $c_t$. The velocity of the longitudinal wave is always larger than that that of the transverse wave, since we can see that

$$c_l > \sqrt{\frac{4}{3}}c_t.$$

It is also possible to consider the propagation of elastic waves in anisotropic, as opposed to isotropic, media. In this case, the stress-strain response does depend on the orientation of the sample. An example of an anisotropic medium would be a crystal. The equations of motion now become

$$\rho\ddot{u}_i = \lambda_{iklm}\frac{\partial^2 u_m}{\partial x_k \partial x_l},$$

where $\lambda_{iklm}$ is a symmetric tensor denoting the adiabatic moduli of elasticity. If we consider a monochromatic elastic wave in an anisotropic medium, we look for a solution of the equations of motion of the form

$$u_i = u_{0i}e^{i(\mathbf{k}\cdot\mathbf{r}-\omega t)}.$$

Substituting a solution of this form into the equations of motion yields

$$(\rho\omega^2\delta_{im} - \lambda_{iklm}k_k k_l)u_m = 0.$$

This is a set of homogeneous equations for the unknown components of the displacements. These equations only have non-trivial solutions when the determinant of the coefficients vanishes, so we must have an equation for the relation between the wave frequency and the wave vector:

$$|\lambda_{iklm}k_k k_l - \rho\omega^2\delta_{im}| = 0.$$

This relation is known as the dispersion relation. If we find the roots of the equation and plug them back into our equations of motion, we can solve to find the directions of the displacement vector in the waves in question. These directions are known as the directions of polarization. One sees that the directions of polarization of the waves with the same wave vector are mutually perpendicular[99].

Up to now, we have considered elastic media which are homogeneous. However, it is also possible to have inhomogeneous media (Oleinik in fact studied strongly inhomogeneous media). In the non-homogeneous case, there is spatial variation of elastic properties. Elastic media in a geophysical setting are usually inhomogeneous, because the elastic properties of the material depend on depth. Recall that Hooke's law states that each component of the stress tensor is linearly related to every component of the strain tensor. This can be written as:

$$\sigma_{ij} = C_{ijkl}e_{kl},$$

where $C_{ijkl}$ is an elasticity tensor. In the inhomogeneous case, this tensor obviously depends on the coordinates, so we instead have

$$\sigma_{ij} = C_{ijkl}(x_m)e_{kl}.$$

This complicates the equations of motion considerably, so one usually formulates problems which can be solved analytically by assuming that the spatial variation of the elastic properties takes a simple form. You might, for example, assume that the inhomogeneity can be captured by a linear variation of the elastic moduli:

$$C_{ij}(x) = C_{ij}^0(1 + ax).$$

Some examples of inhomogeneous elasticity problems which can be solved are uniaxial tension of an inhomogeneous sheet or torsion of an inhomogeneous cylinder[100].

Oleinik was interested in the mathematical problems which arise in elasticity theory, such as the asymptotic behaviour of the eigenfunctions and eigenvalues of boundary value problems which describe inhomogeneous elastic bodies. In particular, she often considered techniques from a body of theory known as homogenization. For example, she found formulas for the homogenization of boundary value problems for a system of linear elasticity with periodic coefficients that oscillate rapidly (including the homogenization of boundary value problems in perforated domains) and proved estimates for the difference between the displacement vector, stress tensor and energy integral of a strongly inhomogeneous inelastic body, as well as the characteristics of the body described by the homogenized system[101]. Maybe the simplest problem to consider is that of homogenization of second order elliptic operators with periodic coefficients. Start with a periodic matrix $\mathbf{A}(x)$, where $x$ is a tuple with real entries. The elements should also be bounded and satisfy an ellipticity condition

$$a_{ij}(x)\xi_i\xi_j \geq \nu_1|\xi|^2, \qquad \forall x, \xi \in \mathbb{R}^m.$$

---

[99] Lev Landau and Evgeny Lifshitz, *Theory of Elasticity* (Oxford: Butterworth-Heinemann, 1999).

[100] Martin Sadd, *Elasticity: Theory, Applications, and Numerics* (USA: Academic Press, 2009).

[101] Olga Oleinik, A. S. Shamaev and G. A. Yosifian, *Mathematical Problems in Elasticity and Homogenization* (The Netherlands: North Holland, 1991).

If we include a parameter $\varepsilon$, then the matrix

$$\mathbf{A}^\varepsilon(x) = \mathbf{A}(\varepsilon^{-1}x)$$

characterizes a medium which in physical terms has inhomogeneous microstructure. Finding the parameters involves solving some differential equations, but the coefficients of these equations will be given by functions which oscillate rapidly. In physical terms, we are looking to reduce the medium to something which is effectively homogeneous: in formal terms, that amounts to homogenizing the matrix and the differential equations.

A constant positive definite matrix $\mathbf{A}^0$ is the homogenized matrix for $\mathbf{A}$ if for a bounded domain $Q \subset \mathbb{R}^m$ and any $f \in H^{-1}(Q)$, the solutions $u^\varepsilon$ of the Dirichlet problem

$$\operatorname{div}(\mathbf{A}^\varepsilon \nabla u^\varepsilon) = f, \qquad u^\varepsilon \in H_0^1(Q),$$

converge as follows:

$$u^\varepsilon \rightharpoonup u^0, \qquad \mathbf{A}^\varepsilon \nabla u^\varepsilon \rightharpoonup \mathbf{A}^0 \nabla u^0,$$

as $\varepsilon \to 0$. The convergence is in the $H_0^1$ norm for the former and the $L^2$ norm for the latter. $u^0$ is the solution of the associated Dirichlet problem:

$$\operatorname{div}(\mathbf{A}^0 \nabla u^0) = f, \qquad u^0 \in H_0^1(Q).$$

This equation is called the homogenized equation, and $\operatorname{div}(\mathbf{A}^\varepsilon \nabla)$ is a homogenized operator. The vector fields

$$\mathbf{A}^\varepsilon \nabla u^\varepsilon = \mathbf{p}^\varepsilon, \qquad \mathbf{A}^0 \nabla u^0 = \mathbf{p}^0.$$

The simplest example would be if we take $Q$ to be the unit interval. In that case, we have

$$\frac{d}{dx}\left(a(\varepsilon^{-1}x)\frac{du^\varepsilon}{dx}\right) = f, \qquad\qquad u^\varepsilon \in H_0^1([0,1]).$$

After some working, it can be shown that $u^0$ solves a Dirichlet problem

$$\frac{d}{dx}\left(\langle a^{-1}\rangle^{-1}\frac{du^0}{dx}\right) = f, \qquad\qquad u^0 \in H_0^1([0,1]).$$

Hence the homogenized coefficient can be computed as

$$a^0 = \langle a^{-1}\rangle^{-1}.$$

We have given an example of homogenization in terms of the Dirichlet problem, but it could obviously be replaced with a Neumann problem (or any other type of boundary value problem) with a suitable description of homogenization. For other types of boundary problem, one gets a result for the convergence of flows which implies a homogenization rule. The rule says that if $\mathbf{v} \in L^2_{\text{pot}}(\cdot)$ and $\mathbf{A}\mathbf{v} \in L^2_{\text{pot}}(\cdot)$, then

$$\mathbf{A}^0\langle \mathbf{v}\rangle = \langle \mathbf{A}\mathbf{v}\rangle.$$

This rule simplifies computation of the homogenized matrix[102].

Oleinik also worked on boundary layer theory. It was observed by Prandtl that when we are considering a flow around a region with a solid boundary, we need to think about the thin boundary later which exists between the surface and the rest of the fluid in the flow: you could think of it as the layer where the surface influences the flow. For example, if one has an inviscid flow around an object governed by the Euler equations, shear stresses due to viscosity will clearly become non-negligible as one gets close to the boundary, so the flow will become viscous in this region. We will not define the Reynolds number here, but you could think of it as a quantity which is higher when the viscosity is lower. The bulk of the flow region where viscosity can be neglected is called the inviscid outer flow and the very thin region close to the wall where viscosity is included is known as the boundary layer.

The simplest case which we can begin with is to consider boundary layers when the flow is laminar (ie. regular and not turbulent). For simplicity, take a steady two-dimensional and incompressible fluid in the $xy$-plane with no body forces due to gravity or similar, such that fluid flows past a fixed flat plate with zero incidence. The visual picture to have in this case is that of a constant velocity distribution perpendicular to the flat plate where the flow hits the leading edge, with a layer of particles being slowed down due to friction. This layer gets larger as one gets further away from the edge at the point of contact. From this picture, it appears that the profile for the boundary layer could be given as a monotonically increasing function of $x$. We should emphasise that the transition from the outer inviscid flow to the boundary layer flow is continuous, so that there is not in reality a boundary which can be identified, and when we talk about the thickness of the layer, we are talking about a layer which has been defined arbitrarily. A reasonable way of defining the boundary would be to set its location at the point where the inner velocity reaches some percentage of the velocity of the outer flow away from the plate: 99%, say[103].

The continuity equation is

$$\frac{\partial u}{\partial x} + \frac{\partial v}{\partial y} = 0.$$

The Navier-Stokes equations become

$$\rho\left(u\frac{\partial u}{\partial x} + v\frac{\partial u}{\partial y}\right) = -\frac{\partial p}{\partial x} + \mu\left(\frac{\partial^2 u}{\partial x^2} + \frac{\partial^2 u}{\partial y^2}\right),$$

$$\rho\left(u\frac{\partial v}{\partial x} + v\frac{\partial v}{\partial y}\right) = -\frac{\partial p}{\partial y} + \mu\left(\frac{\partial^2 v}{\partial x^2} + \frac{\partial^2 v}{\partial y^2}\right).$$

---

[102] S. M. Kozlov, Olga Oleinik and Vasilii Jikov, *Homogenization of Differential Operators and Integral Functionals* (Berlin: Springer-Verlag, 1994).

[103] Hermann Schlichting and Klaus Gersten, *Boundary-Layer Theory* (Berlin: Springer-Verlag, 2017).

The first thing we might like to do is to make these equations dimensionless. This can be done by making some substitutions.

$$x^* = \frac{x}{L}, y^* = \frac{y}{\delta}, u^* = \frac{u}{U}, v^* = \frac{vL}{U\delta}, p^* = \frac{p}{\rho U^2},$$

where $L$ is the length of the plate, $\delta$ is some small distance used to quantify the thickness of the boundary layer and $U$ is the magnitude of the uniform flow in the free stream. Using these we get

$$\frac{\delta^2}{L^2}\left(u^*\frac{\partial u^*}{\partial x^*} + v^*\frac{\partial u^*}{\partial y^*}\right) = -\frac{\delta^2}{L^2}\frac{\partial p^*}{\partial x^*} + \frac{\mu}{\rho L U}\left(\frac{\delta^2}{L^2}\frac{\partial^2 u^*}{\partial x^{*2}} + \frac{\partial^2 u^*}{\partial y^{*2}}\right),$$

$$\frac{\delta^2}{L^2}\left(u^*\frac{\partial v^*}{\partial x^*} + v^*\frac{\partial v^*}{\partial y^*}\right) = -\frac{\partial p^*}{\partial y^*} + \frac{\mu}{\rho L U}\left(\frac{\delta^2}{L^2}\frac{\partial^2 v^*}{\partial x^{*2}} + \frac{\partial^2 v^*}{\partial y^{*2}}\right).$$

If we are defining the thickness of the boundary layer in terms of a percentage of the outer velocity, we might use this number as a subscript so that there is no ambiguity, so in the case of 99% we could write $\delta_{99}$. Since we are assuming that the thickness of the boundary layer is extremely small in comparison to the length of the plate, we must have

$$\frac{\delta^2}{L^2}\frac{\partial^2 u^*}{\partial x^{*2}} \ll \frac{\partial^2 u^*}{\partial y^{*2}}.$$

The other terms will have a similar order of magnitude as long as the Reynolds number is large, which tells us that

$$\frac{\delta}{L} = O\left(\mathrm{Re}^{-\frac{1}{2}}\right).$$

If we substitute this into the second equation and neglect terms of higher order in the Reynolds numbers, we arrive at

$$-\frac{\partial p^*}{\partial y^*} = 0.$$

After dimensionalising this, we recover

$$\frac{\partial p}{\partial y} = 0.$$

This integrates to

$$p = C(x),$$

since the derivative on the left-hand side is only with respect to $y$, so the constant which we get can still depend on the other spatial coordinate. Hence if we remain within the boundary layer, the Navier-Stokes equations reduce to one equation

$$u\frac{\partial u}{\partial x} + v\frac{\partial v}{\partial y} = -\frac{1}{\rho}\frac{dp}{dx} + \frac{\mu}{\rho}\frac{\partial^2 u}{\partial y^2}.$$

We also have the continuity equation as before:

$$\frac{\partial u}{\partial x} + \frac{\partial v}{\partial y} = 0.$$

By considering our original Navier-Stokes equations and taking the limit as $y$ becomes large, we obtain

$$\rho U \frac{dU}{dx} = -\frac{dp}{dx}.$$

Substitute this in and we get the equation of motion for the fluid in the boundary layer (plus the continuity equation)[104].

$$u\frac{\partial u}{\partial x} + v\frac{\partial v}{\partial y} = U\frac{dU}{dx} + \nu\frac{\partial^2 u}{\partial y^2}.$$

We have already written down a rough estimate of the thickness of the laminar plate boundary layer. Qualitatively, our equation might also suggest that the boundary layer will become narrower as we increase the Reynolds number. A better estimate can be obtained by beginning with the fact that the inertial and frictional forces are in equilibrium in the boundary layer. For a plate of length $x$ we have

$$\frac{\partial u}{\partial x} \propto \frac{U_\infty}{x},$$

where $U_\infty$ is the velocity of the outer flow. This means that the inertial force has order of magnitude $\rho U_\infty^2/x$. The friction force per unit volume is equal to $\partial\tau/\partial y$, where $\tau$ is the shear stress. For a laminar flow, this is equal to $\mu\partial^2 u/\partial y^2$. However, $\partial u/\partial y$ is just the velocity gradient perpendicular to the plate and this is of order $U_\infty/\delta$, so for the friction force we must have the relation

$$\frac{\partial\tau}{\partial y} \sim \frac{\mu U_\infty}{\delta^2}.$$

The inertial and frictional forces are in equilibrium, so we have

$$\frac{\mu U_\infty}{\delta^2} \sim \frac{\rho U_\infty^2}{x}.$$

Re-arrange for the boundary layer thickness:

$$\delta \sim \sqrt{\frac{\mu x}{\rho U_\infty}} = \sqrt{\frac{\nu x}{U_\infty}},$$

using the definition of kinematic viscosity.

[104] The Open University, *Mathematical Methods and Fluid Mechanics IV* (Milton Keynes: The Open University, 2009).

We still need to determine the exact numerical value which is multiplying the right hand side, but this can be done by considering the Blasius solution for the boundary layer velocity profile. Assuming the figure of 99% for the velocity as before, we obtain

$$\delta_{99}(x) = 5\sqrt{\frac{\nu x}{U_\infty}}.$$

The ratio of the thickness to the length of the plate is then

$$\frac{\delta_{99}(x)}{l} = \frac{5}{\sqrt{\mathrm{Re}}}\sqrt{\frac{x}{l}}.$$

Also as before, we can see that the estimate predicts a decreasing boundary layer thickness as the Reynolds number is increased. However, we can now see quantitatively that the boundary layer thickness grows in proportion to $\sqrt{x}$, and the square root function is a monotone increasing function as we predicted. It is reasonably obvious to show that the square root is monotone increasing, since if $\Delta x = x' - x > 0$, then

$$\Delta f(x) = \sqrt{x'} - \sqrt{x} = \frac{\left(\sqrt{x'} - \sqrt{x}\right)\left(\sqrt{x'} + \sqrt{x}\right)}{\sqrt{x'} + \sqrt{x}} = \frac{x' - x}{\sqrt{x'} + \sqrt{x}}.$$

But we just said $x' - x$ is strictly positive, so this must also be strictly positive[105]. In fact, the Reynolds number also grows like $\sqrt{x}$.

We have introduced the boundary layer thickness $\delta$ in an arbitrary way, but there is also a more formal way quantifying the thickness, given by what is called the displacement thickness.

$$\delta_1(x) = \frac{1}{U}\int_0^\infty (U - u)\ dy.$$

This is the distance in the $y$-direction by which the boundary of the plate must be displaced to obtain the same mean flow rate as for the inviscid flow: it quantifies the distance by which the streamlines associated with the outer inviscid flow are displaced by the boundary layer. $\delta_1$ can obtained from the momentum integral equation (an integral version of the boundary layer equations):

$$\frac{d\delta_2}{dx} + (\delta_1 + 2\delta_2)\frac{1}{U}\frac{dU}{dx} = \frac{\tau_0}{\rho U^2},$$

where $\tau_0$ is the shear stress at the surface

$$\tau_0 = \nu\rho\left(\frac{\partial u}{\partial y}\right)_{y=0}.$$

[105] Hermann Schlichting and Klaus Gersten, *Boundary-Layer Theory* (Berlin: Springer-Verlag, 2017).

$\delta_2$ is known as the momentum thickness. This is defined to be the thickness of a layer of the external flow which has a momentum flow rate equal to the reduction in momentum flow caused by the presence of the boundary layer.

$$\delta_2 = \int_0^\infty \frac{u}{U}(U - u)\ dy.$$

In the case of laminar flat plate flow, the momentum integral equation reduces to[106]

$$\frac{d\delta_2}{dx} = \frac{\tau_0}{\rho U^2}.$$

For a plate flow, the ratio of $\delta_1$ to the length $l$ becomes

$$\frac{\delta_1(x)}{l} = \frac{1.721}{\sqrt{\mathrm{Re}}}\sqrt{\frac{x}{l}}.$$

Oleinik also did foundational work on weak solutions to nonlinear PDEs. For example, she showed that the method of vanishing viscosity could be used to construct a weak solution to a nonlinear PDE[107]. A vanishing viscosity solution is a type of viscosity solution. Oleinik proved existence, uniqueness and global stability of vanishing viscosity solutions for one spatial dimension and a scalar conservation law[108]. This work on weak solutions was closely allied to her interest in elasticity and PDEs which arise in the setting of mechanics[109]. In the context of boundary layer theory, Oleinik proved global existence of classical solutions to the two-dimensional steady Prandtl equation modelling boundary layer flows for a class of positive data and suggested as an open problem that it would be interesting to study the local structure of the solution to the Prandtl problem close to the separation point, where the boundary layer separates[110].

---

[106] The Open University, *Mathematical Methods and Fluid Mechanics IV* (Milton Keynes: The Open University, 2009).

[107] Olga Oleinik, 'Construction of a generalized solution of the Cauchy problem for a quasi-linear equation of first order by the introduction of 'vanishing viscosity'', *Uspekhi Matematicheskikh Nauk*, 14 (2(86)), 159 – 164 (1959).

[108] Ola Oleinik, 'Discontinuous solutions of non-linear differential equations', *Amer. Math. Soc. Transl.* 26, 95 – 172 (1963).

[109] Olga Oleinik, 'On the uniqueness of the generalized solution of the Cauchy problem for a nonlinear system of equations occurring in mechanics', *Usp. Mat. Nauk.* 12, 169 – 176 (1957).

[110] Olga Oleinik and V. N. Samokin, *Mathematical Models in Boundary Layer Theory* (London: Chapman & Hall, 1999)

# *Charlotte Fischer*

Fischer is an applied mathematician and computer scientist known for expertise in atomic structure calculations. She developed a method known as multi-configuration Hartree-Fock theory, having studied under both Dirac and Hartree. The Hartree-Fock method fits into a wider body of theory related to computational quantum physics and quantum chemistry. Many of us are familiar with the Schrödinger equation as the basic equation of motion in quantum mechanics:

$$i\hbar \frac{d}{dt}\Psi = \hat{H}\Psi.$$

This is a differential equation for the dynamics of the quantum mechanical wave function. The problem is that this equation will only be solvable analytically for situations which are extremely simple, and these solutions (if they exist) were likely found a long time ago. An example of a simplification would be to take a free particle with zero potential energy. The equation then becomes

$$i\hbar \frac{\partial}{\partial t}\Psi = -\frac{\hbar^2}{2m}\frac{\partial^2 \Psi}{\partial x^2}.$$

Another example of a simple case would be the hydrogen atom. However, even in this case, the equation which one gets is not really trivial. The Schrödinger equation for the hydrogen atom in spherical polar coordinates is as follows:

$$-\frac{\hbar^2}{2\mu}\left(\frac{1}{r^2}\frac{\partial}{\partial r}\left(r^2\frac{\partial}{\partial r}\right) + \frac{1}{r^2 \sin\theta}\frac{\partial}{\partial \theta}\left(\sin\theta \frac{\partial}{\partial\theta}\right) + \frac{1}{r^2\sin^2\theta}\frac{\partial^2}{\partial^2\phi}\right)\psi(r,\theta,\phi) + U(r,\theta,\phi)\psi(r,\theta,\phi) = E\psi(r,\theta,\phi).$$

(You might recognise the presence of the Laplacian in spherical polar coordinates). The wave function can be decomposed into a product of a radial function with a spherical harmonic,

$$\psi_{nlm}(r,\theta,\phi) = R_{nl}(r)Y_l^m(\theta,\phi),$$

where the spherical harmonic satisfies an eigenvalue equation for the magnitude of the orbital angular momentum

$$\hat{L}^2 Y_l^m(\theta,\phi) = \hbar^2 l(l+1)Y_l^m(\theta,\phi)$$

and the radial function is as follows

$$R_{nl}(r) = -\left(\frac{(n-l-1)!}{2n\big((n+l)!\big)^3}\right)^{\frac{1}{2}}\left(\frac{2Z}{na_0}\right)^{l+\frac{3}{2}} r^l e^{-\frac{Zr}{na_0}} L_{n+l}^{2l+1}\left(\frac{2Zr}{na_0}\right).$$

Mathematically speaking, this is all just a separation of variables from elementary PDE theory. The radial function gives the radial part of the energy eigenfunction and the spherical harmonic gives the angular part. Note that the spherical harmonics are also eigenfunctions of the $z$-component of the orbital angular momentum, with eigenvalues of $m\hbar$ in this case.

The number $n$ is the principal quantum number which determines the energy levels of the hydrogen atom: it can take positive integer values. However, this also restricts the values which $l$ can take, since if we consider a trial solution to the radial equation

$$u(\rho) = C(\rho)\rho^{l+1}e^{-\beta\rho},$$

it can be shown that this is a solution to the radial part of the wave equation, provided that $C(\rho)$ can be expressed as a power series

$$C(\rho) = \sum_{k=0}^{\infty} c_k \rho^k,$$

where the coefficients satisfy a recurrence relation

$$c_{k+1} = \frac{2\beta(k+l+1)-2}{(k+1)(k+2l+2)} c_k.$$

This relation only generates a coefficient $c_{k_{max}+1}$ which vanishes to zero if

$$\frac{1}{\beta} = k_{max} + l + 1.$$

It is necessary to assume the existence of a coefficient of this form: otherwise, the function will diverge at infinity. If you assume that the series terminates, then the function does not diverge because the exponential function decays faster than polynomial growth.

$k_{max}, l \geq 0$, which implies that $1/\beta$ is an integer greater than or equal to 1. This implies that the bound-state energies of the hydrogen atom are given by

$$E = -\frac{E_R}{n^2},$$

where $E_R$ is the Rydberg energy

$$E_R = \left(\frac{e^2}{4\pi\varepsilon_0}\right)^2 \frac{\mu}{2\hbar^2}.$$

This is exactly what is predicted by the Bohr model of the atom, so Bohr's model was correct (on this occasion, at least). However, if we replace $1/\beta$ with $n$, then we have

$$n = k_{max} + l + 1.$$

$k_{max} \geq 0$, which implies that $n \geq l + 1$ and therefore that

$$l \leq n - 1.$$

This is the origin of the selection rule which says that the quantum number for the orbital angular momentum can only take values of $0, 1, \ldots, n-1$. The energy levels of the hydrogen atom only depend on the principal quantum number, but there could be different values of $l$ for the same value of $n$, so the energy levels are degenerate. Similar considerations tell us that the magnetic quantum number can only take the following values: $0, \mp 1, \ldots, \mp l$. If you have studied chemistry, you may be more familiar with the spectroscopic notation, where $l = 0,1,2,3,4, \ldots$, is replaced by the letters $s, p, d, f, g, \ldots$, and the principal quantum number is placed in front of the letter, but the meaning is the same[111]. In the case where $n = 1$, $l = 0$ and $m = 0$, we have a very simple wave function:

$$\psi_{100} = \frac{1}{\sqrt{\pi}}\left(\frac{Z}{a_0}\right)^{3/2} e^{-\frac{Zr}{a_0}},$$

where $a_0$ is the Bohr radius.

The other molecules which we study in Nature are of much higher complexity than the hydrogen atom, and it becomes necessary to introduce some kind of approximation. The Hamiltonian operator for a free particle with no potential energy is simple, but consider instead the Hamiltonian for a set of interacting nuclei and electrons.

$$\widehat{H} = -\sum_{I=1}^{M}\frac{\nabla_I^2}{2M_I} - \sum_{i=1}^{N}\frac{\nabla_i^2}{2} - \sum_{i,I=1}^{N,M}\frac{Z_I}{|\mathbf{R}_I - \mathbf{r}_i|} + \sum_{i>j=1}^{N}\frac{1}{|\mathbf{r}_i - \mathbf{r}_j|} + \sum_{I>J=1}^{M}\frac{Z_I Z_J}{|\mathbf{R}_I - \mathbf{R}_J|},$$

where $\mathbf{r}_j$ is the position vector of electron $j$ and $Z_I$ and $\mathbf{R}_I$ are the charge and position vector for nucleus $I$, respectively. The first term represents the kinetic energy of the nuclei, the second term is the kinetic energy of the electrons, the third term is the Coulomb interaction energy between the nuclei and the electrons, the fourth term is the Coulomb interaction energy which exists solely between the electrons, and the final term is the Coulomb interaction between the nuclei. A historically important way of simplifying this problem was to use the Born-Oppenheimer approximation, whereby one assumes that since the nuclei are moving with negligible velocity compared to the electrons, we can neglect in the first approximation the kinetic energy associated with the nuclear motion and write a total wave function as follows for $N$ electrons and $M$ nuclei:

$$\Psi(\mathbf{r}, \mathbf{R}) = \psi_{elec}(\mathbf{r}, \mathbf{R})\phi_{nuc}(\mathbf{R}).$$

This allows us to separate the Schrödinger equation into two separate but coupled Schrödinger equations: one for the wave function of the electrons with nuclear coordinates as parameters $\psi_{elec}(\mathbf{r}, \mathbf{R})$ and another for the nuclear wave function $\phi_{nuc}(\mathbf{R})$ for the nuclei as they move in the potential energy whose values should have been obtained by solving the eigenvalue problem for the electronic Schrödinger equation.

$$H_{elec}(\mathbf{r}, \mathbf{R})\psi_{elec}(\mathbf{r}, \mathbf{R}) = E_{elec}(\mathbf{R})\psi_{elec}(\mathbf{r}, \mathbf{R}),$$

$$H_{nuc}(\mathbf{R})\phi_{nuc}(\mathbf{R}) = \left(T_{nuc} + E_{elec}(\mathbf{R})\right)\phi_{nuc}(\mathbf{R}) = E_{tot}\phi_{nuc}(\mathbf{R}).$$

The potential energy surface for the molecule is determined as $E_{elec}(\mathbf{R})$ by solving the electronic Schrödinger equation when the nuclei are fixed in different positions. The surface will be a $3N-6$ dimensional hypersurface, as the electron energies are a function of all the nuclear coordinates.

[111] John Bolton and Stuart Freake (ed.), *Quantum Mechanics of Matter* (Milton Keynes: The Open University, 2009).

Another principle of approximation which is used to study changes in the vibrational and electronic energy states is a molecule is the Franck-Condon principle. In semi-classical terms, this states that electronic transitions occur so rapidly in comparison to vibrational motions that the inter-nuclear distance remains unchanged as a result of the transition. This is because electrons are much lighter than nuclei, and it implies that the most probable transitions $\nu_i \to \nu_f$ are vertical. This vertical line intersects any number of vibrational levels $\nu_f$ in the upper electronic state, hence transitions to many vibrational states of the excited states will occur with transition probabilities which are proportional to the Franck-Condon factors, which are themselves proportional to the overlap integral of the wave functions of the initial and final vibrational states. A vibrational progression is observed as a result, whose shape is determined by the horizontal positions of the two electronic potential energy curves relative to each other. It follows that the most likely transitions are transitions to excited vibrational states which have wave functions with large amplitude when evaluated at the inter-nuclear position vector $\mathbf{R}_C$.

Transitions between electronically excited states of a molecule in the gas phase could occur due to various causes. For example, in non-radiative decay, collisions cause excess energy to be distributed to other modes of vibration and rotation of the molecule. In fluorescence, there is de-excitation of electrons between two states of the same multiplicity. The difference in wave number between the two energy bands is the Stokes shift. In radiationless intersystem crossing, there is a transition between two electronic states with different spin multiplicity as a result of the spin-orbit interaction. Phosphorescence is emission due to the radiative transition between two states of different spin multiplicity. Phosphorescence reduces the observed fluorescence and is often delayed relative to the radiation causing excitation, persisting for several seconds after the exciting source is removed. More details regarding molecular energy levels can be looked up in a textbook of molecular physics. See Landau and Lifschitz for a more mathematical treatment using symmetry groups and Young diagrams[112].

In fact, the Schrödinger equation is not solvable even in the case of the helium atom. We obtain the equation

$$-\frac{\hbar^2}{2m_e}(\nabla_1^2+\nabla_2^2)\psi(\mathbf{r}_1,\mathbf{r}_2)-\frac{2e^2}{4\pi\varepsilon_0}\left(\frac{1}{r_1}+\frac{1}{r_2}\right)\psi(\mathbf{r}_1,\mathbf{r}_2)+\frac{e^2}{4\pi\varepsilon_0 r_{12}}\psi(\mathbf{r}_1,\mathbf{r}_2)=E\psi(\mathbf{r}_1,\mathbf{r}_2).$$

We cannot separate the variables here because of the term corresponding to electrostatic repulsion between the electrons. The two main approximation methods which one has are perturbation theory or a variational treatment. In the variational approach, one chooses some type of trial function with free variational parameters and then determines the parameters using the variational principle. One can then change the trial function and the free parameters until the lowest possible energy is reached. More generally, one can consider a linear combination of several wave functions in the variational approach:

$$\chi = c_1\chi_1 + c_2\chi_2.$$

[112] Lev Landau and Evgeny Lifschitz, *Quantum Mechanics: Non-Relativistic Theory* (Oxford: Butterworth-Heinemann, 1981).

In the case of the helium atom with two electrons, the two trial functions could be formed from the products of the ground state and first excited state for the atom:

$$\chi_1(\mathbf{r}_1, \mathbf{r}_2) = \phi_{1s}(\mathbf{r}_1)\phi_{1s}(\mathbf{r}_2),$$

$$\chi_2(\mathbf{r}_1, \mathbf{r}_2) = \phi_{2s}(\mathbf{r}_1)\phi_{2s}(\mathbf{r}_2).$$

If we include the effect of spin and take account of the Pauli exclusion principle, you may know that it is a corollary of this principle that particles with half-integer spin (fermions) are described by antisymmetric wave functions, whereas particles with zero or integer spin (bosons) are described by symmetric wave functions. In the case of a single-valued, many-particle wave function, this wave function must be antisymmetric.

It follows that with spin, the total electronic wave function for the helium atom can be written as

$$\Psi(1,2) = \frac{1}{\sqrt{2}}\big(1s\alpha(1)1s\beta(2) - 1s\beta(1)1s\alpha(2)\big) = -\Psi(2,1).$$

This can also be written in the form of a Slater determinant:

$$\Psi(1,2) = \frac{1}{\sqrt{2}}\begin{vmatrix} 1s\alpha(1) & 1s\beta(1) \\ 1s\alpha(2) & 1s\beta(2) \end{vmatrix}.$$

Note that we have assumed that the spin and spatial parts of the wave function are separable: this is not true for the case of more than two electrons.

The goal of Hartree-Fock theory is to find the best possible single-determinant approximation to the true wave function for the electronic ground state of an $N$-electron atom.

$$\Psi_{\mathrm{HF}}(\mathbf{x}_1, \mathbf{x}_2, \dots, \mathbf{x}_N) = \frac{1}{\sqrt{N!}}\det|\psi_1(\mathbf{x}_1) \dots \psi_N(\mathbf{x}_N)|.$$

Each of the $\psi_i$ is an orthonormal spin orbital for a single electron which can be written as a separable product of a spatial orbital and a spin function. Although we gave an example for the helium atom, the general approach is to construct the ground state electronic wave function as an antisymmetric product of one-electron orbitals and then evaluate the energy of the atom using that wave function. Next, we use the variational principle to find the energy minimizer, where the equation of constraint is the fact that the wave function has to be normalized to unity. This yields the Hartree-Fock equation, a non-linear one-electron equation for the optimal spin-orbitals. Solving these equations yields the spin-orbitals, the electron energies, and the central potential $U(\mathbf{r})$ which approximates the electron-electron interaction. This is all done in a self-consistent way, hence the terminology of 'self-consistent fields'.

The general Hartree-Fock equation for a one-electron spin-orbital can be written as

$$\hat{F}_i \phi_i = \varepsilon_i \phi_i,$$

where $\varepsilon_i$ is the one-electron energy, $\phi_i$ is the one-electron function, and $\hat{F}_i$ is the Fock operator

$$\hat{F}_i = \hat{f}_i + \sum_{j=1}^{N} \left(2\hat{J}_i - \hat{K}_i\right)$$

such that

$$\hat{f}_i = -\frac{\nabla_i^2}{2} - \frac{Z}{r_i},$$

$$\hat{J}_j(\mathbf{r}_1)\phi_i(\mathbf{r}_1) = \phi_i(\mathbf{r}_1) \int \phi_j^*(\mathbf{r}_2) \frac{1}{r_{12}} \phi_j(\mathbf{r}_2)\, d\mathbf{r}_2,$$

$$\hat{K}_j(\mathbf{r}_1)\phi_i(\mathbf{r}_1) = \phi_j(\mathbf{r}_1) \int \phi_j^*(\mathbf{r}_2) \frac{1}{r_{12}} \phi_i(\mathbf{r}_2)\, d\mathbf{r}_2.$$

As we said, the equation is non-linear, since the Fock operator also depends on the spin-orbitals. The last operator is non-local and purely quantum mechanical in nature, since there is not a function which can give the action of that operator at a point specified by a position vector. This operator is known as the exchange operator in analogy with the exchange interaction energy which is used to define the energy of an electron in a molecular orbital

$$E = 2\sum_{j=1}^{N} I_j + \sum_{i=1}^{N}\sum_{j=1}^{N} \left(2J_{ij} - K_{ij}\right),$$

where

$$I_j = \int \phi_j^*(\mathbf{r}_j) \left(-\frac{\nabla_j^2}{2} - \frac{Z}{r_j}\right) \phi_j(\mathbf{r}_j)\, d\mathbf{r}_j,$$

$$J_{ij} = \iint \phi_i^*(\mathbf{r}_1)\, \phi_j^*(\mathbf{r}_2) \frac{1}{r_{12}} \phi_i(\mathbf{r}_1)\phi_j(\mathbf{r}_2)\, d\mathbf{r}_1\, d\mathbf{r}_2,$$

$$K_{ij} = \iint \phi_i^*(\mathbf{r}_1)\, \phi_j^*(\mathbf{r}_2) \frac{1}{r_{12}} \phi_i(\mathbf{r}_2)\phi_j(\mathbf{r}_1)\, d\mathbf{r}_1\, d\mathbf{r}_2.$$

$I_j$ is the energy of the electron in a spin-orbital in the nuclear field, $J_{ij}$ represents the classical Coulomb interaction between the smoothed-out charge distributions of two electrons, and $K_{ij}$ is the exchange integral representing quantum mechanical reduction of energy of interaction between electrons with parallel spins in different orbitals.

The Hartree-Fock equation can be written down for an $N$-electron closed-shell atom or an atom with one valence electron, but things start to become complicated, so we will go back to the helium atom for now. We begin by writing out the independent particle Hamiltonian which approximates a two-electron Hamiltonian:

$$H_0(\mathbf{r}_1, \mathbf{r}_2) = h(\mathbf{r}_1) + h(\mathbf{r}_1),$$

where

$$h(\mathbf{r}) = h_0(\mathbf{r}) + U(r) = -\frac{1}{2}\nabla^2 + V(r).$$

The independent particle model is also known as the one-electron approximation. This Hamiltonian describes two particles moving independently of each other in a potential of the form

$$V(r) = -\frac{Z}{r} + U(r),$$

where $U(r)$ approximates the effect of a standard Coulomb repulsion. The full Hamiltonian is then

$$H = H_0 + V(\mathbf{r}_1, \mathbf{r}_2),$$

where

$$V(\mathbf{r}_1, \mathbf{r}_2) = \frac{1}{r_{12}} - U(r_1) - U(r_2).$$

If we take an orbital $\psi_\alpha(\mathbf{r})$ which solves the one-electron Schrödinger equation

$$h(\mathbf{r})\psi_\alpha(\mathbf{r}) = \epsilon_\alpha \psi_\alpha(\mathbf{r}),$$

the two-electron problem is solved by a product wave function:

$$H_0 \Psi_{ab}(\mathbf{r}_1, \mathbf{r}_2) = E_{ab} \Psi_{ab}(\mathbf{r}_1, \mathbf{r}_2),$$

where $E_{ab} = \epsilon_\alpha + \epsilon_b$ to zero order.

The lowest energy eigenstate of $H_0$ is the product of the lowest energy spin-orbitals. Recall that we are looking for an antisymmetric product of one-electron orbitals which approximates the atomic ground state of a 2-electron wave function. (Actually, the result we quote is for a two-electron ion, but we can perform some modifications to get the two-electron atom). This is

$$\Psi_{1s,1s}(\mathbf{r}_1, \mathbf{r}_2) = \frac{1}{\sqrt{2}}(\chi_{1/2}(1)\chi_{-1/2}(2) - \chi_{-1/2}(1)\chi_{1/2}(2).$$

This product is an eigenstate of the magnitude of the spin angular momentum $S^2$ and the $z$-component of the spin angular momentum $S_z$, where the eigenvalue is trivial in both cases. To get the Hartree-Fock equation for the helium atom, we approximate the electron interaction to obtain the two-electron energy eigenvalue and note that the energy can be written in terms of the radial part of the wave function and integrated to obtain

$$E_{1s,1s} = \int_0^\infty \left( \left( \frac{dP_{1s}}{dr} \right)^2 - 2\frac{Z}{r} P_{1s}^2(r) + v_0(1s,r) P_{1s}^2(r) \right) dr,$$

where

$$P_{1s} = 2r\zeta^{3/2} e^{-\zeta r}.$$

There is also an equation of constraint due to the fact that the wave function needs to be normalized:

$$N_{1s} = \int_0^\infty P_{1s}(r)^2 \, dr = 1.$$

One then simply uses the variational principle with a Lagrange multiplier

$$\delta\left(E_{1s,1s} - \lambda N_{1s}\right) = 0$$

and requires that the expression solve the boundary conditions and vanish for arbitrary variations. This leads to the following Hartree-Fock equation which can be solved iteratively:

$$-\frac{1}{2}\frac{d^2 P_{1s}}{dr^2} - \frac{Z}{r} P_{1s}(r) + v_0(1s,r) P_{1s}(r) = \epsilon_{1s} P_{1s}(r).$$

To get a more accurate solution, you would need to drop the independent particle assumption and model the fact that the two electrons are influencing each other's motion and becoming correlated[113]. It is worth mentioning that the iterative method may not always converge, but other schemes can be chosen if necessary[114].

There is also a body of literature on numerical solution of the Hartree-Fock equations which we will not go into here. One can also include relativistic effects in the theory by taking the one-electron Hamiltonian to be a Dirac Hamiltonian, but this starts to lead us into many of the difficulties which are associated with QFT. The crucial fact about Hartree-Fock theory is that there is only one determinant. This differs from the multi-configurational Hartree-Fock theory pioneered by Fisher. Whilst the neglect of electron correlation which we mentioned in the previous paragraph does not stop us from obtaining a result for the helium atom, the deviation from experimental data can be large when we do not include these effects. In multi-configurational Hartree-Fock, one expands the $N$-electron wave function as a linear combination of Slater determinants. (Another method which is widely used is many-body perturbation theory).

To be more specific, we now use CSFs (configuration state functions), where a CSF is a linear combination of Slater determinants which is adapted to the wave function of the system which is being studied. A CSF is distinct from a configuration, which is just an assignment of electrons to orbitals familiar even from elementary chemistry. In the original method known as 'configuration interaction', one attempts to naively represent the linear combination of Slater determinants for the Hartree-Fock wave function as the lowest CSF and the first excited CSF of single excitations, double excitations, and so on.

[113] Walter Johnson, *Atomic Structure Theory: Lectures on Atomic Physics* (Berlin: Springer, 2007).
[114] Charlotte Froese Fischer, 'General Hartree-Fock program', *Computer Physics Communications,* 43(3), 355 – 365 (1987).

$$\Psi(\{\mathbf{x}\}) = C_{\mathrm{HF}}\Phi_{\mathrm{HF}}(\{\mathbf{x}\}) + \sum_{i}^{n_{\mathrm{occ}}}\sum_{a}^{n_{\mathrm{vir}}} C_{i\to a}\Phi_{i\to a}(\{\mathbf{x}\}) + \sum_{i,j}^{n_{\mathrm{occ}}}\sum_{a,b}^{n_{\mathrm{vir}}} C_{ij\to ab}\Phi_{ij\to ab}(\{\mathbf{x}\}) + \cdots,$$

where $\{\mathbf{x}\}$ is a vector in a vector space with spatial and spin part. This configuration interaction wave function is used to compute the matrix elements of the configuration interaction Hamiltonian

$$H_{IJ} = \int \Phi_I^*(\{\mathbf{x}\})\hat{H}\,\Phi_J(\{\mathbf{x}\}).$$

One can then solve the configuration interaction determinant:

$$|\mathbf{H} - E\mathbf{I}| = 0.$$

However, the accuracy of this method is limited by having a finite amount of basis functions and the number of CSFs increases exponentially as one increases the number of electrons. In multi-configuration Hartree-Fock, the number of CSFs is limited by restricting them only to exchanges of valence orbitals. The Hartree-Fock computations are then performed with the configuration interaction wave functions for the electronic structures. The MCHF wave function is represented by the linear combination of the CSFs generated from the space which is spanned by the ground state configuration and configurations of excitations from occupied to virtual orbitals (known as the active space).

$$\Psi = \sum_{I}^{n_{\mathrm{active}}\ \mathrm{CSF}} C_I\Phi_I.$$

We use this wave function and solve the multi-configuration Hartree-Fock instead of the standard Hartree-Fock equation:

$$\sum_{s}\hat{F}_{rs}\,\phi_s = \sum_{s}\epsilon_{sr}\,\phi_s.$$

The corresponding Fock operator is

$$\hat{F}_{rs} = \left(\sum_{I,J}^{n_{\mathrm{active}}\ \mathrm{CSF}} C_I^*\,C_J A_{sr}^{IJ}\right)\hat{h} + \sum_{t,u}^{n_{\mathrm{orb}}}\left(\sum_{I,J}^{n_{\mathrm{active}}\ \mathrm{CSF}} C_I^*\,C_J B_{su,rt}^{IJ}\right)V_{tu}^{ee},$$

where we have the corresponding version of the exchange operator, and so on.

$$A_{sr}^{IJ} = \frac{1}{N!}\int \Phi_I^*(\{\mathbf{x}\})\Phi_j^{s\to r}(\{\mathbf{x}\})\,d^3\{\mathbf{x}\},$$

$$B_{su,rt}^{IJ} = \frac{1}{N!}\int \Phi_I^*(\{\mathbf{x}\})\Phi_j^{sr\to ut}(\{\mathbf{x}\})\,d^3\{\mathbf{x}\},$$

$$V_{tu}^{ee} = \int \phi_t^*(\mathbf{r}_2)\frac{1}{r_{12}}\phi_u\,(\mathbf{r}_2)d^3\mathbf{r}_2.$$

It is possible to solve this equation using standard manipulations of the Hamiltonian matrix elements and so on, but we will leave out the details[115].

On the more theoretical side, Fisher made the important discovery of the negative calcium ion. Fisher, Lagowski and Vosko derived this result using multi-configuration Hartree-Fock and density functional theory to show that the negative calcium and scandium ions have stable ground states, where the extra electron sits in the $4p$ orbital (recall that this means that the electron is in a $n = 4, l = 1$ state). A computation with multi-configuration Hartree-Fock predicts an electron affinity of $0.045$ electron volts in very good agreement with experiment[116]. Recall that a positive ion results when you try to put energy into an atom beyond its state of maximum energy. The limiting energy in the case of the hydrogen atom is $13.6$ electron volts. The atom is ionized and an electron is ejected, leaving behind something with overall positive charge. Although we tend to think of energies as being 'discrete' in quantum mechanics, the ion/electron system can take values within a continuous range of energies. As you may know from chemistry, ionic bonding can occur between ions of opposite charge, forming a solid crystal via electrostatic attraction (salt is the classic example). Fisher has also applied her methods to completely calculate the lower energy levels of the sodium-like up to the argon-like sequences[117].

[115] Takao Tsuneda, *Density Functional Theory in Quantum Chemistry* (Japan: Springer, 2014).
[116] Charlotte Froese Fischer, Jolanta Lagowski and S.H. Vosko, 'Ground states of $Ca^-$ and $Sc^-$ from two theoretical points of view', *Physical Review Letters*, 59, 2263 – 2266 (1987).
[117] Charlotte Froese Fischer, Georgio Tachiev and Andrei Irimia, 'Relativistic energy levels, lifetimes, and transition probabilities of the sodium-like to argon-like sequences', *Atomic Data and Nuclear Data Tables*, 92(5), 607 – 812 (2006).

# *Karen Uhlenbeck*

Uhlenbeck is known for her fundamental work in differential geometry, gauge theory and theory of PDEs. In particular, she contributed to the founding of modern geometric analysis (a field where one studies problems or properties of objects in differential geometry using tools and ideas from the theory of differential equations). Uhlenbeck is one of the most celebrated living female mathematicians and recently became the first woman to win the prestigious Abel Prize. Besides her academic achievements, she has also done work to widen minority and female participation in mathematics and co-founded the Women and Mathematics Program at the Institute for Advanced Study.

The physical origin of a gauge theory is that it is a theory whose Langrangian is invariant under a particular Lie group of local gauge transformations known as the gauge group. (You should think of the word 'gauge' as implying change of scale or 'locality'). The simplest possible example of a gauge theory is classical electrodynamics, where the Lagrangian is invariant under local $U(1)$-transformations. $U(1)$ is a Lie group known as the circle group, the most basic example of a unitary group of degree $n$. A gauge theory requires a gauge field, which is the electromagnetic four-potential for the case of electromagnetism. Things become inherently geometric at this point, since the gauge field is a connection. One also has a definition for curvature which is constructed from the connection. In order to be able to discuss electrodynamics as a gauge theory, we must express the usual objects of electromagnetism as differential forms on a spacetime (Minkowski spacetime in our case). General electromagnetic fields can be combined into two 2-forms $F$ and $G$.

In coordinate-free notation, the Maxwell equations can be written as:

$$\mathrm{d}F = 0,$$

$$\mathrm{d}G = j,$$

where $j$ is the current density 3-form. We will write one of the computations out in detail for once. Start with the first two Maxwell equations:

$$\operatorname{curl} E + \partial_t B = 0,$$

$$\operatorname{div} B = 0.$$

These can be written in index notation as:

$$\epsilon^{ijk}\partial_i E_j + \partial_0 B^k = 0,$$

$$\partial_j B^j = 0.$$

We use the definition for the 2-form $F$.

$$\mathrm{d}F = \mathrm{d}\left(\frac{1}{2}\epsilon_{ijk}B^i \mathrm{d}x^j \wedge \mathrm{d}x^k - E_i \mathrm{d}x^0 \wedge \mathrm{d}x^j\right),$$

$$= \frac{1}{2}\epsilon_{ijk}\partial_\mu B^i \mathrm{d}x^\mu \wedge \mathrm{d}x^j \wedge \mathrm{d}x^k - \partial_\mu E_i \mathrm{d}x^\mu \wedge \mathrm{d}x^0 \wedge \mathrm{d}x^j,$$
$$= \frac{1}{2}\epsilon_{ijk}\partial_0 B^i \mathrm{d}x^0 \wedge \mathrm{d}x^j \wedge \mathrm{d}x^k + \frac{1}{2}\epsilon_{ijk}\partial_l B^i \mathrm{d}x^l \wedge \mathrm{d}x^j \wedge \mathrm{d}x^k - \partial_j E_i \mathrm{d}x^j \wedge \mathrm{d}x^0 \wedge \mathrm{d}x^i,$$
$$= \left(\frac{1}{2}\epsilon_{ijk}\partial_0 B^i - \partial_k E_j\right) B^i \mathrm{d}x^0 \wedge \mathrm{d}x^j \wedge \mathrm{d}x^k + \frac{1}{2}\epsilon_{ijk}\partial_l B^i \mathrm{d}x^l \wedge \mathrm{d}x^j \wedge \mathrm{d}x^k = 0.$$

Each of these has to equal 0 in order for the right hand side to equal 0, which implies that

$$\frac{1}{2}\epsilon_{ijk}\partial_0 B^i - \partial_{[k}E_{j]} = 0,$$
$$\frac{1}{2}\epsilon_{i[jk}\partial_{l]} B^j = 0.$$

Contract with the Levi-Civita permutation symbol:

$$\frac{1}{2}\epsilon^{jkl}\epsilon_{ijk}\partial_0 B^i - \epsilon^{jkl}\partial_{[k}E_{j]} = 0,$$
$$\frac{1}{2}\epsilon^{jkl}\epsilon_{i[jk}\partial_{l]} B^i = 0.$$

This becomes

$$\delta_i^l \partial_0 B^i - \epsilon^{jkl}\partial_k E_j = 0,$$
$$\delta_i^l \partial_l B^i = 0,$$

which is

$$\epsilon^{jkl}\partial_j E_k + \partial_0 B^l = 0,$$
$$\partial_i B^i = 0.$$

These are the first two Maxwell equations. In fact, one immediately knows that the first equation is true, since $F$ is a closed form and the exterior derivative of a closed form is zero: however, showing that you can recover the original Maxwell equations takes some computation. The second equation is left as an exercise (slightly harder, as you will have to expand the 3-form $j$ in components to match with the second set of Maxwell equations).

Given that electromagnetism is described successfully by an abelian gauge group, the question arises as to whether it might be possible to describe other interactions in Nature using non-abelian gauge theories. This turns out to be the case, and such a theory is known as a Yang-Mills theory. In the Yang-Mills case, the gauge group must be compact and the definition of curvature is modified by adding an extra term:

$$F = \mathrm{d}A + A \wedge A.$$

The Maxwell equations become the classical Yang-Mills equations, which can be derived by varying a Yang-Mills action of the type

$$\frac{1}{g^2}\int \mathrm{tr}\,(F \wedge * F) + \int (D\Phi)^+ \wedge * D\Phi - m^2 \int \Phi^+ \wedge * \Phi,$$

where $g$ is the coupling constant, $m$ is the scalar mass, and tr is a quadratic form on the Lie algebra of the compact gauge group[118].

Gauge theories are the most important examples of quantum field theories. The Standard Model is a gauge theory, for example, but the Lie group in this case is a direct product of several other Lie groups due to the complexity of the theory and the number of gauge bosons which must be accounted for. However, it is fairly common in mathematics for an area to be initially inspired or created from discoveries in theoretical physics which then develops into its own branch to the point where the physical origin becomes quite distant. This is the case for gauge theory in mathematics, which is a branch of geometry, although physical references might still appear every now and then. Essentially, one is usually studying the theory of connections when one studies gauge theory from the mathematical point of view.

Start with a smooth principal bundle $\pi: P \to X$ with a Lie group $G$. A principal $G$-bundle over a smooth manifold $X$ is a manifold with a smooth right $G$-action and an orbit space $P/G = X$. In general, if a Lie group acts on a manifold, one can define an equivalence relation $p \sim q$ if there exists a group element $g$ such that $g \cdot p = q$. The equivalence classes for this relation are the orbits of the Lie group in the manifold, where the set of orbits is denoted by the quotient $M/G$: if we equip the resulting space with the quotient toplogy, it is known as the orbit space of the action. It can be proved that if the Lie group acts smoothly, freely and properly, the orbit space is a topological manifold of dimension $\dim M - \dim \mathrm{G}$ equipped with a unique smooth structure such that the quotient map from $M$ to $M/G$ is a smooth submersion. If the action is locally equivalent to the trivial action on $U \times G$, then we say that the bundle $P$ has a structure group $G$. There are several ways of defining a connection on this type of bundle: for example, as a Lie algebra-valued 1-form called the connection 1-form (ie. as a 1-form on the bundle which takes its values in the Lie algebra $\mathfrak{g}$ of $G$). As hinted at above, the curvature of the connection is a $\mathfrak{g}$-valued 2-form $\Omega$ given by the equation

$$\Omega = \mathrm{d}\omega + [\omega, \omega].$$

Another possibility is to define it as a $G$-invariant field of tangent $n$-planes $\tau$ on $P$ such that the linear map $\pi_*: \tau_p \to T_{\pi p}(X)$ is an isomorphism for all $p \in P$. Probably the most common way of defining connections is to work with a vector bundle $E$ and make it into a principal bundle. Given a vector bundle, a connection on the frame bundle is defined via a map called the covariant derivative:

$$\nabla: \Omega_X^0(E) \to \Omega_X^1(E).$$

$\Omega_X^0(E)$ denotes a $p$-form in $E$. You might notice the proliferation of linear maps in differential geometry, hence why many people say the subject is just advanced linear algebra. Bundle automorphisms are gauge transformations whose local representations are matrix-valued functions. You will recognise the idea of a connection if you have studied differential geometry. In that setting, one studies connections on the tangent bundle of a manifold (usually the Levi-Civita connection). The novelty in general connection theory (known as Yang-Mills theory) is that one studies connections on auxiliary bundles. This is slightly different to what a physicist means when they are talking about Yang-Mills theory, but they *are* linked. Historically, it was realised that the correct mathematical formulation for electromagnetism came by considering a connection on a principal $U(1)$-bundle, and this gave the impetus for physical Yang-Mills theories, which then provided fertile ground for new geometric ideas[119].

[118] M. Göckeler and T. Schücker, *Differential Geometry, Gauge Theories, and Gravity* (Cambridge: Cambridge University Press, 1997).

An important theorem in this area is Uhlenbeck's theorem for the existence of local Coulomb gauges: this result is significant because one needs to know that it is possible to construct Coulomb gauges for different connections. In formal terminology, the result is as follows. There exist positive non-zero constants $M$ and $\varepsilon_1$ such that any connection $A$ on the trivial bundle over the closed unit ball $\bar{B}^4$ in $\mathbb{R}^4$ with the $L^2$ norm of the curvature $F_A$ less than $\varepsilon_1$ is gauge equivalent to a different connection $A'$ over the open ball $B^4$ such that

$$\mathrm{d}^* A' = 0,$$

$$\lim_{|x| \to 1} A'_r = 0,$$

$$\|A'\|_{L_1^2} \leq M \|F_{A'}\|_{L^2}.$$

Also, for correctly chosen constants, the connection $A'$ is determined in a unique way by these properties, modulo a local gauge transformation

$$A' \to u_0 A' u_0^{-1},$$

for some constant $u_0$ drawn from the group $U(n)$. Essentially, if we assume that the curvature is not too large and work in a nice space $L^2$, one can always satisfy the Coulomb gauge condition, and this condition provides a matrix-valued connection with small entries. The proof of this theorem is quite involved and uses a clever continuity method. One has to show that there exists a positive constant $\zeta$ such that if $B'_t$ is a one-parameter family of rescaled connections for times in the unit interval on the trivial bundle over $S^4$ with a bound on the rescaled curvature $\|F'_{B_t}\| < \zeta$, and $B'_0$ the product connection, then for every time $t$ there exists a nice gauge transformation such that an estimate is satisfied on the $L^2$ norm of $B_t$. The continuity idea comes into play here, as one takes the set of times in the unit interval for which such a gauge transformation exists and shows that the set is both open and closed, and that it contains the points $t = 0$ and $t = 1$, which proves the proposition[120].

To explain the terminology of a classical Coulomb gauge, one needs to explain the idea of 'fixing the gauge'. In physical terms, this corresponds to the fact that you might have loose degrees of freedom in the system which you are trying to describe, and you need a choice of gauge and a reasonable choice of ancillary equation so that these can be picked up and dealt with. The Lorenz gauge is another classical example from electromagnetism (not named after Lorentz, although it is Lorentz-invariant). The Coulomb gauge condition is the classical condition for the case where the base manifold is a domain in $\mathbb{R}^3$ where one insists that the vector potential be normalised by

$$\mathrm{d}^* A = 0.$$

[119] Simon Donaldson and Peter Kronheimer, *The Geometry of Four-Manifolds* (Oxford: Oxford University Press, 1990).

[120] Karen Uhlenbeck, 'Connections with $L^p$ bounds on curvature', *Communications in Mathematical Physics,* 83, 31 – 42 (1982).

Geometrically speaking, a choice of gauge is a choice for a bundle trivialisation. If we start with an open set in $X$ over which $E$ is trivial and fix a diffeomorphism $\tau$ called a local trivialisation which maps from the open set to the fibres of $E$, the open set on which this trivialisation can be performed is called a trivialising neighbourhood. In the case where we are working with a vector bundle, the fibres of $E$ are copies of $\mathbb{C}^n$ or $\mathbb{R}^n$ (depending on whether it is a real or a complex vector bundle, respectively). In particular, if we choose to represent the connection via connection matrices, these matrices will depend on the choice of bundle trivialisation (or choice of gauge). If we take an automorphism $u$ of the bundle $\mathbb{C}^n$, then $u\tau$ is also a local trivialisation and one has a transformation formula for the connection which can be viewed as the effect of a change of local trivialisation on the connection matrices for a fixed connection:

$$A \to uAu^{-1} - (\mathrm{d}u)u^{-1}.$$

In the physics literature, you might see this as

$$A \to \gamma A\gamma^{-1} - (\mathrm{d}\gamma)\gamma^{-1},$$

or alternatively

$$A \to \gamma A\gamma^{-1} + \gamma(\mathrm{d}\gamma^{-1}),$$

where $\gamma$ is a local gauge transformation which maps from an open subset of $\mathbb{R}^4$ into the structure group.

One can also show that the connection 1-form transforming in this way implies that the curvature 2-form transforms as

$$F \to \gamma F\gamma^{-1}.$$

This can be proved with another nitty-gritty computation, which I will reproduce here (but try it yourself before looking, and see if there is someone you can discuss it with if you get stuck). Start with the definition for the curvature:

$$F = \mathrm{d}A + A \wedge A.$$

Perform the transformation and use the graded Leibnitz rule:

$$\begin{aligned} F &\to \mathrm{d}(\gamma A\gamma^{-1} + \gamma\mathrm{d}\gamma^{-1}) + (\gamma A\gamma^{-1} + \gamma\mathrm{d}\gamma^{-1}) \wedge (\gamma A\gamma^{-1} + \gamma\mathrm{d}\gamma^{-1}), \\ &= \mathrm{d}\gamma \wedge A\gamma^{-1} + \gamma\mathrm{d}A\gamma^{-1} - \gamma A\mathrm{d}\gamma^{-1} + \mathrm{d}\gamma \wedge \mathrm{d}\gamma^{-1} + \gamma\mathrm{d}^2\gamma^{-1} + \gamma A\gamma^{-1} \wedge \gamma A\gamma^{-1} + \gamma\mathrm{d}\gamma^{-1} \\ &\qquad \wedge \gamma A\gamma^{-1} + \gamma A\gamma^{-1} \wedge \gamma\mathrm{d}\gamma^{-1} + \gamma\mathrm{d}\gamma^{-1} \wedge \gamma\mathrm{d}\gamma^{-1}. \end{aligned}$$

The fifth term vanishes by definition of the ordinary exterior derivative, and there are some cancellations of $\gamma$ wedged with $\gamma^{-1}$ (in the sixth term, for example).

$$F' = \gamma(\mathrm{d}A + A \wedge A)\gamma^{-1} + (\mathrm{d}\gamma) \wedge A\gamma^{-1} - \gamma A\mathrm{d}\gamma^{-1} + \mathrm{d}\gamma \wedge \mathrm{d}\gamma^{-1} + \gamma\mathrm{d}\gamma^{-1} \wedge \gamma A\gamma^{-1}$$
$$+ \gamma A\mathrm{d}\gamma^{-1} + \gamma\mathrm{d}\gamma^{-1} \wedge \gamma\mathrm{d}\gamma^{-1}.$$

We know that

$$\mathrm{d}(\gamma\gamma^{-1}) = \mathrm{d}(1) = 0.$$

However, we also have

$$\mathrm{d}(\gamma\gamma^{-1}) = (\mathrm{d}\gamma)\gamma^{-1} + \gamma(\mathrm{d}\gamma^{-1}).$$

This implies that

$$(\mathrm{d}\gamma)\gamma^{-1} = -\gamma(\mathrm{d}\gamma^{-1}),$$

which implies that

$$\mathrm{d}\gamma^{-1} = -\gamma^{-1}(\mathrm{d}\gamma)\gamma^{-1}.$$

Plug this back in and use the definition of the curvature in the first term, we have

$$F' = \gamma F\gamma^{-1} + \mathrm{d}\gamma \wedge A\gamma^{-1} + \gamma A \wedge \gamma^{-1}(\mathrm{d}\gamma)\gamma^{-1} - \mathrm{d}\gamma \wedge \gamma^{-1}(\mathrm{d}\gamma)\gamma^{-1} - \mathrm{d}\gamma(\gamma^{-1}) \wedge \gamma A\gamma^{-1}$$
$$- \gamma A \wedge \gamma^{-1}(\mathrm{d}\gamma)\gamma^{-1} + (\mathrm{d}\gamma)\gamma^{-1} \wedge (\mathrm{d}\gamma)\gamma^{-1},$$
$$= \gamma F\gamma^{-1},$$

after some cancellations. This is the required result and it shows that the curvature 2-form transforms in the adjoint representation under a gauge transformation, since the above is the transformation for an adjoint field.

This ability to find Coulomb gauges is used in a crucial way in Uhlenbeck's paper on removal of singularities in four-dimensional Yang-Mills connections[121]. Uhlenbeck's result was extended to higher dimensions by Tao and Tian[122]. The theorem states that if we take a unitary connection $A$ over a ball punctured at the origin $B^4\backslash\{0\}$ which is ASD with respect to a smooth metric on $B^4$ and if

$$\int_{B^4\backslash\{0\}} |F(A)|^2 < \infty,$$

then there exists a smooth ASD connection over $B^4$ which is gauge equivalent to $A$ over the punctured ball. ASD stands for anti-self-dual. In physical terms, a self-dual gauge field is one such that

$$F = * F,$$

and an ASD gauge field is such that

$$F = - * F,$$

[121] Karen Uhlenbeck, 'Removable singularities in Yang-Mills fields', *Communications in Mathematical Physics*, 83, 11 – 29 (1982).
[122] Terence Tao and Gang Tian, 'A singularity removal theorem for Yang-Mills fields in higher dimensions', *J. Amer. Math. Soc.,* 17, 557 – 593 (2004).

where the $*$ denotes the Hodge duality operator. If we equip a manifold with an orientation and a metric, this is a linear map

$$*: \Lambda^p V \to \Lambda^{n-p} V,$$

given in local coordinates by

$$*\left(e^{i_1} \wedge \dots \wedge e^{i_p}\right) \coloneqq \varepsilon_{i_1 \dots i_n} \eta^{i_1 i_1} \dots \eta^{i_p i_p} e^{i_{p+1}} \wedge \dots \wedge e^{i_n},$$

where the $e^i$ form an orientated orthonormal basis of the dual space of $V$. On a four-manifold, the Hodge operator maps 2-forms to 2-forms. The Maxwell equations can be rewritten using the Hodge star as:

$$\mathrm{d}F = 0,$$

$$\delta F = * j,$$

where $\delta$ is the co-derivative, another differential operator on forms defined by

$$\delta\varphi \coloneqq (-1)^{np+n+1+s} * \mathrm{d} * \varphi.$$

In more geometric terms, the self-dual and ASD forms are defined to be the $1$ and $-1$ eigenspaces of the Hodge dual. In terms of sections, they are sections of bundles $\Lambda^+$ and $\Lambda^-$, respectively, where $\Lambda^2$ is the direct sum of $\Lambda^+$ and $\Lambda^-$, and the 2-form $\alpha \wedge \alpha$ is equal to

$$\alpha \wedge \alpha = \mp|\alpha|^2 \, \mathrm{d}\mu.$$

The subspaces of self-dual and ASD forms only depend on the conformal class of the Riemannian metric. Take a Riemannian four-manifold equipped with an orientation: as we said, 2-forms on that manifold can be split into parts which are self-dual and ASD and this splitting extends to the curvature 2-form of a connection on a bundle $E$ over the manifold, so we can write $F_A$ as the direct sum of $F_A^+$ and $F_A^-$. A connection is ASD if $F_A^+ = 0$ and self-dual if $F_A^- = 0$: generally speaking, ASD connections are more convenient[123]. It can be proved fairly easily that every ASD metric connection solves the Yang-Mills equation. The Yang-Mills equation is

$$D^* F = 0.$$

An identity shows that this is the same as

$$D * F = 0.$$

This is a non-linear second-order PDE for $A$ which can also be obtained by varying the Yang-Mills action and assuming that there are no matter fields. One can show as an exercise that the square of the Hodge map is always $1$ or $-1$. If we choose an orthonormal frame and use local coordinates, a $p$-form can be written as

$$\varphi = \frac{1}{p!} \varphi_{i_1 \dots i_p} e^{i_1} \wedge \dots \wedge e^{i_p},$$

[123] Simon Donaldson and Peter Kronheimer, *The Geometry of Four-Manifolds* (Oxford: Oxford University Press, 1990).

hence

$$* \ \varphi = \frac{1}{p!} \frac{1}{(n-p)!} \varphi_{i_1 \dots i_p} \varepsilon_{j_1 \dots j_p j_{p+1}} |\det g_{kl}|^{\frac{1}{2}} \eta^{i_1 j_1} \dots \eta^{i_p j_p} e^{i_{p+1}} \wedge \dots \wedge e^{i_n}.$$

Use the definition for the transformation of the basis by the Hodge dual and then collect everything together. To contract with some of the metric components at the end, you will have to reverse the indices of the remaining Levi-Civita symbol which will introduce a factor of $p!$. You should end with

$$** \ \varphi = (-1)^{p(n-1)+s} \varphi,$$

which proves the claim. In particular, on Euclidean space $\mathbb{R}^4$, $** = 1$. If we let $F$ be ASD (or self-dual, in fact), the Yang-Mills equation becomes

$$D ** F = 0,$$

which is

$$DF = 0.$$

This is the geometric Bianchi identity, so one can find solutions to the Yang-Mills equation by solving either of the first-order equations for $A$:

$$F = \mp * F,$$

meaning that the equations are satisfied by self-dual or ASD connections[124].

A self-dual or ASD connection can also be called an instanton or an anti-instanton. Generally speaking, you could say that an instanton is a special type of solution to the field equations for gauge theories. More specifically, an instanton is a classical solution in a classical Euclidean field theory with finite non-zero action. The name is due to the fact that they happen for an 'instant' (a point) of Euclidean time and so they are important as critical points in the path-integral formulation of a theory which uses Euclidean signature. The instanton solution of the Euclidean Yang-Mills equation leads to an $SU(2)$ bundle over $S^4$, but it can also be proved that any finite action solution of the Euclidean Yang-Mills equations leads to a fibre bundle over $S^4$ (another theorem of Uhlenbeck)[125].

If we consider the best-known example, we take a pure Yang-Mills theory with symmetry group $SU(2)$ in $\mathbb{R}^4$ with Euclidean signature. The equations of motion (ie. the Yang-Mills equations) are $D * F = 0$ and $DF = 0$. When we introduce the ASD condition, these equations reduce to ODEs for the gauge potential *A*. If we make an ansatz for the solution to the ASD equation which only differs from a pure gauge by a function of $r$ at infinite radius, we guarantee that our solutions have *finite, non-zero* action. Our ansatz for the gauge field is such that it becomes a pure gauge as $r$ tends to infinity and the associated field strength disappears, meaning that the action is finite. If we choose the appropriate gauge transformation and use this to evaluate the field strength, we then obtain an equation whose full solution is the 'instanton potential' which is regular on all of $\mathbb{R}^4$. The action is equal to $-8\pi^2/g^2$, so it is obviously finite. One can also consider related bundle structures such as the Dirac monopole: the bundle structure of the Dirac monopole is extremely similar to that of the Yang-Mills instanton bundle. Following this logic we can might attempt to construct a 'gravitational instanton': we look for a metric with Euclidean signature described locally by an orthonormal frame and solve the Einstein equations without matter. We end up with a solution such that $g$ and $f$ tend to 1 as $r$ tends to infinity. This is similar to the Yang-Mills case, where the Yang-Mills instanton potential becomes a pure gauge as $r$ tends to infinity. Gravity ends up being 'analogous' to other gauge theories, rather than directly comparable, since gravity does not quantize well and gauge theories are quantum field theories.

---

[124] Jürgen Jost, *Riemannian Geometry and Geometric Analysis* (Heidelberg: Springer, 2008).

[125] Karen Uhlenbeck, 'Removable singularities in Yang-Mills fields', *Bulletin of the American Mathematical Society*, 1, 579 – 581 (1979).

In more variational terms, the Yang-Mills equations are the Euler-Lagrange equations for the Yang-Mills functional $\|F_A\|^2$ on the space of connections, considered as the square of the $L^2$ norm of the curvature:

$$\|F_A\|^2 = \int_X |F_A|^2 \ \mathrm{d}\mu = \int_X |F_A^-|^2 \ \mathrm{d}\mu + \int_X |F_A^+|^2 \ \mathrm{d}\mu.$$

We also have the ASD condition that a connection is only ASD if

$$\mathrm{Tr}(F_A^2) = |F_A|^2 \ \mathrm{d}\mu.$$

This is just one basic point that can be mentioned about the relation between the Yang-Mills density and the 4-form $\mathrm{Tr}(F^2)$. If we take a bundle $E$ over a compact, orientated Riemannian four-manifold which can be equipped with an ASD connection, then $\kappa(E) \geq 0$, and if $\kappa(E) = 0$ any ASD connection is flat. $\kappa(E)$ is a characteristic number whose definition depends on the group associated with the bundle, but in the case of $SU(r)$ bundles, it is just $c_2(E)$, the second Chern class. For an $SU(2)$ connection, this is

$$c_2 = \frac{1}{8\pi^2}\mathrm{Tr}(F \wedge F).$$

Interestingly, in the case we mentioned previously of $SU(2)$ bundles over $S^4$, the second Chern class takes an explicit topological meaning, since the bundle can be trivialised over the upper and lower hemispheres of $S^4$ up to isomorphisms by the homotopy class of the transition function (similar to the way that the two-sphere can be trivialised over its two hemispheres explicitly with suitable transition functions). $SU(2)$ itself can be identified with $S^3$, so the transition function gives a self-map of the 3-sphere and the second Chern class of the bundle is simply the degree of this map.

Uhlenbeck has also contributed to the theory of integrable systems and we mention some important integrability theorems here. As an example, take the unit hypercube

$$H = \{x \in \mathbb{R}^d \text{ such that } |x_i < 1|\}.$$

If $E$ is a bundle over the hypercube and $A$ is a flat connection on $E$, then there is a bundle isomorphism taking $E$ to the trivial bundle over the cube and $A$ to the product connection. Another example is the integrability theorem for holomorphic structures which tells us when a partial connection comes from a holomorphic structure on the bundle (the connection is said to be integrable in that case). Explicitly, a partial connection $\bar{\partial}_\alpha$ on a $C^\infty$ complex vector bundle over a complex manifold $X$ is integrable if and only if[126]

$$\bar{\partial}_\alpha{}^2 = \Phi_\alpha.$$

In general, Uhlenbeck's work fits in with that of a small group of differential geometers who contributed to making geometric analysis a leading subject. Probably the most obvious examples of equations which fit into geometric analysis are geometric flows (evolution equations for a functional on a manifold with geometric meaning). A famous example is the Ricci flow, used by Perelman to prove Thurston's geometrization conjecture and by implication the Poincaré conjecture. In this flow, the metric involves according to the Ricci curvature.

$$\frac{\partial g}{\partial t} = -2\mathrm{Ric}(g).$$

Classification of the singularities which this flow can develop is a rich and fascinating subject in its own right. Another example is the mean curvature flow. If $F$ is a map from $M^n \times [0, T)$ to a manifold $(N^{n+1}, \bar{g})$, then $F$ is said to solve mean curvature flow if

$$\frac{\mathrm{d}}{\mathrm{d}t} F(p,t) = \vec{H}\big(F(p,t)\big).$$

This is a second-order, quasi-linear, parabolic system. A geometric flow is a gradient flow associated to functional on a manifold: mean curvature flow can be thought of as the gradient flow of area with respect to the $L^2$ norm on the surface.

A trick of Uhlenbeck is often used in this area when one is considering the evolution of geometric quantities such as the curvature tensor under the flow. We might for example have an evolution equation for the curvature and wish to apply some kind of maximum principle to find curvature conditions preserved under Ricci flow. If we have an evolution equation for the Einstein tensor, we need to see this tensor not only as a section of $T^*M \otimes TM$, but also as a section of the sub-bundle of sections which are symmetric with respect to $g(t)$. This sub-bundle inherits the time dependence of the metric, which needs to be dealt with somehow[127]. The trick is to take an arbitrary vector bundle $V$ which is isomorphic to the tangent bundle via a bundle isomorphism and fix a pullback metric $h$ on the fibres. We choose a bundle isometry $u_0$ at time 0 between $V$ and $TM$ and allow a one-parameter family of bundle endomorphisms to evolve under the equation

$$\frac{\partial u(v)}{\partial t} = \mathrm{Ric}_{g(t)}\big(u(v)\big).$$

In this case the pullback metric $h = u_0^*(g_0)$ is constant in time as it has a zero time derivative, and the evolution equation for the metric can be used to show that

[126] Simon Donaldson and Peter Kronheimer, *The Geometry of Four-Manifolds* (Oxford: Oxford University Press, 1990).
[127] Peter Topping, *Lectures on the Ricci Flow* (Cambridge: Cambridge University Press, 2010).

$h = u_t^*\big(g(t)\big)$ for all $t$, hence $u$ remains a bundle isometry between the time-dependent metric on the tangent bundle and the fixed metric on the arbitrary bundle $V$. From there we can pull back with $u$, pull the curvature tensor back to a tensor on $V$ and define a connection on $V$ which is the pull-back of the Levi-Civita connection on the tangent bundle under $u$ at every time $t$[128].

[128] Richard Hamilton, 'Four-manifolds with positive curvature operator', *Journal of Differential Geometry*, 24, 2, 153 – 179 (1986).

## *Krystyna Kuperberg*

Kuperberg has researched topology and dynamical systems, including work related to the Seifert conjecture. The Seifert conjecture is a statement regarding the 3-sphere. The $n$-sphere is the most intensively studied of all manifolds, from the basic level to the advanced. A standard exercise with the $n$-sphere is to compute the transition map given by $\alpha_1 \circ \alpha_2^{-1}$ where $\alpha_1$ is stereographic projection from the south pole of the sphere and $\alpha_2$ is stereographic projection from the north pole. As I say, this is a problem which you will likely see in a homework sheet or an exam if you take a course on smooth manifolds. I will outline some of the solution if it is helpful, but try it yourself and discuss it with someone before looking up the answer. I will study the case of $S^2$, but all the arguments can be generalised.

We begin by parametrising the line which passes through the south pole and some other point $P$. From coordinate geometry, the parametric representation for the line passing through two points $P$ and $P'$ is given by

$$r(t) = P + tD,$$

where $D = P' - P$ and $t \in [0,1]$. Applying this to the situation we are studying, we have

$$\begin{aligned} r(t) &= S + t(P - S), \\ &= (0,0,-1) + t\big((x^1, x^2, x^3) - (0,0,-1)\big), \\ &= (0,0,-1) + t(x^1, x^2, x^3 + 1), \\ &= \big(tx^1, tx^2, -1 + t(x^3 + 1)\big). \end{aligned}$$

This is equivalent to a new set of coordinates:

$$z^1 = tx^1, z^2 = tx^2, z^3 = -1 + t(x^3 + 1).$$

Points on this pencil of lines will intersect the $x^3 = 0$ plane, in which case, doing some re-arranging we have:

$$z^3 = -1 + t(x^3 + 1) = 0.$$

$$t(x^3 + 1) = 1.$$

$$t = \frac{1}{x^3 + 1}.$$

Plug this back into the new set of coordinates:

$$(z^1, z^2, z^3) = \left(\frac{x^1}{1 + x^3}, \frac{x^2}{1 + x^3}, 0\right).$$

This defines a stereographic projection map $\alpha_1: S^2 - \{S\} \longrightarrow \mathbb{R}^2$ by

$$\alpha_1(x^1, x^2, x^3) = \frac{(x^1, x^2)}{1 + x^3}.$$

We repeat the process to define a stereographic projection from the north pole:

$$P = (0,0,1), \qquad P' = (x^1, x^2, x^3).$$

$$r(t) = (0,0,1) + t(x^1, x^2, x^3 - 1).$$

$$y^1 = tx^1, y^2 = tx^2, y^3 = t(x^3 - 1) + 1.$$

$$y^3 = \frac{1}{1 - x^3}.$$

Plug this back into $y^1$, $y^2$ and $y^3$:

$$(y^1, y^2, y^3) = \left(\frac{x^1}{1 - x^3}, \frac{x^2}{1 - x^3}, 0\right).$$

This gives us the corresponding map for projection from the opposite pole:

$$\alpha_2(x^1, x^2, x^3) = \frac{(x^1, x^2)}{1 - x^3}.$$

The map $\alpha_1$ projects the north pole to zero and sends the south pole to infinity, whilst $\alpha_2$ does the opposite. If we define the set $u_1 \cap u_2$ to be the sphere minus the north and south poles, then we must have

$$\alpha_1(u_1 \cap u_2) = \mathbb{R}^2 - \{0\},$$

$$\alpha_2(u_1 \cap u_2) = \mathbb{R}^2 - \{0\}.$$

We already have

$$\alpha_2(x^1, x^2, x^3) = \left(\frac{x^1}{1 - x^3}, \frac{x^2}{1 - x^3}\right) \coloneqq (y^1, y^2).$$

Consider

$$\begin{aligned}
(y^1)^2 + (y^2)^2 + 1 &= \frac{(x^1)^2 + (x^2)^2}{(1 - x^3)^2} + 1, \\
&= \frac{(x^1)^2 + (x^2)^2 + (x^3)^2 + 1 - 2x^3}{(1 - x^3)^2}, \\
&= \frac{2}{1 - x^3}.
\end{aligned}$$

Also

$$y^1 = \frac{x^1}{1 - x^3} \Longrightarrow x^1 = (1 - x^3)y^1 = \frac{2y^1}{1 + (y^1)^2 + (y^2)^2}.$$

$$y^2 = \frac{x^2}{1 - x^3} \Longrightarrow x^2 = \frac{2y^2}{1 + (y^1)^2 + (y^2)^2}.$$

$$1 - x^3 = \frac{2}{1 + (y^1)^2 + (y^2)^2}.$$

This implies that

$$x^3 = 1 - \frac{2}{1+(y^1)^2+(y^2)^2} = \frac{-1+(y^1)^2+(y^2)^2}{1+(y^1)^2+(y^2)^2}.$$

This tells us that

$$\alpha_2^{-1}(y^1, y^2) = (x^1, x^2, x^3).$$

We can use this to explicitly compute the transition map $\alpha_1 \circ \alpha_2^{-1}$:

$$\alpha_1 \circ \alpha_2^{-1}(y^1, y^2) = \left(\frac{x^1}{1+x^3}, \frac{x^2}{1+x^3}\right),$$

$$= \left(\frac{\frac{2y^1}{1+(y^1)^2+(y^2)^2}}{1+\frac{-1+(y^1)^2+(y^2)^2}{1+(y^1)^2+(y^2)^2}}, \frac{\frac{2y^2}{1+(y^1)^2+(y^2)^2}}{1+\frac{-1+(y^1)^2+(y^2)^2}{1+(y^1)^2+(y^2)^2}}\right)$$

$$= \left(\frac{\frac{2y^1}{1+(y^1)^2+(y^2)^2}}{\frac{-1+(y^1)^2+(y^2)^2+1+(y^1)^2+(y^2)^2}{1+(y^1)^2+(y^2)^2}}, \frac{\frac{2y^2}{1+(y^1)^2+(y^2)^2}}{\frac{-1+(y^1)^2+(y^2)^2+1+(y^1)^2+(y^2)^2}{1+(y^1)^2+(y^2)^2}}\right),$$

$$= \left(\frac{2y^1}{2(y^1)^2+2(y^2)^2}, \frac{2y^2}{2(y^1)^2+2(y^2)^2}\right),$$

$$= \frac{(y^1, y^2)}{|y|^2}.$$

Obviously, the amount of working I have shown here is not typical and shows a lot of messy, aesthetically unappealing things which are normally kept 'behind the scenes'. You should appreciate when studying mathematics that a phrase like 'after some algebra' or 'it is easy to show' could hide many pages of working and computation.

In order to verify that the atlas consisting of the two charts $(S^2 - \{S\}, \alpha_1)$ and $(S^1 - \{N\}, \alpha_2)$ defines a smooth structure on $S^2$, we need to confirm that $\alpha_1 \circ \alpha_2^{-1}(y^1, y^2)$ is infinitely differentiable on $\mathbb{R}^2 - \{0\}$. We can see that any derivative of

$$\left(\frac{y^1}{(y^1)^2+(y^2)^2}\right)$$

is always going to involve ratios of monomials in $(y^1)^m(y^2)^n$ over $((y^1)^2+(y^2)^2)^k$. These are continuous and infinitely differentiable on $\mathbb{R}^2 - \{0\}$, so the transition map is also infinitely differentiable on this domain. Bear in mind that two charts can provide an atlas for a manifold and an associated differentiable structure without being compatible for each other. As an example, take $\mathbb{R}$ to be a topological manifold. For general manifolds, we can only talk about existence of coordinate representations of points once we have a chart. However, in this case, we can write the maps down explicitly, since points on the real line can be denoted by a coordinate $x$. Say we have two charts given by $(\mathbb{R}, \alpha)$ and $(\mathbb{R}, \beta)$ where

$$\alpha(x) = x, \qquad \beta(x) = x^{\frac{1}{3}}.$$

The two charts can only be part of the same differentiable structure if $(\beta \circ \alpha^{-1})(x)$ is a diffeomorphism. Doing the computation we find

$$\alpha^{-1}(x) = x, \qquad (\beta \circ \alpha^{-1})(x) = \beta(x) = x^{\frac{1}{3}}.$$

The function $(\beta \circ \alpha^{-1})(x)$ is clearly not smooth at every point on the real line, since the derivative

$$\frac{\partial\left(x^{\frac{1}{3}}\right)}{\partial x} = \frac{1}{3}x^{-\frac{2}{3}}$$

blows up at $x = 0$. This implies that the two charts are not compatible.

The $n$-sphere also appears a lot in minimal surface theory. If $M$ is an $n$-dimensional closed minimal hypersurface in $S^{n+1}$, then one has

$$R = n(n-1) - S,$$

where $R$ is the scalar curvature and $S$ is the length squared of the second fundamental form of $M$. A well-known rigidity theorem of Chern and others states that if $S \leq n$, then either $S = 0$, or $S = n$. This implies that $M$ is the great sphere or the Clifford torus defined by

$$S^k\left(\sqrt{\frac{k}{n}}\right) \times S^{n-k}\left(\sqrt{\frac{n-k}{n}}\right).$$

An interesting problem in this area is whether there might be several scalar curvature pinching results for closed minimal hypersurfaces in a unit sphere. For example, it has been proved that if the scalar curvature of $M$ is a constant, then these exists a positive constant $\alpha(n)$ such that if

$$n \leq S \leq n + \alpha(n),$$

then $S = n$. Another related theorem states that in the case of $n \leq 5$, if $M$ is a closed minimal $n$-hypersurface in $S^{n+1}$, then there is a positive constant $\delta(n)$ such that if

$$n \leq S \leq n + \delta(n),$$

then $S = n$ identically. Take $M$ is an $n$-dimensional closed minimal hypersurface in $S^{n+1}$ where $n \geq 6$. An open problem is whether a positive constant $\delta(n)$ such that if

$$n \leq S \leq n + \delta(n),$$

then $S = n$. A solution to this problem was proposed by Ogiue and Sun for arbitrary $n$, but a fatal flaw was found in their proof[129].

When seeking to find relationships between curvature and topology, one usually focusses on the effect of positive or negative Ricci curvature, sectional curvature, or scalar curvature.

---

[129] Si-Ming Wei and Hong-Wei Xu, 'Scalar curvature of minimal hypersurfaces in a sphere', *Math. Res. Lett.*, 14(3), 423 – 432 (2007).

Other curvature assumptions are possible, however. One can assume positive curvature on totally isotropic two-planes. Under this assumption, and using the Sacks-Uhlenbeck theory of minimal 2-spheres, Micallef and Moore showed that if $M$ is a compact simply connected Riemannian $n$-manifold with $n \geq 4$ and positive curvature on totally isotropic two-planes, then $M$ is homeomorphic to a sphere. The curvature of a Riemannian manifold on totally isotropic two-planes is said to be positive if $K(\sigma) > 0$ whenever $\sigma \subseteq T_pM \otimes \mathbb{C}$ is a totally isotropic two-plane situated at any point in the manifold[130]. Homeomorphisms to a sphere also show up in one of the most famous conjectures in mathematics, the Poincaré conjecture. This conjecture states that every closed, smooth, simply connected 3-manifold is diffeomorphic to $S^3$. There is no need to specify that the manifold is homeomorphic to $S^3$, as classifications of topological 3-manifolds up to diffeomorphism or up to homeomorphism are equivalent. The conjecture can be written very simply in equation form:

$$\pi_1(M) = 0 \Longrightarrow M \cong S^3.$$

The conjecture is incredibly difficult to prove, and there are various abortive attempts in the literature[131]. It was finally proved by Perelman using techniques and ideas of Hamilton[132].

Interestingly, what is called an outermost minimal surface is always a sphere, so you may actually seem them referred to as outermost minimal spheres. (It could be a sphere or a torus technically, but a stability argument removes the latter possibility). Outermost minimal just means a minimal surface which is not contained completely inside another minimal surface. Define $\mathcal{S}$ to be the collection of surfaces which are smooth compact boundaries of open sets $U$ in $M^3$, where $U$ contains all the points at infinity $\{\infty_k\}$ which are added to compactify each of the ends $M_k$. Given a $\Sigma \in \mathcal{S}$, define $\tilde{\Sigma} \in \mathcal{S}$ to be the outermost minimal area enclosure of $\Sigma$. In the mathematical physics literature, one defines $\Sigma \in \mathcal{S}$ in a Cauchy initial data $(M^3, g, k)$ to be a future apparent horizon if

$$H_\Sigma + \mathrm{Tr}_\Sigma(k) = 0 \text{ on } \Sigma.$$

If $\Sigma \in \mathcal{S}$ has zero mean curvature (ie. it is a minimal surface), then the Penrose inequality is only satisfied when $A$ is the area of the outermost minimal area enclosure of $\Sigma$, with equality occurring only if $(M^3 - U, g)$ is isometric to Schwarzchild spacetime.

The Seifert conjecture (now known to be false) states that every nonsingular continuous vector field on $S^3$ has a fixed point or a closed orbit. A counterexample was found to this conjecture fairly early on, but that was not the end of the story, as one would like to improve the differentiability (smoothness) of these examples. Seifert was led to this conjecture when he showed that a flow-generating, nonsingular, continuous vector field on $S^3$ whose vectors are close to vectors which are tangent to the fibres of the Hopf fibration has a closed integral curve. Schweitzer gave a counterexample in dimension 3 by exhibiting an example of an aperiodic $C^1$ vector field on $S^3$. Kuperberg significant improved the differentiability of this example by finding a nonsingular $C^\infty$ vector field with no closed circular orbits on $S^3$: this was done independently of Schweitzer's example. The starting point of Kuperberg's construction is a mirror-image, Wilson-type plug which is modified by an operation called self-insertion[133].

---

[130] Mario Micallef and John Moore, 'Minimal two-spheres and the topology of manifolds with positive curvature on totally isotropic two-planes', *Annals of Mathematics*, 127, 199 – 227 (1988).

[131] David Gabai, 'Valentin Poenaru's program for the Poincaré conjecture', in Shing-Tung Yau (ed.) *Geometry, Topology, and Physics for Raoul Bott* (Cambridge, Massachusetts: International Press, 1994).

[132] John Morgan and Gang Tian, *Ricci Flow and the Poincaré Conjecture* (Providence, Rhode Island: American Mathematical Society, 2007).

Inspired by this real analytic counterexample, Greg Kuperberg modified the Schweitzer counterexample and proved that every 3-manifold admits a $C^1$ volume-preserving flow without fixed points or closed orbits. The main construction is a volume-preserving version of Schweitzer's plug (specifically, a Dehn-twisted Wilson-type plug). He also proved that every 3-manifold admits a $C^\infty$ volume-preserving flow with no fixed points or closed orbits, along with a transversely measured $PL$ foliation with discrete closed leaves. The proof of the theorem constructs a plug and inserts it into the flow, but there are technical obstructions due to the fact that a volume-preserving plug cannot plug up an open set. These are overcome by using Dehn-twisted plugs: these plugs can change the topology of a manifold during insertion. One important aspect of a volume-preserving flow in dimension 3 with no fixed points is that the parallel 1-dimensional foliation is transversely symplectic. Such a volume-preserving flow on a 3-manifold can be understood as a Hamiltonian flow coming from a symplectic structure on $M \times \mathbb{R}$. This links back to the Weinstein conjecture which we mentioned in the context of symplectic topology. In that case, the result is that the flow on a closed $(2n + 1)$-dimensional contact manifold which is transversely symplectic has a closed orbit if

$$H^1(M, \mathbb{Z}) = 0.$$

The Weinstein conjecture has been shown to be true for $S^3$, so the volume-preserving flow we are discussing here cannot be contact[134]. Recall that in slightly more old-fashioned terminology, to say that a manifold is closed implies that it is compact and that it has no boundary ie. the boundary is the empty set.

The Kuperbergs collaborated to construct generalized counterexamples to the conjecture. Using the theory of plugs and the self-insertion construction due to Krystyna, they proved that a foliation of any codimension of any manifold can always be modified in a real analytic or piecewise-linear way so that all minimal sets have codimension 1. A corollary of this theorem is that $S^3$ has an analytic dynamical system such that all limit sets are 2-dimensional. In particular, it has no circular trajectories. The generalized counterexamples are based on construction of aperiodic plugs. An insertible, untwisted, attachable plug is an orientated, 1-dimensional foliation $\mathcal{F}$ of a manifold with boundary given by the Cartesian product $F \times I$ of an $(n - 1)$-dimensional manifold and the interval $I$. The foliation agrees with the trivial foliation in the $I$ direction on a neighbourhood of the boundary of $F \times I$, if a leaf of the foliation connecting $(p, 0)$ with $(q, 1)$ implies that $p = q$ and if there is a non-compact leaf containing $(p, 0)$. An important result in foliation theory is the global Frobenius theorem, what states that if $D$ is an involutive distribution on a smooth manifold $M$, then the collection of all maximal connected integral manifolds of $D$ forms a foliation of $M$.

---

[133] Krystyna Kuperberg, 'A smooth counterexample to the Seifert conjecture', *Annals of Mathematics*, 2(140), 723 – 732 (1994).

[134] Greg Kuperberg, 'A volume-preserving counterexample to the Seifert conjecture', *Comment. Math. Helvetici*, 71, 70 – 97 (1996).

A distribution on $M$ is a choice of $k$-dimensional linear subspace of the tangent space at every point $p \in M$. A distribution is said to be smooth if the union of all of these linear subspaces forms a smooth subbundle

$$D = \coprod_{p \in M} D_p \subset TM.$$

The local frame criterion for subbundles immediately provides with a criterion for smoothness of a distribution. Essentially, every point of a manifold has a neighbourhood on which there are vector fields such that one can form a basis for the corresponding $k$-dimensional subspace in the distribution. In this case, one can say that the distribution is spanned locally by those vector fields. Now if we start with a distribution, an immersed submanifold $N$ contained in $M$ is called an integral manifold if we have the following

$$T_p N = D_p$$

at every point $p$ in $N$. As an example, if $V$ is a smooth nowhere-vanishing vector field on a manifold, then $V$ trivially spans a smooth 1-dimensional distribution on that manifold. The image of any integral curve of $V$ is an integral manifold of $D$.

Integral manifolds might fail to exist, however. Take $D$ a smooth tangent distribution on a manifold. $D$ is said to be an involutive distribution if for any pair of smooth local sections of $D$ (essentially smooth vector fields defined on the appropriate open subset), the Lie bracket is also a local section of the distribution. $D$ is said to be integrable if each point of $M$ is contained in an integral manifold of the distribution. It can be proved that every integrable distribution is involutive. Involutivity can also be described in terms of differential forms: in general, tangent distributions can be described locally by smooth 1-forms, as well as by vector fields. Recall that the graded algebra of smooth differential forms on a smooth manifold $M$ is defined by

$$\mathcal{A}^*(M) = \mathcal{A}^0(M) \oplus \dots \oplus \mathcal{A}^n(M),$$

where $\mathcal{A}^k(M)$ is the vector space of smooth sections of $\Lambda^k(M)$. $\Lambda^k(M)$ is defined as

$$\Lambda^k(M) = \coprod_{p \in M} \Lambda^k(T_p M).$$

$\Lambda^k V$ is the vector space of $k$-linear alternating forms defined as follows:

$$\Lambda^k V \coloneqq \{\varphi : V \times \dots \times V \longrightarrow \mathbb{R}\}.$$

The notation and terminology is slightly confusing: just think about everything in terms of maps if you are not sure what is going on, and it should make sense.

Recall that the map being alternative means that

$$\varphi(\ \ldots, v_i, \ldots, v_j, \ldots\ ) = -\varphi(\ \ldots, v_j, \ldots, v_i, \ldots\ ).$$

This immediately implies that an alternating multilinear form evaluated for a linearly dependent set of vectors vanishes to zero. For example, in the case of a bilinear form, if either of the vectors is zero, bilinearity immediately implies that the form is zero:

$$\varphi(v_1, v_2) = 0.$$

From linear algebra, if the two vectors are linearly dependent, then there must be some scalar $c$ such that

$$v_1 = cv_2.$$

The alternating property implies that

$$\varphi(v_1, v_2) = -\varphi(v_2, v_1).$$

Plugging in:

$$\varphi(cv_2, v_2) = -\varphi(v_2, cv_2).$$

Bilinearity implies that

$$c\varphi(v_2, v_2) = -c\varphi(v_2, v_2).$$

The only way this is satisfied is if

$$\varphi(v_2, v_2) = 0,$$

which establishes the result. Note that the above is not always zero, as it depends on the base field (if the base field is characteristic 2, then the above is not necessarily zero). This argument can be generalized to an arbitrary number of vectors, using multilinearity instead of bilinearity.

Involutivity can be expressed by saying that a smooth distribution is involutive if and only if $\mathcal{I}(D)$ is a differential ideal. $\mathcal{I}(D)$ is defined as

$$\mathcal{I}(D) = \mathcal{I}^0(D) \oplus \ldots \oplus \mathcal{I}^n(D) \subset \mathcal{A}^*(M),$$

where $\mathcal{I}^p(D) \subset \mathcal{A}^p(M)$ is the space of smooth $p$-forms which annihilate $D$. A central result of smooth manifold theory which comes from all this is that every involutive distribution is completely integrable, where a distribution is completely integrable if there always exists a flat chart for the distribution in a neighbourhood of every point in the manifold. This is known as the Frobenius theorem. The global Frobenius theorem follows once we prove that for an involutive distribution $D$ and $\{N_\alpha\}_{\alpha \in A}$ a collection of connected integral manifolds of the distribution with a point in common, the following union

$$N = \bigcup_{\alpha} N_\alpha$$

has a unique smooth structure as a smooth manifold which makes it into a connected integral manifold of $D$ in which each of the $N_\alpha$ is an open submanifold[135].

To return to our original foliation, the base $F$ is an $(n-1)$-manifold which admits a bridge immersion in $\mathbb{R}^{n-1}$, an immersion which lifts to an embedding $\mathbb{R}^n$. A plug is said to be aperiodic if it has no closed leaves. The idea which Schweitzer proposed is to use plugs to modify the foliation and break circular leaves, and he observed that the base of a plug only needs to admit a bridge immersion, and not an embedding. The theorem which we stated on the foliation of any codimension directly implies that the 3-sphere admits an analytic dynamic system such that all limit sets are 2-dimesional. The Schweitzer-Harrison plug $C^2$ aperiodic plug has two 1-dimesional minimal sets. In particular, it can be proved that a 1-foliation of a manifold of dimension at least 3 can be modified in the $PL$ category so that are no closed leaves but all minimal sets are 1-dimensional. If the manifold is closed, then there is an aperiodic $PL$ modification with only one minimal set, and that minimal set is 1-dimensional.

As we have said, the main purpose of a plug is to be able to carry out a geometric operation called insertion. The idea is to start with a foliation $\chi$ of an $n$-manifold $X$ and an $n$-dimensional plug $P$ with a base $F$. We then look for a leaf-preserving embedding of the trivially foliated $F \times I$ in the manifold, and see if it can be replaced with the plug. The technical obstructions involved when one tries to do this are insertibility, attachability and twistedness. An insertion map for a plug into a foliation is defined as an embedding

$$\sigma: F \longrightarrow X$$

of the base of $P$ which is transverse to $\chi$. This insertion map can be extended to the product to get an embedding

$$\sigma: F \times I \longrightarrow X$$

which takes the fibre foliation of $F \times I$ to $\chi$. An $n$-dimensional plug $P$ is said to be insertible if $F$ admits an embedding in $\mathbb{R}^n$ which is transverse to vertical lines. An embedding of this kind is equivalent to a bridge immersion of $F$ in $\mathbb{R}^{n-1}$. To be clear, it is equivalent to an immersion which lifts to an embedding of $F \times I$ in $\mathbb{R}^n$. As an example, if $F$ is 2-dimensional, orientable, and with a boundary which is not equal to the empty set, then $F$ must have a bridge immersion.

A plug is said to be attachable if every leaf in the parallel boundary is finite. If $M$ is a manifold with boundary, take $N_M$ to be an open neighbourhood of the boundary of $M$ in $M$. The next step in insertion of the plug is to remove $\sigma\big((F \times I) - N_{F\times I}\big)$ from the manifold and glue the lip $\sigma(N_{F\times I})$ to the support of an attachable plug via a leaf-preserving homeomorphism defined as

$$\alpha: N_{F\times I} \longrightarrow N_{\mathrm{supp}\,P},$$

[135] John M. Lee, *Introduction to Smooth Manifolds* (New York: Springer, 2003).I

where $N_{\mathrm{supp}\,P}$ is a neighbourhood of the boundary of the support of $P$. The identification also satisfies the following

$$\alpha(p,0) = \alpha_-(p), \qquad \alpha(p,1) = \alpha_+(p).$$

A map $\alpha$ with these properties is an attaching map for the plug. We can let $G$ be the parallel boundary for the support $\mathcal{P}$ of the plug, so that

$$\partial\mathcal{P} = G \cup F_- \cup F_+.$$

All of the leaves in this parallel boundary have two endpoints: a lemma then shows that there is a leaf-preserving homeomorphism $N_G \longrightarrow N_F \times I$. The equivalence ensures the existence of an attaching map for any attachable plug or un-plug[136]. By this point, you are starting to see how sophisticated the geometry is which is going into this: in this case, just the fact that someone found one particular counterexample to the conjecture was not the end of the story at all.

Aside from work on the Seifert conjecture, Kuperberg has also done other work on the application of topology to dynamical systems theory. In particular, with Reed, she studied trajectories of dynamical systems on $\mathbb{R}^3$. She constructed an example of a rest point free dynamical system on $\mathbb{R}^3$ with uniformly bounded trajectories and no circular trajectories. The construction was based on an example of Schweitzer's[137]. She also solved a problem related to the question of whether or not a homogeneous space has to be bihomogeneous. A counterexample to this conjecture was found many years ago, but whether or not it was true for continua remained an open problem for some time, until Kuperberg found a counterexample in the continuum case. A continuum is a topological space, often taken to be a nonempty compact connected metric space[138]. Textbook examples of continua which you might already have come across in topology are the Warsaw circle and the topologist's sine curve.

[136] Greg Kuperberg and Krystyna Kuperberg, 'Generalized counterexamples to the Seifert conjecture', *Annals of Mathematics*, 143(3), 547 – 576 (1996).
[137] Coke Reed and Krystyna Kuperberg, 'A dynamical system on $\mathbb{R}^3$ with uniformly bounded trajectories and no compact trajectories', *Proceedings of the American Mathematical Society*, 106(4), 1095 – 1097 (1989).
[138] Krystyna Kuperberg, 'A homogeneous nonbihomogeneous continuum', *Topology Proceedings*, 13, 399 – 401 (1988).

# *Nicole Tomczak-Jaegermann*

Tomczak-Jaegermann researches geometric functional analysis, especially the theory of Banach spaces. She has received the CRM-Fields-PIMS prize for her research, previously awarded to Coxeter and Tutte. The basic definition for a Banach space is that it is a complete normed vector space (where by complete, we mean complete in the metric induced by the norm). A norm on a real or complex vector space $X$ is a real valued function on $X$ whose value at a point in the space is denoted by $\|x\|$. The norm has to satisfy the following four properties:

$$\|x\| \geq 0,$$

$$\|x\| = 0 \Leftrightarrow x = 0,$$

$$\|\alpha x\| = |\alpha|\, \|x\|,$$

$$\|x + y\| \leq \|x\| + \|y\|,$$

where $x$ and $y$ are any two vectors in $X$ and $\alpha$ is any scalar. A norm on $X$ defines a metric $d$ on $X$ given by

$$d(x, y) = \|x - y\|.$$

This metric is called the induced metric. It can be checked easily that $d$ does define a metric space, and so normed spaces and Banach spaces are metric spaces.

An example of a norm on $C[0,1]$ is the supremum norm:

$$\|f\|_\infty = \sup_{t \in [0,1]} |f(t)|.$$

Other common examples are the $L^1$ and $L^2$ norms.

$$\|f\|_{L^1} = \int_0^1 |f(t)|\, dt,$$

$$\|f\|_{L^2} = \left( \int_0^1 |f(t)|^2\, dt \right)^{1/2}.$$

Different norms can be introduced on the same vector space, so we would like an equivalence relation which enables us to compare two norms. We say that two norms $\|\cdot\|_1$ and $\|\cdot\|_2$ on a vector space $V$ are equivalent if there are positive constants such that

$$c_1 \|x\|_1 \leq \|x\|_2 \leq c_2 \|x\|_1,$$

for all $x \in V$. We can then write

$$\|\cdot\|_1 \sim \|\cdot\|_2.$$

Any two norms on a Euclidean space are equivalent. As an example of two norms which are not equivalent, one can consider the $L^1$ and the supremum norm on $C[0,1]$.

If $V$ and $W$ are two normed spaces, a linear map $L$ from $V$ to $W$ such that

$$\|L(x)\| = \|x\|$$

is called a linear isometry. Two normed spaces are isometrically isomorphic if there exists an invertible linear isometry between them. A norm can be pulled back by a linear invertible map, since if $V$ is a space equipped with a norm, $W$ is a vector space, and $L$ is a linear isomorphism, then one can define a norm on $W$ as

$$\|x\|_W \coloneqq \|L(x)\|_V.$$

You can check that this does satisfy the requirements for a norm, including the triangle inequality. If $V$ is a finite dimensional vector space and $n = \dim V$, then $V$ is linearly isomorphic to the space of tuples with entries from $\mathbb{R}$ or $\mathbb{C}$. However, since a norm can be pulled back by a linear invertible map, this implies that any finite dimensional vector space can be equipped with a norm and any normed space of dimension $n$ is isometrically isomorphic to $\mathbb{R}^n$ with a suitable norm. Since two norms on $\mathbb{R}^n$ are equivalent, we also have that all norms on a finite dimensional vector space are equivalent.

We would now like to consider complete normed spaces. Recall that a sequence $(x_n)_{n=1}^{\infty}$ in a normed space $V$ is Cauchy if for any positive non-zero $\varepsilon$ there exists an $N$ such that

$$\|x_n - x_m\| < \varepsilon$$

for all $m, n > N$. A sequence of real numbers converges if and only if it as Cauchy sequence. One can also prove that any convergent sequence is Cauchy, and that any Cauchy sequence is bounded. To prove the former for a metric space, realise that if $x_n$ converges to $x$ in the limit, then for every $\varepsilon > 0$ there must be some $N$ depending only on $\varepsilon$ such that

$$d(x_n, x) < \frac{\varepsilon}{2}$$

for $n > N$. By the triangle inequality we find that for $m, n > N$

$$d(x_m, x_n) \leq d(x_m, x) + d(x, x_n) < \frac{\varepsilon}{2} + \frac{\varepsilon}{2} = \varepsilon.$$

However, you can see from the definition that this implies that $(x_n)$ is a Cauchy sequence. Since a sequence of real numbers converges only if it is Cauchy, one has that $\mathbb{R}$ is complete, and by extension $\mathbb{C}$. Every finite dimensional normed space can be seen to be complete by considering basis vectors.

A good example of a Banach space is $\ell^p(\mathbb{R})$ equipped with the usual $\ell^p$ norm. To begin with the proof of completeness, just take a tuple of sequences in $\ell^p(\mathbb{R})$, assume it to be Cauchy, and show that the sum which one obtains implies that

$$\left| x_j^{(m)} - x_j^{(n)} \right| < \varepsilon,$$

since all the terms in the sum will be positive. From there you will need to show that $x^{(m)}$ converges to a tuple of real constants as $m$ tends to infinity, which establishes that $\ell^p(\mathbb{R})$ is complete. The argument also goes through the same if we replace real with complex numbers. As another example, $C[0,1]$ equipped with the supremum norm is complete. The proof uses a basic fact from analysis that the uniform limit of a sequence of continuous functions is also continuous. It is possible to have a space with is complete when equipped with one norm, but not complete when equipped with another one. For example, $C[0,1]$ equipped with the $L^1$ norm is not complete. This can be proved by taking a sequence of continuous functions as follows:

$$f_n(t) = \begin{cases} (2t)^n, & \text{for } 0 \leq t \leq 1/2, \\ 1, & \text{for } 1/2 \leq t \leq 1. \end{cases}$$

It can then be shown that $f_n$ does not converge to a continuous function in the $L^1$ norm via an argument by contradiction whereby one assumes that the limit does exist and calls it $f$. Another example of an incomplete metric space is $\mathbb{Q}$ equipped with the usual metric. A normed space $X$ can be viewed as a dense subset of a larger Banach space $\mathcal{X}$ called a completion of $X$. A natural question to ask is whether an incomplete normed vector space can always be completed in this way. (In fact, we have given the game away with this question, as we stated near the beginning of the book that a normed space can always be completed by equipping with an inner product and taking the closure) [139].

A subset $X \subset V$ is dense in $V$ if for any $v \in V$ and any $\varepsilon > 0$ there is $x \in X$ such that

$$\|x - v\|_V < \varepsilon.$$

An obvious example is $\mathbb{Q}$, which is dense as a subset of $\mathbb{R}$. The real numbers can be obtained from this subset by completing $\mathbb{Q}$ with Dedekind cuts[140]. Similarly, the space of tuples with rational entries $\mathbb{Q}^n$ is a dense subset of $\mathbb{R}^n$. If we take two normed spaces and a linear map $L : V \to W$, this map is called a linear isometry if $L$ is linear, bijective, and

$$\|L(x)\|_W = \|x\|_V$$

for all $x \in V$. We say that the two spaces are linearly isometric. Any normed vector space $X$ is linearly isometric to a dense subset of a Banach space $\mathcal{X}$ equipped with a norm $\|\cdot\|_{\mathcal{X}}$. If $X$ is already Banach, then the linear isometry is just the identity. The proof of the completion theorem is quite long, but essentially you are just looking to construct $\mathcal{X}$ using a suitable equivalence relation between Cauchy sequences, show that is a normed vector space, and then show that it is complete and that $X$ is always isometric to a subset of $\mathcal{X}$ such that the subset is dense.

Obviously, at this point, we do not know if the completion is unique, but it can be shown with less work than the proof of the previous theorem that if $\mathcal{X}$ and $\tilde{\mathcal{X}}$ are Banach spaces where the isometries $i : X \to \mathcal{X}$ and $\tilde{\imath} : X \to \tilde{\mathcal{X}}$ are such that

$$\|x\|_X = \|i(x)\|_{\mathcal{X}} = \|\tilde{\imath}(x)\|_{\tilde{\mathcal{X}}},$$

then if $i(x)$ is a dense subset of $\mathcal{X}$ and $\tilde{\imath}(x)$ is a dense subset of $\tilde{\mathcal{X}}$, $\mathcal{X}$ and $\tilde{\mathcal{X}}$ are linearly isometric ie. all completions of a space are the same modulo linear isometries. There are also more advanced theorems relating to linear maps between Banach spaces. We will start by recalling the Baire category theorem. This starts that if $\{G_i\}_{i=1}^{\infty}$ is a countable family of dense open subsets of a complete metric space, then the intersection

$$G = \bigcap_{i=1}^{\infty} G_i$$

is dense in $X$. In words, every countable collection of dense open subsets has a dense intersection. It is also possible to formulate this theorem with reference to a locally compact Hausdorff space, and this is one of the main properties that such spaces share with complete metric spaces. A topological space for which every countable union of dense open subsets is dense is said to be a Baire space. The Baire property is topological so any space which is homeomorphic to a complete metric space must be a Baire space. As a corollary, every meagre subset of a Baire space has a dense complement, where a subset of a topological space is said to be meagre if it can be expressed in the form of a union of a countable number of nowhere dense subsets.

[139] Erwin Kreyszig, *Introductory Functional Analysis with Applications* (USA: John Wiley & Sons, 1978).
[140] G. H. Hardy, *A Course of Pure Mathematics* (Cambridge: Cambridge University Press, 1952).

The first big theorem we would like to state for linear maps between Banach spaces is the principle of uniform boundedness. Start with a Banach space $X$ and a normed space $Y$. Take $S$ to be a subset of bounded linear maps from $X$ to $Y$ such that

$$\sup_{T \in s} \|T(x)\|_Y < \infty$$

for every $x \in X$. If these assumptions hold, then

$$\sup_{T \in s} \|T\| < \infty.$$

This a very useful theorem in applications. As a corollary, take $T_n$ in the space of bounded linear maps and suppose that the norm $\|T_n\|$ is unbounded. There must exist $x \in X$ such that $\|T_n x\|$ is unbounded. This can be proved by a contradiction: assuming that this $x$ does not exist implies that $\|T_n\|$ is bounded. As another corollary, suppose that $X$ is a Banach space, $Y$ a normed space and $T_n$ a bounded linear map as before. Suppose that the following limit exists for every $x \in X$:

$$T(x) \coloneqq \lim_{n \to \infty} T_n(x).$$

If this is the case, then $T \in B(X, Y)$. This can be proved by checking that $T$ is both linear and bounded. The latter follows from the principle of uniform boundedness.

Our next big theorem is the open mapping theorem. Be warned that the statement and proof of this result in some textbooks of functional analysis can be difficult to understand. We can state the theorem as follows: if $T: X \to Y$ a bounded surjective linear map from a Banach space $X$ to another Banach space $Y$, then $T$ maps open sets in $X$ to open sets in $Y$. The proof uses the Baire category theorem. As a corollary, one obtains the inverse mapping theorem. This says that if $T: X \to Y$ is a bounded bijective linear map from a Banach space $X$ onto a Banach space $Y$, then the inverse $T^{-1}$ is bounded. We already know that $T$ has an inverse. This inverse is linear because

$$T\big(\alpha T^{-1}(y_1) + \beta T^{-1}(y_2)\big) = \alpha y_1 + \beta y_2 = T\big(T^{-1}(\alpha y_1 + \beta y_2)\big)$$

which implies that

$$T^{-1}(\alpha y_1 + \beta y_2) = \alpha T^{-1}(y_1) + \beta T^{-1}(y_2).$$

The inverse is also unique, and the result then follows from the open mapping theorem. From standard spectral theory, we know that in infinite dimension, the spectrum of a bounded linear self-map $T$ consists of the set of all $\lambda$ such that the operator $T - \lambda I$ does not have a bounded inverse. The inverse mapping theorem then tells us that if $T$ is a bounded linear self-map and $T - \lambda I$ is a bijection, then the inverse $(T - \lambda I\,)^{-1}$ is bounded. This implies $\lambda$ does not belong to the spectrum of the operator. A corollary of the inverse mapping theorem states that if $X$ is a Banach space which is complete with respect to different norms $\|\cdot\|_1$ and $\|\cdot\|_2$ such that

$$\|x\|_2 \leq C\|x\|_1,$$

then the two norms are equivalent.

Our final theorem is the closed graph theorem. This states that if $T: X \to Y$ is a linear map between Banach spaces and if the graph of $T$

$$G := \big\{\big(x, T(x)\big) \in X \times Y : x \in X\big\}$$

is a closed subset of the product $X \times Y$ with an appropriately defined norm, then $T$ is a bounded operator. The graph $G$ is a closed linear subspace of $X \times Y$, so it must be a Banach space which equipped with the product norm which we just mentioned. Consider the projector map which takes the graph to one of the halves of the product, defined as

$$\pi_X(x, y) = x.$$

This map is a linear, bounded bijection. By the inverse mapping theorem, we know that the inverse is also bounded, so

$$\|\pi_X^{-1}(x)\|_{X\times Y} = \|x\|_X + \|T(x)\|_Y \leq M\|x\|.$$

This implies that

$$\|T(x)\|_Y \leq M\|x\|,$$

which completes the proof. The theorem can be used to prove an important result known as the Hellinger-Toeplitz theorem, which states that if $H$ is a Hilbert space and $T : H \to H$ is a linear operator which satisfies

$$(T(x), y) = \big(x, T(y)\big),$$

then $T$ is bounded. In words, a symmetric operator which is defined everywhere on a Hilbert space is bounded. This has important consequences, as it says that an unbounded linear operator which satisfies the above cannot be defined everywhere on the Hilbert space. However, unbounded operators arise naturally in quantum mechanics (the operator corresponding to the energy observable, for example), but we will not go any further into the unbounded theory here.

There is also a body of more advanced theory relating to Banach spaces. The starting point for this theory is to start thinking about bases in Banach spaces. In common with the usual linear algebra theory, suppose that a Banach space has a countable basis $\{e_n\}$ such that any element of the space can be written uniquely as a decomposition:

$$x = \sum_{j=1}^{\infty} x_j e_j,$$

where the sum converges in the Banach space. If we consider the following finite expansions up to $n$ terms

$$P_n(x) = \sum_{j=1}^{n} x_j e_j,$$

a question can we ask is whether there exists a constant such that

$$\|P_n(x)\| \leq C\|x\|,$$

for every $n$?

To prove this we need to show that $P_n \in B(X, X)$ for every $n$. For this, we need some fairly involved arguments which show that if $(e_n)$ is a countable sequence in a Banach space such that $\|e_n\| = 1$ and such that the closed linear span is the whole of $X$, then

$$|\|x\|| = \sup_{n} \left\| \sum_{j=1}^{n} a_j e_j \right\|_X$$

defines a norm on the Banach space such that the space is complete with respect to this norm. A corollary is that if a sequence $(e_n)$ satisfies the previous conditions, then there is a constant $C$ such that

$$\left\| \sum_{j=1}^{n} x_j e_j \right\| \leq C\|x\|,$$

for $n \in \mathbb{N}$[141].

It is of course not clear that infinite dimensional Banach spaces can have a basis, since the concept of a basis originates in finite dimensional linear algebra. If a normed space contains a sequence of vectors $(e_n)$ with the property that there always exists a unique sequence of scalars $(\alpha_n)$ such that

$$\|x - (\alpha_1 e_1 + \cdots + \alpha_n e_n)\| \longrightarrow 0 \text{ as } n \longrightarrow \infty,$$

then $(e_n)$ is called a Schauder basis for the Banach space. We then have as before the expansion of $x$ with respect to $(e_n)$:

[141] Bryan Rynne and Martin Youngson, *Linear Functional Analysis* (London: Springer, 2008).

$$x = \sum_{k=1}^{\infty} \alpha_k e_k.$$

In particular, this means that the following sequence

$$\left( \sum_{k=1}^{N} \alpha_k e_k \right)_{N=1}^{\infty}$$

must converge to $x$ in the norm topology of $X$. In the case of the sequence space $\ell^p$, a Schauder basis is given by the sequence of tuples whose $n$th term is $1$ and all other terms are zero:

$$e_1 = (1,0,0, \dots),$$

$$e_2 = (0,1,0, \dots),$$

and so on. It can be shown that a normed space with a Schauder basis is separable. Banach asked whether the converse is true: that is, whether every separable Banach space has to have a Schauder basis. Enflo proved that the answer is negative by constructing a counterexample[142].

The concept of a basis in an infinite dimensional Banach space is distinct from the concept of a Hamel basis, which is a vector space basis, being defined merely as a linearly independent subset of the space which spans $V$. More concretely, a Hamel basis for a vector space is a collection of linearly independent vectors such that every element of the space can be represented uniquely as a finite linear combination of the members of the basis. One can prove with the axiom of choice or Zorn's lemma that every vector space has a Hamel basis. Also, if a vector space has a finite Hamel basis, then every Hamel basis for the space has the same number of elements. This is essentially just a linear algebra argument. The Baire category theorem can be used to deduce that if $(e_i)$ is a Hamel basis for an infinite dimensional Banach space, then it must be uncountable. The closed graph theorem can be used to show that a sequence in a separable Banach space is only a Schauder basis if it is equivalent to the type of basis we are already familiar with. The proof is elegant, but the result is not as useful as it sounds, as it is rare that you would know that something was a basis without already knowing that it was a Schauder basis.

If $(e_n)_{n=1}^{\infty}$ is a basis for a Banach space, then the number

$$K_b = \sup_n \|S_n\|$$

is the basis constant, where $S_n$ form a sequence of projections associated with the elements of the Schauder basis. One also has the corresponding notion of a basic sequence. A sequence $(e_k)_{k=1}^{\infty}$ in a Banach space is called a basic sequence if it is a basis for the closed linear span of $(e_k)_{k=1}^{\infty}$. There is also some very interesting theory on the nonlinear geometry of Banach spaces. Since a Banach space is also a metric space by definition, the question of nonlinear geometry is to what extent does the metric structure of a Banach space determine its linear structure. All the theory which we have discussed so far has been regarding linear operators between Banach spaces, but there are many natural maps between these spaces which are non-linear (though they might still have good regularity properties). Investigating these maps is a relatively new field of research in functional analysis.

---

[142] Erwin Kreyszig, *Introductory Functional Analysis with Applications* (USA: John Wiley & Sons, 1978).

If we have very little information about the metric structure of a Banach space, we might ask how much we can know about the space if we only know the homeomorphic class of $X$ as a topological space. In effect, we would like to know when two separable Banach spaces are homeomorphic as topological spaces. This problem was resolved a few decades after it was posed by Fréchet: all separable infinite dimensional Banach spaces are homeomorphic. In the opposite case we assume that the map between the Banach spaces is an isometry, meaning that the map essentially carries all the information we might wish to know about the metric structure of $X$. In this situation it was shown that one can obtain complete information about the real linear structure. Specifically, if $X$ and $Y$ are real normed spaces and if $\Phi$ is an isometry which maps $X$ to $Y$ and $0$ to $0$, then $\Phi$ is a bounded linear operator. If two real Banach spaces are isometric as Banach spaces, then they are also linearly isometric.

In general, we will find ourselves between these two extremes, so the goal is to find out how much of the information about the structures on the first Banach space is transferred to the second space under some map $f$. For example, if it can be proved that if $f$ is a uniform homeomorphism between two normed spaces, then $f$ is a coarse surjective Lipschitz embedding. A basic question in the area is whether the fact that there is a Lipschitz map between Banach spaces with some desirable property such as being a linear isometry establishes the existence of a corresponding linear map which also has that same property[143]. There is no room to go any further into the theory of Banach spaces, but hopefully this has been good motivation. We might also mention that the work of Tomczak-Jaegermann was used in an essential way by Gowers in his solution to Banach's homogeneous space problem. A Banach space is said to be homogeneous if it is isomorphic to all of its infinite dimensional Banach spaces. The homogeneous space problem asks if every homogeneous Banach space isomorphic to a separable Hilbert space. The specific result of Tomczak-Jaegermann was the fact that a Banach space $X$ with cotype $q$ for some finitie $q$ either has a subspace without an unconditional basis or it has a subspace which is isomorphic to $\ell_2$. In addition, if $X$ is a homogeneous Banach space, then either it is isomorphic to $\ell_2$ or it fails to have an unconditional basis[144].

Tomczak-Jaegermann is also interested in probability theory, especially the theory of random matrices and random graphs. Random matrices are often more useful than regular matrices because in the limit of large matrix dimension, the statistical correlations for the spectra of an ensemble of random matrices do not depend on the probability distribution which defines the ensemble, and only depend on the invariant properties of the distribution. This is known as the universal property of random matrices. Techniques from random matrix theory also enable analytic computations which is often impossible to achieve otherwise. This concept of universality means that a matrix does not need to have independent Gaussian random variables as matrix entries. One could theoretically create a random matrix by filling each slot with $1$ or $-1$ with a $50/50$ probability for each. Depending on the overall symmetry of the matrix, one can end up with the common eigenvalue statistics of the Gaussian random matrix ensembles even though we have different rules for creating the matrix elements.

---

[143] Fernando Albiac and Nigel Kalton, *Topics in Banach Space Theory* (Switzerland: Springer, 2016).
[144] R. A. Komorowski and Nicole Tomczak-Jaegermann, 'Banach spaces without local unconditional structure', *Israel J. Math.* 89, 205 – 226 (1995).

However, we need to be slightly careful with statements such as these, as it is only the local eigenvalue statistics (those that operate over a few mean spacings, perhaps) which are universal in this way, and there are other global properties which can depend on the distribution from which the elements are drawn. If you compute the density of the eigenphases on the unit circle, the density will only be uniform for the Gaussian elements, and increasing the size of the matrix will not change this, whereas local properties (nearest neighbour statistics) will become universal in the limit of a large matrix. Filling a matrix with random elements will not give universality for these global properties unless the random number generator samples from a Gaussian distribution. As we are interested in eigenvalues, we consider random unitary matrices. The concept of a random unitary matrix is usually referring to the Haar measure on the group of unitary matrices. One takes the Lie group of matrices $U(n)$ and then makes it into a probability space by assigning the Haar measure, which is the unique measure invariant under multiplication in the Lie group. There are numerical methods for generating random unitary matrices which sample from either of the classical Lie groups $U(n)$ or $O(n)$, where the sampling is uniform since it is done with respect to the uniform Haar measure[145].

After generating these matrices, one finds that for each $n \times n$ matrix there is a list of eigenvalues

$$e^{i\theta_1}, \dots, e^{i\theta_n}.$$

These eigenvalues lie in the range

$$0 < \theta_1 \leq \cdots \leq \theta_n \leq 2\pi.$$

These eigenvalues have a density of $n/2\pi$, such that there are $n$ eigenangles for every interval with length $2\pi$. When considering the eigenangle statistics, it is common to scale them so that the average spacing between consecutive eigenangles has unit size:

$$\phi_1 = \theta_1 \frac{N}{2\pi}, \dots, \phi_n = \theta_n \frac{n}{2\pi}.$$

If we compute a list of differences between neighbouring scaled eigenangles and make a histogram of this list (normalized to unity), the histogram will approximate the probability density of the eigenangle spacings. Integrating the spacing distribution from $0$ to $a$ gives the probability that two consecutive eigenangles will have a distance of separation which is between $0$ and $a$ (same as for any probability distribution like the Maxwell distribution). If we average this spacing distribution over many random matrices, we find that in the limit the distribution will approach the theoretical curve for the relevant ensemble: in the case of $U(n)$ this is the Circular Unitary Ensemble.

[145] Francesco Mezzadri, 'How to generate random matrices from the classical compact groups', *Notices of the American Mathematical Society*, 54, 592 – 604 (2007).

Tomczak-Jaegermann showed that for $A$ an $n \times n$ random matrix with identically independently distributed entries of zero mean, unit variance and a bounded subgaussian moment, the condition number $s_{\max}(A)/s_{\min}(A)$ satisfies a small ball probability estimate:

$$\mathbb{P}\left\{\frac{s_{\max}(A)}{s_{\min}(A)} \leq \frac{n}{t}\right\} \leq 2\exp(-ct^2), \qquad t \geq 1,$$

where $c > 0$ might only depend on the subgaussian moment. A random variable $\xi$ has subgaussian moment bounded above by $K > 0$ if

$$\mathbb{P}\{|\xi| \geq t\} \leq \exp\left(1 - \frac{t^2}{2K^2}\right), \ \ t \geq 0.$$

The values $s_i(A)$ are the singular values of the matrix we have defined and the condition number is defined as

$$\kappa(A) = \frac{s_{\max}(A)}{s_{\min}(A)},$$

where $s_{\max}$ can be viewed as $s_1$ and $s_{\min}$ can be viewed as $s_n$ when the singular values are arranged in non-increasing order[146].

In another work, Tomczak-Jaegermann shows that for a large integer $d$. If $n \geq 2d$ and $A_n$ is the adjacency matrix of a random directed $d$-regular graph on $n$ vertices with the uniform distribution, then the rank of the adjacency matrix is at least $n - 1$ with probability which tends to $1$ as as $n$ tends to infinity[147]. She also showed that if $d$ satisfies an inequality

$$d \leq \log^C n$$

for every large integer $n$, then the empirical spectral distribution of an appropriately rescaled adjacency matrix converges weakly in probability to the circular law. This result gives a complete solution to the problem of weak convergence of the empirical distribution in the setting of directed $d$-regular graphs, developing in the process a technique for bounding intermediate singular values of the adjacency matrix[148].

---

[146] Alexander Litvak, Konstantin Tikhomirov and Nicole Tomczak-Jaegermann, 'Small ball probability for the condition number of random matrices', arXiv:190-08655vs (2019).
[147] Alexander Litvak, Anna Lytova, Konstantin Tihhomirov, Nicole Tomczak-Jaegermann and Pierre Youssef, 'The rank of random regular digraphs of constant degree', *Journal of Complexity*, 48, 103 – 110 (2018).
[148] Alexander Litvak, Anna Lytova, Konstantin Tihhomirov, Nicole Tomczak-Jaegermann and Pierre Youssef, 'Circular law for sparse random regular digraphs', arXiv:1801-05576v2 (2018).

# *Dusa McDuff*

McDuff has received many prizes for her contributions to symplectic geometry and topology. She is also a Fellow of the Royal Society. We have already mentioned symplectic geometry a few times, so we should probably explain the topological side of things a bit more. Symplectic topology studies the global aspects of symplectic geometry, whereas the local structures of a symplectic manifold are equivalent to Euclidean structure. This follows from Darboux's theorem, which states that every symplectic form $\omega$ on $M$ is locally diffeomorphic to the standard form $\omega_0$ on Euclidean $\mathbb{R}^{2n}$. To prove this, one first needs a lemma. Take a smooth manifold $M^{2n}$ and $Q \subset M$ a compact submanifold. Suppose that $\omega_1$ and $\omega_2$ are closed 2-forms such that at every point $q \in Q$ the $\omega_1$ and $\omega_2$ are equal and non-degenerate on the tangent space $T_qM$. If this is the case, then there exist open neighbourhoods $U_1$ and $U_2$ of $Q$ and a diffeomorphism $\psi: U_1 \to U_2$ such that

$$\psi = \mathrm{Id}, \qquad \psi^* \omega_2 = \omega_1,$$

where for the first equation we have restricted to $Q$. To prove this, it is sufficient to show that there exists a 1-form $\sigma \in \Omega^1(U_1)$ such that

$$\sigma = 0, \qquad d\sigma = \omega_2 - \omega_1.$$

One proves this by restricting the exponential map to the normal bundle $TQ^{\perp}$ of the submanifold $Q$ with respect to any Riemannian metric $g$[149].

The normal bundle is an example of a vector bundle which is relatively easy to visualise. To be clear, a smooth $n$-dimensional vector bundle is a pair of smooth manifolds $E$ (known as the total space) and a base space $M$, along with a surjective map $\pi: E \to M$ called the projection map, such that every fibre $\pi^{-1}(p)$ of $E$ over a point $p$ is equipped with the structure of a vector space. In addition for every point in $M$ there is a neighbourhood $U$ of $p$ and a diffeomorphism

$$\varphi: \pi^{-1}(U) \to U \times \mathbb{R}^n$$

called a local trivialisation of $E$ such that one has a commutative diagram when one includes the projection onto the first factor in the product

$$\pi_1: U \times \mathbb{R}^n \to U.$$

The restriction of $\varphi$ to a fibre $\varphi: \pi^{-1}(p) \to \{p\} \times \mathbb{R}^n$ is a linear isomorphism. You see that we are essentially back in linear algebra again, as this construction has essentially glued together a set of vector spaces into a union, but remember that this bundle is also a manifold.

The prototype example is the tangent bundle $TM$, which is the disjoint union of the tangent spaces which one gets at every point. If you take the cotangent spaces instead, you construct a different bundle called the cotangent bundle $T^*M$. Also, consider a bundle such that the fibre at every point is the orthogonal complement of $T_pM$ for a submanifold of Euclidean space: this is called the normal bundle to the submanifold[150]. It is a result of Whitney (the embedding theorem) that every smooth $n$-manifold has a proper smooth embedding into $\mathbb{R}^{2n+1}$. It is a corollary of this theorem that every smooth $n$-manifold is diffeomorphic to a closed embedded submanifold of $\mathbb{R}^{2n+1}$. This follows because we already know that there is a proper smooth embedding, and the image of that has to be an embedded submanifold which is closed in the ambient space because proper continuous maps are closed. As a corollary of the previous statement, every smooth manifold has a distance function whose metric topology is equivalent to the induced topology. This follows because any subspace of a metric space has to inherit a metric which controls the topology of the subspace.

[149] Dusa McDuff and Dietmar Salamon, *Introduction to Symplectic Topology* (New York: Oxford University Press, 1998).

The Whitney embedding theorem answered one of those obvious-sounding questions which no-one seemed able to answer for a long time: are there smooth manifolds which are not diffeomorphic to submanifolds of Euclidean space? A similar result called the Whitney immersion theorem states that every smooth $n$-manifold has an immersion into $\mathbb{R}^{2n}$. There are also stronger versions of both these theorems, which state that if $n > 1$, then every smooth $n$-manifold has an immersion into $\mathbb{R}^{2n-1}$, and if $n > 0$, then every smooth $n$-manifold has a smooth embedding into $\mathbb{R}^{2n}$. Nash's embedding theorem is a similar result for Riemannian manifolds[151]. For a point $p$ in a Riemannian manifold $M$, we define the exponential map

$$\exp_{\mathrm{p}}: V_p \longrightarrow M,$$

$$v \longmapsto \gamma_v(1),$$

where

$$\gamma_v(t): [0, T] \longrightarrow M$$

is the maximal geodesic such that

$$\gamma_v(0) = p, \qquad \partial_t \gamma_v(0) = v.$$

You may see some authors write $\exp_p^\nabla$ to emphasise the dependence on the connection, but we will drop this, as it is implicit. $V_p$ is the set of elements of the tangent space at $p$ such that $\gamma_v$ is defined on the unit interval. The exponential map takes a neighbourhood of $0 \in T_pM$ and maps it diffeomorphically onto a neighbourhood of $p \in M$, so you can view it as a map which starts at the point 0 at the tangent space to a point and then pushes down to the corresponding point on the manifold. This can be proved with a short computation.

$T_pM$ is a vector space, so it can be identified with $\mathbb{R}^n$. For the same reason, $T_0(T_pM)$ is isomorphic to $T_pM$ itself. The derivative of $\exp_{\mathrm{p}}(0)$ is then a linear self-map for the tangent space:

---

[150] John M. Lee, *Riemannian Manifolds: An Introduction to Curvature* (New York: Springer, 1997).
[151] John M. Lee, *Introduction to Smooth Manifolds* (New York: Springer, 2003).

$$d\exp_{\mathrm{p}}(0)\colon T_pM \longrightarrow T_pM.$$

By the identification with Euclidean space, the map $\exp_p^{-1}$ defines a local coordinate chart for neighbourhood of $p$, where $p$ is mapped to $\exp_p^{-1}(p) = 0$. The local coordinates defined by the chart $\left(\exp_p^{-1}, U\right)$ are called Riemannian normal coordinates with centre $p$. In normal coordinates, $x = \left\{x^i\right\}$ around the point $\exp_p^{-1}(p) = 0$, one has the following for a Riemannian metric $g$:

$$g_{ij}(x) = \delta_{ij} - \frac{1}{3}R_{ipjq}x^p x^q + O(|x|^3),$$

$$g_{ij}(0) = \delta_{ij}, \quad \partial_k g_{ij}(0) = 0, \qquad \Gamma_{ij}^k(0) = 0,$$

$$\Delta\left(g_{ij}\right) = \left(\sum_{k=1}^{n} \frac{\partial^2}{\partial(x^k)^2} g_{ij}\right) = -\frac{2}{3}R_{ij}.$$

Computations are usually significantly simpler with this choice of local coordinates, so you should check that you can reproduce the rest of the results in this lemma, and look things up in a standard textbook if you need some pointers. The fact that $g_{ij}(0) = \delta_{ij}$ follows immediately from identifying the tangent space with Euclidean space, since this identification must map an orthonormal basis of $T_pM$ with respect to the Riemannian metric onto the corresponding Euclidean orthonormal basis of the Euclidean space. The expansion for the metric follows from a Taylor expansion:

$$g_{ij}(x) = \delta_{ij} + \frac{1}{2}\partial_p\partial_q g_{ij}x^p x^q + O(|x|^3).$$

A computation shows that[152]

$$\frac{1}{2}\partial_p\partial_q g_{ij}x^p x^q = -\frac{1}{3}R_{ipjq}x^p x^q.$$

Another way of looking at it (a different choice of notation, really) is to take $p \in M$ and $v \in T_pM$. The map

$$\sigma_{p,v}\left(-T_{-}(p,v), T(p,v)\right) \longrightarrow M$$

is the maximal solution to

$$\nabla_{\dot{\sigma}}\dot{\sigma} = 0.$$

The initial conditions in this case are

$$\sigma_{p,v}(0) = p,$$

$$\dot{\sigma}_{p,v}(0) = v.$$

A unique solution exists by standard theorems on solutions for ODEs. We also have

[152] Reto Müller, *Differential Harnack Inequalities and the Ricci Flow* (Germany: European Mathematical Society, 2004).

$$\sigma_{p,\lambda v}(t) = \sigma_{p,v}(\lambda t), \qquad \forall \lambda \in \mathbb{R}, \forall t \in \left(-\frac{T_-}{\lambda}, \frac{T}{\lambda}\right).$$

$\exp_{\mathrm{p}}(v)$ is diffeomorphic with respect to both $p$ and $v$ by smooth dependence of solutions to ODEs on initial conditions.

To prove that the exponential map is a diffeomorphism onto a neighbourhood of $p \in M$, we need to calculate the derivative at the zero vector of the tangent space. This derivative can be computed by considering a curve inside the tangent space which starts at the origin. Start with a curve in $T_pM$ given by

$$\xi(t) = tw,$$

where $w \in T_0(T_pM) \cong T_pM$.

$$\begin{aligned}\left(d\exp_p(0)\right)(w) &= \frac{d}{dt}\exp_p\big(\xi(t)\big),\\ &= \frac{d}{dt}\sigma_{p,\xi(t)}(1),\\ &= \frac{d}{dt}\sigma_{p,w}(t).\end{aligned}$$

Everything on the right hand side is evaluated for $t = 0$. Define a new exponential map acting on pairs as follows:

$$\mathrm{Exp}(p, v) \coloneqq \left(p, \exp_p(v)\right).$$

If we evaluate the derivative at zero for a pair of vectors we have

$$\big(d\mathrm{Exp}(p,0)\big)(v,w) = \big(d\mathrm{Exp}(p,0)\big)(v,0) + \big(d\mathrm{Exp}(p,0)\big)(0,w).$$

Some computations and considerations of the initial conditions above for a curve in the manifold show us that:

$$\begin{aligned}\big(d\mathrm{Exp}(p,0)\big)(v,0) &= \frac{d}{dt}\mathrm{Exp}(\gamma(t),0),\\ &= \frac{d}{dt}\big(\gamma(t),\gamma(t)\big),\\ &= (v,v).\end{aligned}$$

Again, the right-hand sides of the first two equations are evaluated for $t = 0$. Also we have that

$$\begin{aligned}\big(d\mathrm{Exp}(p,0)\big)(0,w) &= (0, \left(d\exp_p(0)\right)(w),\\ &= (0,w).\end{aligned}$$

This implies that

$$d\mathrm{Exp}(p,0) = (0,w).$$

By the inverse function theorem, there exists a neighbourhood of the zero section of $TM$ such that

$$\text{Exp}\colon \mathcal{U} \longrightarrow \text{V}$$

is a diffeomorphism, where V is a neighbourhood of $M \times M$. The inverse function theorem states that if $M$ and $N$ are smooth manifolds, $p \in M$, and $F$ a smooth map between the two manifolds such that the pushforward $F_*$ is a bijection, then there exist connected neighbourhoods $U_0$ of $p$ and $V_0$ of $F(p)$ such that the restriction of $F$ from $U_0$ to $V_0$ is a diffeomorphism. Obviously, these criteria are met in this case if you consider what is being mapped where. The version for maps between open subsets of Euclidean spaces is similar. Essentially, once you have geodesics you can describe the exponential map[153].

One can often use the desirable properties of the exponential map to obtain proofs for other lemmas and theorems. For example, in a proof of a weaker version of Perelman's no local volume collapse theorem, one sees the assumption that if $\omega_n$ is the volume of the unit ball in Euclidean $n$-space and $V(p,r)$ is the volume of the geodesic ball at a point $p$ with radius $r$ given a complete Riemannian manifold, then the volume ratio always tends to the volume of the Euclidean unit ball[154]. In equation form:

$$K(p,r) = \frac{V(p,r)}{r^n} \longrightarrow \omega_n.$$

One way of seeing that is true is to do a short computation with the exponential map. Start with a unit ball in the tangent space at a point and then (assuming the radius $r$ to be small), you can map that ball down to the manifold $M$ with the exponential map (call the particular map we define $\varphi_r$). This is what the map does to a tangent vector.

$$v \longmapsto \exp_p(rv).$$

We said that $r$ was small, so one can define a metric on the ball via the pullback of the metric on the manifold via the exponential map. This implies that the volume for the geodesic ball is equal to the volume of the unit ball in the tangent space with respect to the pullback metric, but then the volume ratio must now be equal to the volume of that unit ball with respect to a modified metric which is

$$h_r = \frac{g_r}{r^{\frac{n2}{n}}} = \frac{g_r}{r^2}.$$

The modification of the metric follows from the fact that multiplying a Riemannian metric by a function $f$ multiplies the volume form by $f^{n/2}$.

Now just use well-known properties of the exponential map which can be looked up in any textbook of differential geometry and evaluate for a pair of vectors and a point $q$ in the unit ball[155].

$$(g_r)_q(v,w) = g_{\exp_p(rq)}\left(d\left(\exp_{\text{p}}\right)_{rq}(rv), d\left(\exp_{\text{p}}\right)_{rq}(rw)\right),$$

[153] Jürgen Jost, *Riemannian Geometry and Geometric Analysis* (Heidelberg: Springer, 2008).
[154] Peter Topping, *Lectures on the Ricci Flow* (Cambridge: Cambridge University Press, 2006).
[155] Michael Spivak, *Comprehensive Introduction to Differential Geometry* (USA: Publish or Perish, 1979).

$$= r^2 g_{\exp_p(rq)}\left(d\left(\exp_{\mathrm{p}}\right)_{rq}(v), d\left(\exp_{\mathrm{p}}\right)_{rq}(w)\right),$$

where we pulled the factors out of the bilinear form. This immediately implies that

$$(h_r)_q(v,w) = g_{\exp_p(rq)}\left(d\left(\exp_{\mathrm{p}}\right)_{rq}(v), d\left(\exp_{\mathrm{p}}\right)_{rq}(w)\right).$$

Take limits.

$$\lim_{r\to 0}(h_r)_q(v,w) = g_{\exp_p(0)}\left(d\left(\exp_{\mathrm{p}}\right)_0(v), d\left(\exp_{\mathrm{p}}\right)_0(w)\right),$$
$$= g_p(u,v).$$

In the limit, the metric associated to the volume form converges to the Euclidean metric, but this implies that the volume form converges to the Euclidean volume form. The diffeomorphic properties of the exponential map are also used to prove the lemma which we mentioned earlier in connection with the Darboux theorem. The main result then follows by applying the lemma in the case where $Q$ is a single point. The Darboux theorem is a local result, so taking it and attempting to smear it over the entire manifold will not work because of topological obstructions. This proof of the Darboux theorem was originally suggested by Moser.

In symplectic topology, we study the global structure of a symplectic manifold and the non-local behaviour of symplectomorphisms which are far from the identity. A number of new constructions have appeared in the development of these global studies. For example, Lalonde and McDuff discovered a method for constructing symplectic embeddings of balls, with the method being use to prove the non-squeezing theorem for arbitrary symplectic manifolds. This theorem states that there is no symplectic embedding which takes a standard $(2n+2)$-ball with unit radius into a cylinder $(M \times D^2(a), \omega \oplus \sigma)$ whose base space $D^2(a)$ is a closed disc with $\sigma$-area $a < \pi$[156]. This is one of the most important theorems in symplectic geometry and was first proved by Gromov[157]. Another key question in the transition from local to global is the question of the existence of periodic orbits for Hamiltonian flows. Start with $c \in \mathbb{R}$ a regular value of a smooth Hamiltonian $H: \mathbb{R}^{2n} \to \mathbb{R}$. If $\Phi$ is a smooth map between two manifolds $M$ and $N$, a point $p$ in $M$ is said to be a regular point of $\Phi$ if the associated pushforward

$$\Phi_*: T_pM \longrightarrow T_{\Phi(p)}N$$

is surjective. It is a critical point if this not the case.

A point in the other manifold $N$ is said to be a regular value of $\Phi$ if every point of the set $\Phi^{-1}(c)$ is a regular point. $\Phi^{-1}(c)$ is known as a level set. For technical reasons, it is common to think of a submanifold as a level set of some smooth map, although one needs a criterion to distinguish from level sets which are not embedded submanifolds. For example, if you happen to have a level set which is a curve containing a cusp, then this curve cannot be an embedded submanifold in the plane. In the case where we map to $\mathbb{R}^k$, the level set $\Phi^{-1}(0)$ is called the zero set of the map $\Phi$. (A level sets approach was crucial in the Huisken-Ilmanen proof of the Riemannian Penrose inequality which we mentioned earlier). There are various theorems relating to level sets: for example, the regular level set theorem says that every regular level set of a smooth map is a closed embedded submanifold whose codimension is equal to the dimension of the range. This is a corollary of the constant-rank level set theorem.

---

[156] François Lalonde and Dusa McDuff, 'Local non-squeezing theorems and stability', *Geometric and Functional Analysis* 5(2), 364 – 386 (1995).

[157] Mikhail Gromov, 'Pseudo-holomorphic curves in symplectic manifolds', *Invent. Math.* 82, 307 – 347 (1985).

Let us suppose given by the level set of the smooth Hamiltonian $S = H^{-1}(c)$ is compact and non-empty. It follows that $S$ is invariant under the flow of the Hamilton equations which are familiar from classical mechanics:

$$\dot{x} = \frac{\partial H}{\partial y}, \qquad \dot{y} = -\frac{\partial H}{\partial x}.$$

We already said earlier that at every regular point of $H$ the Hamiltonian vector field is tangent to the level sets of $H$. Furthermore, the vector field $XH$ has no zeroes on the surface $S$. One simple thing which we would like to know is whether or not the flow has a periodic orbit. In some cases it is not hard to see what the periodic orbits would be (we already mentioned characteristics on the sphere, which are all closed circular, and so periodic). Similarly, for an ellipsoid viewed as a hypersurface of $\mathbb{R}^{2n}$, there are no less than $n$ distinct closed characteristics. In this situation, a characteristic is a leaf of a foliation of the surface $S$ called the characteristic foliation. This in turn is the foliation determined by the integral curves of a 1-dimesional tangent distribution constructed by choosing a subspace

$$L_z = \{J_0 v \text{ such that } v \perp T_z S\} \subset T_z S$$

of the tangent space $T_z S$ for every $z \in S$. If $H$ is a Hamiltonian with $S$ the associated regular level set, then one can see that

$$XH \in L_z$$

for every $z \in S$, so the solutions of the Hamilton equations above are the characteristics, modulo reparametrisations of the time.

The Weinstein conjecture relates to a particular type of hypersurface called a contact hypersurface (the hypersurface must also be compact). The specific statement of the conjecture is that the Reeb vector field of such a hypersurface must have at least one periodic orbit. Contact geometry has been described as the version of symplectic geometry for odd dimension. Start with a manifold of dimension $2n+1$ and $\xi \subset TM$ a transversally orientable hyperplane field. Take a 1-form $\alpha$ such that $\xi = \ker \alpha$, then $d\alpha$ is nondegenerate on $\xi$ if and only if

$$\alpha \wedge (d\alpha)^n \neq 0.$$

$\xi$ is a contact structure and $\alpha$ a contact form for $\xi$. Informally speaking, a contact structure can be viewed as a field of hyperplanes which could not be further from integrability. One can also view the contact structure as an equivalence class of forms partitioned up by an equivalence relation given by the equation above: explicitly, the equivalence relation is that the form $\alpha$ is equivalent to $\alpha'$ if and only if

$$\alpha' = f\alpha,$$

where $\alpha'$ and $\alpha$ are 1-forms with

$$\xi = \ker\alpha = \ker\alpha'.$$

This implies that specifying a contact form is equivalent to specify some positive non-zero function $f$, but a function on a symplectic manifold generates a flow. Given a contact form $\alpha$, there is a unique vector field $Y_\alpha: M \to TM$ such that

$$\iota(Y)d\alpha = 0, \qquad \alpha(Y) = 1.$$

We call this vector field the Reeb vector field determined by the contact form $\alpha$. An example of a contact structure would be the standard contact structure on the odd-dimensional Euclidean space $\mathbb{R}^{2n+1}$ given by the 1-form:

$$\alpha_0 = dz - \sum_j y_j dx_j.$$

One can specify another contact form on $\mathbb{R}^{2n+1}$ via

$$\alpha_1 = dz + \frac{1}{2}\sum_j (x_j dy_j - y_j dx_j.$$

In this case, the two contact forms are diffeomorphic, but this will not always be the case. Not all of the diffeomorphisms of $S^n$ preserve contact forms, for example, even up to a multiplying factor.

Hypersurfaces of contact type are characterised by the existence of a vector field $X$ close to the level surface of the Hamiltonian such that

$$\mathcal{L}_X\omega_0 = \omega_0.$$

Weinstein proposed his conjecture for all hypersurfaces of contact type and it was subsequently proved by Viterbo that every hypersurface of contact type in $\mathbb{R}^{2n+1}$ has a closed characteristic. The general conjecture for contact manifolds which are not embedded contact hypersurfaces in Euclidean space is still open, however. The conjecture was proved for all closed 3-dimensional manifolds by Taubes using a version of Seiberg-Witten Floer homology[158]. Another basic problem in symplectic topology which marks the transition from local to global is that of rigidity. Suppose we take a sequence of symplectomorphisms which converge on a compact set to a diffeomorphism. This limit preserves volume, but it is not clear that the symplectric structure is preserved, and if this structure is not preserved, then we can no longer guarantee that we are even studying symplectic topology. Put another way, we would like to know the $C^0$-closure of the group of symplectic diffeomorphisms. It turns out that the key to characterizing symplectomorphisms is to use the non-squeezing theorem to define a symplectic invariant which is continuous with respect to the $C^0$ topology. Symplectic diffeomorphisms are then simply the diffeomorphisms which preserve this symplectic invariant, hence we see that the four basic strands of the subject influence each other.

---

[158] Clifford Taubes, 'The Seiberg-Witten equations and the Weinstein conjecture', *Geometry and Topology*, 11(4), 2117 – 2202 (2007).

The final basic problem in symplectic topology is that of counting the fixed points of symplectic diffeomorphisms. Again, there is a conjecture which encapsulates this problem, known as Arnold's conjecture. The conjecture relates the minimum number of fixed points of a Hamiltonian symplectomorphism on a closed manifold to Morse theory[159]. We will not go into Morse theory here, but the starting point is to take a compact Riemannian manifold along with a function $F \in C^\infty(M)$. Morse theory begins by studying the behaviour of $F$ near every critical point (assuming that there is a finite number of these points) and studies the possible topological changes which the level sets of $F$ can undergo[160]. The minimum number of critical points for a function defined on a compact manifold is actually a topological invariant of the manifold. This minimum number obviously has to be at least 2 from elementary calculus, because any function will have a distinct local maximum and local minimum. One way of phrasing the conjecture is that if $\phi$ is the time-1 map of a time-dependent Hamiltonian flow on a compact symplectic manifold, then $\phi$ has a number of distinct fixed points which is at least the minimum number of critical points on a compact manifold, taken over all smooth functions $f$. The Arnold conjecture can be viewed as a global version of the local statement that a function defined on a manifold has to have more critical points then there are zeroes for a vector field or convector field defined on that manifold. The exact details of the interplay between local symplectic structure and global properties (and how the former manifests itself in the latter) is at the core of symplectic topology[161].

McDuff has made many contributions to this still rapidly growing subject. McDuff was the first person to construct simply connected, closed symplectic manifolds which are not Kähler manifolds. An example she provided is as follows. Take Thurston's symplectic Kähler 4-manifold $(M, \omega)$. This manifold is a quotient $\mathbb{R}^4/\Gamma$, where $\Gamma$ is the discrete affine group generated by the unit translations along the $x_1, x_2, x_3$-axes along with a transformation

$$(x_1, x_2, x_3, x_4) \longmapsto (x_1 + x_2, x_2, x_3, 1 + x_4).$$

$M$ is a $T^2$-bundle over $T^2$ and the symplectic form $\omega$ lifts to a form on $\mathbb{R}^4$:

$$\omega_1 = dx_1 \wedge dx_2 + dx_3 \wedge dx_4.$$

[159] Dusa McDuff and Dietmar Salamon, *Introduction to Symplectic Topology* (New York: Oxford University Press, 1998).
[160] John Milnor, *Morse Theory* (Princeton: Princeton University Press, 1963).
[161] Dusa McDuff and Dietmar Salamon, *Introduction to Symplectic Topology* (New York: Oxford University Press, 1998).

The first and third Betti numbers of the manifold are odd, which implies that $M$ has no Kähler structure[162]. Recall that the alternating sum of the Betti numbers is equal to the Euler characteristic.

$$\chi(U) = \sum_{p=0}^{n} (-1)^p \, \beta_p(U).$$

In the context of almost complex structures, McDuff proved that for $(M, \omega)$ a compact connected symplectic 4-manifold which contains a symplectically embedded 2-sphere $S$ with an intersection number defined by the positive definite inner product $S \cdot S$, we have that $(M, \omega)$ is the blow-up of a rational or ruled symplectic 4-manifold. If $(M, \omega)$ is minimal, then it is either the complex projective plane or an $S^2$-bundle over a Riemannian surface. The fibres of the ruling may also be assumed to be symplectic[163]. In a paper on circle actions, McDuff proved that a symplectic circle on a closed 4-manifold is Hamiltonian if and only if it has fixed points, also providing an example of a symplectic, non-Hamiltonian circle action on a compact symplectic 6-manifold with fixed points[164]. There are also various results related to blow-ups and blow-downs. For example, the fact that in dimension 4, the blow-ups $(\widetilde{M}, \widetilde{\omega}_\psi)$ and $(\widetilde{M}, \widetilde{\omega}_\phi)$ are isotopic if and only if the normalized symplectic embeddings $\psi$ and $\phi$ are isotopic[165].

---

[162] Dusa McDuff, 'Examples of simply connected symplectic non-Kählerian manifolds', *Journal of Differential Geometry*, 20, 267 – 277 (1984).

[163] Dusa McDuff, 'Rational and ruled symplectic 4-manifolds', *Journal of the American Mathematical Society*, 3, 679 – 712 (1990).

[164] Dusa McDuff, 'The moment map for circle actions on symplectic manifolds', *Journal of Geometrical Physics*, 5, 149 – 160 (1988).

[165] Dusa McDuff, 'Remarks on the uniqueness of symplectic blowing-up'. In *Symplectic Geometry* (ed. D. Salamon), 157 – 168. (Cambridge: Cambridge University Press, 1993).

# *Karen Vogtmann*

Vogtmann is a mathematician who works primarily in geometric group theory. She is known for the introduction of Culler-Vogtmann Outer space, often shortened simply to Outer space. The idea that group theory can be applied to geometry is fairly embedded at this point, both in geometry and in mathematical physics. This began with Klein's so-called Erlangen programme in the nineteenth century, which planned to use group theoretic methods to study geometry. As an example, the symmetries of spherical geometry in one dimension form a group $O(2)$ formed of rotation and reflection matrices:

$$\begin{bmatrix} \cos\theta & -\sin\theta \\ \sin\theta & \cos\theta \end{bmatrix}, \begin{bmatrix} \cos\theta & \sin\theta \\ \sin\theta & -\cos\theta \end{bmatrix}$$

Is it possible to do the reverse and apply geometric methods in order to prove things and study objects in group theory? It is, and this is what is known as geometric group theory. Gromov is usually credited with the introduction of geometric group theory as a distinct field, although there are many isolated examples in the twentieth and even the nineteenth century where mathematicians used methods which could loosely be labelled as 'geometric group methods': for example, when Dehn applied hyperbolic geometry backwards to solve the word problem for a surface group[166].

The basic tool which one starts with is a Cayley graph. This is a graph in every sense of the word (a geometric object), but it also encodes information about a finitely generated group and so allows us to view a finitely generated group as a geometric object. A group $G$ is finitely generated if it has a finite generating set $A$ and a group is generated by a finite subset within it if $G = \langle A \rangle$, where $\langle A \rangle$ is the intersection of all the subgroups of $G$ which contain $A$. It is possible to view a graph $\Gamma$ as a metric space, where the graph is thought of as a set of vertices $V$ along with a set of edges $E$ (by metric space, I just mean that there is a notion of distance, roughly speaking). The way of doing this is to realise $K$ as a 1-complex, by viewing each edge of the graph as a copy of the unit interval where the endpoints are vertices. You can think of a cell complex as a space which is constructed by beginning with a set of points (the vertices, in this case) and then attaching cells of increasing dimension, where each cell is a space which is homeomorphic to a ball in Euclidean space. In our case, we stop with cells of dimension 1 when we construct our spaces, as we are just adding edges to vertices[167].

The 1-complex can now be viewed as a set, on which we define a notion of distance by parametrising each edge and so defining length for every interval in an edge. This means that a path (that is, a combinatorial sequence of edges or vertices) has a well-defined finite, non-negative length $l$. Given two points on the graph, we can define the distance between them to be the smallest length for any path which connects them together. The set of lengths of all such paths is a discrete set, so this minimum has to be achieved at some point.

---

[166] Brian Bowditch, *A Course on Geometric Group Theory* (Japan: Mathematical Society of Japan, 2006).
[167] John M. Lee, *Introduction to Topological Manifolds* (New York: Springer, 2011).

If $K$ is a connected graph, then it follows that the set derived from the realisation of $K$ as a 1-complex equipped with the metric which we have defined together make up a metric space, which induces a topology on the set if we assume that the graph is locally finite. However, one can also associate a graph to a group. Start with a group $G$ and a subset of that group $S$: we are looking to construct a graph $\Delta(G; S)$. Let the set of vertices of this graph $V(\Delta)$ be equal to $G$. For every $g$ in the group and every $a$ in the subset, there is a directed edge from $g$ to the composition of elements $ga$ which can be labelled by $a$ (or $a^{-1}$ if you go in the opposite direction). The group $G$ acts on the vertex set by left multiplication, which extends to a free action on the graph $\Delta$. It is also possible to start with an element of the group and get a path from $g$ to the composition $gp(w)$, where $p$ maps a word from the set of words in $A$ to the group $G$. A word is simply a finite sequence of elements of $A$ (think of the way that you can take an alphabet for a language and form any word by making a sequence from some of the elements in that alphabet).

This path is the image of a path from $1$ to $p(w)$ acted on by the element $g$ in the group action we have already specified. Since $\langle S \rangle$ is the set of elements which can be expressed as a word in the alphabet $A$, it follows by considering the above construction that the group $G$ equals $\langle S \rangle$ if and only if the graph $\Delta(G; S)$ is connected. It means that a finitely generated group acts freely on a locally connected finite graph. In particular, if $S$ is a generating set for $G$, then $\Delta$ is called the Cayley graph of the group with respect to that generating set. As a simple example, the integers $\mathbb{Z}$ can be generated by one single generator such that $\mathbb{Z} = \langle a \rangle$. In this case, the Cayley graph is just the line of real numbers. $\mathbb{Z}$ acts freely on the Cayley graph via translation. The quotient graph for the action of the group on the set of vertices always consists of a single vertex and a number of loops equal to the number of generators in the generating set, so in this case the quotient graph is one vertex and one loop ie. a circle[168].

The main group which we are interested in modelling and studying is the free group with $n$ generators and no relations:

$$F_n = \langle a_1, \ldots, a_n \rangle.$$

We then consider the group of automorphisms of the free group $\mathrm{Aut}(F_n)$, where a group automorphism is just an isomorphism which takes a group 'back to itself'. In linear algebra, the corresponding linear transformation is an endomorphism or 'linear self-map' from a vector space back to that same vector space. In the case of automorphisms, the transformation must be invertible, but this is clear from the fact that the automorphism group must have an inverse to be defined as a group. An example of an automorphism of the free group might be inversion, where $a_i$ is taken to $a_i^{-1}$ and $a_j$ is left untouched, such that $a_i \neq a_j$. Another example might be the map which multiplies the $i$-th generator by the $j$-th generator from the right or the left and leaves the other group elements untouched: this is an isomorphism, as one can write down the inverse. (We can multiply from the left or the right as we are in a free group). The automorphisms of the free group can be divided further into inner and outer automorphisms. The inner automorphisms of the free group

[168] Brian Bowditch, *A Course on Geometric Group Theory* (Japan: Mathematical Society of Japan, 2006).

$\mathrm{Inn}(F_n)$ form a normal subgroup of the automorphism group, where an inner automorphism is a group automorphism given by conjugation of the group by one of its elements. More importantly, we can also form the group of outer automorphisms by taking the quotient of the automorphism group with the group of inner automorphisms.

$$\mathrm{Out}(F_n) = {\mathrm{Aut}(F_n)}\big/{\mathrm{Inn}(F_n)}.$$

The groups form what is called a short exact sequence.

$$1 \to \mathrm{Inn}(F_n) \to \mathrm{Aut}(F_n) \to \mathrm{Out}(F_n) \to 1.$$

Bear in mind that the group of inner automorphisms of the free group is isomorphic to the free group itself. It is also necessary to study a group which you can think of as combining elements from the free group and its automorphism group. For any $k$, we obtain a group $A_{n,k+1}$

$$A_{n,k+1} = \{((w_1, .., w_k), g) \text{ such that } w_i \in F_n, g \in \mathrm{Aut}(F_n)\}.$$

This allows us to form another short exact sequence:

$$1 \to F_n \to A_{n,2} \to \mathrm{Aut}(F_n) \to 1.$$

The goal in general is to try to learn more about $\mathrm{Aut}(F_n)$ and $\mathrm{Out}(F_n)$.

Generally speaking, less is known about automorphism groups than other groups and you can learn a lot about a mathematical object by studying its automorphisms, whatever they might be. Automorphism groups are usually complicated and can be studied at several levels. For example, to study the group of automorphisms $G$ of a Riemann surface $X$ of genus $g$ requires information about the character $\chi$ of the action of $G$ on the space of holomorphic abelian differentials on $X$, the signature of the Fuchsian group $\Gamma$ of the surface, and a surface kernel epimorphism $\Phi: \Gamma \to G$. Recall that a Riemann surface is just a connected topological 2-manifold equipped with a smooth structure such that the relevant transition maps between charts are holomorphic, conformal homeomorphisms. The different types of information are related: for example, the numbers $|\mathrm{Fix}_{X,u}(h)|$ determine the character of the Riemann surface by the Eichler trace formula. Under some mild assumptions on the relevant groups and taking a compact Riemann surface with genus $g \geq 2$, where $h$ in $G^{\mathrm{X}}$ is of order $m$, these numbers are defined as

$$|\mathrm{Fix}_{X,u}(h)| = |C_G(h)| \sum_{1 \leq i \leq r,\ \ m | m_i, h \sim \Phi(c_i)^{m_i u/m}} \frac{1}{m_i}.$$

The Eichler trace formula then says that if $\sigma$ is an automorphism of order greater than 1 of a compact Riemann surface with genus $g \geq 2$ and $\chi$ the character of the action of the automoprhism group $\mathrm{Aut}(X)$ on the space of holomorphic differentials on $X$, then the following holds[169]:

---

[169] Thomas Breuer, *Characters and Automorphism Groups of Compact Riemann Surfaces* (Cambridge: Cambridge University Press, 2000).

$$\chi(\sigma) = 1 + \sum_{u \in I(m)} \left|\mathrm{Fix}_{X,u}(\sigma)\right| \frac{\zeta_m^u}{1-\zeta_m^u}.$$

An important geometric model for $\mathrm{Out}(F_n)$ introduced by Culler and Vogtmann is known as Outer space[170]. The development of Outer space was inspired partly by Teichmüller spaces which turn up in the theory of surfaces. They also have some applications in theoretical physics, appearing in string theory scattering and in some superconformal field theories, where they are associated with spaces of exactly marginal deformations. In more technical terms, it is analogous to the Teichmüller space associated with the mapping class group of a surface. The mapping class group of a topological surface $F$ is the quotient group

$$\frac{\mathrm{Out}_t^+\big(\pi_1(F)\big)}{\mathrm{Inn}\big(\pi_1(F)\big)},$$

where $\pi_1(F)$ is the fundamental group for the surface. The subgroup which fixes the parabolic conjugacy classes of the group and fixes the punctures in the surface is known as the pure mapping class group. A moduli space of Riemann surfaces can be formed by quotienting the Teichmüller space through by the associated pure mapping class group

$$\mathcal{M}_{g,n} = \frac{\mathcal{T}_{g,n}}{PMCG},$$

where $g$ is the genus of the surface and $n$ is the number of punctures. An important subset of this moduli space is the set of Riemann surfaces such that the length of the shortest closed geodesic is at least some positive number $\epsilon$: this subset is denoted by $\mathcal{M}_{g,n}(\epsilon)$[171]. Mumford's compactness theorem states that the space of Riemann surfaces $\mathcal{M}_{g,n}(\epsilon)$ with length of closed geodesics at least $\epsilon$ is compact. Furthermore, a sequence of representations with minimum lengths of closed geodesics bounded from below has a convergent subsequence of $PMCG$-translates[172]. From Mumford compactness and so-called thick-thin decomposition, it follows that studying the compactification of the moduli space $\mathcal{M}_{g,n}$ is equivalent to studying the way that lengths of geodesics tend to zero.

The mapping class group acts properly discontinuously on $\mathcal{T}_{g,n}$, and by analogy Outer space is a contractible space on which the group of outer automorphisms of the free group $\mathrm{Out}(F_n)$ acts properly discontinuously. An action of a group $\Gamma$ on a proper length space $X$ is said to be properly discontinuous if for all $r \geq 0$ and all $x \in X$, the set of group elements $g$ such that $d(x, gx) \leq r$ is finite. The definition can also be stated without a metric by saying that it is the set of group elements $g$ such that $gK \cap K \neq \emptyset$ is finite for all compact $K$ (since a compact set can only meet finitely many subsets of itself). As mentioned previously, the natural metric on a finite graph is as an assignment of lengths to edges. A marking of a graph is a homotopy equivalence $f$ which maps the rose $R_n$ (the graph with 1 vertex and $n$ edges) to the graph in question. This allows for an identification between the fundamental group of the graph $\pi_1(\Gamma)$ and the free group $F_n$ and two markings $f$ and $f'$ are equivalent are equivalent if that is a homeomorphism $\phi$ mapping $\Gamma$ to $\Gamma'$ such that there $\phi f$ is homotopic to $f'$. (If you are not familiar with some of these topological definitions, just think of a homotopy as taking a function and deforming it continuously into another one). Consider two triples $(\Gamma, l, f)$ and $(\Gamma', l', f')$ such that the valence of every vertex of the graphs is at least 3 and $l$ is a metric with unit volume, an isometry $\phi$ from $\Gamma$ to $\Gamma'$ such that $\phi f$ is homotopic to $f'$ is an equivalence. This allows for an alternative definition for Outer space

$$\chi_n = \frac{\{(\Gamma, l, f)\}}{\sim}$$

as the set of equivalence classes of finite marked metric graphs with vertices of valence at least 3 and unit volume. It is more common to leave out the equivalence class notation and to simply think of Outer space as a space of marked metric graphs[173]. Probably the most succinct definition is to define Outer space as the space of free minimal actions of the free group by isometries of metric simplicial trees.

[170] Marc Culler and Karen Vogtmann, 'Moduli of graphs and automorphisms of free groups', *Invent. Math.* 86, no. 1, 91 – 119 (1987).
[171] Scott Wolpert, *Families of Riemann Surfaces and Weil-Petersson Geometry* (USA: American Mathematical Society, 2009).
[172] David Mumford, 'A remark on Mahler's compactness theorem', *Proceedings of the American Mathematical Society*. 28, 289 – 294 (1971).

Outer space $\chi_n$ has its own rich geometry and topology and $\chi_n$ decomposes into open simplices. If $\Gamma$ is a graph and $f$ is a marking, the set of possible metrics on $\Gamma$ is an open simplex:

$$\left\{(l_1, \dots, l_E) : \; l_i > 0, \sum_i l_i = 1\right\}.$$

The dimension of the simplex is one below the number of edges $E$. If $T$ is a forest (a disjoint union of trees) in the graph and another graph $\Gamma'$ is obtained by identifying all the edges of $T$ with points, then the set of metrics on $\Gamma'$ can be identified with the open face of the set of metrics on $\Gamma$ for which the coordinates of edges in $T$ are $0$. This can be stated equivalently by saying that one has collapsed a forest to get $\Gamma'$ from the graph $\Gamma$, or that one has blown up a forest to get $\Gamma$ from the graph $\Gamma'$. The union of the set of metrics with all of these open faces as $T$ cycles through all the possible forests in $\Gamma$ is a simplex with missing faces which can be obtained from the corresponding closed simplex by removing the appropriate open faces: this union is denoted by $\Sigma(n)$, and its smallest possible dimension is always $n - 1$.

In topological terms, $\chi_n$ is a simplicial complex of simplices with missing faces and a simplicial topology can be defined on the space analogous to that of any other simplicial complex. If one restores the missing faces, one gets a simplicial complex which is isomorphic to a sphere complex. Without going into too many details, a sphere complex $S_{n,k}$ is a simplicial complex with vertices which are isotopy classes of spheres and $k$-simplices which are compatible collections of $k + 1$ isotopy classes $[\Delta_0], \dots, [\Delta_k]$. A compatible collection has distinct isotopy classes and a sphere system has isotopy classes

[173] Mladen Bestvina, 'Geometry of Outer space', in *Geometric Group Theory*, ed. by Mladen Bestvina, Michah Sageev and Karen Vogtmann. (USA: American Mathematical Society, 2014).

$[s_0], \ldots, [s_k]$. Outer space is a subspace of a sphere complex $S(M_{n,s})$ consisting of open simplices $\sigma(s_0, \ldots, s_k)$, where the $s_i$ form a complete sphere system and $M_{n,s}$ is defined as

$$M_{n,s} = \#S^1 \times S^1 \setminus \coprod_s B^3.$$

Another basic fact about Outer space is that it is contractible. One can prove this using sphere systems. The idea of the proof is to fix a maximal sphere system and then construct paths from an arbitrary point to a simplex $\Sigma$ in the sphere complex and so contract $\chi_n$ along those paths[174]. During the proof, one can employ Hatcher's normal form theorem. This states that a sphere system $S$ is isotopic to a system which intersects $\Sigma$ transversally in circles and that in each piece $p$ of $M_{n,s}$ (technically a 3-punctured sphere), this intersection is always a union of discs, cylinders and pairs of pants such that different boundary components are joined by cylinders or such that there is a boundary piece of each of the 3-spheres. This in turn can be proved using Laudenbach's lightbulb trick: homotopic embedded spheres are isotopic[175].

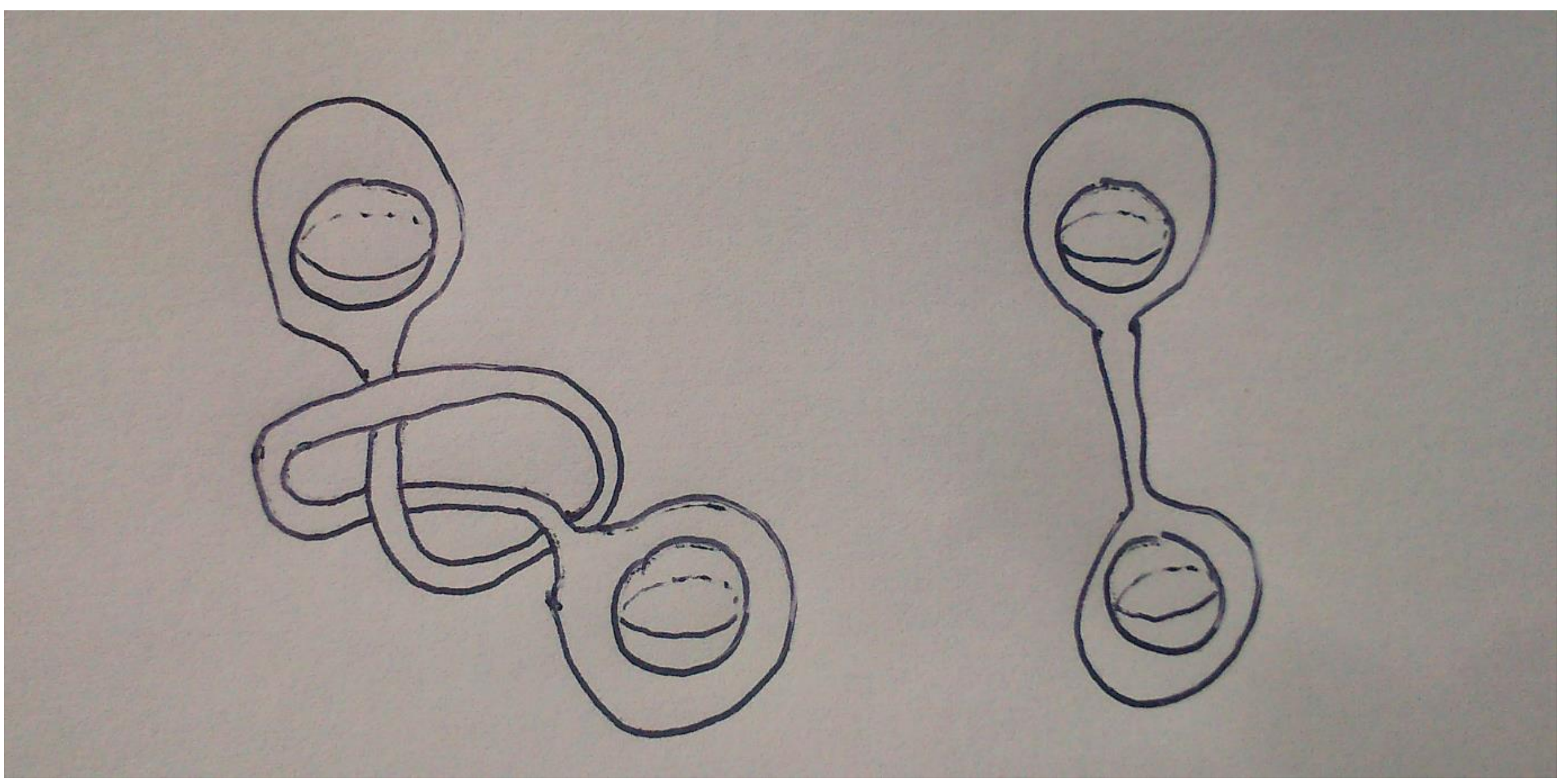

Figure 4: Visual demonstration of Laudenbach's lightbulb trick.

Finally, we will mention that there is a natural right action of $\mathrm{Out}(F_n)$ on $\chi_n$ which is obtained by precomposition of the marking. An element $\Phi$ of $\mathrm{Out}(F_n)$ can be considered as a homotopy equivalence from the $n$-rose to the $n$-rose. In terms of the equivalence classes we defined before, the action is then

$$[(\Gamma, l, f)] \cdot \Phi = [(\Gamma, l, f\Phi)].$$

This action is proper and it acts properly and cocompactly (ie. the quotient by the action is compact) on the thick part and on the spine $K_{n,s}$, where the spine is a specific deformation retract of Outer space. In fact, $K_{n,s}$ is a contractible cube complex whose vertices are complete sphere systems in $M_{n,s}$[176].

[174] Karen Vogtmann, 'Automorphisms of free groups and outer space', *Proceedings of the Conference on Geometric and Combinatorial Group Theory, Part I* (Haifa, 2000*)*. Vol. 94, pp. 1 – 31 (2002).
[175] Allen Hatcher, *Algebraic Topology* (Cambridge: Cambridge University Press, 2009).

We have mentioned surface groups and free groups. Another important class of groups are the right-angled Artin groups. As with surface groups and many other groups, the right-angled Artin groups have finite presentations. Whether a free group has a finite presentation depends on whether or not the group is finitely generated. If $\Gamma$ is a finite graph with a vertex set $V$ and an edge set $E$ consisting of edges, $\{[x_i, x_j]\}$, the right-angled Artin group is defined by the following finite presentation

$$A_\Gamma \coloneqq \langle V\colon [x_i, x_j], \;\; \text{when } [x_i, x_j] \in \mathrm{E} \,\rangle.$$

The first set of square brackets in the presentation is meant to indicate that all the relations are commutator relations, whereas the second set denotes an edge between two vertices. As an example, if the graph $\Gamma$ contains no edges, then the right-angled Artin group is the free group on $n$ generators[177]. The outer automorphisms of the right-angled Artin groups have also been studied, and it has been shown that for a connected, triangle-free graph $\Gamma$, $\mathrm{Out}(A_\Gamma)$ can be defined in terms of maximal join subgraphs in $\Gamma$. In fact, it is also possible to form the corresponding analogue of an Outer space for $\mathrm{Out}(A_\Gamma)$ and show that $\mathrm{Out}(A_\Gamma)$ acts properly on this outer space and that the outer space is contractible with a spine. Using the fact that the Tits alternative is known for outer automorphism groups of free groups, one can also show that it holds for $\mathrm{Out}(A_\Gamma)$[178]. The Tits alternative states that for a finitely generated linear group over a field, either the group contains a solvable subgroup of finite index, or it contains a non-abelian free group (given the name in a Russian article by Roman'kov)[179].

An intriguing recent possibility is that of applying Outer space to QFT, although these ideas are quite far from publication at the time of writing and I will not go into too much detail as this is primarily a book about mathematics, rather than physics. For example, it has been observed that the so-called Cutkosky rules can be related to Outer space[180]. The Cutkosky rules are rules for dealing with certain Feynman diagrams, so I will say a few words about QFT to set the scene. If we start from non-relativistic quantum mechanics (we might be considering an electron which is moving at speed well below light speed, for example) and it is being scattered by something, we would like to study the scattering amplitude of the outgoing wave compared to the incoming wave which is incident on the thing which is doing the scattering (everything is in terms of wave functions which solve the Schrodinger 'wave' equation, hence viewing things as waves). The amplitude could be imaginary, but you always square the modulus and so the probability ends up being a real number. The collision can be elastic or inelastic (kinetic energy may or may not be conserved), but what you can do is consider the scattering amplitude as a function of energy where energy is a complex variable. Since it is a function of a complex variable, it can have branch points in the complex plane.

---

[176] Mladen Bestvina, 'Geometry of Outer space', in *Geometric Group Theory*, ed. by Mladen Bestvina, Michah Sageev and Karen Vogtmann. (USA: American Mathematical Society, 2014).
[177] Cornelia Drutu and Michael Kapovich, *Geometric Group Theory* (USA: American Mathematical Society, 2018).
[178] Ruth Charney, John Crisp and Karen Vogtmann, 'Automorphisms of 2-dimensional right-angled Artin groups', *Geometry and Topology*. 11 2227 – 2264 (2007).
[179] Vitaly Roman'kov, 'Embedding some interlacings in groups of automorphisms of finitely generated solvable groups', *Algebra and Logic*. Vol. 15, Issue 3, 187 – 192 (1976).
[180] Dirk Kreimer, 'Cutkosky rules from Outer space', arXiv:1607.04861 (2016).

This scattering amplitude has a branch point at $E = 0$ which extends into a branch cut along the entire positive real axis, meaning that the imaginary part of the scattering amplitude is a discontinuity across this cut. This essentially means that the classical optical theorem can be stated for quantum mechanics as well as for classical physics. The optical theorem states that the imaginary part of the scattering amplitude for scattering straight forwards with angle $0$ is proportional to the cross section of the scattering object (not necessarily the same as its geometric area)[181]. The optical theorem can be extended to field theory as you might expect. The easiest way to prove the optical theorem is to use an operator called the $S$-matrix, or scattering matrix, which gives the amplitude for $n$ particles all with their own wave functions to scatter and turn into $m$ particles with their own new wave functions (defined in an asymptotic way by assuming they are coming in from a long way out and going back out a long way to infinity). The theorem is then a consequence of the fact that the $S$-matrix is a unitary operator and it can be stated equivalently as an identity for the imaginary part of an S-matrix element. In QFT, the number of particles going in does not have to equal the number going out because of Einstein's equation $E = mc^2$.

The question is, how does this imaginary part of an $S$-matrix element arise when you do a Feynman diagram expansion? If you have a quantum field theory with allows for interaction between fields (unlike a free field theory where things travel along and never 'see' each other), then you might be able to do a perturbation expansion. A field theory has to have a functional integral (ie. a Feynman path integral), so if you expand the interaction part as a Taylor series you end up with a number of terms which are basically the same and which can be collected together with a Feynman diagram[182]. Consider, for example, if we take a functional defined as

$$Z_1(J) \propto \exp\left(\frac{i}{6} Z_g g \int \mathrm{d}^4 x \, \left(\frac{1}{i}\frac{\delta}{\delta J(x)}\right)^3\right) Z_0(J),$$

where $Z_0(J)$ is the path integral for the corresponding free-field theory, then the Taylor expansion in powers of $g$ and $J$ is

$$Z_1(J) \propto \sum_{V=0}^{\infty} \frac{1}{V!}\left(\frac{iZ_g g}{6} \int \mathrm{d}^4 x \, \left(\frac{1}{i}\frac{\delta}{\delta J(x)}\right)^3\right)^V \sum_{P=0}^{\infty} \frac{1}{P!}\left(\frac{i}{2} \int \mathrm{d}^4 y \, \mathrm{d}^4 z \, J(y) \, \Delta(y - z) J(z)\right)^P .$$

Each term in this dual expansion can be seen to produce many expressions which are algebraically identical and which are organised with Feynman diagrams.

A tree level diagram with no loops corresponds to zero order, one loop is first order, and so on. The Feynman rules are rules for perturbation theory and the ordinary assumptions which are made for applications of perturbation theory are exactly the same ie. the parameter multiplying the non-quadratic part of the potential (a physical 'coupling constant' in this case) has to be sufficiently small. The coupling constant in quantum electrodynamics is a constant equal to $1/137$, so QED is very amenable to perturbation theory. As you might know already, if you go too far in the expansion, you will quickly find that the integrals always diverge and so the diagrams are infinite, hence the need for renormalization to somehow 'get rid' of the infinities and 'subtract them out of existence'.

---

[181] Lev Landau and Evgeny Lifschitz, *Quantum Mechanics: Non-Relativistic Theory* (Oxford: Butterworth-Heinemann:, 1981).
[182] Mark Srednicki, *Quantum Field Theory* (Cambridge: Cambridge University Press, 2016).

It turns out that a Feynman diagram only contributes an imaginary part to a matrix element when virtual particles are on-shell. On-shell means that they obey the relativistic energy-momentum relation from special relativity: since they are virtual and 'flitter' in and out of existence during a process such as a collision, they don't have to obey this relation. In mathematical terms, as you might guess, there will not be an imaginary part to the matrix element unless there is a branch cut. The branch cut is along the real axis as before but now it starts at a threshold energy $s_0$ (again, we are considering the matrix element as a function of energy as a complex variable, or as a function of momentum where there is a threshold momentum). It might not be obvious, but if you write out a Feynman diagram and then construct the expression for it, there could be a branch cut along the real axis (this might well be from a $\log$ term, which has a well-known branch cut). This means you can compute the imaginary part which a Feynman diagram contributes to a matrix element by computing the discontinuity of the Feynman diagram across its branch cut. What Cutkosky showed is that you can do this for any order and that you can always do it with a simple set of rules[183].

For the first step, you draw the Feynman diagram and cut through it in every way such that the cut propagators can be put on-shell at the same time (for example, if you have a diagram with a loop in the middle). Both sides of the loop contribute their own propagators, since a line corresponds to a propagator in a Feynman diagram (a propagator being the amplitude for a particle to travel from one place to another, but you can also view it as a Green's function if preferable). Each of those contributes a propagator term when you write out the expression for the diagram using the Feynman rules, so it's a question of whether both of the integrals can be written in such a way that the momenta in the expressions obey the on-shell condition: it might depend on the region of momentum integration as to whether or not this is possible.

In the second step, for each cut propagator you replace the propagator term with a Dirac delta function multiplied by $2\pi i$. This essentially comes from the Cauchy residue theorem, as you are somehow closing the contour of integration around the poles which are present and picking up the residues. You integrate over the delta functions for each cut. Finally, you sum all the contributions from every cut and you have computed the discontinuity of any Feynman diagram, and that in turn gives you its imaginary contribution to an S-matrix element, the S-matrix elements being the things with direct physical interpretation which can be measured. Because we can do this, we can prove the optical theorem in QFT to all orders and generalise the theorem to all Feynman diagrams, which means that you can use them to evaluate other amplitudes which can be defined in terms of Feynman diagrams, even if they are not $S$-matrix elements: this is significant for studying unstable particles. Unstable particles and 'resonances' are studied in non-relativistic quantum mechanics, but we would also like to be able to deal with them in field theory, especially in high-energy physics where most things you might discover in a collider experiment are unstable.

[183] Michael Peskin and Daniel Schroeder, *An Introduction to Quantum Field Theory* (USA: Westview Press, 1995).

Another related idea is that the analytic structure of Feynman amplitudes viewed as functions into the complex plane can be connected to Outer space[184]. One can study functions which select a Feynman diagram from a Hopf algebra of diagrams and evaluate them with renormalized Feynman rules to map them into $\mathbb{C}$, depending on the rules which in turn depend on the field theory. One then has a multi-valued amplitude whose properties can be studied using techniques from complex analysis and algebraic geometry, an approach going back to Landau or earlier. The evaluation of a Feynman diagram is a loop integral and it is possible to always regard this amplitude as an iterated integral built up from amplitudes for one-loop graphs such as tadpoles. It is not possible to define the function which evaluates a diagram with rules as a unique function if you define it as an iterated integral, but you can define an equivalence relation on principal sheets. The iteration of the integral then gives the function a structure which is closely linked to that of the correct coloured Outer space.

As a Feynman amplitude is a many-valued function with branch points, there exists a group generated by the branch cycles known as the 'monodromy' of the graph, in the sense that the monodromy group encodes the extent to which an amplitude fails to be single-valued as it goes around a branch point with non-trivial monodromy. The generators of the monodromy group for an amplitude can be found by computing the variation of the amplitude with Cutkosky's theorem. The wider importance of this is that the generators of the fundamental group for the graph (a simple object) can be mapped to the generators of the monodromy group of the amplitude, so, for example, a single loop like the triangle graph has one generator for its fundamental group, meaning that there is just one generator for the monodromy group, which then can be used to define the multi-valued functions which are assigned to the amplitude.

These different multi-valued functions are needed to study the sheet structure of the amplitude. The functions have representations in terms of coloured graphs (for example, there will be 3 coloured graphs in the case of 3 functions) and the functions are needed for the simplex or cell which is being integrated over. In the case of what is called the 3-edge banana graph $b_3$, the simplex is a triangle. As the corners of the triangle are not in Outer space, they must be changed to arcs, so there are 3 functions associated to the 3 arcs and then equivalence relations to change between the arcs[185]. For the example of $b_3$, one can take the 3 equivalent markings for the graph and each marking determines a multi-loop integral for the amplitude. However, there are two generators for the monodromy group of $b_3$, so essentially there is a connection between the number of markings for the appropriate subgraph (the two-edge banana graph $b_2$ which has 2 markings) and the generators of the monodromy group of the graph itself, with the choice of edges determining the integral. The ultimate aim is to get information on the sheet structure of the associated amplitude and to evaluate a parametric Feynman integral by assigning an open simplex from Outer space to the Feynman diagram (for example, a tetrahedron for the dunce's cap graph). One then computes the Feynman integral by integrating over the volume of that simplex.

---

[184] Dirk Kreimer, 'Multi-valued Feynman graphs and scattering theory', arXiv:1807:00288v3 (2018).
[185] Paolo Aluffi and Matilde Marcolli, 'Feynman motives of banana graphs', *Communications in Number Theory and Physics*, 3(1) (2008).

In a more general setting, Vogtmann and other mathematicians have constructed a geometric space of trees which models the set of all phylogenetic trees. These trees are diagrams familiar from evolutionary biology which detail the relationships between species: the trees generally have a distinguished vertex called the root which indicates the common descendant of all the species from which the diagram grows[186]. The diagrams will also have a number $n$ of 'final' nodes of degree 1 which are called leaves in the formal terminology of the space of trees. Each space of trees will consist only of trees with the same number of leaves. The space has a metric of non-positive curvature, giving a formal way of measuring distances between different phylogenetic trees or finding trees which are a fixed distance from other trees. This model also allows rigorous justification for the practice of discarding sections of a set of trees which agree with each other under comparison[187].

[186] Richard Dawkins and Yan Wong, *The Ancestor's Tale: A Pilgrimage to the Dawn of Life* (London: W&N, 2017).
[187] Louis Billera, Susan Holmes and Karen Vogtmann, 'Geometry of the space of phylogenetic trees', *Advances in Applied Mathematics* 27, 733 – 767 (2001).

# Carolyn Gordon

Gordon is a fellow of the American Mathematical Society and the American Association for the Advancement of Science. She specializes in the study of homogeneous spaces and isospectral geometry, and is best known for providing a negative answer (with Webb and Wolpert) to the question 'Can you hear the shape of the drum?' (first asked by Kac). This was done by constructing two subsets of $\mathbb{R}^2$ which have different shapes but the same eigenvalues[188]. To be more technical, they found two domains which have different shapes and yet the spectrum of their Laplacians has the same frequency. In fact, the theory of the Laplacian on a Riemannian manifold is very deep. If $(M, g)$ is a compact Riemannian manifold, there has to be a Laplace operator $\Delta$ on $M$ which acts on smooth functions (0-forms) on $M$:

$$\Delta = -\mathrm{div}(\mathrm{grad}\, f).$$

In local coordinates this takes the form

$$\begin{aligned}\Delta f &= \frac{1}{|g|^{\frac{1}{2}}} \sum_{k,j=1}^{n} \frac{\partial}{\partial x^k}\left(|g|^{\frac{1}{2}}\left(\frac{\partial f}{\partial x^j}\right) g^{kj}\right),\\ &= \sum_{k,j=1}^{n} g^{kj} \frac{\partial^2 f}{\partial x^k \partial x^j} + \text{lower order terms.}\end{aligned}$$

I will recall some spectral theory for finite dimensional normed spaces. If $T: X \to X$ is a symmetric linear operator and $X$ a finite dimensional normed space, then there exists an orthonormal basis of eigenvectors of $X$ with eigenvalues

$$0 \leq \lambda_1 \leq \cdots \leq \lambda_n.$$

The set of all the $\{\lambda_i\}$ is the spectrum of the operator $\sigma(T)$ and we also have that

$$\begin{aligned}\lambda \in \sigma(T) &\Leftrightarrow \exists (T - \lambda I)^{-1}\\ &\Leftrightarrow \mathrm{Ker}(T - \lambda I) = 0.\end{aligned}$$

The spectrum of a linear operator is a generalization of the concept of a set of eigenvalues of a matrix which you are probably familiar with from linear algebra. Depending on the choice of basis for $X$, $T$ can be represented by a matrix. In that case, the spectral theory just becomes the elementary theory of eigenvalues for matrices. For a given square matrix $A$, eigenvalues and eigenvectors are defined in terms of an equation

$$Ax = \lambda x.$$

An eigenvalue of such a square matrix is a number $\lambda$ such the above equation has a non-trivial solution and $x$ is an eigenvector of the matrix corresponding to that eigenvalue. I will not describe the spectral theory for square matrices any further as it can be looked up in any textbook on finite dimensional linear algebra.

---

[188] Carolyn Gordon, David Webb and Scott Wolpert, 'One cannot hear the shape of a drum', *Bulletin of the American Mathematical Society* 27, 134 – 138 (1992).

In infinite dimensions, things become more complicated. It is important to allow infinite dimensional spaces for physical applications, since the state spaces of quantum mechanics must be allowed to have infinite dimension. Historically, quantum mechanics provided much of the impetus for the development of functional analysis. If $X$ is infinite dimensional, it is possible to have a linear operator $T$ which has values in its spectrum which are not eigenvalues of that operator. For example, take the space $l^2$ and define a linear operator called the right-shift operator which takes us back into $l^2$.

$$T: l^2 \rightarrow \; l^2,$$

$$(\xi_1, \xi_2, \dots) \longmapsto (0, \xi_1, \xi_2, \dots),$$

where $(\xi_j)$ is a sequence in $l^2$. The right shift is a bounded operator. One can also invert this operator to get a left-shift operator defined by

$$(\xi_1, \xi_2, \dots) \longmapsto (\xi_2, \xi_3, \dots).$$

One can show with the definition of this operator that $T$ contains the value $0$ in its spectrum, and yet $0$ is not an eigenvalue of $T$ because the zero vector is not an eigenvector. However, a lot of the relevant theory does carry over if we stick to compact operators: the spectral properties of compact operators on Hilbert spaces are very close to those of operators on finite dimensional spaces. If $X$ and $Y$ are two normed spaces, an operator $T: X \rightarrow Y$ is compact if for every bounded subset $M$ of $X$ the closure of $T(M)$ is also compact[189].

In particular, the eigenvector decomposition of the normed space $X$ goes through the same if we take $X$ an infinite dimensional Hilbert space and $T$ a compact operator. Again, the Hilbert space has an orthonormal basis of eigenvectors for the operator and the eigenspaces are finite. The linear operator is also bounded by definition. A linear operator on an infinite dimensional space might be unbounded even if the range is finite dimensional. Maybe the most important unbounded linear operator is the differentiation operator which might be familiar from quantum mechanics:

$$D: \mathfrak{D}(D) \rightarrow L^2(-\infty, \infty),$$

$$x \longmapsto i\frac{dx}{dt}.$$

Differential equations often give rise to unbounded linear transformations. Some of the properties of the spectra of bounded self-adjoint linear operators still hold for unbounded self-adjoint linear operators. The eigenvalues are all real and the spectrum is real and closed, for example.

We can go further and try to show that the eigenvector decomposition also holds for particular unbounded differential operators on compact Riemannian manifolds, where the space $X$ is the Hilbert space of all the forms on the manifold (this is the same as the space of functions on the manifold in the case of 0-forms). An example is the circle, where $X = L^2(S^1, \mathbb{C})$. The simplest differential operator on the circle is $d/d\theta$, but this generalizes to the exterior derivative, which does not have a spectrum. Instead, we take the second derivative:

$$\Delta = -\frac{d^2}{d\theta^2}.$$

[189] Erwin Kreyszig, *Introductory Functional Analysis with Applications* (USA: John Wiley & Sons, 1978).

The eigenfunction decomposition of $L^2(S^1)$ is well known from Fourier analysis. An orthonormal basis is given by the trigonometric polynomials $\{e^{in\theta}\}$ and the eigenfunction decomposition is given by

$$f = \sum_n \langle f, e^{in\theta} \rangle e^{in\theta},$$

with

$$\langle f, g \rangle = \frac{1}{2\pi} \int_{S^1} f(\theta)\overline{g(\theta)}\, d\theta.$$

By considering the limits, it is easy to see that $\Delta$ is unbounded[190].

Two Riemannian manifolds are isospectral if the spectra of their Laplacians coincide, so the question arises as to whether two isospectral Riemannian manifolds must also be isometric. If $M$ has nonempty boundary, one can consider the Dirichlet spectrum, which is the spectrum of the Laplacian acting on smooth functions which vanish at the boundary (so-named because of the corresponding Dirichlet boundary condition). For a region of $\mathbb{R}^2$, the Dirichlet eigenvalues of the Laplacian are basically the frequencies which a drumhead under tension would produce if it were shaped like that region, hence the formulation in terms of vibrating drums. Gordon, Webb and Wolpert answered the question in the negative by constructing two simply connected plane regions which are non-isometric and whose Dirichlet and Neumann spectra coincide. In fact, Milnor disposed of the question for higher dimensional manifolds quite early on by constructing a pair of isospectral 16-dimensional tori which were non-isometric, but the corresponding problem in the plane remained open for several decades afterwards. One notes that simply being in lower dimension does not necessarily make a conjecture easier to prove or disprove, as it always depends on the nature of the problem at hand.

A key tool in the proof was Sunada's theorem. This states that if $M$ is a Riemannian manifold on which a finite group $G$ acts by isometries, $H$ and $K$ are subgroups of $G$ which act freely and there is a bijection $f: H \to K$ taking every element of $H$ to an element $f(h)$ which is conjugate to $h$ in $G$, then the quotient manifolds $H\backslash M$ and $K\backslash M$ are isospectral. Gordon, Webb and Wolpert used Sunada's method to construct a pair of isospectral of 2-orbifolds with boundary. However, the Neumann spectra of the orbifolds are the same as the Neumann spectra of the underlying manifolds, but the underlying spaces are simply connected plane regions. A slightly more complicated argument can be used to deduce the Dirichlet isospectrality of the underlying spaces[191].

[190] Steven Rosenberg, *The Laplacian on a Riemannian Manifold: An Introduction to Analysis on Manifolds* (Cambridge: Cambridge University Press, 1997).

Gordon has continued to do wide-ranging research on various aspects of isospectral Riemannian manifolds (much of this research being motivated by the idea that geometric and topological information about compact manifolds can be encoded in the spectra of natural operators on those manifolds). In particular, the overarching question is: to what extent does the spectrum of the Laplacian on a Riemannian manifold actually determine the geometry of the manifold? Low-dimensional round spheres are certainly determined by their spectra in a unique way, but there are other cases where geometry is not controlled by the spectrum. For dimension $n \geq 8$ there are isospectral deformations of Riemannian metrics on $S^n$ which can be chosen to be arbitrarily close to the round metric, and for dimension $n \geq 9$ there are isospectral deformations of Riemannian metrics on the ball in $\mathbb{R}^n$ which can be chosen to be arbitrarily close to the flat metric[192].

Some of Gordon's earlier work focussed on isometry groups of homogeneous manifolds. Homogenous spaces in some sense look the same wherever you are on the space, and so they play a key role in modelling the Universe in cosmology, where it is assumed that the Universe looks the same everywhere on the large scale. A Riemannian manifold $(M, g)$ will have a group of Riemannian isometries $\mathrm{Iso}(M, g)$ which map the manifold back to itself. A Riemannian isometry between manifolds $(M, g_M)$ and $(N, g_N)$ is a diffeomorphism from the first manifold to the second such that the metric on the first manifold is equal to the metric on the second manifold pulled back by the diffeomorphism:

$$F^* g_N = g_M.$$

This can also be written as

$$g_M(v, w) = g_N\big(DF(v), DF(w)\big),$$

where $v$ and $w$ are tangent vectors in $T_pM$. This follows from the definition of the pullback. If $F: M \to N$ is more generally a smooth map between smooth manifolds and $\omega$ is a smooth $k$-form on $N$, then the pullback of that form is a smooth form on $M$ defined by

$$(F^* \omega)_p(X_1, \dots, X_k) = \omega_{F(p)}(F_* X_1, \dots, F_* X_k),$$

where $X_1, \dots, X_k \in T_pM$ and $F_* X_p \in T_{F(p)}N$ is a vector obtained by pushing forward a vector field at a point on $M$ by the smooth map $F$. Although vector fields do not necessarily push forward to new vector fields, smooth covector fields do pull back to smooth covector fields. If the map $F$ is not an injection (ie. different elements in the domain can be mapped to the same element in the codomain), then there might be different vectors at a point of $N$ which are obtained by pushing forward a vector field on $M$ from different points in the manifold[193]. Informally speaking, the metric measures lengths of tangent vectors and so it can be used to measure lengths of curves.

[191] Carolyn Gordon, David Webb and Scott Wolpert, 'One cannot hear the shape of a drum', *Bulletin of the American Mathematical Society* 27, 134 – 138 (1992).

[192] Carolyn Gordon, 'Isospectral deformations of metrics on spheres', Inventiones Mathematicae 145, 317 – 331 (2001).

The manifold is homogeneous if the isometry group acts transitively such that there is always an $F \in \mathrm{Iso}(M, g)$ with the property that $F(p) = q$ for two points $p, q \in M$. The isotropy group $\mathrm{Iso}_p(M, g)$ at the point $p$ is the set of isometries in the group such that $F(p) = p$. For a homogeneous space such as Euclidean space, one has that the space is isomorphic to $\mathrm{Iso}/\mathrm{Iso}_p$ for any point $p$. As an example, if we take the $n$-sphere equipped with the canonical metric, the Riemannian isometry group is as follows:

$$\mathrm{Iso}(S^n(r), \mathrm{can}) = O(n+1),$$

which is also the isotropy group of $(\mathbb{R}^{n+1}, \mathrm{can})$ at the point $0$. A smooth manifold equipped with a transitive smooth action by a Lie group is often referred to as a homogeneous manifold. One can take a Lie group $G$ (a smooth manifold which is also a group such that the multiplication and inverse maps given by

$$m(g, h) = gh, i(g) = g^{-1}$$

are both smooth) and trivialise the tangent space using left or right translations on $G$:

$$TG \simeq G \times T_e G.$$

The left and right translation maps are self-maps of the Lie group defined by

$$L_g(h) = gh, R_g(h) = hg.$$

$T_e G$ is the tangent space of the Lie group at the identity element. A tangent vector in $T_e G$ can be extended to a left (or right) invariant vector field by left (or right) translation to other parts of the manifold, hence $T_e G$ can be identified with the space of left invariant vector fields. A vector field on $G$ is left-invariant if for any two elements in the group

$$(L_g)_* X_{g'} = X_{gg'}.$$

However, the space of all smooth left-invariant vector fields on a Lie group is an example of a Lie algebra, which is referred to as the Lie algebra $\mathfrak{g}$ of the group $G$. A Lie algebra in general is simply a real vector space equipped with a map from $\mathfrak{g} \times \mathfrak{g}$ to $\mathfrak{g}$ called the bracket which satisfies the Jacobi identity along with bilinearity and antisymmetry. One can also show that the connected component of the Lie group which contains the identity element is the only connected open subgroup of $G$ and that every connected component of the group is diffeomorphic to this component.

Because of the above, the canonical inner product on $T_e G$ induces a Riemannian metric on the Lie group, so left translations are Riemannian isometries. If $H$ is a closed normal Lie subgroup of $G$, the quotient $G/H$ is a smooth manifold by the quotient manifold theorem, and it follows that it also a Lie group and that the quotient map is a Lie group homomorphism. In fact, if we equip $G$ with a metric for which all right translations are Riemannian isometries, then there is a unique Riemannian metric on the quotient group such that the quotient map is a submersion. If the metric is left invariant, then $G$ acts as an isometry group on the quotient, making $G/H$ into a homogeneous manifold[194]. If we start with the complex space $\mathbb{C}^{n+1}$, this has $S^{2n+1}$ as a subset. The Lie group acts as an isometry group on $\mathbb{C}^{n+1}$ and $S^{2n+1}$ and the quotient map

$$\pi: \mathbb{C}^{n+1} \backslash \{0\} \to \mathbb{CP}^n$$

is a surjective submersion when restricted to $S^{2n+1}$, so we have

$$\frac{S^{2n+1}}{S^1} \simeq \mathbb{CP}^n.$$

[193] John M. Lee, *Introduction to Smooth Manifolds* (New York: Springer, 2003).

Since $S^1$ acts smoothly as a Lie group by isometries of the round metric, the metric on $S^{2n+1}$ induces a homogeneous, isotropic metric on $\mathbb{CP}^n$ such that the map $S^{2n+1} \to \mathbb{CP}^n$ is a Riemannian submersion: this metric is called the Fubini-Study metric[195].

A very important Lie group is the special unitary group $SU(2,\mathbb{C})$, which is defined as the following set of matrices with complex entries:

$$SU(2,\mathbb{C}) = \left\{ \begin{bmatrix} \alpha & \beta \\ -\bar{\beta} & \bar{\alpha} \end{bmatrix} : |\alpha|^2 + |\beta|^2 = 1 \right\}.$$

The corresponding Lie algebra

$$\mathfrak{su}(2,\mathbb{C}) = \left\{ \begin{bmatrix} i\alpha & \beta + i\gamma \\ -\beta + i\gamma & -i\alpha \end{bmatrix} : \alpha, \beta, \gamma \in \mathbb{R} \right\}$$

has the following basis:

$$u_1 = \frac{1}{2}\begin{bmatrix} i & 0 \\ 0 & -i \end{bmatrix}, u_2 = \frac{1}{2}\begin{bmatrix} 0 & 1 \\ -1 & 0 \end{bmatrix}, u_3 = \frac{1}{2}\begin{bmatrix} 0 & i \\ i & 0 \end{bmatrix},$$

which can be visualised as left invariant vector fields on $SU(2,\mathbb{C})$. If the basis is orthonormal, then we obtain a metric which is left invariant. Another Lie group which is important in theoretical physics is the special linear group $SL(2,\mathbb{C})$: this is the subgroup of the general linear group such that the matrices have unit determinant. The fundamental group of a connected Lie group with identity element as base point is always abelian. Although $SL(2,\mathbb{C})$ is non-compact, it is of a particular form such that its fundamental group can be reduced to computation of a compact classical group: in this case, $SU(2)$. The fundamental group of $SU(2)$ is isomorphic to the fundamental group of the 3-sphere, but the fundamental group for a sphere apart from the circle has to be the trivial group (because any loop can be shrunk down to a point). One can show that

$$\pi_1\big(SL(2,\mathbb{C})\big) \simeq \pi_1\big(SU(2)\big),$$

which shows that $\pi_1\big(SL(2,\mathbb{C})\big)$ is the trivial group[196].

---

194 Peter Petersen, *Riemannian Geometry* (New York: Springer, 2006).
195 John M. Lee, *Riemannian Manifolds: An Introduction to Curvature* (New York: Springer, 1997).
196 Brian C. Hall, *Lie Groups, Lie Algebras, and Representations: An Elementary Introduction* (Switzerland: Springer, 2015).

Spinors are defined via representations of $SL(2,\mathbb{C})$. They can also be defined in terms of representations of $SU(2,\mathbb{C})$, but these are space spinors, rather than the usual spacetime spinors. Space spinors could be used for problems where one is working on Riemannian 3-manifolds (assuming existence of a spinor structure equipped with a Hermitian product)[197]. $SL(2,\mathbb{C})$ is the universal covering group of the Lorentz group $SO(3,1)$. (If $G$ is a connected Lie group, then there exists a simply connected Lie group $G'$ called the universal covering group and a smooth covering map $\pi: G \to G'$ which is a homomorphism of Lie groups: see a textbook on algebraic topology or smooth manifolds if you need to revise elementary covering space theory). If $V$ is the fundamental representation of $SL(2,\mathbb{C})$ on $\mathbb{C}^2$, this representation is equipped with an invariant symplectic form, but there is no hermitian structure on $\mathbb{C}^2$, as the size of vectors does not change.

Take a ball of radius $\pi$ and consider the map $B^3(\pi) \to SO(3)$ which takes the pair $(r,\theta)$ to a rotation by the angle $r$ with the axis along $\theta$. This is a surjective map, as one can always find a point by picking an axis and taking any rotation. Rotations can only be the same if we rotate by $\pi$ and $-\pi$, so only antipodal points can be identified. This means that

$$SO(3) \simeq {}^{B^3}\!/_{\sim} = \mathbb{RP}^3.$$

There exists a covering map $\pi: S^3 \to \mathbb{RP}^3$, but this gives a covering map by projection $\mathrm{Spin}(3) \to SO(3)$, because $\mathbb{RP}^3 \simeq SO(3)$ and by one of the accidental isomorphisms amongst the classical Lie groups in low dimension, $\mathrm{Spin}(3) = SU(2)$, which is diffeomorphic to $S^3$. One can take the inclusion maps $\mathrm{Spin}(3) \to SL(2,\mathbb{C})$ and $SO(3) \to SO(3,1)$. We already know that there is a covering map $SL(2,\mathbb{C}) \to SO(3,1)$, so the maps all commute with each other to form a commutative diagram. This means that $V$ can be regarded as a $\mathrm{Spin}(3)$ representation with a Hermitian structure, which can be thought of as an isomorphism given by Clifford multiplication or as a new Hermitian inner product[198]. Covering maps are themselves local Riemannian isometries.

An important and useful class of homogeneous spaces are the symmetric spaces (introduced and studied by Cartan). A Riemannian manifold $(M,g)$ is a symmetric space if for every point in the manifold, the isotropy group for that point $\mathrm{Iso}_p$ contains an isometry $A_p$ such that the differential $DA_p: T_pM \to T_pM$ is the antipodal map. The fact that isometries preserve geodesics can be used to show that such a space is homogeneous. Symmetric spaces are geodesically complete, meaning that any geodesic can be extended indefinitely in both directions. This can be proved using the fact that an involution of a symmetric space switches geodesics through $p$. A Riemannian manifold is symmetric if for every point in the manifold there is an isometry $\sigma_p$ called an involution such that

$$\sigma_p(p) = p,$$

$$D\sigma_p(p) = -Id.$$

[197] Thomas Bäckdahl and Juan Valiente Kroon, 'A formalism for the calculus of variations with spinors', *Journal of Mathematical Physics* 57, 022502 (2016).
[198] Thomas Parker and Clifford Henry Taubes, 'On Witten's Proof of the Positive Energy Theorem', *Commun. Math. Phys.* 84, 223 – 238 (1982).

$\mathbb{CP}^n$ equipped with the Fubini-Study metric is an example of a symmetric space. In a symmetric space, any two points can be connected by a geodesic. Any Lie group with a bi-invariant metric is a symmetric space, since the inverse operation forms the required symmetry around the identity element. These include some of the well-known Lie groups: $SU(n+1)$, $SO(2n+1)$, $Sp(n)$ and $SO(2n)$. Interestingly, the symmetric spaces come in pairs composed of a compact and a non-compact space. In more geometric terms, symmetric spaces have parallel curvature tensor: that is[199],

$$\nabla R = 0.$$

Gordon has recently worked on a particular type of complete Riemannian manifold known as a geodesic orbit manifold: $(M, g)$ is a geodesic orbit manifold if every geodesic is an orbit of a one-parameter group of isometries. A manifold of this kind must also be homogeneous. With Nikonorov, Gordon showed that if $(M, g)$ is an $n$-dimensional, simply connected, geodesically orbit manifold, then $(M, g)$ is isometric to a non-compact symmetric space, a simply connected, geodesically orbit Riemannian nilmanifold, or the total space of a Riemannian submersion with totally geodesic fibres, where the base space is a non-compact symmetric space and the fibres are isometric to a simply connected, geodesically orbit Riemannian nilmanifold. They also showed that if $(M, g)$ is a geodesically orbit manifold which is diffeomorphic to $\mathbb{R}^n$, then $M$ admits a simply transitive solvable isometry group of the form $S \times N$, where $S$ is an Iwasawa subgroup of a semisimple Lie group and $N$ is the group in a geodesically orbit nilmanifold $(N, g)$ of step size at most two[200].

In connection with the spectral theory of the Laplacian and bi-invariant metrics on Lie groups, Gordon, Schueth and Sutton showed that a bi-invariant metric on a compact connected Lie group is spectrally isolated in the class of left-invariant metrics. The term 'spectral isolation' refers to the question of whether special types of Riemannian manifold such as symmetric spaces have spectra which can be distinguished from those of other Riemannian manifolds. Aside from manifolds of constant curvature, there are not any examples of metrics which are spectrally isolated from other arbitrary Riemannian metrics, but one can restrict to classes of metrics: for example, any set of mutually isospectral compact symmetric spaces is finite. Gordon considered the question as to whether symmetric spaces can be spectrally distinguished from other homogeneous Riemannian manifolds: more specifically, if we restrict to compact symmetric spaces which are given by bi-invariant Riemannian metrics on compact Lie groups, is such a metric spectrally isolated if we stay in the class of left-invariant Riemannian metrics? The answer is yes, and Gordon obtained a stronger quantitative result: given a bi-invariant metric $g$ on $G$ there is a positive integer $n$ such that in a neighbourhood of the metric in the class of left-invariant metrics with the same volume or less, $g$ is uniquely determined by the first $n$ distinct non-zero eigenvalues of the Laplacian[201].

[199] Jürgen Jost, *Riemannian Geometry and Geometric Analysis* (Heidelberg: Springer, 2008).
[200] Carolyn Gordon and Yurii Nikonorov, 'Geodesic orbit Riemannian structures on $\mathbb{R}^n$', *J. Geom. Phys.,* 134, 235 – 243 (2018).

In a slightly more applied setting, Gordon studied the case of a real $2m$-torus $M$ with a translation-invariant metric $g$ and a translation-invariant symplectic form $\omega$, where the latter is considered to be a magnetic field on $M$. To this form and metric, one can associate a Hamiltonian system $(T^*M, \Omega, H)$, where $\Omega$ is the twisted symplectic form given by $\omega_0 + \pi^*\omega$, $\pi$ is the projection mapping $T^*M$ to $M$, and $H$ is the Hamiltonian function given by

$$H(q, \xi) = \frac{1}{2} h_q(\xi, \xi).$$

This system describes the dynamics of a charged particle moving in the magnetic field, and could presumably be used to model the dynamics of plasma inside a tokamak[202]. The idea of applying Hamiltonian dynamics to plasma confinement is not restricted to tokamaks, and is also being considered for devices known as stellarators. In this case, one is looking to have dynamics for the guiding-centre which is integrable or almost integrable to maximise confinement and no or small toroidal current. The condition on the toroidal current requires deviation from axisymmetry. Integrability is guaranteed in the case of quasisymmetry, but it is believed that the only true quasisymmetries are axisymmetries, so scientists will have to develop ways of testing the deviation from integrability which will necessarily result in the design of stellarator devices.

In general, magnetic confinement uses electromagnetic fields to confine a hot plasma of deuterium and tritium as it circulates in a region. A major practical difficulty involved in controlled fusion is that the fuel must be confined over a long time at extreme temperatures, with the fields being necessary to keep the fuel away from the walls of the tokamak. The fields can only confine the plasma for a short period of time and the plasma is prone to instabilities due to the temperatures and the strong magnetic fields. The tokamak itself must also be able to withstand the electromagnetic forces which are being generated. To assess these, one assumes the fusion device to be a three-dimensional conducting structure $V_c$ which surrounds a plasma torus. The magnetic vector potential is derived using the Biot-Savart law:

$$\mathbf{A}(\mathbf{r}) = \int_{V_c} \frac{\mathbf{J}(\mathbf{Q})}{|\mathbf{r} - \mathbf{Q}|} \, dV_c + \mathbf{A}_{\text{plasma}}(\mathbf{r}).$$

The electric vector potential $\mathbf{T}$ is expanded in terms of edge elements $\mathbf{N}_i$ once we have a finite element discretization of $V_c$:

$$\mathbf{T} = \sum_i I_i \mathbf{N}_i.$$

[201] Carolyn Gordon, Dorothee Schueth and Craig Sutton, 'Spectral isolation of bi-invariant metrics on compact Lie groups', *Ann. Inst. Fourier* 60(5), 1617 – 1628 (2010).
[202] Carolyn Gordon, William Kirwin, Dorothee Schueth and David Webb, 'Classical equivalence and quantum equivalence of magnetic fields on flat tori', arXiv:1108.5113v1 (2011).

The model for the conducting structure can be combined with a model for the plasma, which would typically be modelled in the form of single-fluid magnetohydrodynamics equations defined between the plasma and the walls of the device. The plasma model would need to be discretized using similar methods and then combined with the differential equation for the potentials. As mentioned earlier, there is always the possibility of instabilities developing in the plasma, even when it is in its equilibrium configuration: these perturbations can die down or they can grow. In particular, there is interest in so-called resistive wall modes. These are unstable modes which develop when the plasma requires a perfectly conducting wall to stabilize its external kink modes. A vacuum chamber surrounding a plasma is a good approximation to a perfectly conducting wall, but its conductivity is finite. Working through the eigenvalue problem, it can be shown that finite wall resistivity destroys kink stability and so resistive wall modes are inherently unstable. The modes can be computed using inverse iteration algorithms, but this becomes very computationally intense, and there would probably need to be some development of parallel computing techniques before monitoring of instabilities in tokamaks becomes feasible[203]. Another possibility for a magnetic confinement device is the Reversed Field Pinch, which has similar magnetic field strength in both the toroidal and the polar direction, whereas the tokamak has a stronger field in the toroidal direction.

[203] Bruno Carpentieri, 'Fast computational techniques for modelling RFX-mod fusion devices on hybrid CPU-GPU architectures', URSI Commission B International Symposium on Electromagnetic Theory (Espoo, Finland: EMTS, 2016).

# *Frances Kirwan*

Kirwan mostly works in algebraic and symplectic geometry. She was the second-youngest person to be elected as President of the London Mathematical Society and the first female mathematician to be elected to the historic post of Savilian Professor of Geometry. We have discussed algebraic geometry earlier in the book, but did not mention algebraic curves. These curves form possibly the richest part of algebraic geometry and have a unique status in terms of their beauty, classical prestige and importance in current mathematical research. An algebraic curve is an algebraic variety of dimension 1. If we take an affine or projective subset (a subset of affine or projective space) together with its Zariski topology, one can define a corresponding sheaf of functions $\mathcal{O}_X$ with the appropriate restriction. The triple of a subset with the topology and the sheaf is known as an affine or a projective variety, respectively (this definition absorbs the one which we suggested earlier). One can have an affine or a projective curve accordingly, where a projective curve is a zero set of a polynomial in a projective plane.

A very important family of algebraic curves are the elliptic curves: these were used by Wiles in his famous proof of Fermat's Last Theorem. An elliptic curve is by definition smooth and non-singular. It has no self-intersections or cusps (it cannot loop around on itself or curve up and back down in a way that you have a spike). Especially in applications and more computational settings, it is common to write an explicit equation for the curve known as the Weierstrass form.

$$y^2 = x^3 + Ax + B,$$

where $A$ and $B$ are constants. If $A$ and $B$ are taken to be elements of some field $k$, then we say that the elliptic curve $E$ is defined over that field. Drawing or visualising an elliptic curve over an arbitrary field is usually impossible, but one can visualise more familiar graphs over $\mathbb{R}$ to aid intuition[204]. For example, the plane curve in $\mathbb{R}^2$

$$y^2 = x^3$$

is a cuspidal cubic curve and

$$y^2 = x^3 + x^2$$

is a nodal cubic curve with a self-intersection. These cubics can be parametrised much as you might be used to parametrising a circle or an ellipse[205]. It is also useful to add the point at infinity to an elliptic curve, so one usually includes this point when defining the curve over a field:

$$E(L) = \{\infty\} \cup \{(x, y) \in L \times L \text{ such that } y^2 = x^3 + Ax + B\}.$$

---

[204] Lawrence Washington, *Elliptic Curves: Number Theory and Cryptography* (New York: Chapman & Hall, 2008).

[205] Miles Reid, *Undergraduate Algebraic Geometry* (Cambridge: Cambridge University Press, 2010).

It is an interesting exercise to consider how one might take two points on an elliptic curve and obtain a new point. This leads to the group law for elliptic curves. Take an elliptic curve $E$ in Weierstrass form. Let $P_1 = (x_1, y_1)$ and $P_2 = (x_2, y_2)$ be two points on $E$ (not the point at infinity). We define a third point $P_3$ given by addition of these two points as follows: if $x_1$ is different to $x_2$, then

$$x_3 = \left(\frac{y_2 - y_1}{x_2 - x_1}\right)^2 - x_1 - x_2, \quad y_3 = \frac{y_2 - y_1}{x_2 - x_1}(x_1 - x_3) - y_1.$$

If $x_1 = x_2$ but $y_1$ is different to $y_2$, then $P_1 + P_2 = \infty$. If $P_1 = P_2$ and $y_1$ is non-zero, then

$$x_3 = \left(\frac{3x_1^2 + A}{2y_1}\right)^2 - 2x_1, \quad y_3 = \frac{3x_1^2 + A}{2y_1}(x_1 - x_3) - y_1.$$

If $P_1 = P_2$ and $y_1$ is zero, then $P_1 + P_2 = \infty$ and $P + \infty = P$ for all $P$. It can be shown that addition of points on an elliptic curve satisfies the group axioms and hence the points of $E$ form an abelian group where the formal point at infinity $\infty$ is the identity. Obviously, over a finite field there will only be a finite number of points with coordinates drawn from that field, so the group will be finite as well as abelian. In the important case of an elliptic curve defined over $\mathbb{Q}$, $E(\mathbb{Q})$ is a finitely generated abelian group. This a key result known as the Mordell-Weil theorem and for such a group we have

$$E(\mathbb{Q}) \cong \mathbb{Z}^r \oplus F,$$

for some finite group $F$. Another key theorem is the Lutz-Nagell theorem, which says that if $E$ is an elliptic curve in Weierstrass form with $A$ and $B$ in $\mathbb{Z}$ and $P$ in $E(\mathbb{Q})$, then $x$ and $y$ are also in $\mathbb{Z}$, and if $y$ is non-zero, then $y^2$ divides $4A^3 + 27B^2$. This theorem can often be used to find possible torsion points for an elliptic curve defined over $\mathbb{Q}$, where a torsion point is a point of finite order (ie. a member of the torsion group). All points in an elliptic curve over a finite field are torsion[206].

Kirwan has worked in particular on moduli spaces of complex algebraic curves. As you might expect, a complex algebraic curve in $\mathbb{C}^2$ is a subset of $\mathbb{C}^2$ of the form

$$C = \{(x, y) \in \mathbb{C}^2 : P(x, y) = 0\},$$

where $P$ is a polynomial in two variables with coefficients which are complex numbers. Obviously, one can take real coefficients for the polynomial to obtain a corresponding real algebraic curve. This type of curve has been studied intensively since ancient times (the cubic curves which Newton studied and classified were real). However, complex algebraic curves can often be easier to use than their real versions. The same polynomial might factorise completely when the coefficients are complex compared to the version with real coefficients, for example. You have probably realised that a complex algebraic curve cannot be drawn in $\mathbb{C}^2$ because it has four real dimensions (the complex plane $\mathbb{C}$ already has two real dimensions and $\mathbb{C}^2 = \mathbb{C} \times \mathbb{C}$). If we add the point at infinity, we can have pictures of curves which give us a good idea of what is happening topologically, but in this case when the point $\infty$ is added to a complex line it becomes a sphere in topological terms. This is effectively what happens with the Riemann sphere, which can be considered as the complex projective line. This all links back in with complex analysis and the theory of Riemann surfaces, since the Riemann sphere (and every compact Riemann surface in general) is a projective algebraic curve. The symbol for a sheaf $\mathcal{O}_X$ actually originates in the complex analysis setting (from the Italian word 'olomorfico', meaning 'holomorphic').

[206] Lawrence Washington, *Elliptic Curves: Number Theory and Cryptography* (New York: Chapman & Hall, 2008).

We will mention some of the topological properties of complex algebraic curves. A non-singular projective curve of genus $g$ in the projective plane $\mathbb{P}^2$ is topologically speaking a sphere with $g$ handles (for a handle, imagine that you have cut through a torus to disconnect it and make it into a cylinder and then sewn both sides onto the sphere). The genus of this curve always satisfies the degree-genus formula:

$$g = \frac{1}{2}(d-1)(d-2).$$

If the curve is an irreducible projective curve with a set of singular points $p_i$, one can assign an integer $\delta$ to each singularity such that Noether's formula holds:

$$g = \frac{1}{2}(d-1)(d-2) - \sum_{j=1}^{r} \delta(p_i).$$

One can prove the degree-genus formula by using the fact that every non-singular projective curve in the projective plane can be viewed as a branched cover of the projective line. Every non-singular projective curve in the plane (ie. every compact Riemann surface) can be triangulated and the Euler characteristic of a curve does not depend on this choice of triangulation. For the case of a compact orientable surface with no boundary and genus $g$, the Euler characteristic is simply $2-2g$. One can prove this by starting with a surface of genus zero (a sphere). This has Euler characteristic 2 and one inducts on this surface[207].

A very significant result in this area is the Riemann-Roch theorem: this theorem directly relates the complex analysis of a Riemann surface to the genus (a topological invariant). The theorem says that if $D$ is a divisor on an algebraic curve $X$, then

$$\dim L(D) - \dim H^1(D) = \deg \mathrm{D} + 1 - \dim H^1(0),$$

where $L(D)$ is the space of meromorphic functions on $X$ with poles bounded by $D$. A problem where one attempts to compute the dimension of this space is known as a Riemann-Roch problem. Recall from complex analysis that a function is meromorphic on a domain if it is holomorphic (infinitely differentiable and locally equal to its own Taylor series within a neighbourhood of a point) everywhere except for a set of points. If a holomorphic function of a complex variable $f(z)$ can be written in the form

$$f(z) = (z - z_0)^n g(z)$$

[207] Frances Kirwan, *Complex Algebraic Curves* (Cambridge: Cambridge University Press, 1992).

for integer $n$, the function has a zero of multiplicity $n$ at $z_0$. If $n = 1$, we say that $z_0$ is a simple zero. A point where a meromorphic function is no longer holomorphic is called a singularity. If we can expand a function about a point $a$ via a Laurent series (like a Taylor series, but terms of negative degree are allowed),

$$f(z) = \cdots + C_2(z-a)^2 + C_1(z-a) + C_0 + \frac{C_{-1}}{z-a} + \frac{C_{-2}}{(z-a)^2} + \cdots,$$

then that point is a singularity. If there is only a finite number of terms containing negative powers of $z - a$, then the point is a pole[208].

If $X$ is a Riemann surface, a divisor on $X$ is a function $D : X \to \mathbb{Z}$ such that set of points in $X$ where $D(p) \neq 0$ is a discrete subset of that curve. The divisors of the Riemann surface form a group when equipped with pointwise addition as a group operation. In the case of a compact Riemann surface, the degree of a divisor is just the sum of the values of that divisor:

$$\deg(D) = \sum_{p \in X} D(p).$$

A divisor which can be defined in terms of the order function is known as a principal divisor on $X$:

$$\operatorname{div}(f) = \sum_{p} \operatorname{ord}_p(f) \cdot p.$$

In the case where we have a meromorphic 1-form rather than a function, the divisor is canonical.

$$\operatorname{div}(\omega) = \sum_{p} \operatorname{ord}_p(\omega) \cdot p.$$

The set of canonical divisors form a coset of the principal divisors of $X$, which themselves form a subgroup of the divisors. As it is, our statement of the Riemann-Roch theorem is not very helpful from a computational point of view, as it turns the problem of computing the dimension of $L(D)$ into a problem of computing first cohomology groups, which must be identified. This can be done using the Serre duality theorem, which states that for any divisor $D$ on an algebraic curve $X$ (a compact Riemann surface if we are working over $\mathbb{C}$), the residue map

$$\operatorname{Res} : L^{(1)}(-D) \to (H^1(D))^*$$

is an isomorphism between complex vector spaces. For a canonical divisor $K$ and a general divisor $D$, one has

$$\dim H^1(D) = \dim L^{(1)}(-D) = \dim L(K - D).$$

[208] Louis Milne-Thomson, *Theoretical Hydrodynamics* (London: Macmillan & Co Ltd, 1968).

This eventually leads to the improved version of the Riemann-Roch theorem. Start with an algebraic curve $X$ of genus $g$. For any general divisor $D$ and any canonical divisor $K$, one has the beautiful result

$$\dim L(D) - \dim L(K - D) = \deg(D) + 1 - g.$$

The theorem is especially useful for computing the degree of large divisors[209].

A moduli space of algebraic curves is a space whose points can be viewed as isomorphism classes of algebraic curves. In other words, it is a quotient space such that the equivalence relation on the set of (non-singular) curves is isomorphism, or a space which parametrises all curves of genus $g$. (The definition also holds if we take some other type of algebraic geometric object apart from a curve). We might typically denote this space with $\mathcal{M}_g$ where $g$ is the genus or algebraic genus of the curves, bearing in mind that the space will become more complicated as the genus gets larger. A moduli space is often what is called an algebraic stack, to the point where we might refer to the more fundamental moduli stack $\mathcal{M}_g$. I will not give the definition of a stack as it is technical and you can look into it elsewhere. In some cases, it is easier to consider moduli spaces of curves with several points which are marked or 'labelled': this is denoted by $\mathcal{M}_{g;n}$, where $n$ is the number of marked points. As a trivial example, the moduli space $\mathcal{M}_{0;3}$ is a single point, because two points on the complex projective line given be mapped to each other by an automorphism of the line. In the case of elliptic curves with one marked point, things are already more complicated. The corresponding moduli space $\mathcal{M}_{1;1}$ is the quotient of the upper half-plane by a discrete $\mathrm{SL}(2, \mathbb{Z})$ group action, but it can also be represented as the quotient of a different moduli space $\mathcal{M}_{0;4}$ by an action of a permutation group.

The general theory behind construction of moduli spaces of curves is difficult, but we can say a few general things about such spaces. The space must have a topology and a complex structure. There is some sense which can be made precise in which a projective curve which has an equation which is a small perturbation of the equation for another curve must give a point in the moduli space which is a small perturbation away from the other point. In the case with no marked points and genus greater than or equal to 3, the moduli space of curves has to have a universal curve. More explicitly, the moduli space has to come with a complex variety $\mathcal{C}_{g;0}$ and a map $\pi\colon \mathcal{C}_{g;0} \to \mathcal{M}_{g;0}$ such that the fibre over a point $\pi^{-1}(c)$ is a complex curve whose equivalence class coincides with the corresponding point in the moduli space. A moduli space has to contain all possible families of curves. If $p$ is a holomorphic map between two complex varieties $E$ and $B$ such that the fibre $p^{-1}(b)$ over every point $b$ in $B$ is the quotient of a complex curve modded out by its automorphism group, then there must exist maps from $E$ to $\mathcal{C}_{g;0}$ and from $B$ to $\mathcal{M}_{g;0}$ such that one can form a commutative diagram.

In general, one has to distinguish between the fine and the coarse moduli problem. A coarse moduli space of algebraic curves with genus $g$ is a triple consisting of $\mathcal{C}_{g;0}$, $\mathcal{M}_{g;0}$ and a map $\pi$ between them such that every curve of genus $g$ appears once as a fibre of $\pi$ and such that for every holomorphic family $p: E \to B$ whose fibres are quotients as above there exist holomorphic maps which allow us to form a commutative diagram as above[210]. One of the simplest non-trivial examples of a moduli space is a Hilbert scheme. This a moduli space of subvarieties (subschemes, more generally) with a Hilbert polynomial in some projective space. If we take a subvariety in a complex projective space and $I_X^{(n)} \subset \mathbb{C}^{(n)}[x_1, \dots, x_n]$ the space of polynomials which vanish on the variety $X$, then for large $n$ the dimension of the quotient is a polynomial in $n$ called the Hilbert polynomial. In the special case of points in the complex plane, the Hilbert scheme is a smooth irreducible variety:

$$\mathrm{Hilb}(\mathbb{C}^2, n) = \{\text{ideals } I \subset \mathbb{C}[x_1, x_2] \text{ of codimension } n\}.$$

[209] Rick Miranda, *Algebraic Curves and Riemann Surfaces* (USA: American Mathematical Society, 1995).

A collection of points in the complex plane is specified by an ideal in the coordinate ring

$$I_P = \{f(p_1) = \cdots = f(p_n) = 0\} \subset \mathbb{C}[x_1, x_2],$$

and the codimension of this ideal is always $n$, so this is the type of ideal which is being considered in the definition. The map from $I_P$ to $P$ extends to a map

$$\pi_{\mathrm{Hilb}}: \mathrm{Hilb}(\mathbb{C}^2, n) \to \frac{(\mathbb{C}^2)^n}{S(n)},$$

which is what is known as a resolution of singularities of $(\mathbb{C}^2)^n/S(n)$. This makes this particular Hilbert scheme into a type of variety known as an equivariant symplectic resolution[211].

Symplectic geometry is a branch of differential geometry which studies the properties of symplectic manifolds which are invariant under symplectomorphisms, where a symplectic manifold is a smooth manifold equipped with a smooth, closed, nondegenerate 2-form $\omega$. A choice of symplectic form on a manifold is also known as a symplectric structure and a diffeomorphism between two symplectic manifolds is called a symplectomorphism if the symplectic form on one manifold is equal to the form pulled back to the other manifold by that diffeomorphism:

$$F^* \widetilde{\omega} = \omega.$$

In the linear algebraic setting (known as linear symplectic geometry), a nondegenerate alternating 2-tensor is called a symplectic tensor and a vector space with a symplectic structure is called a symplectic vector space. A 2-tensor is said to be nondegenerate if the equation $\omega(X, Y) = 0$ implies that $X = 0$ for any $X$ and $Y$ in the vector space. If $\omega$ is a symplectic tensor on a vector space, it can be shown by an induction argument that the vector space must have even dimension and that there must exist a basis called the symplectic basis such that $\omega$ can be written as

[210] Maxim Kazaryan, Sergei Lando and Victor Prasolov, *Algebraic Curves: Towards Moduli Spaces* (Switzerland: Springer, 2018).
[211] Andrei Okounkov, 'On the crossroads of enumerative geometry and geometric representation theory' (2018), arXiv:1801.09818.

$$\omega = \sum_{i=1}^{n} \alpha^i \wedge \beta^i,$$

where the $\alpha^i$ and $\beta^i$ are elements of the dual basis for $V^*$. This implies that every symplectic manifold has even dimension. The most obvious example of a symplectic vector space is the Euclidean space $\mathbb{R}^{2n}$ equipped with a bilinear form[212]

$$\omega = \sum_{i=1}^{n} \mathrm{d}x_i \wedge \mathrm{d}y_i.$$

We have mentioned that moduli spaces are often constructed by quotienting out varieties by an action of a linear group $G$: for simplicity, we can take a linear group to be a semi-direct product of a unipotent group with a reductive group. A very important way of developing quotients in this way for moduli spaces is Mumford's geometric invariant theory. The basic aim of general invariant theory is to study group actions on a variety or scheme $X$. In geometric invariant theory, one studies these actions with a view to forming a geometric quotient of $X$ by $G$ which is itself a nice scheme. To be more technical, if $X$ is a normal quasi-projective $G$-variety, then there exists a $G$-linearization for the $G$-action on that variety. The goal is then to look for a natural variety which parametrises $G$-orbits in an affine or projective variety using invariant functions in the coordinate ring $k[X]$[213]. As a specific example, there is a theorem which states that if $X$ is an affine scheme over a field $k$, $G$ is a reductive algebraic group and $\sigma: G \times X \to X$ is a group action of $G$ on $X$, then a uniform categorical quotient $(Y, \phi)$ of $X$ by $G$ exists and $Y$ is an affine scheme. If $X$ is algebraic, then $Y$ is algebraic over $k$. If the field $k$ has characteristic zero, then $(Y, \phi)$ is a universal categorical quotient, and if $X$ is Noetherian, then $Y$ is Noetherian[214].

We did not fully define a scheme previously, so we will make a quick digression and do this here. Start with a topological space $X$. A sheaf of abelian groups on the base space $X$ is a topological space $\mathcal{O}$ along with a map $\pi$ from $\mathcal{O}$ to $X$ such that $\pi$ is a local homeomorphism, the fibres $\pi^{-1}(p)$ are abelian groups (the definition can easily be adjusted to allow for fibres which are rings or modules), and the group operations on each fibre are continuous in the relative topology from the space $\mathcal{O}$. The fibres of the sheaf are known as stalks in the jargon. Now, take a ringed space (ie. a sheaf of rings). If $\mathcal{O}_X$ and $\mathcal{O}_Y$ are two ringed spaces, a morphism of ringed spaces is defined to be a continuous function $\sigma: X \to Y$ and a collection $\psi$ of homomorphisms $\psi_U$ from $\mathcal{O}_Y(U)$ to $\mathcal{O}_X\left(\sigma^{-1}(U)\right)$ which respect composition with the relevant restriction maps for the ringed spaces. An affine scheme can be defined as a ringed space which is isomorphic to the structure sheaf $(\operatorname{Spec} A, \mathcal{O})$ of some ring. One can define morphisms of affine schemes which are locally morphisms of the schemes as ringed spaces. A general scheme is a ringed space $(X, \mathcal{O})$ such that every point of $X$ has an open neighbourhood so that the restriction to that neighbourhood is an affine scheme modulo appropriate isomorphisms. Schemes along with their morphisms form a category, as you might expect. The definition for a scheme may seem strange, but you will have to believe me that they are the only tool which allows mathematicians to formulate problems and prove results in algebraic geometry without being in some way tied down to their geometric intuition (this was the problem which ultimately derailed the great Italian school of classical algebraic geometry). The other algebraic objects such as varieties which you might be more familiar with can be recovered from schemes[215].

[212] John M. Lee, *Introduction to Smooth Manifolds* (New York: Springer, 2003).
[213] Brent Doran and Frances Kirwan, 'Towards non-reductive geometric invariant theory', *Pure and Applied Mathematics Quarterly* 1 (3) 61 – 105 (2007).
[214] David Mumford, John Fogarty and Frances Kirwan, *Geometric Invariant Theory*, Ergebnisse der Mathematik und ihrer Grenzgebiete (2), 23 (New York: Springer, 1994).

We have seen that much of the motivation for geometric invariant theory comes from the construction of moduli spaces, but it is also a developed theory in its own right, which in more algebraic terms studies the actions of reductive algebraic groups. An algebraic group $G$ over a field $k$ is a group pre-scheme $G/k$ which is also an algebraic scheme. Kirwan has recently studied the problem of generalizing geometric invariant theory to non-reductive group actions. This is a necessary step, since many of the interesting problems in the theories of moduli spaces and classical invariants to which we would we like to apply geometric invariant theory involve this type of action. There are some key aspects of reductive actions which make them much nicer to work with. For example, take an affine variety $X$ with a reductive $G$-action. If there are two disjoint $G$-invariant closed subvarieties $Y_1$ and $Y_2$ contained in $X$, then there exists an invariant $f$ in the ring $k[X]^G$ such that $f(Y_1) = 0$ and $f(Y_2) = 1$. As another example, if $X$ is a $G$-invariant affine subvariety of $\mathbb{A}^n$ such that $G$ acts reductively on $\mathbb{A}^n$, then given $f \in k[X]^G$, this $f$ can be extended to an invariant $F$ on $\mathbb{A}^n$ and if $I_X$ is the defining ideal of $X$, then

$$k[X]^G = {k[\mathbb{A}^n]^G}\Big/{I_X \cap \mathbb{A}^n}.$$

As another example, if $G$ is a reductive algebraic group acting on an affine variety $X = \mathrm{Spec}(A)$, then $A^G$ is a finitely generated $k$-algebra. For a final example, take $G$ to be a reductive group as above acting on an affine or projective variety. In this case, the quotient map $q: \mathrm{Spec}(A) \to \mathrm{Spec}(A^G)$ is a surjection, but it is possible to construct a counter-example where one has the quotient of an affine variety by a non-reductive group action which is not itself an affine variety (the above statement is effectively saying that in the reductive case the quotient would also have to be affine).

Most importantly, reductive geometric invariant theory often allows one to compute the open subsets made up of stable and semi-stable points, since in the case of a reductive group action on an affine variety, there always exists an invariant function which separates two disjoint closed $G$-invariant subvarieties of $X$. Although many properties fail for non-reductive actions, unipotent groups do have some nice features which give motivation for non-reductive geometric invariant theory apart from greater generality and applicability to a wider range of problems (unipotent can be read here as non-reductive). As examples, every orbit for a unipotent group action on a quasi-affine variety is closed, which is useful since the theory we discuss here does not distinguish between an orbit contained inside another one. Also, every homogeneous space of a unipotent group is isomorphic to an affine space and a connected unipotent group over a field of characteristic zero has no proper finite subgroups[216].

[215] Anthony Knapp, *Advanced Algebra* (New York: Birkhäuser, 2007).

Kirwan has also studied the links between geometric invariant theory and so-called moment maps in symplectic geometry. The theory of classical geometric invariants in the case of complex algebraic geometry is related to reduction in symplectic geometry. Start with a compact, connected Lie group $G$ with a Lie algebra $\mathfrak{g}$ which acts smoothly on a symplectic manifold $(M, \omega)$ by symplectomorphisms. The self-map $\psi_g: M \rightarrow M$ is a symplectomorphism for every $g \in G$ and one has that

$$\psi_{gh} = \psi_g \circ \psi_h, \;\; \psi_1 = \mathrm{Id}.$$

The infinitesimal action gives us a Lie algebra homomorphism defined as

$$X_\xi = \frac{d}{dt}\psi_{\exp(t\xi)},$$

for $\xi \in \mathfrak{g}$ when the derivative is evaluated at $t = 0$. As $\psi_g$ is a symplectomorphism, it follows that $X_\xi$ as defined above is always a symplectic vector field. The action of $G$ on $M$ is called weakly Hamiltonian if each of the vector fields $X_\xi$ is Hamiltonian. If we have a symplectic manifold and a smooth function $f \in C^\infty(M)$, the Hamiltonian vector field of $f$ is the vector field defined as

$$Xf = \widetilde{\omega}^{-1}(df),$$

where $\widetilde{\omega}$ is the bundle isomorphism from the tangent bundle to its dual. We could also write the function as $H$ to indicate that it is a Hamiltonian function. If $M$ is a closed manifold, then the vector field $XH$ generates a one-parameter family of diffeomorphisms $\phi_H^t$ such that

$$\frac{d}{dt}\phi_H^t = XH \circ \phi_H^t, \;\; \phi_H^0 = Id.$$

This is called the Hamiltonian flow associated with $H$ and it can be shown that $H$ is constant along the flow of $XH$ and that at every regular point of $H$ the Hamiltonian vector field is tangent to the level sets of $H$. For an example which can be visualised, take $H$ to be the height function on $S^2$. The level sets of this function (ie. the values where it is constant) are circles of constant height and the family of diffeomorphisms rotates every circle with constant speed, so the Hamiltonian flow simply rotates the sphere about its vertical axis (this is an example of a symplectic circle action). The action of $G$ is Hamiltonian if the map from $\xi$ to $H\xi$ can be chosen to be a homomorphism of Lie algebras with respect to the Lie algebra structure of $\mathfrak{g}$ and the Poisson structure on $C^\infty(M)$[217].

---

[216] Brent Doran and Frances Kirwan, 'Towards non-reductive geometric invariant theory', *Pure and Applied Mathematics Quarterly* 1 (3) 61 – 105 (2007).

[217] Dusa McDuff and Dietmar Salamon, *Introduction to Symplectic Topology* (New York: Oxford University Press, 1998).

A moment map for the action of $G$ on $M$ is a smooth map $\mu$ from $M$ to $\mathfrak{g}^*$ which is equivariant with respect to this action and the coadjoint action of $G$ on $\mathfrak{g}^*$. The map must also satisfy the relation

$$d\mu(x)(v)\xi = \omega_x(v, X_\xi),$$

for $x$ a point in the manifold and $v \in T_x M$. The quotient $\mu^{-1}(0)/G$ has a symplectic structure and is known as the symplectic reduction of $M$ by the action of $G$ at $0$ (or alternatively, the symplectic quotient). There are certain situations where the type of quotient we have been considering in geometric invariant theory $X//K$ can be canonically identified with the symplectic quotient (this can happen in the case of a non-singular complex projective variety $X$ and $K$ a complex reductive group acting via a complex linear representation) [218]. A theorem of Kirwan says that if we take a symplectic manifold acted on by a Lie group with a moment map $\mu$, plus the symplectic quotient $\mu^{-1}(p)/G$ and the equivariant cohomology ring of $M$ denoted by $H_G^*(M)$, then the Kirwan map

$$H_G^*(M) \rightarrow H^*\left(\frac{\mu^{-1}(p)}{G}\right)$$

is a surjection on rational coefficients[219]. If you are not familiar with equivariant cohomology, think of cohomology where you can have group actions on the relevant spaces.

We will mention a few other results. If we take a compact Lie group $G$ and fix an inner product on the Lie algebra $\mathfrak{g}$, then given a Hamiltonian action of the group on a symplectic manifold with an associated moment map, the norm squared of that map defines a certain Morse stratification of the manifold by locally closed symplectic submanifolds such that the stratum to which any point in $M$ belongs is determined by the limiting behaviour of its trajectory under the gradient flow with respect to a compatible metric. Kirwan showed that one can construct symplectic quotients in a natural way for the action of the group on the unstable strata[220]. With Jeffrey, Kirwan proved a residue formula which enables one to evaluate elements of the cohomology ring $H^*\left(\frac{\mu^{-1}(0)}{G}\right)$ whose degree is the dimension of $\mu^{-1}(0)/G$, provided that $0$ is a regular value of the moment map. Techniques used in the proof of this formula enabled a new proof of Witten's nonabelian localization formula for Hamiltonian actions of compact groups on symplectic manifolds[221].

---

[218] Frances Kirwan, 'Quotients by non-reductive algebraic group actions' (2008), arXiv:0801.4607v4.

[219] Frances Kirwan, *Cohomology of Quotients in Complex and Algebraic Geometry* (Princeton: Princeton University Press, 1984).

[220] Frances Kirwan, 'Symplectic quotients of unstable Morse strata for normsquares of moment maps' (2018), arXiv:1802.09237v1.

[221] Lisa Jeffrey and Frances Kirwan, 'Localization for nonabelian group actions'. (Balliol College, Oxford: Technical Report, 1993).

# *Leila Schneps*

Schneps is a mathematician working on various aspects of analytic number theory and Galois theory. Galois theory is a body of theory which allows us to say that the following polynomial equation

$$x^5 - 6x + 3 = 0,$$

cannot be solved with radicals, a statement which is not immediately obvious. Given a polynomial $f$ with coefficients in $\mathbb{Q}$, Galois theory provides a splitting field of $f$ which is the smallest possible field which contains all the zeroes of the polynomial. The splitting field has an associated finite Galois group $G$. Galois theory connects field theory and group theory in a profound way, and translates the problem of solubility of $f$ into the problem of solubility of the associated Galois group, and as most of us know, there are many techniques in group theory which would enable us to answer this question with relative ease.

This raises the question as to whether all finite groups are in fact Galois groups of an extension of the radicals: this is known as the inverse Galois problem. Noether proposed that mathematicians could make progress on the inverse problem by embedding $G$ in the permutation group $S_n$, one obtains a $G$-action on the field $\mathbb{Q}(X_1, \dots, X_n)$. If $E$ is the fixed field under this action, then it follows that $\mathbb{Q}(X_1, \dots, X_n)$ is a Galois extension of $E$ with Galois group $G$. Geometrically, the Galois extension of $E$ corresponds to a projection of varieties

$$\pi \colon \mathbb{A}^n \longrightarrow {}^{\mathbb{A}^n}\!/_{G},$$

for an affine $n$-space over $\mathbb{Q}$. If $P$ is a $\mathbb{Q}$-rational point of $\mathbb{A}^n/G$ which is lifted to $Q \in \mathbb{A}^n(\overline{\mathbb{Q}})$, then the conjugates of $\mathbb{Q}$ under the action of the absolute Galois group $\mathrm{Gal}(\overline{\mathbb{Q}}/\mathbb{Q})$ are the $sQ$ such that $s \in H_Q \subset G$, and $H_Q$ is the decomposition group at $Q$. If $H_Q$ is equal to the group $G$, then $Q$ generates a field extension of $\mathbb{Q}$ where the Galois group is $G$. A variety is rational over $\mathbb{Q}$ if it is birationally isomorphic over $\mathbb{Q}$ to $\mathbb{A}^n$, or if the function field is isomorphic to $\mathbb{Q}(T_1, \dots, T_n)$. It is a theorem of Hilbert's that if $\mathbb{A}^n/G$ is rational over $\mathbb{Q}$, then there is an infinite number of points $P$ such that $H_Q = G$. This is a corollary of Hilbert's irreducibility theorem.

We should mention that $\overline{\mathbb{Q}}$ is an infinite algebraic extension of $\mathbb{Q}$ and that it is, in a sense which can be directly quantified, very large. It is defined as

$$\overline{\mathbb{Q}} = \{\alpha \in \mathbb{C} \text{ such that } \alpha \text{ is algebraic over } \mathbb{Q}\}.$$

This is a very large object, because we have that

$$\mathbb{Q}\left(\sqrt{2}, \sqrt{3}, \sqrt{5}, \dots\right) \subseteq \overline{\mathbb{Q}},$$

but we must also have

$$\mathbb{Q}\left(\sqrt{2}, \sqrt[2]{2}, \sqrt[4]{2}, \sqrt[5]{2}, \dots\right) \subseteq \overline{\mathbb{Q}},$$

and so on. It is possible to prove that $\overline{\mathbb{Q}}$ is countably large. You probably know already that $\mathbb{Q}$ is countable. We know that

$$\overline{\mathbb{Q}} = \bigcup_{f \in \mathbb{Q}[x]} \{\text{roots of } f\},$$

because any algebraic number is the zero of some polynomial. Now, we know that $\mathbb{Q}[x]$ is countable.

$$\mathbb{Q}[x] = \bigcup_{n=0}^{\infty} \mathbb{Q}[x]_{\leq n},$$

where $\mathbb{Q}[x]_{\leq n}$ is the polynomial of degree $\leq n$. We also have that

$$\mathbb{Q}[x]_{\leq n} \simeq \mathbb{Q}^{n+1},$$

under the isomorphism given explicitly by

$$a_0 + a_1 x + \cdots + a_n x^n \longmapsto [a_0, a_1, \ldots, a_n].$$

We know that $\mathbb{Q}^N$ is countable for any $N$, which implies that $\mathbb{Q}[x]$ is finite. By the fundamental theorem of arithmetic, the set of zeroes of $f$ is finite with size equal to the degree of $f$ unless $f = 0$, so $\overline{\mathbb{Q}}$ is a countable union of finite sets, hence it is countable. One can also show that $[\overline{\mathbb{Q}} : \mathbb{Q}]$ is countably infinite. To start with, it must be the size of a basis for $\overline{\mathbb{Q}}/\mathbb{Q}$ and any basis is countable in $\overline{\mathbb{Q}}$, which implies that $[\overline{\mathbb{Q}} : \mathbb{Q}]$ is countable. All that remains is to show that $[\overline{\mathbb{Q}} : \mathbb{Q}]$ is not finite. For another exercise, one could prove that the linear combination $\sqrt{2} + \sqrt{3}$ is algebraic.

As a similar exercise, take the polynomial in $\mathbb{Q}[x]$ given by

$$f_n(x) = x^3 + x + n,$$

for any integer $n$. We would like to prove that $f_n(x)$ has only one zero over $\mathbb{R}$. We already know that $f_n(x)$ has at most 3 zeroes in $\mathbb{C}$. By the fundamental theorem of arithmetic, if $\alpha \in \mathbb{C}$, then $f_n(\alpha) = 0$ implies that $f_n(\bar{\alpha}) = 0$. Since $f_n$ is of degree 3, that implies that at least one solution is real. This follows from elementary analysis. If we think of $f_n$ as a function over $\mathbb{R}$, then

$$f_n'(x) = 3x^2 + 1$$

implies that $f_n(x)$ is monotonically increasing. This in turn implies that there exists one intersection with the $y = 0$ axis and the result follows.

We might also like to show that the polynomial is irreducible over $\mathbb{Q}$. By the Gauss lemma, this is equivalent to proving that the polynomial is irreducible over $\mathbb{Z}$. Start by writing out the polynomial.

$$f_n(x) = (x - c)(x^2 + ax + b).$$

The polynomial only splits in this way, as there is only one real zero. If $c$ is integer, then the polynomial is irreducible.

$$f_n(x) = x^3 + (a-c)x^2 + (b-ac)x - bc = x^3 + x + m.$$

We now compare the coefficients.

$$a - c = 0 \implies a = c,$$

$$b - ac = 1 \implies b = 1 + c^2,$$

$$bc = m \implies m = -(c^3 + c).$$

We conclude that if $m = -(c^3 + c)$ for some $c \in \mathbb{Z}$, then

$$f_n(x) = (x-c)\big(x^2 + cx + (1 + c^2)\big).$$

It follows that the polynomial is irreducible over $\mathbb{Z}$, and hence over $\mathbb{Q}$.

The irreducibility theorem states if $K$ is a number field, then for every $n$, the affine space $\mathbb{A}^n$ (or the projective space $\mathbb{P}^n$) has the Hilbert property over $K$. A (quasi-projective) variety $X$ over a field $K$ is said to satisfy the Hilbert property if $X$ is not thin, where a subset $A$ of $X(K)$ is 'thin' if it is contained in a finite union of sets of type $(C_1)$ or $(C_2)$. You can also say that the set is thin if there exists a morphism

$$\pi: V \longrightarrow W,$$

where the dimension of $W$ is less than or equal to the dimension of $V$ and the morphism has no rational cross-sections. Furthermore, we must have[222]

$$A \subset \pi\big(W(K)\big).$$

We earlier mentioned the definition of a field extension in terms of it being a homomorphism (a monomorphism, to be more specific). For example, $\mathbb{R}/\mathbb{Q}$ is a field extension where the homomorphism is just the trivial inclusion map. $\mathbb{C}/\mathbb{R}$ and $\mathbb{C}/\mathbb{Q}$ are also field extensions with the inclusion map as the defining homomorphism. Now if we take a field $K$ and a non-empty subset of that field denoted by $S$. The subfield of $K$ generated by $S$ is defined to be the intersection of all the subfields of $K$ which contain $S$. For example, the subfield of $\mathbb{R}$ generated by the single element $\{1\}$ is $\mathbb{Q}$ and the subfield of $\mathbb{C}$ generated by $\{i\}$ is $\mathbb{Q}(i)$. $\mathbb{Q}(i)$ denotes an adjoining of a root. If $L/K$ is a field extension and $A$ is contained in $L$, then we generally write $K(A)$ for the subfield of $L$ generated by the union $K \cup A$. This is known as the field obtained by adjoining $A$ to $K$. The adjoining of one root is the building-block for field extensions. It can be shown that

$$\mathbb{Q}\big(\sqrt{2}\big) = \big\{a + b\sqrt{2}\big\},$$

$$\mathbb{Q}\big(\sqrt[3]{2}\big) = \Big\{a + b\sqrt[3]{2} + c\big(\sqrt[3]{2}\big)^2\Big\},$$

where all the coefficients are rational. A field extension $L/K$ is said to be a simple extension if we can write

$$L = K(\alpha),$$

[222] Jean-Pierre Serre, *Topics in Galois Theory* (Boston: Jones and Bartlett Publishers, 1992).

for $\alpha \in L$.

It can proved that if $K$ is a field and $f$ in $K[x]$ is an irreducible polynomial, where $L \coloneqq K[x]/f(x)$, then one has that the following are equivalent: $L$ is a field, the map from $K$ to $L$ given by $a \longmapsto a + (f)$ is a field extension, the element $x + (f)$ in $L$ is a root of $f$, and $L = K(\alpha)$ where $\alpha = x + (f)$. If $K$ is a field and $f \in K[x]$ is a polynomial of degree $n$ where $\alpha_1, \dots, \alpha_n$ are the roots of $f$ in some extension of $K$, then $K(\alpha_1, \dots, \alpha_n)$ is the splitting field of $f$ over $K$. The splitting field can be viewed as the smallest field over which $f$ can be split as a product of linear factors. As an example, take the polynomial

$$f = x^4 - 3x^3 + 2x^2.$$

This can be factored as

$$f = x^2(x-1)(x-2).$$

This implies that the splitting field for $f$ over $\mathbb{Q}$ is $\mathbb{Q}$ adjoined with 0, 1 and 2, which is just $\mathbb{Q}$. It can also be proved that splitting fields always exist and that they are unique. If $K$ is a field and $f \in K[x]$ is a polynomial, then a splitting field $L/K$ for $f$ exists. Also, if $L_1/K$ and $L_2/K$ are both splitting fields for $f$, then there is a field isomorphism $\iota\colon L_1 \longrightarrow L_2$ such that

$$\iota(a) = a,$$

for all $a \in K$.

We should also define the notion of the degree of a field extension. Start by noting that for $L/K$ a field extension, $L$ is a vector space over the field $K$ (recall that in linear algebra, vector spaces are always defined over fields). You might wish to show to yourself that all the usual axioms for a vector space are satisfied. This kind of strange fusion between linear algebra and field theory is typical of Galois theory. If $L/K$ is an extension, then the degree of $L/K$ (written as $[L : K]$) is defined to be the dimension of $L$ as a vector space over $K$. The extension is said to be finite if the degree is finite, which obviously implies that $L$ is a finite-dimensional $K$-vector space. If not, then the extension must be infinite. For example, if $L/K$ is an extension such that $[L : K] = p$ for some prime $p$, then the extension is simple if and only if $L = K(\alpha)$ for some $\alpha \in L$. To prove this, start by considering the degree. We already know that

$$[L : K] = p > 1,$$

which implies that $L \neq K$, so there must exist an $\alpha$ in $L \backslash K$. Consider the chain of inclusions:

$$K \subseteq K(\alpha) \subseteq L.$$

From the Tower Law (which we will get to shortly), we must have

$$p = [L : K] = [L : K(\alpha)][K(\alpha) : K].$$

This is a multiplication of two factors, so one of the degrees has to be unity and the other has to be $p$. If $[L : K(\alpha)] = p$, then

$$[K(\alpha) : K] = 1.$$

This is a contradiction, as that would imply that $K = K(\alpha)$, which implies that $\alpha \in K$, but we have already said that $\alpha \in L$. If, on the other hand, $[L : K(\alpha)] = 1$, then we have

$$L = K(\alpha),$$

which is what we wanted to prove. As an aside, note that every finite extension of field extensions of $\mathbb{Q}$ is simple.

For some more definitions, start once more with a field extension $L/K$ with $\alpha \in L$. $\alpha$ is said to be algebraic over $K$ if there is some non-zero polynomial $f$ in $K[x]$ such that

$$f(\alpha) = 0.$$

$\alpha$ is transcendental over $K$ if it is not algebraic. Specifically, a complex number is said to be an algebraic number if it is algebraic over $\mathbb{Q}$. In other words, it is the zero of some polynomial with rational coefficients. A complex number is called a transcendental number if it is not an algebraic number. Although it can be very hard to prove that a number is transcendental (think about $\pi$, for example), the complex plane is almost entirely composed of transcendental numbers, and it is relatively rare for a number to be algebraic. If you do not believe this, you could try looking up some pictures which provide a visualisation of the algebraic numbers in the complex plane compared to the transcendental ones. Note that this is distinct from saying that $\pi$ is irrational. The proof in either case is very long, but Cartwright noticed that Hermite's proof for the irrationality of $\pi$ can be greatly simplified. The first example of an algebraic number which will probably occur to you is $\sqrt{2}$. Again, we emphasise that this is separate from rationality or irrationality, as $\sqrt{2}$ is certainly not a rational number. This should also make it immediately clear to you that the real line also contains algebraic numbers.

It can be proved quite easily that finite field extensions are also algebraic, but the converse is not true. As an example, consider the extension given by

$$\mathbb{Q}(\sqrt{2}, \sqrt{3}, \sqrt{5}, \dots)/\mathbb{Q}.$$

This is an algebraic extension with infinite degree. Now if we start with a field extension $L/K$ such that $\alpha \in L$ is algebraic over $K$, then the minimal polynomial of $\alpha$ over $K$ is the monic polynomial $m \in K[x]$ of smallest degree such that

$$m(\alpha) = 0.$$

We would like to actually do something with these definitions, so let us make some statements for a field extension $L/K$. If $\alpha \in L$ is algebraic over $K$, then the minimal polynomial exists and is unique. The minimal polynomial over $\alpha$ is the unique monic irreducible polynomial $m \in K[x]$ which satisfies $m(\alpha) = 0$. If $f \in K[x]$ satisfies $f(\alpha) = 0$, then $m$ divides $f$. This lemma means that we can determine whether or not a polynomial is minimal, but to do this we have to check that it is irreducible. There are some nice tests for irreducibility of polynomials. We mentioned Gauss's lemma: in fact, this lemma is true whenever you have a polynomial over a unique factorisation domain, not just the integers.

Another test is Eisenstein's criterion. This states that if $p$ is a prime number and

$$f = a_n x^n + \cdots + a_1 x + a_0 \in \mathbb{Z}[x]$$

a polynomial which satisfies

$$p \nmid a_n,$$

$$p \text{ divides } a_i \text{ for } i = 0, 1, \dots, n-1,$$

$$p^2 \nmid a_0.$$

If these criteria are met, then the polynomial is irreducible over $\mathbb{Q}$. We are now moving towards a basic result known as the Tower Law. Before this, we will state a few other results. If $\alpha$ is algebraic over $K$ with minimal polynomial $m \in K[x]$ and $(m)$ is the principal ideal generated by $m$, then the map given by

$$\tilde{\phi} : {}^{K[x]}\!\big/_{(m)} \longrightarrow K(\alpha),$$

$$\tilde{\phi}\big(f + (m)\big) = f(\alpha),$$

is an isomorphism. Also, if $L/K$ is a field extension and let $\alpha \in L$ be algebraic over $K$. Suppose that the minimal polynomial $m$ of $\alpha$ over $K$ has degree $d$, then

$$K(\alpha) = \{a_0 + a_1\alpha + \cdots + a_{d-1}\alpha^{d-1} \text{ for } a_0, \dots, a_{d-1} \in K].$$

Furthermore, a basis for $K(\alpha)$ over $K$ is given by $1, \alpha, \dots, \alpha^{d-1}$. Also,

$$[K(\alpha) : K] = d.$$

To finish with the Tower Law, let us say we start with a chain of field extensions of finite degree $K \subseteq L \subseteq M$. Let $l_1, l_2, \dots, l_r$ be a basis for $L/K$ and $m_1, m_2, \dots, m_r$ a basis for $M/L$. It follows that $\{l_i m_j\}$ is a basis for $M/K$, where $i$ runs from 1 to $r$ and $j$ runs from 1 to $s$. Also, the degree of the extensions can be decomposed as a product:

$$[M : K] = [M : L][L : K].$$

The proof is essentially linear algebra. There is much more to say about Galois theory, but hopefully this has been good motivation. I should say that the modern (as opposed to classical) approach to Galois theory is due to Grothendieck and is based on abstract category theoretic properties, rather than linear algebra[223]. You can also go much further with the inverse Galois problem and link it to some of the other areas which we have discussed in the book (moduli spaces and mapping class groups, for example)[224].

Schneps has also researched Lie algebras, especially an elliptic version of a Lie algebra called the Kashiwara-Vergne Lie algebra[225]. We earlier gave the basic definition of a Lie algebra in terms of vector spaces. One can define notions of a Lie subalgebra as one would expect. A linear map between two Lie algebras $\phi: \mathfrak{g} \to \mathfrak{h}$ is called a Lie algebra homomorphism if

$$\phi[X, Y] = [\phi(X), \phi(Y)].$$

[223] Ian Stewart, *Galois Theory* (USA: Chapman & Hall, 2004).
[224] Leila Schneps and P. Lochak, *The Inverse Galois Problem, Moduli Spaces and Mapping Class Groups* (Cambridge: Cambridge University Press, 1997).
[225] Leila Schneps, 'Double shuffle and Kashiwara-Vergne Lie algebras', *J. Algebra* 367, 54 – 74 (2012).

If the map is a bijection, then $\phi$ is a Lie algebra isomorphism, and if it is a self-map, then it is a Lie algebra automorphism. If $\mathfrak{g}$ is a Lie algebra and $X$ is an element of $\mathfrak{g}$, then the linear self-map of Lie algebras defined by

$$\mathrm{ad}_X(Y) = [X, Y]$$

is known as the adjoint map or the adjoint representation. This is obviously just the Lie bracket again, so if one is composing $n$ Lie brackets, one can write $(\mathrm{ad}_X)^n(Y)$. The adjoint map can also be viewed as the linear map from $\mathfrak{g}$ to the space of linear operators on $\mathfrak{g}$. A computation shows that

$$\mathrm{ad}_{[X,Y]} = \mathrm{ad}_X \mathrm{ad}_Y - \mathrm{ad}_Y \mathrm{ad}_X = [\mathrm{ad}_X, \mathrm{ad}_X].$$

This means that the adjoint map is a Lie algebra homomorphism.

If $\mathfrak{g}_1$ and $\mathfrak{g}_2$ are Lie algebras, then the direct sum of $\mathfrak{g}_1$ and $\mathfrak{g}_2$ is the direct sum in terms of vector spaces, where the bracket is given as you might expect:

$$[(X_1, X_2), (Y_1, Y_2)] = ([X_1, Y_1], [X_2, Y_2]).$$

It can be checked easily that this direct sum is another Lie algebra with respect this Lie bracket. If $\mathfrak{g}$ is a finite-dimensional real or complex Lie algebra and $X_1, \ldots, X_N$ is a basis for $\mathfrak{g}$ viewed as a vector space, then one has unique structure constants such that

$$[X_j, X_i] = \sum_{k=1}^{N} c_{jik} X_k.$$

These are known as structure constants and are used in theoretical physics. The structure constants satisfy two conditions which are derived from the skew symmetry of the Lie bracket and the Jacobi identity.

As with any mathematical object, one is usually interested in some special types which have nice properties. For example, a Lie algebra is irreducible if the only ideals in $\mathfrak{g}$ and $\{0\}$. In Lie algebra theory, a Lie subalgebra $\mathfrak{h}$ of a Lie algebra $\mathfrak{g}$ is an ideal in $\mathfrak{g}$ if $[X, H] \in \mathfrak{h}$ for all $X \in \mathfrak{g}$ and $H \in \mathfrak{h}$. A Lie algebra is simple if it is irreducible and $\dim \mathfrak{g} \geq 2$. As an example, a one-dimensional Lie algebra cannot have any ideals apart from itself or the trivial ideal. It is an irreducible Lie algebra, but it cannot be a simple Lie algebra because of the dimensionality. Less trivially, we can show that $\mathfrak{sl}(2, \mathbb{C})$ is a simple Lie algebra. If we take the Lie algebra to be the space of $3 \times 3$ upper triangular matrices with zeroes on the leading diagonal, then this Lie algebra is nilpotent. If $\mathfrak{g}$ is the space of $2 \times 2$ matrices of the form

$$\begin{pmatrix} a & b \\ c & 0 \end{pmatrix}$$

with $a$, $b$ and $c$ complex, then the Lie algebra is solvable but not nilpotent.

As we have already said, a Lie algebra can always be associated with a Lie group, and vice versa. A good starting point for this would be to consider the Lie algebra associated to a matrix Lie group. It is fairly typical to focus on the matrix Lie groups when beginning the study of Lie theory, as these groups are much simpler to get a handle on. A problem in Lie group theory can be transferred to the corresponding Lie algebra, and as we have emphasised, Lie algebras can be studied with elementary linear algebra techniques which we all know. For example, the Lie algebra of Euclidean space is just Euclidean space again, because it turns out to be abelian and is isomorphic to $\mathbb{R}^n$ with the trivial Lie bracket. The Lie algebra of the circle group is isomorphic to the real line, and the Lie algebra of the torus $\mathbb{T}^n$ is also isomorphic to $\mathbb{R}^n$, essentially because the torus is equal to a product of circles.

Take $G$ to be a matrix Lie group. The Lie algebra $\mathfrak{g}$ of $G$ is the set of all matrices $X$ such that $e^{tX} \in G$ for all real $t$. If we take the general linear group, there is some ambiguity because there is a difference between $\mathrm{Lie}\big(\mathrm{GL}(n, \mathbb{R})\big)$ and $\mathfrak{gl}(n, \mathbb{R})$, since they have different Lie algebras coming from a Lie bracket of vector fields in the first case and from a commutator bracket of matrices in the second case (ie. it is a matrix Lie algebra). However, one can prove that there is an explicit isomorphism between the Lie algebras given by a composition of natural maps

$$\mathrm{Lie}\big(\mathrm{GL}(n, \mathbb{R})\big) \longrightarrow T_{I_n}\mathrm{GL}(n, \mathbb{R}) \longrightarrow \mathfrak{gl}(n, \mathbb{R}).$$

The first isomorphism is a vector space between the Lie algebra and the tangent space to the Lie group at the identity matrix. However, the general linear group is an open subset of $\mathfrak{gl}(n, \mathbb{R})$, so there is an isomorphism from the tangent space to the matrix Lie algebra. If $G$ is a matrix Lie group with a corresponding matrix Lie algebra $\mathfrak{g}$ and $X$ and $Y$ are elements of $\mathfrak{g}$, then the following are true:

$$AXA^{-1} \text{ for all } A \in G,$$

$$sX \in \mathfrak{g} \text{ for all } s \in \mathbb{R},$$

$$X + Y \in \mathfrak{g},$$

$$XY - YX \in \mathfrak{g}.$$

It follows that the matrix Lie algebra of a matrix Lie group is a real Lie algebra where the bracket is given by what we could call the commutator:

$$[X, Y] = XY - YX.$$

A matrix Lie group is complex if its associated Lie algebra is a complex subspace of $M_n(\mathbb{C})$. The general and special linear groups over the complex numbers are obvious examples of complex matrix Lie groups.

To give some examples, the Lie algebra of $\mathrm{GL}(n, \mathbb{C})$ is the space of all $n \times n$ matrices with complex entries, the Lie algebra of $\mathrm{GL}(n, \mathbb{C})$ is the same with real entries, the Lie algebra of $\mathrm{SL}(n, \mathbb{C})$ is the space of all $n \times n$ complex matrices with trace zero, and can guess what the Lie algebra of $\mathrm{SL}(n, \mathbb{R})$ is. The Lie algebras of these groups are denoted with the usual gothic letters. Continuing in this way, the Lie algebra of the unitary group $\mathrm{U}(n)$ consists of the space of complex matrices such that one has for the adjoint

$$X^* = -X$$

and the Lie algebra of $\mathrm{SU}(n)$ is the subspace of the previous space where one takes the traceless matrices. The Lie algebra of the orthogonal group $\mathrm{O}(n)$is the space of real matrices such that

$$X^{\mathrm{Tr}} = -X$$

and the Lie algebra of the special orthogonal group is the same again. An interesting finitely generated group is the discrete Heisenberg group $H$. It can be defined as the following array for integer $x$, $y$ and $z$:

$$\begin{array}{ccc} 1 & x & z \\ 0 & 1 & y. \\ 0 & 0 & 1 \end{array}$$

The Lie algebra of this group is the space of all matrices of the form

$$\begin{array}{ccc} 0 & a & b \\ 0 & 0 & c, \\ 0 & 0 & 0 \end{array}$$

where a, $b$ and $c$ are real numbers[226].

We might also wonder if the link between Lie groups and Lie algebras transfers over to the maps between them, and this turns out to be the case: a Lie group homomorphisms between two Lie groups induces a map between the corresponding Lie algebras. It can be shown as a consequence that isomorphic Lie groups must have isomorphic Lie algebras. Specifically, if $G$ and $H$ are matrix Lie groups with respective Lie algebras $\mathfrak{g}$ and $\mathfrak{h}$. If $\Phi$ is a Lie group homomorphism, then there exists a real-linear unique map of Lie algebras such that

$$\Phi(e^X) = e^{\phi(X)}.$$

This map has several properties:

$$\phi(AXA^{-1}) = \Phi(A)\phi(X)\Phi(A)^{-1},$$

$$\phi([X, Y]) = [\phi(X), \phi(Y)],$$

$$\phi(X) = \frac{d}{dt}\Phi(e^{tX}).$$

Linked with her studies of Kashiwara-Vergne Lie algebras, Schneps has studied the Grothendieck-Teichmüller group[227]. The elliptic Kashiwara-Vergne Lie algebra has some properties which are similar to the Grothendieck-Teichmüller Lie algebra $\mathfrak{grt}$, where $\mathfrak{grt}$ is the space of polynomials $b \in \mathfrak{lie}_2$ which satisfies the pentagon relation, equipped with the Poisson (or Ihara) bracket

$$\{b, b'\} = [b, b'] + d_b(b') - d_{b'}(b).$$

[226] Brian C. Hall, *Lie Groups, Lie Algebras, and Representations: An Elementary Introduction* (Switzerland: Springer, 2015).

[227] Pierre Lochak and Leila Schneps, 'Open problems in Grothendieck-Teichmüller theory', *Proceedings of Symposia in Pure Mathematics*, 75, 165 – 186 (2006).

The Grothendieck-Teichmüller group $GT(\mathbb{Q})$ can be defined as the group formed by the automorphisms of the unitary operad in Malev complete groupoids $PaB_+^{\wedge}$ which reduce to the identity map on the object sets of the operad. Grothendick-Teichmüller theory was created the aim of studying the automorphism groups of the profinite mapping class groups, where these groups are the fundamental groups of moduli spaces of Riemann surfaces of all genera and any number of marked points. The aim is to find new properties of $\mathrm{Gal}(\overline{\mathbb{Q}}/\mathbb{Q})$, or even ultimately to show $GT(\mathbb{Q})$ is equal to $\mathrm{Gal}(\overline{\mathbb{Q}}/\mathbb{Q})$[228]

The set of Drinfeld associators, denoted by $Ass(\mathbb{Q})$, is in bijection with the set of categorical equivalences of unitary operads in Malcev complete groupoids

$$v^{\wedge}: PaB_+^{\wedge} \longrightarrow CD_+^{\wedge}$$

from the completed parenthesized braid operad to the operad of chord diagrams. A Drinfeld associator can be interpreted as a morphism of operads in groupoids. There is a natural action of $GT(\mathbb{Q})$ on the set of Drinfeld associators $Ass(\mathbb{Q})$ by translation on the right. Explicitly, there is a morphism

$$v: PaB \longrightarrow CD^{\wedge}$$

which represents an element of $Ass(\mathbb{Q})$. The morphism

$$v \circ \gamma: PaB \longrightarrow CD^{\wedge}$$

represents the action of an element of $GT(\mathbb{Q})$ on $v$ and can be defined as a composite where the morphism $v$ is extended to the completion of the parenthesized braid operad[229].

[228] Leila Schneps, 'The Grothendieck-Teichmüller group: a survey', in Pierre Lochak and Leila Scneps (eds.), *Geometric Galois Actions, 1* (Cambridge: Cambridge University Press, 1997).
[229] Benoit Fresse, *Homotopy of Operads and Grothendieck-Teichmüller Groups: The Applications of Rational Homotopy Theory Methods* (Providence. Rhode Island: American Mathematical Society, 2017).

# *Claire Voisin*

Voisin has received many awards for her mathematical achievements and is known for her work in algebraic geometry, especially relating to Hodge theory. This body of theory has many applications. For example, given the moduli space of rank $k$ stable holomorphic vector bundles over a complex Riemann surface, one knows via Hodge theory that the moduli space is equipped with a canonical Kähler metric. The Riemann-Roch theorem gives us a formula for the dimension of this moduli space, but as mentioned earlier, it might not be possible to actually use that formula in practice[230]. A Kähler metric is a Riemannian metric on a complex manifold ie. one has a Kähler manifold with a Hermitian metric

$$h = g - i\omega,$$

where $\omega$ is a Kähler form defined by

$$g(X, Y) = \omega(X, JY).$$

$J$ is an almost-complex structure on the manifold. A structure of this kind is a smooth choice of linear maps on each tangent space

$$J_X: T_X M \to T_X M$$

such that

$$J_X^2 = \mathrm{Id}.$$

In other words, the tangent spaces have the structure of complex vector spaces where $J$ acts as multiplication by $i$. It can be shown that if a manifold $M$ has an almost-complex structure, then it must have a $\mathrm{Spin}^c$ structure. Similar constructions are possible when one works with almost-quaternionic structures, where the tangent spaces have the structure of quaternionic vector spaces.

For $r \geq 2$, fix an orthonormal frame $e_1, \ldots, e_r \in \mathbb{R}^r$. The Clifford algebra $\mathrm{Cl}_r$ is the algebra generated by $e_1, \ldots, e_r$ subject to the following relations

$$e_i^2 = -1, \qquad e_i e_j = -e_j e_i.$$

The even Clifford algebra is the subalgebra generated by the even products:

$$\mathrm{Cl}_r^0 = \langle e_i e_j \big| i < j \rangle.$$

The concept of an almost-complex structure can be generalised by saying that a rank $r$ almost-even Clifford structure on a manifold $M$ is a smooth choice of algebra representations:

$$\Phi: \mathrm{Cl}_r^0 \to \mathrm{End}(T_x M),$$

[230] Edward Witten, 'The Verlinde algebra and the cohomology of the Grassmannian', in Shing-Tung Yau (ed.) *Geometry, Topology, and Physics for Raoul Bott* (Cambridge, Massachusetts: International Press, 1994).

$$e_i \cdot e_j \mapsto \Phi(e_i, e_j) = J_{ij}.$$

The endomorphisms in this case are skew symmetric. A natural question to ask is whether $M$ having a rank $r$ almost-even Clifford structure implies that $M$ must have a $\mathrm{Spin}^r(n)$ structure. It can be shown that a manifold having this type of structure plus a large automorphism group implies that the manifold is always a spin manifold when $r = 5, 7$ and that around half of all manifolds are spin when $r = 3, 4, 6, 8$, though the proof is technical[231]. A $\mathrm{Spin}^r$ structure is just an appropriate generalisation of a $\mathrm{Spin}^c$ structure, given that not every manifold will have a $\mathrm{Spin}^c$ structure and so the group needs to be enlarged. This can be done by considering the twisted spin group:

$$\mathrm{Spin}^r(n) = \frac{\mathrm{Spin}(n) \times \mathrm{Spin}(r)}{\langle(-1, -1)\rangle}.$$

The corresponding theorems are analogous: just replace $SO(2)$ with $SO(r)$, and so on. Every orientable Riemannian manifold has one of these structures.

Hodge theory is a way of studying cohomology groups on smooth manifolds. If one wishes to summarise it briefly, one can see it is a theory of harmonic forms. Like many things, the theory actually originated in mathematical physics. The de Rham theorem in smooth manifold theory states if $M$ is a smooth manifold and $p$ a non-negative integer, then the de Rham cohomology group is isomorphic to the corresponding singular cohomology group with coefficients in $\mathbb{R}$ (see Bredon for a nice proof)[232].

$$H^p_{\mathrm{dR}} \cong H^p(M, \mathbb{R}).$$

This map is given by

$$\omega \longmapsto \left[\sigma \longmapsto \int_\sigma \omega\right].$$

If you do not know what a de Rham cohomology group is, take a form $\varphi$ which vanishes when the exterior derivative is applied to it: such a form is said to be closed. Consider a $p$-form $\varphi$ which is the exterior derivative of a $(p-1)$-form $\eta$

$$\varphi = \mathrm{d}\eta.$$

It can be proved that the square of the exterior derivative is always zero. Obviously it is equivalent to show that the result is 0 when acting on an arbitrary $p$-form, so put everything in local coordinates.

$$\begin{aligned}\mathrm{d}\varphi &= \frac{1}{p!}\left(\mathrm{d}\varphi_{i_1 \dots i_p}\right) \wedge \mathrm{d}x^{i_1} \wedge \dots \wedge \mathrm{d}x^{i_p}, \\ &= \frac{1}{p!}\frac{\partial \varphi_{i_1 \dots i_p}}{\partial x^j} \mathrm{d}x^j \wedge \mathrm{d}x^{i_1} \wedge \dots \wedge \mathrm{d}x^{i_p}.\end{aligned}$$

[231] Gerardo Arizmendi, Rafael Herrera and Noemi Santana, 'Almost even-Clifford hermitian manifolds with large automorphism group', *Proceedings of the American Mathematical Society* 144(9), 4009 – 4020 (2016).
[232] Glen Bredon, *Topology and Geometry* (New York: Springer-Verlag, 1993).

Act on the resulting form with the exterior derivative operator:

$$\mathrm{dd}\varphi = \frac{1}{p!}\frac{\partial^2 \varphi_{i_1 \dots i_p}}{\partial x^j \partial x^k} \mathrm{d}x^k \wedge \mathrm{d}x^j \wedge \mathrm{d}x^{i_1} \wedge \dots \wedge \mathrm{d}x^{i_p}.$$

Mixed partial derivatives commute.

$$\frac{\partial^2 \varphi_{i_1 \dots i_p}}{\partial x^j \partial x^k} = \frac{\partial^2 \varphi_{i_1 \dots i_p}}{\partial x^k \partial x^j}.$$

However, the exterior product anti-commutes:

$$\mathrm{d}x^k \wedge \mathrm{d}x^j = -\mathrm{d}x^j \wedge \mathrm{d}x^k.$$

This means that the contraction over the $j$ and $k$ indices will vanish, but that implies that the whole thing vanishes.

$$\mathrm{dd}\varphi = 0.$$

This in turn implies that every exact form is closed. Quotient the set of closed $p$-forms on $M$ through by the set of exact $p$-forms and you have the $p$th de Rham cohomology group. For a closed form $\omega$ on $M$, the equivalence class of $\omega$ in this quotient space is called the cohomology class of $\omega$.

In words, every real cohomology class on a smooth manifold can be represented by a closed form on that manifold, and two closed forms represent the same cohomology class if they differ only by an exact form. This is still considered to be a remarkable result many years after it was proved. We mentioned earlier that the Maxwell equations can be written in terms of a 2-form. In the homogeneous case, the equations can actually be written as

$$\mathrm{d}F = 0,$$

$$\mathrm{d} * \mathrm{F} = 0.$$

The relation $F$ and its Hodge operator is invariant under the group of Lorentz transformations. It was by considering these equations that Hodge first defined the Hodge duality operator and came upon the Hodge theorem. As you probably know, an equivalence class can have many representatives, and there is a question as to whether there is some canonical choice for the form $\omega$ which represents a cohomology class. We might imagine the set $\{\omega + \mathrm{d}\tau\}$ for some exact form $\mathrm{d}\tau$ and consider the 'smallest' member of this set as $\tau$ varies over all $p$-forms. This set could be defined pointwise and integrated over the entire manifold[233].

A smooth $p$-form $\omega$ on a compact Riemannian manifold is said to be harmonic if it vanishes when acted on by the exterior derivative operator and if

$$\|\omega\|_M \leq \|\omega + \mathrm{d}\tau\|_M.$$

[233] Mark Green, Jacob Murre and Claire Voisin, *Algebraic Cycles and Hodge Theory* (Berlin: Springer-Verlag, 1994).

The set of harmonic $p$-forms on the manifold is denoted with $\mathcal{H}^p(M)$. Define the Hodge Laplacian on forms as

$$\Delta: \Omega^p(M) \to \Omega^p(M),$$

$$\Delta\omega = (\mathrm{d}\delta + \delta\mathrm{d})\omega.$$

For $\omega$ a $p$-form, $\Delta\omega = 0$ if and only if $\mathrm{d}\omega = 0$ and $\delta\omega = 0$. This can be proved trivially using the definition of the operator, since

$$(\Delta\omega, \omega) = (\mathrm{d}\delta\omega, \omega) + (\delta\mathrm{d}\omega, \omega),$$

$$= (\delta\omega, \delta\omega) + (\mathrm{d}\omega, \mathrm{d}\omega).$$

The second step follows from the fact that

$$(\delta\omega, \phi) = (\omega, \mathrm{d}\phi),$$

where $\omega$ is a $p$-form and $\phi$ is a $(p-1)$-form: this can be proved using Stokes's theorem for manifolds. If we set $\Delta\omega = 0$, then this clearly implies that both $\delta\omega = 0$ and $\mathrm{d}\omega = 0$, which is one direction of the proof[234]. With some more work, one can show that both of these are equivalent to $\omega$ being a harmonic form.

In essence, we are looking to show that a sequence of smooth forms in a de Rham equivalence class whose sizes converge to the $\inf$ of the sizes of all the forms in the class also converges to a smooth form. This can be done with some elliptic PDE theory, and one arrives at a central result known as the Hodge theorem. This states that for a compact Riemannian manifold $M$ with orientation, the canonical map which sends a harmonic $p$-form to the corresponding de Rham equivalence class is an isomorphism of vector spaces.

$$\mathcal{H}^p(M) \to H^p_{\mathrm{dR}}(M).$$

This implies that every de Rham cohomology class has a unique harmonic form as a representative. This theorem can lead to some powerful geometric results, especially for Kähler manifolds. A related corollary is the Hodge decomposition theorem that there exists an isomorphism as follows:

$$H^k(M, \mathbb{C}) \cong \oplus_{p+q=k} H^{p,q}(M).$$

There are several special results associated with the cohomology of compact Kähler manifolds. For example, if $M$ is a compact Kähler manifold of dimension $n$, then $(-1)^q\langle , \rangle$ is positive definite on the space $L^j\left(P^{p,q}(M)\right)$ and the spaces are orthogonal. $\langle , \rangle$ is a Hermitian inner product on $\Omega^k(M)$ defined by

$$\langle \alpha,\beta \rangle = i^{k^2} \int_M \omega^{n-k} \wedge \alpha \wedge \bar{\beta}.$$

[234] Jürgen Jost, *Riemannian Geometry and Geometric Analysis* (Heidelberg: Springer, 2008).

$P^{p,q}(M)$ denotes a space of primitive differential forms defined as

$$P^{p,q}(M) = \{\Omega^{p,q}(M) \text{ such that } \Lambda\omega = 0\}.$$

Another example is the principle of two types. If we start with a compact Kähler manifold as before and let a form $\alpha \in \Omega^{p,q}(M)$ be such that $\partial\alpha = 0$, $\bar{\partial}\alpha = 0$ and let it be exact either under the differential operator $\partial$ or $\bar{\partial}$, then there exists $\lambda \in \Omega^{p-1,q-1}(M)$ such that

$$\alpha = \partial\bar{\partial}\lambda.$$

$\partial$ and $\bar{\partial}$ are a particular pair of operators which are defined uniquely by the equation[235]

$$\mathrm{d} = \partial + \bar{\partial}.$$

The proof of this cohomology result requires the Hodge theorem, plus some results from the theory of elliptic operators, specifically the existence of a unique Green's function $G$ on a compact orientated Riemannian manifold mapping $\Omega^k(M)$ to $\Omega^k(M)$ such that $G$ commutes with $\mathrm{d}$ and $\delta$, $G\left(\mathcal{H}^k(M)\right) = 0$, and

$$\mathrm{Id} = \pi_{\mathcal{H}} + \Delta G,$$

where $\pi_{\mathcal{H}}$ is the orthogonal projection.

One of the Millennium Prize problems is a problem in Hodge theory, known as the Hodge conjecture. This conjecture states that if $X$ is a non-singular complex projective manifold (ie. a smooth projective algebraic varietiy over $\mathbb{C}$), then every Hodge class (or Hodge cycle) is a linear combination of classes of algebraic cycles with coefficients in $\mathbb{Q}$. By a theorem of Chow, algebraic cycles on a complex projective variety are equivalent to closed analytic subspaces. There are also more generalized versions of the Hodge conjecture: for example, there is a version of the conjecture which takes $X$ to be a complex Kähler variety and states that every Hodge cycle on $X$ should be a rational linear combination of Chern classes of vector bundles on $X$. A different version is exactly the same, but with vector bundles replaced by coherent sheaves. Voisin proved this to be false by constructing a 4-dimensional complex torus with a non-trivial Hodge cycle of degree 4 such that any analytic coherent sheaf $\mathcal{F}$ on the variety has a trivial second Chern class. In the course of finding this counterexample, Voisin found a non-trivial criterion for the vanishing of the second Chern class of an analytic coherent sheaf on $X$. If $X$ is a compact Kähler variety of dimension $n \geq 3$ such that the Néron-Severi group generated by the first Chern classes of holomorphic line bundles on $X$ is the zero group, if $X$ does not contain any proper closed analytic subset of positive dimension, and if for Kähler equivalence class $[\omega] \in H^2(X, \mathbb{R}) \cap H^{1,1}(X)$, the set of Hodge cycles $H\ dg^4$ is perpendicular to $[\omega]^{n-2}$ for the intersection pairing

$$H^4(X, \mathbb{R}) \otimes H^{2n-4}(X, \mathbb{R}) \to \mathbb{R},$$

[235] Mark Green, Jacob Murre and Claire Voisin, *Algebraic Cycles and Hodge Theory* (Berlin: Springer-Verlag, 1994).

then any analytic coherent sheaf $\mathcal{F}$ on $X$ has a trivial second Chern class[236]. Hodge originally proposed a stronger version of the conjecture, but Grothendieck realised that this version was automatically false for trivial reasons[237].

We should take the opportunity at this point to define Hodge structures. A Hodge structure of weight $k$ is a real vector space $V$ together with a lattice $\Gamma \subseteq V$, together with a decreasing filtration $F^{\cdot} V_{\mathbb{C}}$ on $V_{\mathbb{C}} = V \otimes \mathbb{C}$ with $F^0 V_{\mathbb{C}} = V_{\mathbb{C}}$ and $F^{k+1} V_{\mathbb{C}} = 0$ such that

$$V_{\mathbb{C}} = F^p V_{\mathbb{C}} \oplus \overline{F^{k-p+1} V_{\mathbb{C}}},$$

for all $p = 0, \dots, k$. It is also possible to define a sub-Hodge structure. Bear in mind that not every Hodge structure arises from a compact Kähler manifold. An integral polarised Hodge structure of weight $k$ is given by a Hodge structure $(V_{\mathbb{Z}}, F^p V_{\mathbb{C}})$ of weight $k$, along with an intersection form $Q$ on $V_{\mathbb{Z}}$ which we will define, which is symmetric if $k$ is even, alternating if $k$ is odd, and such that the Hodge decomposition is orthogonal for the induced Hermitian form and

$$(-1)^{\frac{k(k-1)}{2}} i^{p-q-k} H(\alpha) > 0,$$

for $\alpha$ non-zero of type $(p, q)$.

We have tried to give a brief idea of how Hodge theory can be applied to learn things about cohomology. It can also be used to draw conclusions about Poincaré duality and duality theorems, though we will not go into this here. One important theorem which we should mention is the Hodge index theorem. Start with $X$ a compact Kähler manifold of dimension $n$, equipped with a Kähler form $\omega$. Take a pairing $\langle , \rangle$ defined by

$$H^k(X, \mathbb{R}) \otimes H^{2n-k}(X, \mathbb{R}) \to \mathbb{R}.$$

If $L$ is the Lefschetz operator which acts on the cohomology, define the intersection form on $Q$ on $H^k(X, \mathbb{R})$:

$$Q(\alpha, \beta) = \langle L^{n-k} \alpha, \beta \rangle = \int_X \omega^{n-k} \wedge \alpha \wedge \beta.$$

One can prove that the subspaces $H^{p,q}(X) \subset H^k(X, \mathbb{C})$ form an orthogonal direct sum for $H_k$. Also, the form defined by

$$(-1)^{\frac{k(k-1)}{2}} i^{p-q-k} H_k$$

is positive definite on the complex subspace defined by the intersection of $H^k(X, \mathbb{C})_{\text{prim}}$ with $H^{p,q}(X)$. The Hodge index theorem is a consequence of this theorem: it describes the index (or rather, the signature) of the intersection form. The signature of the form is defined by

[236] Claire Voisin, 'A counterexample to the Hodge conjecture extended to Kähler varieties', *International Mathematics Research Notices*, 2002 (20), 1057 – 1075 (2002).
[237] Alexander Grothendieck, 'Hodge's general conjecture is false for trivial reasons', *Topology*, 8 (3), 299 – 303 (1969).

$$Q(\alpha, \beta) = \int_X \alpha \wedge \beta$$

on $H^n(X, \mathbb{R})$ is equal to

$$\sum_{a,b} (-1)^a \dim\Big(H^{a,b}(X)\Big).$$

Most of these isolated results (and many others) can be tied up to prove that the rational cohomology of a compact polarised complex manifold can be decomposed into a direct sum of polarised Hodge structures. We can also define integral and rational Hodge structures via the Hodge decomposition of the cohomology. Lefschetz decomposition and the Hodge index theorem allow for a direct sum expression of the cohomology of a compact Kähler manifold, where the elements of the sum are primitive components which are compatible with the Hodge decomposition. If $\alpha$ is a real closed form of type $(1,1)$ on a compact Kähler manifold with integral cohomology class, then $\alpha$ is the Chern form for a holomorphic line bundle endowed with a Hermitian metric. Also, if $L$ is a holomorphic line bundle over a compact complex manifold, and $h$ is a Hermitian bundle metric whose Chern form is Kähler, then for large enough $N$, the holomorphic sections of $L^{\otimes N}$ gives a holomorphic embedding of the manifold into a projective space[238].

There is also a notion of formal Hodge theory in the literature, where one studies formal Hodge structures and extends the Hodge realization $T_{\text{Hodge}}$ of Deligne's 1-motives to a realization from the category of Laumon's 1-motives to the category of torsion free formal Hodge structures $\text{FHS}_1^{\text{fr}}$ of level $\leq 1$. A motive is an object which was introduced to unify different cohomology theories (we have mentioned de Rham and singular cohomology, but there are many others). These formal Hodge structures form a subcategory of the abelian category $\text{FHS}_1$ of all formal mixed Hodge structures:

$$\text{FHS}_1^{\text{fr}} \subset \text{FHS}_1.$$

The category of torsion free graded polarizable mixed Hodge structures is a subcategory of $\text{FHS}_1^{\text{fr}}$ and one can show that there is an equivalence of categories

$$T_{\oint} : \mathcal{M}_1^{\text{a,fr}} \longrightarrow \text{FHS}_1^{\text{fr}},$$

where $\mathcal{M}_1^{\text{a,fr}}$ is the category of Laumon 1-motives and where one obtains a commutative diagram such that the functor from the category of Deligne 1-motives $\mathcal{M}_1^{\text{fr}}$ to the Laumon 1-motives is a canonical inclusion, the functor from $\mathcal{M}_1^{\text{a,fr}}$ to $\mathcal{M}_1^{\text{fr}}$ is a forgetful functor (the left inverse of the inclusion), and

$$T_{\oint}\ (M)_{\text{ét}} = T_{\text{Hodge}}(M_{\text{ét}}).$$

This is equivalent to saying that $T_{\oint}\ = T_{\text{Hodge}}$ if $M = M_{\text{ét}}$, where $M \coloneqq [F \to G]$ is a Laumon 1-motive (defined to be a commutative formal group with torsion free étale part, a commutative connected algebraic group, and a map of abelian sheaves $u$ which are faithfully flat with finite presentations. The 'fppf' topology is a type of topology which is relevant in algebraic geometry (flat topologies in general are examples of Grothendieck topologies)[239].

[238] Claire Voisin, *Hodge Theory and Complex Algebraic Geometry* (Cambridge: Cambridge University Press, 2002).

This talk of categories may seem strange if you are unaccustomed to category theory, but it is by this point embedded in the mathematical landscape. Homology, for example, is a functor from the category of topological spaces to the category of sequences of abelian groups, and does not necessarily make much sense if you try to view it any other way. Start for a class of objects and for every ordered pair of objects in the class, specify a set of morphisms which take object $A$ to object $B$ (denoted by $[A, B]$). For every triple of objects $A, B, C$, specify a rule which assigns to every pair of morphisms another morphism which belongs to set of morphisms from object $A$ to object $C$, called a composition of morphisms. The class of objects along with the morphisms and the composition maps forms a category if composition of morphisms is associative (if $(hg)f = h(gf)$ for any triple of morphisms, where $f \in [A, B]$, $g \in [B, C]$ and $h \in [C, D]$) and if for any object $B$ there is an identity morphism taking $B$ to $B$ such that for any two morphisms $f \in [A, B]$ and , $g \in [B, C]$, one has the following equalities:

$$\mathrm{Id}_B f = f, \qquad g\mathrm{Id}_B = g.$$

There are many examples of categories. Take the category of all sets and their maps, for example. Observe that we never talk about taking the set of all sets, or we get involved in Russell's paradox. The most common layperson's statement of this paradox asks that one imagine a barber in a village who shaves everyone that does not shave themselves and no-one else. The question is whether the barber will have to shave himself, or not: this is the type of paradox which is being described, since you are ending up with something which is and is not a member of itself. In formal terms, the paradox is

$$\text{If } R = \{x \text{ such that } x \notin x\}, \text{then } \mathrm{R} \in \mathrm{R} \iff \mathrm{R} \notin \mathrm{R}.$$

Obviously this is complete nonsense.

The objects of the category are all sets and the set of morphisms is composed of all the possible maps from one set into another. The next most obvious example is the category of groups and group homomorphisms, where the objects are the groups and the morphism are the homomorphisms. The category of topological spaces along with the set of continuous maps is another example. One can also take the category formed from all topological spaces and all classes of homotopic maps: in this case, the morphisms are classes of a particular type of map, and not maps themselves. In general, morphisms do not need to be maps, although we may denote them as if they were. These definitions also absorb the notion of an isomorphism in a more general way. Two objects in a category are isomorphic if there exist two morphisms $f \in [A, B]$ and $f \in [B, A]$ such that

$$fg = \mathrm{Id}_B, \qquad gf = \mathrm{Id}_A.$$

[239] Luca Barbieri-Viale, 'Formal Hodge theory', *Math. Res. Lett.* 14(3), 385 – 394 (2007).

An isomorphism in the category of sets is a bijection, an isomorphism in the category of topological spaces is a homeomorphism, and so on.

Now take two categories and for an object $A$ in the first category, assign an object in the second category $F(A)$. Also, for a morphism $f \in [A, B]$ in the set of morphisms of the first category, assign a morphism $f_* \in [F(A), F(B)]$ in the set of morphisms of the second category. This defines a covariant functor from the first category to the second category, as long as $f$ being an identity morphism implies that $f_*$ is an identity morphism and as long as

$$(fg)_* = f_* g_*,$$

for some well-defined composition. If you define everything to act in the opposite direction and modify the axioms trivially, you have instead a contravariant functor[240]. For an example of a functor, take the category of finite-dimensional vector spaces $\mathrm{V}$ with linear isomorphisms as morphisms. One can define a functor $\mathcal{F}\colon \mathrm{V} \to \mathrm{V}$ by requiring that

$\mathcal{F}(V) = V^*$ for some vector space $V$ living in the category, and that $\mathcal{F}(A) = (A^{-1})^*$ for some isomorphism $A$ living in the set of morphisms.

It can be shown that $\mathcal{F}$ is a smooth covariant functor (ie. it is smooth for every $V$) and that $\mathcal{F}(TM)$ and $T^*M$ are smoothly isomorphic vector bundles via a canonical choice of bundle isomorphism[241]. Also, the image of the Hodge structure on the cohomology of Kähler manifolds relative to the holomorphic maps between any two Kähler manifolds is functorial. It should be clear to you that this kind of theory is very useful for proving things, and that it can be extended to multiple settings in a flexible way. However, as always you have to think about what you are working with can be computed, and then make restrictions if it cannot. For example, a functor for a homology theory will take the category of all topological spaces to the category of all *sequences* of abelian groups, which is going to make practical computation of homology groups extremely difficult unless restrictions are made to the type of space we are studying. In fact, it is rare that you can consider what you might call the completely 'arbitrary' case in pure mathematics, as there is too much going on: one usually has to make a restriction or assumption of some kind before any progress can be made.

One of Voisin's most significant results was the proof of Green's conjecture on the syzygies of the canonical embedding of an algebraic curve in the case of generic curves $C$ of genus $g(C)$ and gonality $G$ in the range

$$\frac{g(C)}{3} + 1 \leq \mathrm{gon}(C) \leq \frac{g(C)}{2} + 1.$$

The gonality is defined to be

$$\mathrm{gon}(C) \coloneqq \mathrm{Min}\{d, \exists L \in \mathrm{Pic}\, C, d^0 L = d, h^0(L) \geq 2\}$$

such that

---

240 Sergey Matveev, *Lectures on Algebraic Topology* (Switzerland: European Mathematical Society, 2006).
241 John M. Lee, *Introduction to Smooth Manifolds* (New York: Springer, 2003).

$$\text{gon}(C) = \begin{cases} \frac{g+3}{2}, & g \text{ odd} \\ \frac{g+2}{2}, & g \text{ even} \end{cases}$$

Specifically, Voisin proved a reformulation of the conjecture: that for a smooth projective curve of genus $g$ in characteristic $0$, the condition that the Clifford index be larger than $l$ is equivalent to the fact that $K_{g-l'-2,1}(C, K_C) = 0$. The conjecture for arbitrary curves is still an open problem[242].

Another important result was the proof that in dimension $\geq 4$ there exist compact Kähler manifolds which are not homotopic to projective complex manifolds. This is interesting, because it was already known that a compact Kähler surface is deformation equivalent to a complex projective surface, but this result says that the same is not true in higher dimensions. Two deformation equivalent complex manifolds are homeomorphic and homotopic. The problem of whether the result for a surface can be extended to higher dimension is known as the Kodaira problem. There is also a more general Kodaira theorem which states that a compact complex manifold is projective if and only if it admits a Kähler form whose cohomology class is integral[243].

[242] Claire Voisin, 'Green's generic syzygy conjecture for curves of even genus lying on a K3 surface', *J. Eur. Math. Soc.* 4, 363 – 404 (2002).
[243] Claire Voisin, 'On the homotopy types of compact Kähler and complex projective manifolds', *Inventiones Math.* 157 (2), 329 – 343 (2003).

# *Olga Holtz*

Holtz is an applied mathematician working in numerical analysis, numerical linear algebra and approximation theory. As the name suggests, numerical analysis is the subject which studies numerical solutions to problems of analysis. The aim is to find accurate approximate solutions to problems which cannot be solved otherwise. We have already mentioned the Newton-Raphson method as an example. Another important branch of numerical analysis is numerical integration, where one seeks to use a numerical method to evaluate a definite integral which cannot be computed in terms of elementary functions (most integrals fall into this category). Start with an integral

$$\int_a^b f(x)\,dx.$$

Assume that the integrand can be evaluated at equally spaced points and that the integral can be approximated with a formula of the form

$$\int_a^b f(x)\,dx \approx \frac{b-a}{N}\sum_{k=1}^{N} c_k f(a+kh).$$

A formula of this type is known as a Newton-Cotes formula. For example, set $N = 2$.

$$\int_a^b f(a+y)\,dx \approx \frac{b-a}{N}\sum_{k=1}^{N} c_k\big(f(a+kh)\big).$$

Use Taylor expansions.

$$f(a+kh) = f(a) + khf'(a) + \frac{1}{2}k^2h^2f''(a) + O(h^3),$$

$$f(a+y) = f(a) + \frac{1}{2}hf'(a) - \frac{1}{2}h^2f''(a) + O(h^3),$$

where $h$ is the step size $(b-a)/N$. This gives us a Taylor representation which must be equivalent to the original function being summed:

$$f(a) + h\left(\frac{3}{2}\big(f'(a) + hf''(a)\big)\right) + O(h^3).$$

By comparison of coefficients, we have $c_1 = c_2 = 3/2$. This gives us a Newton-Cotes formula known as the trapezoid rule:

$$\int_a^b f(x)\,dx \approx \frac{b-a}{3}\left(\frac{3}{2}(f_1 + f_2)\right).$$

One can use a computer to find the third coefficient. This gives us a Newton-Cotes formula known as the Milne rule:

$$\int_a^b f(x)\,dx \approx \frac{b-a}{4}\left(\frac{4}{3}(2f_1 - f_2 + f_3)\right).$$

This can be continued to get further Newton-Cotes formulae.

The idea of approximating an integral with a sum which is closely related is a common one in analysis. The Euler-Maclaurin formula, for example, gives the difference between the integral and the sum. In some cases, one can take a function defined by an integral and write it as a series expansion. As an example, start with the following function:

$$f(x) = \int_0^1 \frac{\sqrt{t}}{(1-xt^2)^\beta}\,dt,$$

where $\beta > 1$. It is known that

$$(1-z)^{-\beta} = \frac{1}{\Gamma(\beta)}\sum_{k=0}^{\infty}\frac{\Gamma(k+\beta)z^k}{k!},$$

where $\Gamma$ is the gamma function. This function is essentially a generalisation of the factorial function to complex analysis. For some time, it was suggested that the Euler characteristic $\chi\big(\mathrm{Out}(F_n)\big)$ might be related to number theory. Borinsky proved that this is the case using the gamma function. In particular:

$$\bar{\chi}\big(\mathrm{Out}(F_{n+1})\big) = -\frac{\Gamma\left(n+\frac{1}{2}\right)}{n\log^2 n} + O\left(\frac{\Gamma\left(n+\frac{1}{2}\right)}{n\log^4 n}\right).$$

$\bar{\chi}\big(\mathrm{Out}(F_{n+1})\big)$ grows at a faster than exponential rate, since the error term is not sufficient to cancel the growth. This solves a conjecture due to Vogtmann that one can take the Euler characteristic of the mapping class groups to get a relation with the Bernoulli numbers. Recall from the definition of the gamma function that

$$\begin{aligned}\Gamma\left(n+\frac{1}{2}\right) &= \frac{2n!}{4n!}\sqrt{\pi},\\ &= \frac{(2n-1)!!}{2^n}\sqrt{\pi}.\end{aligned}$$

This implies the nice identity:

$$\Gamma\left(\frac{1}{2}\right) = \sqrt{\pi}.$$

By comparison with the above series we have

$$\begin{aligned}f(x) &= \int_0^1 dt\,\frac{\sqrt{t}}{\Gamma(\beta)}\sum_{k=0}^{\infty}\frac{\Gamma(k+\beta)}{k!}x^k t^{2k},\\ &= \frac{1}{\Gamma(\beta)}\sum_{k=0}^{\infty}\frac{\Gamma(k+\beta)}{k!}x^k\int_0^1 dt\,t^{2k+\frac{1}{2}},\end{aligned}$$

$$= \frac{1}{\Gamma(\beta)} \sum_{k=0}^{\infty} \frac{\Gamma(k+\beta) x^k}{k! \left(2k + \frac{3}{2}\right)},$$
$$= \frac{2}{\Gamma(\beta)} \sum_{k=0}^{\infty} \frac{\Gamma(k+\beta)}{k! \, (4k+3)} x^k.$$

We would now like to find the radius of convergence for the series. This can be done very easily by finding the limit of the ratio of coefficients:

$$u_k = \frac{\Gamma(k+\beta)}{k! \, (4k+3)}.$$

$$\frac{u_k}{u_{k+1}} = \frac{\Gamma(k+\beta)}{k! \, (4k+3)} \frac{(k+1)! \, (4k+4)}{\Gamma(k+\beta+1)},$$
$$= \frac{(k+1)(4k+4)}{4k+3} \frac{\Gamma(k+\beta)}{(k+\beta)\Gamma(k+\beta)},$$
$$= \frac{(k+1)(4k+4)}{(k+\beta)(4k+3)}.$$

This ratio tends to 1 as $k$ tends to infinity, so the radius of convergence is $1$. We now take the ratio of the coefficient with the previous coefficient and do some algebra which we will show.

$$\frac{u_k}{u_{k-1}} = \frac{\Gamma(k-\beta)}{k! \left(2k + \frac{1}{2} + 1\right)} \frac{(k-1)! \left(2k + \frac{1}{2} - 2\right)}{\Gamma(k+\beta-1)},$$
$$= \frac{(k+\beta-1)\left(2k + \frac{1}{2} - 2\right)}{k\left(2k + \frac{1}{2} + 1\right)},$$
$$= \left(1 + \frac{\beta^{-1}}{k}\right)\left(1 - \frac{3}{4k}\right)\left(1 + \frac{3}{4k}\right)^{-1},$$
$$= \left(1 - \frac{3}{4k} + \frac{\beta^{-1}}{k} + O(k^{-2})\right)\left(1 - \frac{3}{4k} + O(k^{-2})\right),$$
$$= 1 - \frac{3}{4k} + \frac{\beta^{-1}}{k} - \frac{3}{4k} + O(k^{-2}),$$
$$= 1 + \frac{4(\beta-1) - 6}{4k} + O(k^{-2}),$$
$$= 1 + \frac{4\beta - 10}{4k} + O(k^{-2}),$$
$$= 1 - \frac{\frac{10}{4} - \beta}{k} + O(k^{-2}).$$

Comparing with a result of Richards, we see that there is a singularity of order $\alpha = 1 - \beta$, since[244]

$$f(x) \sim A(1-x)^{1-\beta},$$

for some constant $A$ and for $x$ close to the radius of convergence.

Another common method of numerical integration is Gaussian quadrature. One approximates a definite integral representing a function as a weighted sum over $N$ points:

$$\int_{-1}^{1} w(x) f(x)\, dx \approx \sum_{k=0}^{N} c_k f(x_k),$$

where $w(x)$ is a weight function. The weight function is often the unit function, but not always. Be warned that basic quadrature commands in computer algebra programs will likely not look for singularities, as they will probably assume a unit weight function. As we know, the inner product notation denotes

$$(f, g) = \int_{a}^{b} w(x) f(x) g(x)\, dx.$$

If we take the functions in this integral to be the set of orthogonal polynomials $\phi_k(x)$ such that the leading term of $\phi_n$ is $x^n$, we can assume that they satisfy an orthogonality condition

$$(\phi_n, \phi_m) = h_n \delta_{nm}.$$

Written out explicitly, we might choose the weight function so that we have

$$\int_{0}^{1} (1-x)^{3/2}\, \phi_n(x) \phi_m(x)\, dx = h_n \delta_{nm}.$$

These polynomials can be found using Gram-Schmidt orthogonalisation. For example, the first few polynomials would be

$$\phi_0 = 1,$$

$$\phi_1 = x - \frac{2}{7},$$

$$\phi_2 = x^2 - \frac{8}{11}x + \frac{8}{99},$$

and so on. The corresponding coefficients of the Kronecker delta are seen to be:

$$h_0 = \frac{2}{5}, h_1 = \frac{8}{441}, h_2 = \frac{128}{127{,}413}.$$

Now let us say that we wish to approximate the following definite integral:

---

[244] Derek Richards, *Advanced Mathematical Methods with Maple®* (Cambridge: Cambridge University Press, 2001).

$$\int_0^1 (1-x)^{3/2} f(x)\, dx.$$

Write

$$f(x) = \sum_{j=0}^{\infty} a_j \phi_j(x),$$

where

$$a_j = \frac{1}{h_j} \int_0^1 w(x)\, \phi_j(x) f(x)\, dx.$$

You may recognise the similarity to Fourier coefficients. Choose $x_k$ to be the $N$ zeroes of the polynomial $\phi_N(x)$ and the $c_k$ to solve the linear equations

$$\sum_{j=1}^{N} c_j = h_0\, , \sum_{j=1}^{N} c_j\, p_k(x_j),$$

for $k$ up to $N-1$.

Let us choose $N = 4$, then we have an approximation

$$\int_0^1 (1-x)^{3/2} f(x)\, dx \approx \sum_{k=1}^{4} c_k f(x_k),$$

where the numbers $x_k$ are the zeroes of the fourth-order orthogonal polynomial $\phi_4$ and the coefficients are given by a tuple as follows

$$c = [0.122, 0.169, 0.092, 0.017].$$

This quadrature can be implemented explicitly with a computer to get numerical values for certain definite integrals. For example,

$$\int_0^1 (1-x)^{3/2} \cos(\pi \sin x)\;\, dx \approx 0.221.$$

The Gauss quadrature method can be extended to allow for approximations of integrals of the form

$$\int_{-1}^{1} \frac{f(x)}{\sqrt{1-x^2}}\, dx.$$

In this case, the quadrature is Chebyshev-Gauss quadrature. If we take the weight function to be $e^{-x^2}$, we have Gauss-Hermite quadrature for evaluating integrals of the form

$$\int_{-\infty}^{\infty} e^{-x^2} f(x)\;\, dx.$$

Note that this allows for integration over an infinite interval.

Holtz has worked in particular on matrix orthogonal polynomials, extending some classical results from the theory of orthogonal polynomials on the unit circle to the matrix case. Holtz proved what can be called a matrix Szegő theorem: for any matrix probability measure $\sigma \in P_\ell(\mathbb{T})$ and any natural number $n$, one has

$$\int_{\mathbb{T}} \operatorname{Tr}\log\sigma' \frac{d\theta}{2\pi} \le \operatorname{Tr}\log\beta_n = \log\prod_{k=0}^{n-1}\det(1-\alpha_k^*\,\alpha_k).$$

As a corollary, one has that if σ is a Szegő measure, then

$$\int_{\mathbb{T}} \operatorname{Tr}\log\sigma' \frac{d\theta}{2\pi} \le \inf_n \operatorname{Tr}\log\beta_n \le -\sup_n \log\|\beta_n^{-1}\| \le 0.$$

Using these techniques, she obtained an elementary proof of the Helson-Lowdenslager distance formula, which states that for every $\sigma \in P_\ell(\mathbb{T})$,

$$\exp\int_{\mathbb{T}} \frac{1}{\ell}\operatorname{Tr}\log\sigma' \frac{d\theta}{2\pi} = \inf_{A,P}\int_{\mathbb{T}} \frac{1}{\ell}\operatorname{Tr}\big((A+P)^* d\sigma (A+P)\big),$$

where $A$ cycles through all matrices of unit determinant, and $P$ runs over all trigonometric polynomials of the form[245]

$$P\big(e^{i\theta}\big) = \sum_{k>0} A_k e^{ik\theta}.$$

Most of us are familiar with Gram-Schmidt orthogonalisation. More advanced numerical linear algebra considers iterative methods which are suitable for large sparse systems. We mentioned earlier the Gauss-Seidel and conjugate gradient methods, and stated that the Gauss-Seidel method is not as efficient as the conjugate gradient method. In fact, this can be quantified by including a comparison between both methods in terms of the computational time needed to solve a randomly generated positive definite $n \times n$ matrix with density $0.1$. The comparison is conjugate gradient against Gauss-Seidel with two different preconditioners[246].

| $n$ | Gauss-Seidel | Gauss-Seidel (RCM) | Gauss-Seidel (MD) | CG |
|---|---|---|---|---|
| 100 | 0.007928804 | 0.007056294 | 0.006839876 | 0.002512009 |
| 500 | 0.502849338 | 0.405763259 | 0.303914523 | 0.002419189 |
| 1000 | 3.219242909 | 3.281714839 | 2.618037809 | 0.004361086 |
| 3000 | 116.8587276 | 105.7410554 | 99.96394861 | 0.044996817 |

Table 2: Efficiency of the conjugate gradient method in comparison with Gauss-Seidel.

[245] Maxim Derevyagin, Olga Holtz, Segey Krushchev and Mikhail Tyaglov, 'Szegő's theorem for matrix orthogonal polynomials', *Journal of Approximation Theory*, 164(9), 1238 – 1261 (2011).
[246] Jamie Glaves, 'Computational techniques for sparse matrices', Conference Proceedings of Tomorrow's Mathematicians Today, Manchester Metropolitan University (2017).

One can see that the conjugate gradient method remains much more efficient even for very large matrices. Technically, conjugate gradient is an example of a Krylov space method. A crucial part of both methods is the minimization of some measure of error at the $k$th iteration. Another such method is GMRES (generalized minimal residual method). GMRES can be used for non-symmetric systems but one must store a basis for the $k$th Krylov space $K_k$ and the algorithm must be re-started at some point, which increases the convergence time. The $k$th iteration of GMRES is the solution to a least squares problem:

$$\text{minimize}_{x \in x_0 + K_k} \|b - Ax\|_2.$$

Iterative methods are generally effective provided that the problem is preconditioned (where the preconditioner transforms the problem to one more amenable to iterative methods). In practice, preconditioning often involves an attempt to reduce the condition number, where the condition number measures the sensitivity of a function to changes in the input. Reducing the condition number improves the performance of the iteration, as does transforming the problem so that the eigenvalues are clustered near 1.

As an example, when performing an iteration we may wish to replace the linear system

$$Ax = b$$

in the conjugate gradient method with another system which is symmetric positive definite and which has the same solution. We could do this by expressing the preconditioned problem in terms of a new symmetric positive definite matrix $B$ such that

$$A = B^2$$

and then use a two-sided preconditioner which approximates $B^{-1}$. The square of this preconditioner is equal to $M$, where $M$ is a symmetric positive definite matrix close to $A^{-1}$ and so we can now express the preconditioned system in terms of the matrix $SAS$ which has its eigenvalues clustered near 1[247].

It is also possible to study large dense, as opposed to large sparse, systems. The methods for solving these systems are obviously different, since the main operation in the sparse case is that of matrix-vector multiplication, which is only going to have linear computational cost for a sparse system. One way of dealing with a dense system is to replace the original matrix $A$ with an approximate matrix $\tilde{A}$ such that the computational cost of a matrix-vector product is much less than $n^2$, ie. the computational cost is reasonable. To be more specific, the computational cost of a matrix-vector multiplication with the new approximate matrix should be should be bounded by $O(n \log^{\alpha} n)$, the memory needed to store the matrix and the time required to generate the matrix should both be bounded by $O\left(n \log^{\beta} n\right)$, and for some positive $\varepsilon$ and a matrix norm $\|\cdot\|$, the approximate matrix should satisfy

$$\left\|A - \tilde{A}\right\| \leq \varepsilon.$$

[247] C. Kelly, *Iterative Methods for Linear and Nonlinear Equations* (Philadelphia: Society for Industrial and Applied Mathematics, 1995).

In both cases, $\alpha$ and $\beta$ should not depend on $n$.

You might wonder what kind of problem would produce dense matrices. When we are interpolating with radial basis functions, we form an interpolant of the form

$$s(\mathbf{x}) = \sum_{j=1}^{n} \alpha_j K(\mathbf{x}, \mathbf{x}_j),$$

where the kernel is positive definite

$$K(\mathbf{x}, \mathbf{y}) = \Phi(\mathbf{x} - \mathbf{y}) = \phi(\|\mathbf{x} - \mathbf{y}\|_2).$$

The interpolation conditions create a system

$$A\alpha = f$$

with an interpolation matrix $A$ which is symmetric positive definite. This matrix will be positive when the kernel does not have local support. In this case, it so happens that the computational cost of computing the matrix entries is cheap despite the matrix being dense, but this is not a typical example[248].

Another classic of numerical linear algebra which we could mention is Gaussian elimination. This method was known to the Chinese linear algebraists, and rediscovered by Gauss. The result was also known to Newton (as was any method of elementary algebra or linear algebra). Gaussian elimination provides an algorithm for solving systems of simultaneous linear equations. The general procedure is as follows: take a system of $n$ linear equations with a coefficient matrix $A$. Write down the augmented matrix $A|b$ with $n$ rows $R_n$. Subtract multiples of the first row from the other rows until the elements below the leading diagonal are zero in the first column. Repeat the process with multiples of the second row until the elements below the leading diagonal are zero in the second column. Continue doing this until we have an augmented matrix such that the coefficient matrix is upper triangular. Solve the corresponding system of equations by plugging everything in backwards. We will show how this might look in an explicit example. Start with a system of linear equations:

$$x - 2y + 5z = 6,$$

$$x + 3y - 4z = 7,$$

$$2 + 6y - 12z = 12.$$

We need the elements below the leading diagonal in the first column to vanish. This can be done by performing elementary operations with the first row: for example, subtracting the first row from the second row would be the obvious thing to do. For the third row, you have to multiply the first row by 2 before you subtract it. We then have

$$x - 2y + 5z = 6,$$

[248] Holger Wendland, *Numerical Linear Algebra: An Introduction* (Cambridge: Cambridge University Press, 2017).

$$5y - 9z = 1,$$

$$10y - 22z = 0.$$

I will not bother actually writing everything out so that the coefficients are entries in an augmented matrix, but you can do this if you find it useful. Repeat the process by performing elementary operations with the second row until the element below the leading diagonal in the second column vanishes.

$$x - 2y + 5z = 6,$$

$$5y - 9z = 1,$$

$$-4z = -2.$$

The rows of $A$ are linearly independent so there is a unique solution to this system, which can be found trivially. This method can break down in various ways and will not always work but I will not go into this, as it can be found in any book on linear algebra. An application of the method is in polynomial interpolation.

In the field of linear algebra, Holtz has studied the conditions under which a matrix of differentiable functions has to commute with its element-wise derivative. She posed the problem of whether a matrix over a differential field $\mathbb{F}$ with an algebraically closed field of constants such that

$$MM' = M'M$$

has to have eigenvalues which are elements of the field $\mathbb{F}$. A differential field is an algebraic field together with an additional operation called the derivative which satisfies

$$(a + b)' = a' + b',$$

$$(ab)' = ab' + a'b.$$

She also proves various other results relating to triangularizability and diagonalizability of matrices in this context: for example, if $\mathbb{F}$ is a differential field with an algebraically closed field of constants $\mathbb{K}$ and if $M \in \mathbb{F}^{n,n}$ is type 1, then $M$ is $\mathbb{K}$-triangularizable.

This theorem implies that type 1 matrices have $n$ eigenvalues in $\mathbb{F}$ if $\mathbb{K}$ is algebraically closed. A type 1 matrix is one which has the form

$$M_1 = \sum_{\lambda} f_\lambda C_\lambda,$$

where the $C_\lambda$ are pairwise commuting constant matrices. A type 2 matrix is one which has the form

$$M_2 = (f_\alpha g_\beta),$$

where $f_1, \dots, f_n, g_1, \dots, g_n$ are arbitrary functions satisfying the conditions

$$\sum_{\alpha} f_\alpha g_\alpha = \sum_{\alpha} f'_\alpha g_\alpha = 0,$$

which implies that

$$\sum_{\alpha} f_\alpha g'_\alpha = 0.$$

This distinction can be found in a letter of Schur[249].

Closely allied to numerical analysis is a huge body of theory called approximation theory. Consider a curve in the plane which is given by a function and then think about how one might draw a straight line which best fits the curve. This is the prototypical problem in approximation theory. One has a function or a member of a set which needs to be approximated, a set of possible approximations $\mathcal{A}$ (the set of all straight lines in the example we have given), and some method for selecting the best approximation from the set $\mathcal{A}$. How can one measure how good an approximation is? Recall that metric spaces are equipped with a distance function which could be used to do this, and it is often the case when working with an approximation problem that there is a nice metric space which contains both the function to be approximated and the set of approximations $\mathcal{A}$. We could say that $a_1$ is a better approximation than $a_2$ if the inequality

$$d(a_1, f) < d(a_2, f)$$

is satisfied. We can define $a_B$ to be the best approximation if the inequality

$$d(a_B, f) < d(a, f)$$

holds for all $a \in \mathcal{A}$. Bear in mind that the metric space needs to be able to measure the error when we attempt a trial approximation so that we can improve our attempts to get a best fit. Some famous results which you might be familiar with (the Stone-Weierstrass theorem and the Lagrange interpolation formula, for example) are results of approximation theory[250].

[249] Olga Holtz, Volker Mehrmann and Hans Schneider, 'Matrices that commute with their derivative. On a letter from Schur to Wielandt', *Linear Algebra Appl.* 438(5), 2574 – 2590 (2012).
[250] Michael Powell, *Approximation Theory and Methods* (Cambridge: Cambridge University Press, 1981).

# *Maryam Mirzakhani*

Mirzakhani was known for her profound work on symplectic geometry, hyperbolic geometry and moduli spaces of Riemann surfaces. She began publishing mathematical articles on a variety of topics from an early age and was the first woman and the first Iranian to win the prestigious Fields medal for achievements in mathematics. We mentioned Teichmüller spaces earlier in connection with the mapping class group, but did not really define them. A compact Riemann surface with a finite number of punctures is said to be of finite type. If you fix a smooth compact reference surface $F$ with genus $g$, every diffeomorphism $F$ into a Riemann surface $R$ can be pulled back to a conformal structure on the reference surface (denote the set of all the conformal structures on $F$ by $S(F)$). If one denotes the group of orientation-preserving diffeomorphisms as usual by $\mathrm{Diff}^+(F)$, then the quotient $\mathcal{M}(F) = S(F)/\mathrm{Diff}^+(F)$ is the space of conformal structures on $F$ with $n$ labelled, distinguished points (this is simply the moduli space of Riemann surfaces of type $F$. If one quotients through instead by the normal subgroup of diffeomorphisms which are homotopic to the identity map, one obtains the space of marked conformal structures on $F$:

$$\mathcal{T}(F) = {}^{S(F)}\!\big/_{\mathrm{Diff}_0(F)}.$$

This is the Teichmüller space of surfaces of type $F$. If one quotients thorough the first diffeomorphism group by the second group, one obtains the pure mapping class group which we mentioned before.

We will state another definition for something you may have heard mentioned before: a train track for a hyperbolic surface $R$ with cusps is an embedded family of smooth curves with a finite set of vertices called switches and edges called branches. There is a unique tangent at each switch such that curves meet at the switch with the same tangent line and the neighbourhood of a switch is a union of disjoint smoothly embedded arcs. There are some further restrictions on the form that a component of the complement of the train track can take. As with many concepts in modern geometry, the concept of the train track was introduced by Thurston, who needed something which would give a piecewise linear symplectic structure for the space of measured geodesic laminations $\mathcal{MGL}$. A lamination is a foliation of a closed subset[251].

Recall that a foliation of dimension $k$ on an $n$-manifold is a collection of disjoint, connected, immersed $k$-dimensional submanifolds of $M$ such that the union is $M$ and in a neighbourhood of a point in $M$ there is a smooth chart $(U, \varphi)$, where $\varphi(U)$ is a product of connected open subsets $U' \times U'' \subset \mathbb{R}^k \times \mathbb{R}^{n-k}$. Also, each leaf of the foliation intersects $U$ on a countably large union of $k$-dimensional slices. A chart of this kind is called a foliated chart. The immersed submanifolds are actually equivalence classes, known as the leaves of the foliation, where foliation is an equivalence relation. These criteria are met trivially in the case where we take the set of all $k$-dimensional affine subspaces of $\mathbb{R}^n$ parallel to $\mathbb{R}^k$. As another example, if one takes the set of connected components of the curves in the $(y, z)$-plane given by the following equations

$$z = \sec y + c, \qquad c \in \mathbb{R},$$

$$y = \left(k + \frac{1}{2}\right)\pi, \qquad k \in \mathbb{Z},$$

and rotates the curves about the $z$-axis, one obtains a foliation of $\mathbb{R}^3$ with a nice geometric visualisation[252]. A lamination $\mathcal{G}$ is geodesic if its leaves are complete geodesics. A transverse measure $\mu$ for a geodesic lamination is an assignment for each smooth transverse arc $\tau$ to $\mathcal{G}$ a non-negative measure whose support is $\tau \cap \mathcal{G}$. A measured geodesic lamination is the union of a finite number of simple closed geodesics and a perhaps uncountable number of open complete geodesics and the complement of a measured geodesic lamination is a finite union of convex subsurfaces with geodesic boundaries.

[251] Scott Wolpert, *Families of Riemann Surfaces and Weil-Petersson Geometry* (USA: American Mathematical Society, 2009).

A geodesic lamination is carried by a train track $\tau$ provided that there is a differentiable map $f$ of the hyperbolic surface which is homotopic to the identity map and such that $f(\mu)$ is contained in $\tau$ and the differential does not vanish on a leaf. Given a marked hyperbolic surface, the mapping class group acts on the associated space of measured geodesic laminations. A fundamental result is that the action of the mapping class group on $\mathcal{MGL}$ is ergodic[253]. This terminology may strange if you have studied some ergodic theory, but remember that we have defined a measure on the space. In that case, we actually have the group MCG acting on a measure space $(\mathcal{MGL}, \mathcal{A}, \mu)$, where $\mathcal{A}$ is a $\sigma$-algebra. A $\sigma$-algebra on a set $X$ is a family of subsets which contains $X$ and is closed under taking of complements, under formation of countable unions, and under formation of countable intersections.

As a trivial example, take the collection of all the subsets of a set $X$. This is a $\sigma$-algebra, but if we take the collection of all the finite subsets of an infinite set $X$, then this collection cannot contain $X$. It is therefore not closed under taking of complements, and it cannot be a $\sigma$-algebra (or even an algebra) on $X$. As always, we are also interested in constructing new $\sigma$-algebras from old ones, and there are various ways of doing this. The intersection of a collection of $\sigma$-algebras on $X$ is another $\sigma$-algebra, for example. If $X$ is a set and $\mathcal{A}$ is a $\sigma$-algebra on the set, a function $\mu$ whose domain is the $\sigma$-algebra and whose values belong to the extended half-line $[0, +\infty]$ with $\infty$ is included as a point is countably additive if

$$\mu\left(\bigcup_{i=1}^{\infty} A_i\right) = \sum_{i=1}^{\infty} \mu(A_i),$$

[252] John M. Lee, *Introduction to Smooth Manifolds* (New York: Springer, 2003).

[253] Scott Wolpert, *Families of Riemann Surfaces and Weil-Petersson Geometry* (USA: American Mathematical Society, 2009).

where $\{A_i\}$ is an infinite sequence of disjoint sets in $\mathcal{A}$. A measure is a function of this kind such that $\mu(\emptyset) = 0$ and a triple $(X, \mathcal{A}, \mu)$ is called a measure space[254]. A discrete group $G$ acting on a measure space via measure-preserving transformations is said to have an ergodic action if

$$g \cdot a = a, \forall g \in G$$

implies that $\mu(a) = 1$ or $\mu(a) = 0$. In words, the only measurable sets which are invariant under the group action are $X$ and $\emptyset$. As usual, this statement is only valid up to sets of measure zero. In this case, the Thurston volume element is the unique MCG-invariant measure in the measure class.

Mirzakhani generalized an identity due to McShane and showed that certain MCG-sums of geodesic lengths lead to functions which are constant over the Teichmüller space. The identity can be used to obtain partitions of unity on the Teichmüller space which are invariant under the action of the mapping class group. The summand for the length identity is a rational function defined on the plane as follows:

$$H(x,y) = \frac{1}{1+e^{\frac{x+y}{2}}} + \frac{1}{1+e^{\frac{x-y}{2}}},$$

and functions on $\mathbb{R}^3$ defined via the relations

$$\frac{\partial}{\partial x}\mathcal{D}(x,y,z) = H(y+z,x), \qquad \mathcal{D}(0,0,0) = 0,$$

$$2\frac{\partial}{\partial x}\mathcal{R}(x,y,z) = H(z,x+y) + H(z,x-y), \qquad \mathcal{R}(0,0,0) = 0.$$

We take the Teichmüller space $\mathcal{T}(L_1, \dots, L_n)$ of $g$ marked hyperbolic structures, where $L_1, \dots, L_n$ are the lengths of the geodesic boundaries $\beta_1, \dots, \beta_n$of the structures[255]. The McShane-Mirzakhani length identity states that for a hyperbolic surface $R$, we have

$$\sum_{\alpha_1,\alpha_2} \mathcal{D}\left(L_1, l_{\alpha_1}(R), l_{\alpha_2}(R)\right) + \sum_{j=2}^{n}\sum_{\alpha} \mathcal{R}\left(L_1, L_j, l_\alpha(R)\right) = L_1,$$

where the first summation is over pairs of simple closed geodesics such that $\beta_1, \alpha_1, \alpha_2$ bound an embedded pair of pants, and the double summation is over simple closed geodesics such that $\beta_1, \beta_j, \alpha$ bound an embedded pair of pants. This is the version of the length identity for boundaries of positive length, but there are versions for hyperbolic surfaces with cusps, and length identities are a large subject in themselves.

This identity also provides a way of evaluating the WP volume integral $\int_{\mathcal{T}/\mathrm{MCG}} dV$. Mirzakhani realised that the length identity and the Fenchel-Nielsen construction of hyperbolic surfaces lead to a recursion for WP volume integrals known as Mirzakhani volume recursion. Fenchel-Nielsen coordinates are coordinates on Teichmüller space. A collection of $3g-3+n$ disjoint simple closed geodesics can be used to cut a surface up into $2g-2+n$ pairs of pants, such that $n$ of the pants boundaries are cusps with zero length and not geodesics of positive length. The Fenchel-Nielsen coordinates for a point of the Teichmüller space of a surface are given by two sets of parameters which can be used to sew the surface together from the pieces of the decomposition. The first parameter is the boundary length $l_j$ and the second is the Fenchel-Nielsen twist $\tau_j$, which measures the hyperbolic displacement between two boundaries. As you might intuitively guess, a Fenchel-Nielsen twist around a simple geodesic is the action of the cutting the surface open along the geodesic, rotating the boundaries of the two resulting circles relative to each other, and then gluing back together. The lengths of the coordinates are simply the lengths of geodesics which are homotopic to the $3g-3+n$ disjoint simple closed geodesics. It can be shown that the mapping of $\mathcal{T}$ to parameters in this way

$$\mathcal{T} \to \left(l_j, \tau_j\right) \in \prod_{j=1}^{3g-3+n} (\mathbb{R}_+ \times \mathbb{R})$$

is an equivalence. The Fenchel-Nielsen coordinates can be used to show that the variation of a length of a simple geodesic to zero has finite WP length. One can show that these coordinates are symplectic coordinates such that one can define a WP symplectic form as

$$\omega = \frac{1}{2}\sum_j dl_i \wedge d\tau_j \, .$$

[254] David Cohn, *Measure Theory* (New York: Springer, 2013).
[255] Scott Wolpert, *Families of Riemann Surfaces and Weil-Petersson Geometry* (USA: American Mathematical Society, 2009).

This can be proved directly by searching for the coefficients of the 2-form[256].

The Weil-Petersson volume is a volume of a moduli space and volumes of this kind can be computed by using covers for intermediate moduli spaces and considering a decomposition of a surface into different configurations of subsurfaces: this is what eventually leads to the volume recursion theorem. There are finite symmetry considerations to take into account when cutting open a surface in this way, but we will skip ahead slightly to Mirzakhani's volume recursion. Start with the Teichmüller space of genus $g$ marked hyperbolic structures as before $\mathcal{T}(L_1, \dots, L_n)$ of genus $g$ marked hyperbolic structures, where $L_1, \dots, L_n$ are the lengths of the geodesic boundaries $\beta_1, \dots, \beta_n$. As before, one obtains a moduli space by quotienting the Teichmüller space through by the pure mapping class group. The WP volume $V_g(L_1, \dots, L_n)$ of this moduli space is a symmetric function of the boundary lengths. Assuming positive lengths, we set

$$V_0(L_1, L_2, L_3) = 1,$$

$$V_1(L_1) = \frac{\pi^2}{12} + \frac{L_1^2}{48}.$$

The volume also satisfies the following:

[256] Scott Wolpert, 'On the Weil-Petersson geometry of the moduli space of curves', *Amer. J. Math.,* 107(4) 969 – 997 (1985).

$$\frac{\partial}{\partial L_1} V_g(L) = \mathcal{A}_g^{con}(L) + \mathcal{A}_g^{dcon}(L) + \mathcal{B}_g(L),$$

where

$$\mathcal{A}_g^{*}(L) = \frac{1}{2}\int_0^\infty \int_0^\infty \hat{\mathcal{A}}_g^{*}(x, y, L)\, xy\, dx\, dy,$$

$$\mathcal{B}_g(L) = \int_0^\infty \widehat{\mathcal{B}}_g(x, L)\, x dx,$$

$$\widehat{\mathcal{B}}_g(x, L) = \frac{1}{2}\sum_{j=1}^{n} \Big(H\big(x, L_1 + L_j\big) + H\big(x, L_1 - L_j\big)\Big) V_g\big(x, L_2, \dots, \widehat{L_j}, \dots, L_n\big),$$

$$\hat{\mathcal{A}}_g^{con}(x, y, L) = H(x + y, L_1) V_{g-1}\big(x, y, \hat{L}\big),$$

$$\hat{\mathcal{A}}_g^{dcon}(x, y, L) = \sum_{g_1 + g_2 = g,\, I_1 \coprod I_2 = \{2, \dots, n\}} H(x + y, L_1) V_{g_1}\big(x, L_{I_1}\big) V_{g_2}\big(x, L_{I_2}\big),$$

where the summation is over decompositions of pairs of hyperbolic structures and $I_1$ and $I_2$ are unordered sets which provide a partition. The $\hat{L}$ refers to the fact that we are referring to a moduli volume of a subsurface, since $\hat{L} = (L_2, \dots, L_n)$. Without getting too lost in the details, the key fact is that the WP volume of the moduli space can be written as an integral of volumes of surfaces whose decomposition involves one less pair of pants.

One can then derive a general integration recursion. If one starts with a multi curve

$$\gamma = \sum_{j=1}^{m} a_j \gamma_j$$

and some function $f$ which becomes small at infinity, one can write a sum which is invariant under the mapping class group:

$$f_\gamma(R) = \sum_{\mathrm{MCG}/\mathrm{Stab}(\gamma)} f\left(\sum_{j=1}^{m} a_j l_{h(\gamma_j)}(R)\right).$$

The stabilizer in this case is the subgroup of the mapping class group which stabilizes the weighted curves. For a weighted multi curve and a MCG-sum of a function, Mirzakhani showed that one has

$$\int_{\tau(R)/\mathrm{MCG}} f_\gamma\, dV = (|Sym(\gamma)|)^{-1} \int_{\mathbb{R}_{>0}^m} f(|\mathbf{x}|)\, V(R(\gamma); \mathbf{x})\, \mathbf{x} \cdot d\mathbf{x},$$

where $|\mathbf{x}| = \sum_j a_j x_j$, $\mathbf{x} \cdot d\mathbf{x} = x_1 \dots x_m\, dx_1 \dots dx_m$ and $R(\gamma)$ denotes a hyperbolic surface $R$ cut open on the multi curve $\gamma$. The proof of this theorem uses a short exact sequence of mapping class groups.

This theorem can be combined with the McShane-Mirzakhani length identity to derive the recursion formula for the volume of a moduli space, since the length identity has the form of a MCG-sum. Once one has the Mirzakhani volume recursion, the next thing to do is to try to understand the integrals which are involved in the formula. The integrals in the recursion formula actually turn out to be polynomials whose coefficients are in $\pi\mathbb{Q}^+$. Mirzakhani showed that the WP volume $V_g(L)$ is a polynomial of degree $6g-6+2n$ in the boundary lengths which can be written as[257]

$$V_g(L) = \sum_{|j| \leq 3g-3+n} v_j L^{2j},$$

where $L$ is the tuple of boundary lengths, $j$ is a multi-index, $|j|$ is a summation, and the coefficient is in $\pi^{6g-6+2n-|j|}\mathbb{Q}$. The polynomial coefficients here are intersection numbers for tautological classes on the compactified moduli space. Symplectic reduction arguments show that there is an expansion for the WP volume in terms of tautological intersection numbers. The WP volume relative to the symplectic form $2\omega_L$ is a polynomial of degree $3g-3+n$ in the squares of the boundary lengths $L_1^2, \ldots, L_n^2$.

$$\frac{V_g(2\pi L)}{(2\pi^2)^d} = \frac{1}{d!}\int_{\overline{\mathcal{M}}_{g,n}} \left(\kappa_1 + \sum_{j=1}^{n} L_j^2 \psi_j\right)^d$$
$$= \sum_{d_0+\cdots+d_n=d} \prod_{j=0}^{n} \frac{1}{d_j!} \langle \kappa_1^{d_0} \prod_{j=1}^{n} \tau_{d_j} \rangle_{g,n} \prod_{j=1}^{n} L_j^{2d_j},$$

where $\tau_{d_j} = \psi_j^{d_j}$ and the angle brackets denote the integration of a $2d$-form in the Deligne-Mumford compactification of the moduli space[258].

One can consider the following intersection pairing relation:

$$\langle \tau_{d_1} \ldots \tau_{d_n} \rangle_g = \int_{\overline{\mathcal{M}}_{g,n}} \psi_1^{d_1} \ldots \psi_n^{d_n} \;\text{ for } d_1 + \cdots + d_n.$$

The generating function for the $\langle \tau_{d_1} \ldots \tau_{d_n} \rangle_g$ intersection numbers (originally emerging as a partition function when Witten was studying two-dimensional gauge theories in a physical context) is defined as

$$\mathbf{F} = \sum_{g=0}^{\infty} \lambda^{2g-2} F_g,$$

where

[257] Maryam Mirzakhani, 'Simple geodesics and Weil-Petersson volumes of moduli spaces of bordered Riemann surfaces', *Invent. Math.,* 167(1), 179 – 222 (2007).
[258] Maryam Mirzakhani, 'Weil-Petersson volumes and intersection theory on the moduli space of curves', *J. Amer. Math. Soc.,* 20(1), 1 – 23 (2007).

$$F_g(t_0, t_1, \dots) = \sum_{\{d_j\}} \langle \prod_{j=1}^{n} \tau_{d_j} \rangle_g \prod_{r=1}^{\infty} \frac{t_r^{n_r}}{n_r!}.$$

The summation is over all sequences of non-negative integers with a finite number of non-zero terms. You will notice that the generating function only encodes the intersection numbers of the tautological classes $\psi_j$. The Witten conjecture states that the generating function satisfies sets of PDEs known as the Virasoro constraint equations and the Korteweg-de Vries hierarchy. The conjecture was proved by Kontsevich. Mirzakhani took a different approach to the conjecture, applying her volume recursion formula and computing the integrals for integration over one and two boundary lengths

$$\int_0^{\infty} x^{2j+1} H(x,t) dx, \qquad \int_0^{\infty} \int_0^{\infty} x^{2j+1} y^{2k+1} H(x+y,t) \, dx \, dy$$

to derive a result which is equivalent to the Witten-Kontsevich formula[259].

You may be familiar with the prime number theorem from analytic number theory. It was observed many years ago that the density of the primes $p \leq x$ appears to be as $(\log x)^{-1}$. The prime number theorem states that if $\pi(x)$ is the number of primes with magnitude at most $x$, then as $x \to \infty$ we have an asymptotic relation

$$\pi(x) \sim \frac{x}{\log x}.$$

Gauss observed that a better approximation to $\pi(x)$ is given by a singular integral with an asymptotic expansion for $x > 1$.

$$\mathrm{Li}(x) = \int_0^{x} \frac{dy}{\log y} = \int_1^{x} \left(1 - \frac{1}{y}\right) \frac{dy}{\log y} + \log \log \mathrm{x} + \gamma.$$

In his attempts to prove the theorem, Tchebyshev realized that it was easier to work with a function

$$\theta(x) = \sum_{p \leq x} \log p.$$

What is even better is the Tchebyshev function $\psi(x)$ obtained by taking the average

$$\psi(x) = \sum_{n \leq x} \Lambda(n),$$

where $\Lambda(n)$ is the von Mangoldt function defined as

$$\Lambda(n) = \begin{cases} \log p\,, & \text{if } n = p^k \text{ for some prime } p \text{ and integer } k \geq 1 \\ x, & \text{otherwise.} \end{cases}$$

[259] Scott Wolpert, *Families of Riemann Surfaces and Weil-Petersson Geometry* (USA: American Mathematical Society, 2009).

If we use this function, the prime number theorem becomes the asymptotic statement

$$\psi(x) \sim x.$$

The first proofs of the theorem by Hadamard and de la Valleé Poisson were analytic, having now been cleaned up and streamlined into a standard 'proof via complex analysis'. The proof is intricate in its argument, but quite pleasing, as it uses all the basic tools such as contour integrals and the Cauchy residue theorem which you might be familiar with from a course on elementary complex analysis. There is also an elementary proof by Selberg and Erdős, but it is harder to follow than the analytic proof and only elementary in the sense that it does not require complex analysis. Both proofs by Selberg and Erdős rely on the Selberg formula[260]

$$\sum_{p \leq \mathrm{x}} (\log p)^2 + \sum_{pq \leq x} \sum (\log p)(\log q) = 2x\log x + O(x).$$

One can count closed geodesics on compact Riemann surfaces, so is there analogous version of the prime number theorem on this type of surface? The prime geodesic theorem states for a compact hyperbolic surface, the counted number of closed geodesics $\mathbf{c}(L)$ and the count of the closed primitive geodesics $\mathbf{c}_{\text{primitive}}(L)$ with length at most $L$ satisfy the following asymptotic relation as $L \to \infty$

$$\mathbf{c}(L) \sim \mathbf{c}_{\text{primitive}}(L) \sim \frac{e^L}{L}.$$

In the general case of compact negatively curved $n$-manifolds it can be shown that the same count satisfied the following asymptotic relation as $L \to \infty$

$$\mathbf{c}(L) \sim \frac{e^{hL}}{hL},$$

where $h$ is the topological entropy of the geodesic flow. For a multi curve on a hyperbolic surface $R$, one can define the count of the multi curves with length at most $L$ within the mapping class group orbit of the particular curve $\gamma$

$$s_R(L, \gamma) = \#\{\alpha \in \mathrm{MCG} \text{ such that } l_\alpha(R) \leq L\},$$

where $l_\alpha = \sum_{j=1}^{m} a_j l_{\alpha_j}$ for $\alpha = \sum_{j=1}^{m} a_j \alpha_j$. The count of all simple closed geodesics is defined as the following summation[261]:

$$s_R(L) = \sum_{\gamma \in \mathcal{S}/\mathrm{MCG}} s_R(L, \gamma).$$

Mirzakhani used the integrals from her volume recursion formula and the properties of the measured geodesic laminations which we mentioned earlier to arrive at the beautiful prime simple geodesic theorem. This states that for a rational multi curve $\gamma$, there is a positive constant $c(\gamma)$ depending only on $\gamma$ such that

$$\lim_{L\to\infty} \frac{s_R(L,\gamma)}{L^{6g-6+2n}} = \frac{c(\gamma)\mathbf{B}(R)}{\mathbf{b}(R)},$$

where

$$\mathbf{B}(R) = \mu_{\text{Thurston}}(\mathbf{B}_R), \quad \mathbf{b}(R) = \int_{\mathcal{M}} \mathbf{B}(R)\, dV.$$

[260] Henryk Iwaniec and Emmanuel Kowalski, *Analytic Number Theory* (USA: American Mathematical Society, 2004).
[261] Scott Wolpert, *Families of Riemann Surfaces and Weil-Petersson Geometry* (USA: American Mathematical Society, 2009).

The constant is computed using the weights for the multi curve and the coefficients of the volume of the moduli space. As usual, this is for a hyperbolic surface $R(\gamma)$ which is cut open on the curve $\gamma$. For one single curve, the constant can be written as

$$c(\gamma) = \big(|Sym(\gamma)|(6g-6+2n)\big)^{-1} \prod_{R(\gamma)\ \text{components}} V(R';L)_{\text{leading}}.$$

Mirzakhani provided explicit examples of the computation of $c(\gamma)$ in the case of genus 0 and also for separating curves of compact surfaces. In the genus 0 case, the surface can be divided up by a non-trivial simple closed curve $\gamma_j$ into a pair of genus 0 surfaces, one with $j$ cusps and one boundary and the other with $n-j$ cusps and one boundary (assuming that the original surface had greater than 3 cusps). It can be shown that one obtains

$$c(\gamma_j) = \frac{1}{2^{n-4}(j-2)!\,(n-j-2)!\,(2n-6)2^{\delta_{j,n-j}}}.$$

In the case of a separating closed curve, a compact surface of genus $g$ with be divided into a genus $j$ surface with one boundary and another surface of genus $g-j$ with one boundary. The resulting constant is then quite similar:

$$c(\gamma_j) = \frac{1}{2^{3g-4}24^g j!\,(g-j)!\,(3j-2)!\,(3(g-j)-2)!\,(6g-6)2^{\delta_{j,g-j}}}.$$

Things become more complicated when counting multi curves. Mirzakhani showed that for a multi curve $\gamma$, the WP integral of $s_R(L,\gamma)$ is a polynomial

$$P(L,\gamma) = (|Sym(\gamma)|)^{-1} \int_{\left\{0\le\sum_{j=1}^{m} a_j x_j \le L,\, x_j\ge 0\right\}} V(R(\gamma);\mathbf{x})\,\mathbf{x}\cdot d\mathbf{x},$$

where $V(R(\gamma);\mathbf{x})$ is the volume of the moduli space for cut open surfaces as usual and $\mathbf{x}$ is a tuple of boundary lengths. In the case of a rational multi curve, Mirzakhani proved that the scaled orbit measures converge weak-$*$ as follows:

$$\mu_{T,\gamma} \longrightarrow \frac{c(\gamma)}{\mathbf{b}(R)}\mu_{\text{Thurston}},$$

where the constant is the leading coefficient in the polynomial which we just wrote down. We will mention one final remarkable result of Mirzakhani's. For a multi curve $\gamma = \vec{\alpha}\cdot\vec{\gamma}$, the constant is given by a formula

$$c(\gamma) = \sum_{\vec{j}, |\vec{j}| \leq 3g-3+n-m(\gamma)} \frac{v_{\vec{j}}(\gamma)}{|Sym(\gamma)|(6g-6+2n)!} \prod_{i=1}^{m(\gamma)} \frac{(2j_i+1)!}{a_i^{2j_i+2}}.$$

The moduli space integral over Thurston volume is

$$\mathbf{b}(R) = \sum_{\gamma \in \mathcal{MGL}_{\text{unit}}/\text{MCG}} \sum_{\vec{j}, |\vec{j}| \leq 3g-3+n-m(\gamma)} \frac{v_{\vec{j}}(\gamma)}{|Sym(\gamma)|(6g-6+2n)!} \prod_{i=1}^{m(\gamma)} (2j_i+1)!\, \zeta(2j_i+1),$$

where $\zeta$ is the Riemann zeta function[262].

[262] Maryam Mirzakhani, ‘Growth of the number of simple closed geodesics on hyperbolic surfaces’, *Ann. of Math. (2)*, 168(1), 97 – 125 (2008).

# *Maryna Viazovska*

Viazovska is well known for work on sphere packing and combinatorial design theory. She spectacularly solved the sphere-packing problem in dimensions $8$ and $24$. This was something of a surprise, as the solution of the problem for dimension less than or equal to $3$ involved a computer-assisted proof and many arduous computer calculations, whereas Viazovska used the kinds of more traditional arguments which one would associate with a proof. There is still some debate about the necessity of computer-verified proofs and the future role that computers and machine learning will play in the business of formal proof theory. Certainly, it is possible that computers will be used routinely in the future to prove lemmas which would otherwise be very tedious to prove: this type of lemma is common in commutative algebra, for example. Computer scientists have been using formal proof verification tools for decades, but in general, the mathematical community has not paid much attention. The fact that Viazovska used arguments with modular forms to solve cases of a problem which had previously only been thought amenable to computer-assisted proof suggests that there might not necessarily be an *essential* place for this type of proof in the future of mathematics, and that adaptation of these methods might depend on whether mathematicians wish to use them or find them helpful, rather than on any actual necessity.

Maybe the most famous of the computer-assisted proofs are the proofs of the four-colour theorem and the Kepler conjecture. There are many expositions of both these results and the proofs in the popular mathematics literature, but we will go through them anyway[263]. The theorem originates in cartography. Perhaps you yourself have printed off a white copy of a map of some country divided up into states or counties and then seen how few colours you need to colour the whole thing in, assuming that states which are next to each other or have a common border are different colours. Experiment would seem to suggest that the minimum number of colours needed to do this is $4$, but how would one go about proving this? As always, we might begin by formulating the problem in mathematical language. Start with a map $\mathcal{L}$ and call $\mathcal{M}_{\mathcal{L}}$ the set of all the states of the map. One has $n$ colours with which to fill in these states. An $n$-colouring of the map is a mapping $\varphi\colon \mathcal{M}_{\mathcal{L}} \to \{1, \dots, n\}$. An $n$-colouring is said to be admissible if neighbouring countries which share a common border always have distinct function values (these function values corresponding to colours).

We might also want to consider permutations of the colours. Let $\varphi\colon \mathcal{M}_{\mathcal{L}} \to \{1, \dots, n\}$ be an admissible $n$-colouring of a map and $\pi\colon \{1, \dots, n\} \to \{1, \dots, n\}$ be a permutation, then it can be shown that the composition $\pi \circ \varphi$ is also an admissible $n$-colouring. Two colourings of a map are considered to be equivalent if they only differ by a permutation of colours. We can now state the four-colour theorem as follows: every map has an admissible $4$-colouring. Without going into the details, the basic approach to proving this theorem is investigate a kind of counterexample which is often referred to as a 'minimal criminal'. The idea of the counterexample is that if there does exist a map which cannot be coloured with four colours, then there has to be a map of this kind containing the fewest number $m$ of countries. A map with four countries can trivially be coloured in with four colours (or assigned numbers from a set of four numbers, if you prefer), so the fewest number of countries has to be bigger than four. A map which has fewer than $m$ countries admits a 4-colouring. A map with $m$ countries that does not admit a 4-colouring is called a minimal criminal in the literature. The aim is to show that such a minimal criminal cannot exist. The problem can also be formulated as a combinatorial problem with no reference to geometry[264]. A similar problem in combinatorics known as the Hadwiger-Nelson problem asks for the minimum number of colours which are need to colour the plane such that no two points at unit distance from each other have the same colour. Recent progress on this issue was made by de Grey, who showed that the minimum number is at least 5 by presenting a number of finite unit-distance graphs in the plane which are not 4-colourable. The smallest of these has 1581 vertices. This proof was also computer-assisted[265].

[263] Ian Stewart, *The Great Mathematical Problems: Marvels and Mysteries of Mathematics* (London: Profile Books, 2014).

The Kepler conjecture states that the densest possible packing of $\mathbb{R}^3$ by equal spheres is via what is now called the FCC (face-centred cubic) packing, which fills space with a density of $\pi/\sqrt{18} \approx 74\%$. The conjecture originally emerged in the applied setting of efficient stacking of cannonballs and it is part of Hilbert's eighteenth problem. We will not go into the proof here, as it is extremely complicated and would require hundreds of pages even for a sketch. Again, the proof given by Hales and Ferguson is computer-assisted. It was maybe the first proof of its kind at the time and caused a huge amount of controversy and soul-searching about what it really means to say that something has been proved, since the referees of the paper were not able to state certainty about the correctness of the many calculations which had been carried out with a computer[266].

The FCC lattice is one of several lattice types which are used to describe the crystal structures found in Nature. A crystal is a solid structure in which a group of atoms is arranged in a repeating pattern over a regular three-dimensional grid called a lattice. In fact, Gauss had already proved the Kepler conjecture over a regular lattice prior to the modern proofs. The FCC lattice is a simple cubic lattice, but with an additional lattice point at the centre of every face of every cube. A displacement vector which joins two lattice points is called a lattice vector. If we assume that the structure can be modelled as an infinite crystal, then any displacement by a lattice vector will produce no change in the crystal ie. there is translational symmetry. This implies that the potential energy function for an electron in a crystal obeys a periodicity condition for any lattice vector. This implies Bloch's theorem: in an infinite crystal structure, the energy eigenfunctions of electrons can be written in the form

$$\psi(\mathbf{r}) = e^{i\mathbf{k}\cdot\mathbf{r}} u_{\mathbf{k}}(\mathbf{r}),$$

[264] Rudolf Fritsch and Gerda Fritsch, *The Four-Color Theorem: History, Topological Foundations, and Idea of Proof* (New York: Springer-Verlag, 1998).
[265] Aubrey de Grey, 'The chromatic number of the plane is at least 5', arXiv:1804.02385v3 (2018).
[266] Thomas Hales and Samuel Ferguson, *The Kepler Conjecture: The Hales-Ferguson Proof* (New York: Springer, 2010).

where $\mathbf{k}$ is the wave vector and $u_{\mathbf{k}}(\mathbf{r})$ is a function with the periodicity of the lattice

$$u_{\mathbf{k}}(\mathbf{r}+\mathbf{R}) = u_{\mathbf{k}}(\mathbf{r}) \text{ for any lattice vector } \mathbf{R}.$$

The theorem tells us that the energy eigenfunction for an electron in an infinite crystal is the product of the plane wave associated with a free electron multiplied by a periodic function. This explains an old mystery of how it is that electrons appear to move freely through wires even though there must be atoms obstructing them (something which the classical theory cannot explain). In a perfect crystal, the electron eigenfunction extends without stopping as a plane wave (only modulated by the periodic function), so the electron wave is not obstructed by the densely packed atoms around it.

We can apply the Bloch theorem to show trivially that the probability density in a Bloch wave also has the periodicity associated with the lattice. The electron probability density is given by

$$\begin{aligned}|\psi(\mathbf{r})|^2 &= \psi^*(\mathbf{r})\psi(\mathbf{r}),\\ &= \left(e^{i\mathbf{k}\cdot\mathbf{r}}u_{\mathbf{k}}^*(\mathbf{r})\right)\left(e^{i\mathbf{k}\cdot\mathbf{r}}u_{\mathbf{k}}(\mathbf{r})\right),\\ &= |u_{\mathbf{k}}(\mathbf{r})|^2.\end{aligned}$$

The result then follows. The values for the wave vector are restricted to a discrete set using appropriate periodic boundary conditions[267]. When considering the unit cell for the FCC lattice, there are eight lattice points on the corners of the cell and a lattice point in the centre of each of the six faces, so the unit cell contains four lattice points. If we fill Euclidean space with FCC unit cells, one sees that the lattice points of an extended FCC lattice are described as points with the usual coordinates where the coordinates are $a\mathbb{Z}$ for the lattice constant $a$, or two of the three coordinates are half-odd integers multiplying $a$, with the remaining coordinate $a\mathbb{Z}$ as before. The Wigner-Seitz cell for the FCC lattice takes the shape of a rhombic dodecahedron, where every face is the perpendicular bisector between the central point and one of its neighbours. One uses these bisectors to construct the cell, where the Wigner-Seitz cell of a point of a lattice is defined to be the set of all points which are closer to that point than any other lattice point. A simpler lattice is the simple cubic lattice, but most chemical elements cannot form a simple cubic lattice in atomic form, although there is one which does (polonium)[268].

We have mentioned that the Kepler conjecture was proved with an extremely complicated argument which made liberal use of computer calculations. It was therefore rather unexpected when Viazovska announced short proofs of the corresponding conjecture in dimensions $8$ and $24$. The proof in dimension $8$ was done using only modular forms and Fourier analysis and employs an argument which is relatively simple in places, utilizing standard complex analysis and transform theory and quoting results like the Fourier transform of the Gaussian. For $k$ an integer, a function is said to be weakly modular of weight $-2k$ if $f$ is meromorphic on the half-plane $\mathbb{H}$ and if it satisfies the following relation:

[267] John Bolton and Stuart Freake (ed.), *Quantum Mechanics of Matter* (Milton Keynes: The Open University, 2009).

[268] Steven Simon, *The Oxford Solid State Basics* (Oxford: Oxford University Press, 2013).

$$f(z) = (cz+d)^{-2k} f\left(\frac{az+b}{cz+d}\right),$$

for all matrices in $SL(2,\mathbb{Z})$. You will probably recognise the rational function of $z$ as a Möbius transformation of the complex plane. Furthermore, if $f$ is meromorphic on the half-plane, $f$ is a weakly modular function of weight $2k$ if and only it satisfies two relations:

$$f(z+1) = f(z),$$

$$f\left(-\frac{1}{z}\right) = z^{2k} f.$$

A weakly modular function is called modular if it is meromorphic at infinity and a modular function which is holomorphic everywhere is called a modular form. If such a function is zero at infinity, the modular form is called a cusp form. From complex analysis, we know that a modular form of weight $2k$ is given by a series

$$f(z) = \sum_{n=0}^{\infty} a_n q^n = \sum_{n=0}^{\infty} a_n\, e^{2\pi i n z},$$

which has unit radius of convergence. This modular form will be a cusp form if $a_0 = 0$.

If $f_1$ and $f_2$ are two modular forms with weight $2k_1$ and $2k_2$ respectively, then the product $f_1 f_2$ is another modular form of weight $2k_1 + 2k_2$. If we denote by $M$ the set of pairs $(\omega_1, \omega_2)$ of elements of $\mathbb{C}^*$ such that $\mathrm{Im}(\omega_1/\omega_2) > 0$, two elements of $M$ define the same lattice (not related to the lattices we were talking about before) if and only if they are congruent modulo $SL(2,\mathbb{Z})$. This implies that the set of lattices $\mathcal{R}$ of $\mathbb{C}$ can be identified with the quotient of $M$ by the group action of $SL(2,\mathbb{Z})$. We also have that the map given by

$$(\omega_1, \omega_2) \longmapsto \omega_1/\omega_2$$

furnishes us with a bijection of $\mathcal{R}/\mathbb{C}^*$ onto the quotient $H/G$ after passing to the quotient. This implies that an element of $H/G$ can be identified with a lattice of $\mathbb{C}$ modulo homotheties. The quotient $\mathcal{R}/\mathbb{C}^*$ can also be viewed as the set of isomorphism classes of elliptic curves. For a relatively simple example of a modular function, take an integer $k > 1$. If $\Gamma$ is a lattice of $\mathbb{C}$, define the following series:

$$G_k(\Gamma) = \sum_{\gamma\in\Gamma}{}' \frac{1}{\gamma^{2k}}.$$

It can be shown that this is an absolutely convergent series of weight $2k$, known as the Eisenstein series of index $2k$. As this series is a function defined on the space $M$, we have

$$G_k(\omega_1, \omega_2) = \sum_{m,n}{}' \frac{1}{(m\omega_1 + n\omega_2)^{2k}},$$

where the summation is over all pairs of integers which are distinct from the pair $(0,0)$. Evaluating for some $z$, one has

$$G_k(z) = \sum_{m,n}{}' \frac{1}{(mz+n)^{2k}}.$$

One can show that the Eisenstein series $G_k(z)$ is a modular form of weight $2k$, and that evaluated at infinity we have

$$G_k(\infty) = 2\zeta(2k),$$

where $\zeta$ is the Riemann zeta function. To prove this, you will need to show that the Eisenstein series is everywhere holomorphic (including at infinity). There is an interesting link between the Eisenstein series and the Weierstrass $p$-function. Start with a lattice over $\mathbb{C}$ and take the corresponding Weierstrass $p$-function:

$$\mathfrak{P}_\Gamma(u) = \frac{1}{u^2} + \sum_{\gamma\in\Gamma}{}' \Big(\frac{1}{(u-\gamma)^2} - \frac{1}{\gamma^2}\Big).$$

Perform the Laurent expansion.

$$\mathfrak{P}_\Gamma(u) = \frac{1}{u^2} + \sum_{k=2}^{\infty} (2k-1)\, G_k(\Gamma) u^{2k-2}.$$

Setting $x = \mathfrak{P}_\Gamma(u)$ and $y = \mathfrak{P}'_\Gamma\,(u)$, we end up with a cubic curve

$$y^2 = 4x^3 - g_2 x - g_3,$$

where $g_2$ and $g_3$ are multiples of the Eisenstein series of lowest weights. This is a non-singular cubic curve, and one can show that in the projective plane this curve is isomorphic to the elliptic curve $\mathbb{C}/\Gamma$[269].

In more formal terms, Viazovska proved that no packing of unit balls in $\mathbb{R}^8$ can have a density which is greater than that of the $E_8$-lattice packing. We are interested in a quantity called the sphere packing constant which measures how much of Euclidean $n$-space can be covered with unit balls which can touch without overlapping. We start with a standard open ball in Euclidean space denoted by $B_d(x, r)$, where $x$ is the centre and $r$ is the radius. If $X \subset \mathbb{R}^d$ is a discrete set of points such that $\|x - y\| \geq 2$, then the union of all the unit balls

$$\mathcal{P} = \bigcup_{x\in X} B_d(x, 1)$$

is technically what we mean by a sphere packing. If $X$ is a lattice, then we say that $\mathcal{P}$ is a lattice sphere packing (this is the kind of sphere packing which Gauss studied many years ago). The finite density of a packing $\mathcal{P}$ is defined to be the following ratio of volumes:

$$\Delta_{\mathcal{P}}(r) := \frac{\mathrm{Vol}\big(\mathcal{P} \cap B_d(0, r)\big)}{\mathrm{Vol}\big(B_d(0, r)\big)}.$$

[269] Jean-Pierre Serre, *A Course in Arithmetic* (New York: Springer, 1973).

The density of the packing is defined to be the limit supremum:

$$\Delta_{\mathcal{P}} := \limsup_{r\to\infty} \Delta_{\mathcal{P}}(r).$$

The particular quantity which we search for is the sup over the set of all possible packing densities:

$$\Delta_d := \sup_{\mathcal{P}\subset\mathbb{R}^d} \Delta_{\mathcal{P}}.$$

This is known as the sphere packing constant. In the case of $1$ dimension, we trivially have

$$\Delta_1 = 1,$$

so no prizes for that one. In the Euclidean plane, as you might intuitively expect, the most efficient packing is a hexagonal lattice, where a disc touches six other discs. The density of the hexagonal lattice packing is given by

$$\Delta_2 = \frac{\pi}{\sqrt{12}} \approx 91\%.$$

This is much higher than the corresponding figure which we had for $\mathbb{R}^3$ which was proved by Hales and Ferguson:

$$\Delta_3 = \frac{\pi}{\sqrt{18}} \approx 74\%.$$

It is interesting to speculate what this packing constant would be for the case of $\mathbb{R}^8$. Viasovska proved that it is

$$\Delta_8 = \frac{\pi^4}{384} \approx 25\%.$$

This corresponds to the density of the $E_8$-lattice sphere packing. The $E_8$-lattice $\Lambda_8$ is given by

$$\Lambda_8 = \left\{(x_i) \in \mathbb{Z}^8 \cup \left(\mathbb{Z} + \frac{1}{2}\right)^8 \text{ such that } \sum_{i=1}^{8} x_i \equiv 0 \pmod 2\right\}.$$

This is a subset of $\mathbb{R}^8$. This lattice is the unique, positive definite, even, unimodular lattice of rank $8$ (unique modulo isometries). If you are familiar with root systems, the name is due to the fact that it is the root lattice of the $E_8$ root system. The minimum distance between two points in $\Lambda_8$ is $\sqrt{2}$ and the $E_8$-lattice sphere packing is the packing of unit balls whose centres are situated at $\frac{1}{\sqrt{2}}\Lambda_8$.

Viasovska's proof that no packing of unit balls in $\mathbb{R}^8$ has a density greater than the density of the $E_8$-lattice packing involves the technique of linear programming bounds. This technique enables one to obtain upper bounds in many discrete optimization problems, although these bounds are not always particularly sharp. One might hope that there would be linear programming bounds which could be applied specifically to sphere packings, and that turns out to be the case. Viasovska used these techniques (known as Cohn-Elkies linear programming bounds), including a result which is worth quoting. Take $f: \mathbb{R}^d \to \mathbb{R}$ to be some non-zero admissible function. If the function satisfies the following:

$$f(x) \leq 0 \text{ for } \|x\| \geq 1$$

and

$$\hat{f}(x) \geq 0 \text{ for all } x \in \mathbb{R}^d,$$

then the density of $d$-dimensional sphere packings is bounded from above by

$$\frac{f(0)}{\hat{f}(0)} \frac{\pi^{\frac{d}{2}}}{2^d \Gamma\left(\frac{d}{2}+1\right)} = \frac{f(0)}{\hat{f}(0)} \mathrm{Vol}\, B_d\left(0, \frac{1}{2}\right).$$

Combined with Fourier analysis and some theory of modular forms, this yields a relatively short proof for the case of Euclidean $8$-space[270].

Viasovska and her colleagues would later return to this question, but now for $\mathbb{R}^{24}$. Based on the proof in eight dimensions, they studied a lattice called the Leech lattice and proved that the Leech lattice achieves the optimal sphere packing density in $\mathbb{R}^{24}$, and that it is the only periodic packing in $\mathbb{R}^{24}$ with that density, modulo scaling and isometries. It follows that the optimal sphere packing density in $\mathbb{R}^{24}$ is equal to that of the Leech lattice. So, what might this number be? It must be close to zero, surely? In fact, it is[271]

$$\Delta_{24} = \frac{\pi^{12}}{12!} \approx 0.19\%.$$

Linking back to symplectic topology, there is also a symplectic packing problem. Given $k \geq 1$, we would like to find the maximum of all the radii $r$ (denoted by $r_k$) for which there is a symplectic embedding of $k$ disjoint balls of radius $r$ into the symplectic manifold $(M, \omega)$. $(M, \omega)$ can be filled completely by $k$ balls if the maximum radius $r_k$ is such that the volume of $k$ balls of radius $r_k$ is equal to the volume of $M$ with respect to the symplectic volume form $\omega^n/n!$ Informally speaking, the problem asks how much of the volume of a symplectic manifold can be filled up with disjoint embedded open symplectic balls. Significant progress has been made on this problem, but nothing is known beyond dimension six, due to many of the relevant techniques only being applicable in dimension four[272].

One particular case of the symplectic packing problem which has been studied is the torus equivalent case, where the symplectic manifold and the standard symplectic open ball are equipped with a Hamiltonian action of an $n$-torus and the symplectic embeddings of the ball into the manifold are equivariant with respect to these actions. A result which has been found in this case is that if $M$ is a $2n$-dimensional Delzant manifold $M$, then $M$ admits a perfect equivariant symplectic ball packing if and only if there exists some $\lambda > 0$ such that $M$ is equivariantly symplectomorphic to either the complex projective plane or the product of two complex projective lines (where the symplectic form is multiplied by $\lambda$ in both case) when $n = 2$, or such that $M$ is equivariantly symplectomorphic to the complex projective $n$-plane when $n \neq 2$. Another way of stating this is that $(\mathbb{CP}^n, \lambda\omega_{\mathrm{FS}})$ and $(\mathbb{CP}^1 \oplus \mathbb{CP}^1, \lambda\omega_{\mathrm{FS}} \oplus \lambda\omega_{\mathrm{FS}})$ are the only Delzant manifolds which admit a perfect, equivariant, symplectic ball packing[273].

[270] Maryna Viazovska, 'The sphere packing problem in dimension 8', *Annals of Mathematics*, 185(3), 991 – 1015 (2017).
[271] Henry Cohn, Abhinav Kumar, Stephen Miller, Danylo Radchenko and Maryna Viazovska, 'The sphere packing problem in dimension 24', *Annals of Mathematics*, 185(3), 1017 – 1033 (2017).
[272] Dusa McDuff and Dietmar Salamon, *Introduction to Symplectic Topology* (New York: Oxford University Press, 1998).

Viazovska has worked on the theory of spherical designs. These sound esoteric, but have been taken up and applied in approximation theory and quantum mechanics. Start with the unit sphere $S^n$ in $\mathbb{R}^{n+1}$ with the Lebesgue measure $\mu_n$ normalized to unity such that

$$\mu_n(S^n) = 1.$$

A set of points $x_1, \dots, x_N$ in $S^n$ is called a spherical $t$-design if

$$\int_{S^n} P(x)\, d\mu_n(x) = \frac{1}{N} \sum_{i=1}^{N} P(x_i),$$

for all polynomials in $n + 1$ variables over some algebraic field, with total degree at most $t$. Denote by $N(n, t)$ the minimal number of points in a spherical $t$-design in the $n$-sphere. It had already been known since the 1970s that there is a lower bound on this quantity:

$$N(n,t) \geq \begin{cases} \binom{n+k}{n} + \binom{n+k-2}{n}, & t = 2k, \\ 2\binom{n+k}{n}, & t = 2k+1. \end{cases}$$

A spherical $t$-design which attains this lower bound is said to be tight. For a simple example, consider the vertices of a regular $(t + 1)$-gon: these form a tight spherical $t$-design in the circle. Eight tight spherical designs are known to exist for $n \geq 2$ and $t \geq 4$, but in general, tight spherical designs are rare and it has shown that they only exists for certain values of $t$ when $n \geq 2$.

Viazovska and her colleagues studied the asymptotic bounds on $N(n, t)$ for fixed $n \geq 2$ and $t$ tending to infinity in the limit. They used a theorem related to the Brouwer fixed point theorem familiar from elementary algebraic topology: the theorem says that if $f$ is a continuous self-map for Euclidean space and $\Omega$ an open bounded subset with boundary $\partial\Omega$ such that $0 \in \Omega$, then if $\big(x, f(x)\big) > 0$ for all $x \in \partial\Omega$, then there exists some $x \in \Omega$ such that $f(x) = 0$. This theorem is used to prove that for each $N \geq C_n t^n$, there exists a spherical $t$-design in $S^n$ consisting of $N$ points. As well as the previous theorem, the proof also makes use of the concept of an area-regular partition. Take $\mathcal{R} = \{R_1, \dots, R_n\}$ to be a finite collection of closed sets of the $n$-sphere such that

[273] Álvaro Pelayo, 'Toric symplectic ball packing', *Topology and its Applications*, 153, 3633 – 3644 (2006).

$$\bigcup_{i=1}^{N} R_i = S^n$$

and

$$\mu_n(R_i \cap R_j) = 0,$$

for all $1 \leq i < j \leq N$. The partition $\mathcal{R}$ is called area-regular if $\mu_n(R_i) = 1/N$. One also has a partition norm defined by

$$\|\mathcal{R}\| \coloneqq \max_{R \in \mathcal{R}} \operatorname{diam} R,$$

where $\operatorname{diam} R$ is the maximal geodesic distance between two points of $R$. One can show that for every $N$ there is an area-regular partition such that

$$\|\mathcal{R}\| \leq B_n N^{-\frac{1}{n}},$$

for sufficiently large $B_n$[274].

In a related paper, Viazovska pushed this approach further by proving that there exist certain configurations in the $n$-sphere which are spherical $t$-designs with a number of points which is asymptotically minimal in the sense which we just described, plus they have the best separation property in asymptotic terms. A sequence of $N$-point configurations $X_N = \{x_{1N}, \dots, x_{NN}\}$ in $S^n$ is said to be well-separated if

$$\min_{1 \leq i < j \leq N} |x_{iN} - x_{jN}| \geq \lambda_n N^{-\frac{1}{n}},$$

for a constant $\lambda_n$ and $n \geq 2$. Furthermore, there exists another constant $L_n$ such that

$$\min_{1 \leq i < j \leq N} |x_i - x_j| \geq L_n N^{-\frac{1}{n}}.$$

Viazovska and her colleagues proved that for every $n \geq 2$, there exist constants $C_n$ and $\lambda_n$ depending only on $n$ such that for each $N > C_n t^n$ there exists spherical $t$-design in $S^n$ consisting of $N$ points with the inequality

$$|x_i - x_j| \geq \lambda_n N^{-\frac{1}{n}}.$$

This is a direct generalisation of the previous theorem which we mentioned. (As with the proof of that result, the proof uses degree theory and area-regular partitions)[275].

Viazovska has also made contributions to pure solid state and semiconductor physics, although this is maybe not too surprising given the connections between packings, lattices and crystal structures which we have touched upon. For example, along with colleagues, she calculated the temperature-dependent London depth in the crystal $\mathrm{Ba_{1-x}K_xFe_2As_2}$ based on the electronic band structure and momentum-dependent superconducting gap derived from angle-resolved photoemission spectroscopy and muon spin rotation data. The London penetration depth is defined via the following formula:

$$\frac{1}{\lambda^2(T)} = \frac{e^2}{2\pi\varepsilon_0 c^2 h L_c}\int_{\mathrm{FS}} v_{\mathrm{F}}(\mathbf{k})\left(1-\int_{-\infty}^{+\infty}\left(-\frac{\partial f_T(\omega)}{\partial\omega}\right)\left|\mathrm{Re}\frac{\omega+i\Sigma''}{\sqrt{(\omega+i\Sigma'')^2-\Delta_{\mathbf{k}}^2(T)}}\right|d\omega\right)dk,$$

where $v_{\mathrm{F}}$ is the Fermi velocity, $\Delta_{\mathbf{k}}(T)$ is the momentum-dependent superconducting gap, $\Sigma''$ is the scattering rate, $dk$ is the length element for the Fermi surface, $T$ is temperature, $f_T(\omega)$ is the Fermi function, and the other fundamental physical constants are defined as usual[276].

---

[274] Andriy Bondarenko, Danylo Radchenko and Maryna Viazovska, 'Optimal asymptotic bounds for spherical designs', *Annals of Mathematics*, 178(2), 443 – 452 (2013).

[275] Andriy Bondarenko, Danylo Radchenko and Maryna Viazovska, 'Well-separated spherical designs', *Constructive Approximation*, 41(1), 93 – 112 (2015).

[276] D. V. Evtushinsky et al., 'Momentum-resolved superconducting gap in the bulk of $\mathrm{Ba}_{1-x}\mathrm{K}_x\mathrm{Fe}_2\mathrm{As}_2$ from combined ARPES and $\mu$SV', *New Journal of Physics*, 11(5), 55069 (2009).